%% file: thesis.tex
\documentclass[12pt,twoside]{report}
\usepackage{graphicx}
\usepackage{thesis}
\input{preamble}
\begin{document}
\input{TitleStuff}
\addcontentsline{toc}{chapter}{Table of Contents}
\tableofcontents
\addcontentsline{toc}{chapter}{List of Figures}
\listoffigures
\input{Chap1}

\input{Chap2}

\input{Chap3}

\input{AppendA}

\input{AppendB}

\input{AppendC}
\input{biblio}

\end{document}

%% file: preamble.tex
\renewcommand{\baselinestretch}{1.5}
	
\newcommand{\Tn}{{\rm{\bf T}^n}}
\newcommand{\An}{{\rm{\bf A}^n}}
\newcommand{\Rn}{{\rm{\bf R}^n}}

\newcommand{\Zn}{{\rm{\bf Z}^n}}
\newcommand{\T}{{\rm \bf T}}
\newcommand{\A}{{\rm \bf A}}
\newcommand{\Z}{{\rm \bf Z}}
\newcommand{\R}{{\rm \bf R}}

\newcommand{\norm}[1]{{\parallel {#1} \parallel}}
\newcommand{\pderiv}[2]{{\frac{\partial {#1}}{\partial {#2}}}}

\newcommand{\half}{{\frac{\scriptstyle 1}{\scriptstyle 2}}}
\newcommand{\tqrt}{{\frac{\scriptstyle 3}{\scriptstyle 4}}}
\newcommand{\oqrt}{{\frac{\scriptstyle 1}{\scriptstyle 4}}}

\newcommand{\cfStrut}{{\rule{0cm}{2.0ex}}}
\newcommand{\medStrut}{{\rule{0cm}{4.0ex}}}
\newcommand{\bigStrut}{{\rule{0cm}{6.0ex}}}
\newcommand{\hugeStrut}{{\rule{0cm}{7.0ex}}}
\newcommand{\emskip}{{\hskip 1em}}
\newcommand{\delt}[1]{{{\scriptstyle \Delta}{#1}}}
\newcommand{\id}{{\protect \bf I}}
\newcommand{\bfd}{{\protect \bf d}}
\newcommand{\bfm}{{\protect \bf m}}
\newcommand{\bfp}{{\protect \bf p}}
\newcommand{\bfq}{{\protect \bf q}}
\newcommand{\bfw}{{\protect \bf w}}
\newcommand{\bfL}{{\protect \bf L}}
\newcommand{\bfJ}{{\mbox{\protect \boldmath $J$}}}
\newcommand{\bfu}{{\mbox{\protect \boldmath $u$}}}
\newcommand{\bfv}{{\mbox{\protect \boldmath $v$}}}
\newcommand{\bfx}{{\mbox{\protect \boldmath $x$}}}
\newcommand{\bfy}{{\mbox{\protect \boldmath $y$}}}
\newcommand{\bpsi}{{\mbox{\protect \boldmath $\psi$}}}
\newcommand{\bfeta}{{\mbox{\protect \boldmath $\eta$}}}
\newcommand{\bbeta}{{\mbox{\protect \boldmath $\beta$}}}
\newcommand{\bgamma}{{\mbox{\protect \boldmath $\gamma$}}}
\newcommand{\bthet}{{\mbox{\protect \boldmath $\theta$}}}
\newcommand{\bveps}{{\mbox{\protect \boldmath $\varepsilon$}}}
\newcommand{\bfomega}{{\mbox{\protect \boldmath $\omega$}}}
\newcommand{\bLam}{{\mbox{\protect \boldmath $\Lambda$}}}
\newcommand{\lamM}[1]{{\lambda_{-}({#1})}}
\newcommand{\lamP}[1]{{\lambda_{+}({#1})}}
\newcommand{\Tr}[1]{{{\rm Tr}\,[{#1}]}}
\newcommand{\TrP}{{\rm Tr}_{max}}
\newcommand{\TrM}{{\rm Tr}_{min}}
\newcommand{\iprod}[2]{{\langle {#1}, {#2} \rangle}}
\newcommand{\rowSum}[2]{{[{#1}]_{#2\,\star}}}
\newcommand{\matSum}[1]{{[{#1}]_{\star\,\star}}}
\newcommand{\Veps}{{V_{\epsilon}}}
\newcommand{\Vtrg}{{V_{trig}}}
\newcommand{\bfxStar}{{\bfx^{\star}}}
\newcommand{\Tthet}{\tilde{\theta}}
\newcommand{\by}{{$\!\times\,$}}
\newcommand{\eInv}{{\epsilon_{inv}}}
\newcommand{\GT}{{\widetilde{G}}}
\newcommand{\App}{{A_{pp}}}
\newcommand{\Aup}{{A_{up}}}
\newcommand{\Avp}{{A_{vp}}}
\newcommand{\Auu}{{A_{uu}}}
\newcommand{\Auv}{{A_{uv}}}
\newcommand{\Avu}{{A_{vu}}}
\newcommand{\Avv}{{A_{vv}}}
\newcommand{\Ppp}{{P_{pp}}}
\newcommand{\Pup}{{P_{up}}}
\newcommand{\Pvp}{{P_{vp}}}
\newcommand{\Puu}{{P_{uu}}}
\newcommand{\Puv}{{P_{uv}}}
\newcommand{\Pvu}{{P_{vu}}}
\newcommand{\Pvv}{{P_{vv}}}

\newcommand{\Buu}{{B_{uu}}}
\newcommand{\Buv}{{B_{uv}}}
\newcommand{\Bvu}{{B_{vu}}}
\newcommand{\Bvv}{{B_{vv}}}
\newcommand{\ATpp}{{\widetilde{A}_{pp}}}
\newcommand{\ATup}{{\widetilde{A}_{up}}}
\newcommand{\ATvp}{{\widetilde{A}_{vp}}}
\newcommand{\ATuu}{{\widetilde{A}_{uu}}}
\newcommand{\ATuv}{{\widetilde{A}_{uv}}}
\newcommand{\ATvu}{{\widetilde{A}_{vu}}}
\newcommand{\ATvv}{{\widetilde{A}_{vv}}}
\newcommand{\AppInv}{{A^{-1}_{pp}}}

\newcommand{\AuvInv}{{A^{-1}_{uv}}}
\newcommand{\AvuInv}{{A^{-1}_{vu}}}

\newenvironment{typedef}[1]
{
\par
\vskip \abovedisplayskip
\hspace{1\parindent}
\begin{minipage}{0.9\textwidth }
    \singlespace
    \begin{tabbing}
    \tt typedef struct {#1} $\{$  \= \hspace{7em} \= \kill
    \tt typedef struct {#1} $\{$  
}{
    \end{tabbing}
\end{minipage}
\thesisspace \par 
}
\newcommand{\member}[2]{\> {\tt {#1}} \> {\tt {#2}} \\ }
\newcommand{\singlespace}{ \renewcommand{\baselinestretch}{1.0}
			   \protect \footnotesize
			   \protect \normalsize
}
\newcommand{\thesisspace}{ \renewcommand{\baselinestretch}{1.5}
			   \protect \footnotesize
			   \protect \normalsize
}

%% file: TitleStuff.tex
%
%
\pagenumbering{roman}
\thispagestyle{empty}
\begin{center}
{\bf Ghosts of Order on the Frontier of Chaos} \\
\vspace{1.0in}
Thesis by \\
Mark Muldoon \\
\vspace{1.0in}
In Partial Fulfillment of the Requirements \\
for the Degree of \\
Doctor of Philosophy \\
\vspace{1.75in}
California Institute of Technology \\
Pasadena, California \\
\vspace{0.75in}
1989\\
\vspace{0.5in}
(Submitted May, 1989) \\
\end{center}
\newpage
%
%
\singlespace
\begin{quote}
	Then from the heart of the tempest Yahweh spoke and gave Job 
	his answer.  He said:
	\begin{quote}
	    Brace yourself like a fighter; now it is my turn to ask
	    questions and yours to inform me. \newline
	    \newline
	    Where were you when I laid the earth's foundations? \newline
	    Who decided the dimensions of it? Do you know? \newline
	    Who laid its cornerstone when all the stars of morning 
	    were singing with joy? \newline \newline
	    Who pent up the sea when it leapt tumultuous out of
	    the womb, when I wrapped it in a robe of mist and made 
	    black clouds its swaddling bands? \newline \newline
	    Have you ever in your life given orders to the morning
	    or sent the dawn to its post? \newline \newline
	    Have you journeyed all the way to the sources of the
	    sea, or walked where the abyss 
	    is deepest? \newline \newline
	    Have you an inkling of the extent of the earth? \newline
	    Which is the way to the home of the light and where
	    does the darkness dwell?
	    \end{quote}

	    \raggedleft {\em The Jerusalem Bible}
\end{quote}
\vfil
\begin{quote}
    There are seven or eight categories of phenomena in the
    world that are worth talking about, and one of them is
    the weather.  Any time you care to get in your car and
    drive across the country and over the mountains, come
    into our valley, cross Tinker Creek, drive up the road to
    the house, walk across the yard, knock on the door and ask
    to come in and talk about the weather, you'd be welcome.

    \raggedleft {\em Annie Dillard}
\end{quote}
\vfil
\begin{quote}
     Then we would write the beautiful letters of the alphabet,
     invented by smart foreigners long ago to fool time and
     distance.

     \raggedleft {\em Grace Paley}
\end{quote}
\thesisspace
\newpage
%
%
\addcontentsline{toc}{chapter}{Acknowledgments}
\begin{center}
{\bf Acknowledgements} \\
\end{center}

     I offer my thanks to my advisor, Anatole  Katok, to my scientific 
correspondents, Jim Meiss, Robert MacKay, and Rafael de la Lllave,
and to Steve Wiggins of Caltech; without their many intellectual 
gifts I would have written a different, and lesser, thesis.

     More profoundly, I thank my friends, Susan Volman, Dave Wark,
Bette Korber, James Theiler, Paul Stolorz, Brian Warr, Chi-Bin Chien,
Dawn Meredith, Joel Morgan, Morgan Gopnik and Tom Bondy, and especially 
my mother and sister, Lucille and Maureen Muldoon; without
their love and encouragement I could not have written a thesis at all.

     Last, I thank Steve Frautschi for his patience and for providing
me, through {\em The Mechanical Universe}, with the most enjoyable summer
job I have ever had.  I also gratefully acknowledge Caltech's Concurrent
Computation Program, whose machines both performed my calculations and
typeset my thesis.
\newpage
%
%
\addcontentsline{toc}{chapter}{Abstract}
\begin{center}
    {\bf Abstract}\\
\end{center}

What kinds of motion can occur in classical mechanics?
We address this question by looking at the
structures traced out by trajectories in phase space;
the most orderely, completely integrable systems are characterized
by phase trajectories confined to low-dimensional, invariant tori.
The KAM theory examines what happens to the tori when an
integrable system is subjected to a small perturbation and finds that,
for small enough perturbations, most of them survive.

The KAM theory is mute about the disrupted tori, but,
for two dimensional systems, Aubry and Mather discovered an astonishing
picture: the broken tori are replaced by ``cantori,'' 
tattered, Cantor-set remnants of the original invariant curves.
We seek to extend Aubry and Mather's picture to higher 
dimensional systems and report two kinds of studies;
both concern perturbations of a completely integrable, 
four-dimensional symplectic map.  In the first study we compute
some numerical approximations to Birkhoff periodic orbits; sequences
of such orbits should approximate any higher dimensional analogs 
of the cantori.  In the second study we prove converse KAM theorems;
that is, we use a combination of analytic arguments and rigorous,
machine-assisted computations to find perturbations so large that
no KAM tori survive.  We are able to show that the last few of our
Birkhoff orbits exist in a regime where there are no tori.
\newpage

%% file: Chap1.tex
\input{introduction/intro1}

%% file: introduction/intro1.tex
\chapter{Introduction}
\pagenumbering{arabic}
\singlespace
\begin{quotation}

	There is a maxim which is often quoted, that ``The same causes
	will always produce the same effects.'' \dots

	It follows from this, that if an event has occured at a given 
	time and place it is possible for an event exactly similar
	to occur at any other time and place.

	There is another maxim which must not be confused with that
	quoted at the beginning of this article, which asserts ``That
	like causes produce like effects.''

	This is only true when small variations in the intial circumstances
	produce small variations in the final state of the system.  
	In a great many physical phenomena this condition is satisfied; 
	but there are other cases in which a small initial variation may 
	produce a very great change in the final state of the system, 
	as when the displacement of the ``points'' causes a railway
	train to run into another instead of keeping its proper course.

	\raggedleft {\em James Clerk Maxwell, 1877}
\end{quotation}
\thesisspace

\begin{figure}[bht]
	\begin{center}
		\includegraphics[height=6.5cm]{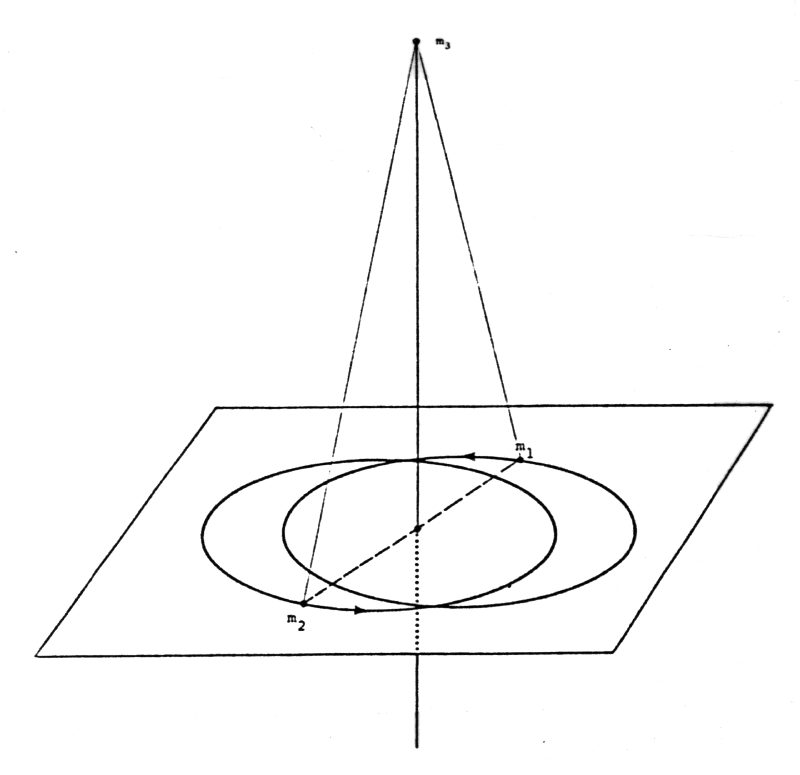}
	\end{center}
	\vspace*{-2\belowdisplayskip}
	\caption{ \em A system of two equally massive stars, 
	$m_1$ and $m_2$, and a test mass, $m_3$, which travels on 
	a line through the center of mass. {\rm [Moser73]} }
\label{fig:Moser}
\end{figure}

    Maxwell's warning, that like causes need not produce like
effects, can apply to even the simplest looking physical
systems.  Consider two equally massive stars bound in a binary system.  
Their orbits both lie in the same plane and, in a suitable coordinate
system, their center of mass is at rest at the origin.  If the
orbits are nearly (but not quite) circular the system will look like the 
one pictured in figure (\ref{fig:Moser}).
Now imagine adding a third body, a test mass so small that it does not 
disturb the motion of the stars.  Place the test mass at the origin 
and give it a velocity $v_0$ normal to the plane of the orbit.
The test mass will bob up and down on the line through the origin and, if the
initial velocity, $v_0$, is near enough to the escape velocity, the
subsequent motion of the test particle will display a fantastically sensitive 
dependence on the value of $v_0$; by suitable choice of $v_0$ one can 
arrange for test mass to begin in the orbital plane, spend $\approx s_1$
periods of the binary system above the plane, pass through to spend 
$\approx s_2$ periods below, then $\approx s_3$ above \dots and so on,
producing a sequence,
\begin{displaymath}
	\cdots s_0, s_1, s_2 \cdots,
\end{displaymath}
where each $s_j$ is an integer counting the number of 
complete periods  of
the binary which pass between visits by the test mass.
The $s_j$ can be chosen completely independently,  subject only to
the restriction $s_j > C$ for a constant $C$.

This system is described by Moser in \cite{Moser:book}.  
He begins his study  by drastically simplifying the problem; when
$t=0$ he notes the phase, $\theta_0$, of the binary orbit and the 
speed, $v_0$, of the test mass, then asks for $\theta_1$ and 
$v_1$, the corresponding phase and speed at the instant when the
test particle first returns to the orbital plane.  Certainly they
depend only on $\theta_0$ and $v_0$, so he constructs some
functions $\theta'(\theta,v)$ and $v'(\theta,v)$ such that 
\begin{displaymath}
     \theta_1 = \theta'(\theta_0, v_0) 
     \qquad {\rm and } \qquad
     v_1 = v'(\theta_0, v_0),
\end{displaymath}
then uses them to find a sequence, 
$ \cdots (\theta_0, v_0), (\theta_1, v_1) \cdots $, which captures 
the essential features of the dynamics.  Moser shows that the wild
behaviour described above occurs because the mapping,
\begin{equation}
(\theta,v) \rightarrow (\theta'(\theta,v), v'(\theta,v)),
\label{eqn:Moser}
\end{equation}
behaves like the 
celebrated horseshoe example of Smale, \cite{Smale:horse}.
Smale constructed the horseshoe by a process of abstraction; he
began by trying to understand the qualitative behaviour of a system
of differential equations\thesisfootnote{
    Smale gives a non-technical account of all this in one of the
    papers collected in \cite{Smale:time}.},
but eventually pared away most of the 
original problem, leaving a simple, illuminating model of the 
dynamics.
A detailed description 
of the horseshoe, along with a host of examples and criteria for
recognizing horseshoe-like behaviour, appear in \cite{Wiggins:book};
for us it will be enough to recognize that complicated dynamics 
arise even in simple classical systems and that these dynamics can
be explained in terms of structures in the phase space.
For the rest of the thesis we will be concerned with a different
relationship between structure and dynamics; we will be examine how the
highly structured phase space of an orderly classical system changes
under perturbation.  

\section{Integrability and the KAM Theorem}

\begin{figure}
	\begin{center}
		\includegraphics[height=5.5cm]{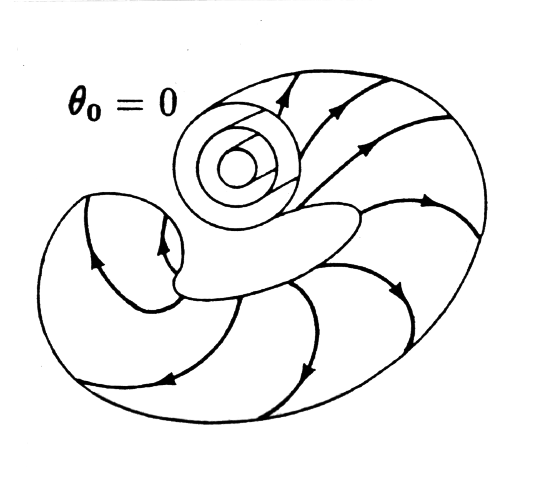}
	\end{center}
	\vspace*{-2\belowdisplayskip}
    \caption{ \em The phase space of a completely integrable system. 
		{\rm [Arn78] }}
    \label{fig:tori}
\end{figure}

    The most orderly of Hamiltonian systems are the {\em completely 
integrable} ones; these systems have so many constants of the motion,
($N$ for an $N$-degree-of-freedom system,) that we can reformulate the 
problem in terms of action-angle variables\thesisfootnote{
	We will use boldface symbols to denote $n\/$-dimensional
	objects, so that $\bthet$ is in $\Tn$, the $n$-dimensional 
	torus, $\bfp$ in $\Rn$.
	We will write 
	$\bthet_j$ for the angular coordinate of the $j\/$th image of 
	some phase point, $(\bthet_0, \bfp_0)$, and  
	$x_j$ (which is in ordinary type) for the real number 
	which is the 
	$j$th component of some $\bfx \in \Rn$.
	Ocassionally we will 
	need to express, ``the $k\/$th coordinate of the $j\/$th 
	image of the phase point $(\bthet_0, \bfp_0)$.'' 
	That will be written $\theta_{j,k}.$ 
}
$(\bthet, \bfJ)$, so that the Hamiltonian, $H(\bfp,\bfq)$, becomes
a funtion of the actions alone.  Then Hamilton's equations are
\begin{eqnarray}
    \dot{J}_i & = & - \pderiv{H}{\theta_i} = 0, \nonumber \\
    \dot{\theta_i} & = & \pderiv{H}{J_i} \equiv \omega_i. 
\label{eqn:ham}
\end{eqnarray}
Figure (\ref{fig:tori}) illustrates the structure of the phase space
for a completely integrable, 2 degree-of-freedom system. Conservation of 
energy restricts the motion to a 3-dimensional energy surface, 
represented here as a solid torus.  A phase trajectory winds around
on a two dimensional torus, covering it densely unless $\omega_1$ and
$\omega_2$ are rationally dependent, that is, unless there are integers
$m_1$ and $m_2$ such that
\begin{equation}
	m_1 \omega_1 = m_2 \omega_2.
\label{eqn:resonance}
\end{equation}
Tori for which (\ref{eqn:resonance}) holds are called {\em resonant}
and they are entirely covered by periodic phase trajectories.  

     Figure (\ref{fig:tori}) also illustrates a construction we 
will use throughout the thesis, the Poincar\'{e} surface of section.
This technique reduces the continuous Hamiltonian flow, (\ref{eqn:ham}),
whose trajectories lie in a 2$n$-1 dimensional energy surface, to a 
discrete-time map, $T$, which acts on a 2$n$-2 dimensional surface. 
In figure (\ref{fig:tori}), the surface of section is given by 
$\theta_1 = 0$ and the map $T$ carries a phase point, $\bfx$, to 
the next point where $\bfx's$ trajectory intersects the surface.
That is,
\begin{displaymath}
    T( \bfJ, \theta_1 =0, \theta_2 ) = 
    ( \bfJ, \theta_1 =0, \theta_2 + 2\pi \frac{\omega_2}{\omega_1} ).
\end{displaymath}
The structures of integrability leave a clear signature on the surface
of section; all the orbits of $T$  are confined to circles, so that the
orbit of a typical point hops around its circle, eventually filling it
densely.  Those circles that are cross sections of resonant tori are
covered by periodic orbits; if a circle arises from a torus obeying a 
relation like (\ref{eqn:resonance}), then the points on it are periodic
with period $m_2$ and hop $m_1$ times around the circle before repeating.

     This extremely regular structure has profound qualitative effects
on the physics of the motion; integrable systems are far from satisfying 
the ergodic hypothesis of statistical mechanics.  A  phase trajectory, 
confined by conservation laws to 
an $n$ dimensional submanifold of the 2$n$-1 dimensional energy surface,
does not even  come close to exploring the whole of energetically 
accessible phase space and so predictions based on the microcannonical
ensemble, which gives equal weight to all points with the same energy,
will certainly be wrong.
These remarks, along with the evident success of statistical mechanics, 
suggest that complete integrability must be rare, that most of the 
structure of integrability cannot survive perturbation.  Indeed,
Fermi believed that the slightest perturbation would 
completely disrupt integrability, \cite{Fermi}.

    The fate of invariant tori is, however, much more complicated and
wonderful; it is the subject of the most spectacular theorem in 
Hamiltonian dynamics.
\newline
{\bf Theorem} (KAM) \newline {\em
If an unperturbed (completely integrable) system is 
non-degenerate\thesisfootnote{
	The non-degeneracy condition is that 
\begin{displaymath}
	{\rm det}\left| \pderiv{\bfomega}{\bfJ} \right| = 
	{\rm det}\left| 
	\frac{\partial^2 H_0}{\partial \bfJ^2}\right| \neq 0,
\end{displaymath}
   where $H_0(\bfJ)$ is the unperturbed Hamiltonian.  It means
   that the $\omega_i(\bfJ)$ are independent as functions.
},
then for sufficiently small conservative Hamiltonian perturbations, most 
non-resonant tori do not vanish, but are only slightly deformed, so that
in the phase space of the perturbed system, too, there are invariant tori
densely filled with phase curves winding around them 
conditionally-periodically, with a number of independent frequencies equal 
to the number of degrees of freedom.  These invariant tori form a majority 
in the sense that the measure of the complement of their union is small 
when the perturbation is small.}
\vskip \belowdisplayskip
That is, most tori survive small perturbations!  The statement above is
taken from Arnold's book, \cite{Arnold:Math-Meth}, but he does not give a proof.
Moser's book, \cite{Moser:book} gives an argument and \cite{Bost:KAM}
gives a thorough review.

\section{The Taylor-Chirikov Standard Map}

     We conclude our introduction with a brief review of an exhaustively
studied example, the Taylor-Chirikov standard map. 
It is a 2-dimensional, area-preserving map acting on the set 
${\rm S}^1 \times \R = \{ (x,p) | x \in [0,1),\; p \in \R \}$.
\begin{eqnarray}
	p' & = & p - \frac{k}{2\pi} \sin(2\pi x), \nonumber \\
	x' & = & x + p' \bmod 1.
\label{eqn:std-def}
\end{eqnarray}
Chirikov 
\cite{Chkv:Osc} describes this example as a periodically-kicked rotor,
sampled at the frequencey of the kicking; $x$ is a 
normalized angle 
variable with  $p$ the corresponding angular momentum.
Chirikov's rotor recieves periodic, impulsive blows whose size and 
direction depend on the rotor's angular position at the moment the 
impulse is delivered.   For $k = 0$, the system is completely 
integrable; $p$ is a constant of the motion and the orbits are
confined to one-dimensional curves.  

\begin{figure}[p]
	\begin{center}
		\includegraphics[height=19.0cm]{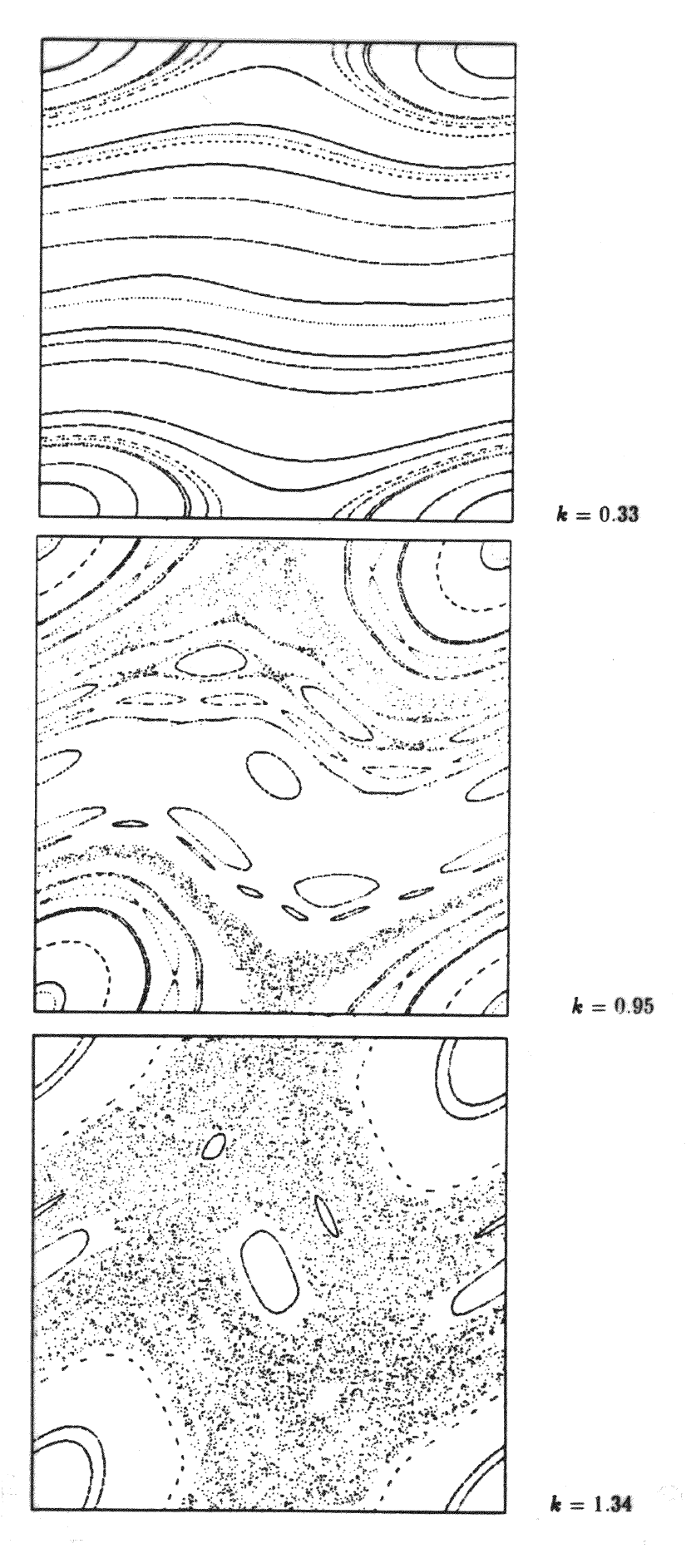}
	\end{center}
	\vspace*{-2\belowdisplayskip}
    \caption{ \em Orbits of the standard map for 
		  several sizes of the perturbation $k$.
		  Each panel shows 200 iterates from the 
		  orbits of 20 different initial conditions.}
\label{fig:std-demo}
\end{figure}

     Figure (\ref{fig:std-demo}) shows the structure of the
phase space for various values of the perturbation. Each panel
shows the orbits of several points from the
the set $\{ (x,p) | x \in [0,1), \; p \in [0,1) \}.$ 
Here we will give a qualitative discussion of these pictures,
at the same time introducing ideas which we will study fully 
in later chapters.
The series begins in the top panel with a small perturbation;
many orbits still seem to lie on or between circles. The 
arcs in the corners of the picture, when associated by 
periodic boundary conditions,  form ovals 
encircling the fixed point $z_0 \equiv (x = 0, p = 0)$.
The ovals 
arise because $z_0$ is an {\em elliptic}
fixed point; that is, the derivative of the map,
\begin{displaymath}
DT = \left[ \begin{array}{cc}
		\displaystyle \pderiv{x'}{x} & 
		\displaystyle \pderiv{x'}{p} \\
		\displaystyle \pderiv{p'}{x} \bigStrut & 
		\displaystyle \pderiv{p'}{p} \bigStrut
	\end{array} \right],
\end{displaymath}
is such that the matrix $DT_{z_0}$ has 
its eigenvalues on the unit circle.  Consequently, points
which start near $z_0$ stay nearby and their orbits form 
the arcs. If we were to restrict our attention to this 
{\em elliptic island} we would find that it
has much the same structure as the whole phase space; the
ovals would play the role of invariant circles and in amongst them
would lie yet smaller  elliptic islands.
If we magnified one of those islands $\dots$
the structure goes on forever.
There is also another fixed point, at $z_1 \equiv (x=\half, p=0)$, 
but it is {\em hyperbolic};   the matrix $DT_{z_1}$ has eigenvalues
off the unit circle, so almost every orbit which begins near it 
eventually
moves away with exponential speed.  Besides the fixed points,
there are always at least two
periodic orbits for every rational rotation number $\frac{p}{q}$.
Chapter \ref{chap:numres} gives a longer and more technical
discussion of periodic orbits and also discusses some special
sets, the {\em cantori}, which are, in a sense, the ghosts 
of disrupted tori. The chapter begins with a review of the two 
dimensional theory then shows some numerical work aimed at 
higher dimensional generalizations.

     In the middle panel, many more elliptic islands are evident,
as is a broad {\em stochastic layer}, a region which no 
longer contains any invariant tori; the orbits in such a region
are quite complicated and chaotic, and are confined to a layer only
because the phase space is two dimensional and thus the
invariant circles divide phase space into two disjoint pieces and 
so pairs of circles can trap even very chaotic orbits.
In higher dimensional systems
the tori have too low a dimension to isolate parts of the phase
space; points not actually contained in tori are free to diffuse
throughout the whole stochastic part of the phase space, though
they do so only very slowly, in a process called {\em Arnold
diffusion} \cite{Arnold:diff,Nek:diff}.
Although we will not have much more to say about  Arnold diffusion,
we will have cause to consider the topological consequences
of higher dimension; in both the remaining chapters we will find
that topology prevents us from proving results as strong as 
those available for two dimensional systems. 

    The final panel shows a  perturbation large enough to 
gaurantee very strong chaos; $k$ is so large that
Mather, \cite{Ma:circ}, has shown analytically that no 
invariant circles (of the type 
which wind all the way around the cylinder) remain.
Numerical experiments by Greene suggest that no circles
exist for $|k| > k_c \approx {\rm 0.971635406}$.
We leave this subject for the moment, but Chapter 
3 is entirely devoted to converse KAM results, theorems that say,
as Mather does, that for large enough perturbations, no tori 
exist at all.  There we will review Mather's  work, as well as the 
computer-assisted arguments of MacKay and Percival  
then discuss higher dimensional generalizations and 
show some new results.

%% file: Chap2.tex
\input{numres/intro2a}

\input{numres/intro2b}
\input{numres/shapes}
\input{numres/regularity}
\input{numres/hedlund}

%% file: numres/intro2a.tex
\chapter{Ghosts of Order}
\label{chap:numres}

In this chapter we ask, ``What becomes of invariant tori?''  
We have seen that the phase space of completely integrable
Hamiltonian systems is filled by such tori and that the  KAM theory
assures us that some of them persist even in the face of small perturbations.
What becomes of the tori for which KAM fails?  In general, one can't say.
But for certain two dimensional, area-preserving maps Mather 
\cite{Ma:top} and, independently, 
Aubry \cite{Aub:ref}, demonstrated the existence 
of some remarkable sets.  They are
reminiscent of invariant tori, but are not complete curves, rather, they 
look like graphs supported above a Cantor set.
Orbits on these ``cantori'' are similar to rotation on an invariant torus;
one may consider Mather's sets the ghosts of destroyed invariant tori.  Here
we review the two dimensional results, then present some numerical 
investigations\thesisfootnote{
Kook and Meiss, \protect \cite{KM:orbs}, have reported
similar studies; J. Meiss has been especially helpful in discussing
this work.}
from on effort to find the higher dimensional
analogs of Mather's sets. 
At the end of the chapter we discuss a 
topological obstacle which prevents simple generalization of the 
Aubry-Mather theory. 

\section{Basic notions and notations}

    In this section we give careful definitions of the maps we will 
study, the spaces they will act on, and the tools we will use to 
understand them.  We will also review the two dimensional theory, 
describing cantori and explaining how to approximate them by 
periodic orbits.  In the course of the review we will introduce a 
variational principle that will be the foundation of all our work.

\subsection{spaces and maps}

We will study maps based on the Poincar\'{e} map of a near-integrable,
action-angle system and so they will act on the {\em n-dimensional 
multi-annulus}, 
$\An = \Tn \times \Rn$, where $\Tn$ is the $n$-torus and $\Rn$ is 
$n$-dimensional Euclidean space.
To avoid having to worry
about factors of $2\pi$, we will always normalize the 
angles, and so write points in $\An$
as $(\bthet, \bfp)$ where $\bthet = (\theta_1, \theta_2 \cdots 
\theta_n)$ and the $\theta_i$ are periodic coordinates with period 1.

\begin{figure}
	\begin{center}
		\includegraphics[height=3.75cm]{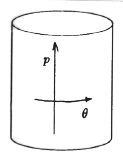}
	\end{center}
	\vspace*{-2\belowdisplayskip}
    \caption{\em The cylinder and its coordinate system. }
    \label{fig:cyl-def}
\end{figure}

The one-dimensional annulus, $\A = \T \times \R$, is conveniently 
represented as a 
cylinder with coordinates as pictured in figure (\ref{fig:cyl-def}).
Maps taking the cylinder to itself will be called $T$, or $T_{\epsilon}$
if they depend on parameters; maps acting on $\An$ for $n > 1$ will be
either $f$ or $f_{\epsilon}$.  In all cases, our maps will be 
{\em symplectic}, that is, they will preserve the standard symplectic
form (see e.g. \cite{Arnold:Math-Meth,KB:Birk}),
\begin{equation}
    \Omega = \sum_{j=1}^{n} d\bthet_i \wedge d\bfp_i.
\label{eqn:symplectic}
\end{equation}
For a map $T$ on the cylinder, preservation of (\ref{eqn:symplectic})
means that $T$ preserves area and orientation and so 
is equivalent to Liouville's 
theorem about the preservation of volume in phase space.  For 
higher dimensional systems, preservation  of (\ref{eqn:symplectic}) 
also implies preservation of volume, but is stronger.

     We will often need to work with a {\em lifting}, 
$F_\epsilon$, of a symplectic map, $f_\epsilon$, to the universal
cover of $\An$. This is essentailly a version of $f_\epsilon$ 
extended periodically so that acts on the whole of $\Rn \times \Rn$.
If $f_\epsilon: \An \rightarrow \An$, 
$f_\epsilon(\bthet,\bfp) = (\bthet'(\bthet, \bfp), \bfp'(\bthet,\bfp))$
then
$F_\epsilon$ acts on $\Rn \times \Rn$
$F_\epsilon(\bfx,\bfp) = (\bfx'(\bfx, \bfp), \bfp'(\bfx,\bfp))$,
and agrees
with $f_\epsilon$ up to an integer translation.  That is, if
$f_\epsilon( \bthet_0, \bfp_0) = ( \bthet_1, \bfp_1)$ and 
$F_\epsilon( \bfx_0 = \bthet_0, \bfp_0) = ( \bfx_1, \bfp_1)$ 
then
\begin{equation}
    \bfx_1 - \bthet_1 = \bfm
\label{eqn:lift}
\end{equation}
for some integer vector $\bfm \in \Zn$.  
Further, 
\begin{displaymath}
    F_\epsilon(\bfx_0 + \bfm, \bfp_0) = 
    F_\epsilon(\bfx_0, \bfp_0) + \bfm.
\end{displaymath}
The choice of a lift, $F_\epsilon$, which comes down to the choice of 
$\bfm$ in (\ref{eqn:lift}) does not affect any qualitative 
features of the dynamics.
For example,
a lift of the standard map is 
\begin{eqnarray*}
	p' & = & p - \frac{k}{2\pi} \sin(2\pi x),  \\
	x' & = & x + p',
\end{eqnarray*}
which is just the same as (\ref{eqn:std-def}) except that 
the position coordinate is no longer taken mod 1.
We will always use the convention that 
$F_\epsilon : \Rn \times \Rn$ is a lift of
$f_\epsilon : \An \rightarrow \An$.

\subsection{a variational principle}

\begin{figure}
	\begin{center}
		\includegraphics[height=3.75cm]{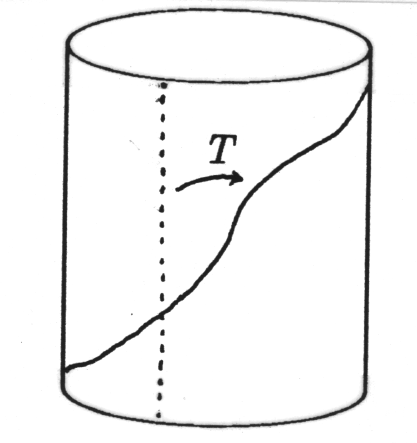}
	\end{center}
	\vspace*{-2\belowdisplayskip}
    \caption{\em A twist map carries vertical lines to 
    monotone curves. }
\label{fig:twist}
\end{figure}

The dynamics of an autonomous Hamiltonian system can be 
characterized with the principle of least action; to specify a
 segment of a phase trajectory, $\gamma(t) = (\bfp(t), \bfq(t))$,
 one need only note the values of the position coordinates
 at the ends of the segment and require that $\gamma$ be an 
 extremal of the ``reduced action'' functional 
 \cite{Arnold:Math-Meth},
 \begin{equation}
	S(\bfq_0, \bfq_1) = \int_{\bfq_0}^{\bfq_1} \bfp d\bfq.
\label{eqn:action}
\end{equation}
In particular, one can get the momenta at the endpoints of the
segment by taking derivatives of $S(\bfq_0, \bfq_1)$;
\begin{displaymath}
    \bfp_1 = \pderiv{S}{\bfq_1}
    \qquad {\rm and} \qquad
    \bfp_0 = - \pderiv{S}{\bfq_0}.
\end{displaymath}
The analogous thing
for a symplectic map $F_\epsilon : \Rn \rightarrow \Rn$
is an {\em action-generating function}, a function, 
\mbox{$H_\epsilon : \Rn \times \Rn \rightarrow \R,$} where
\mbox{$H_\epsilon = H_\epsilon(\bfx, \bfx')$} is such that if
$F_\epsilon(\bfx_0, \bfp_0)  =  (\bfx_1, \bfp_1)$,
then
\begin{equation}
    \bfp_1 = \pderiv{H_\epsilon}{\bfx'}
    \qquad {\rm and} \qquad
    \bfp_0 = -\pderiv{H_\epsilon}{\bfx}
\label{eqn:Hgen}
\end{equation}
The point of constructing a generating function is that it enables us 
to discuss dynamics entirely in terms of the position coordinates. 
In the
next section we will demonstrate the usefulness of variational arguments
by reviewing the theory of area-preserving twist maps of the cylinder.  
These maps get their name because of a geometric property of their action;
a $C^1$ map $T$ is {\em twist} if it carries every vertical line into a 
monotone curve; see figure (\ref{fig:twist}).
More analytically, if $T(\theta,p) = (\theta'(\theta), p'(\theta,p)$ is
a symplectic map of the cylinder, then $T$ is twist if
\begin{displaymath}
    \pderiv{\theta'}{p} \neq 0.
\end{displaymath}

\subsection{area-preserving twist maps}

    Here we will examine the kinds of orbits which can occur for an
area-preserving twist map.   Since we will be wanting to make variational
arguments we require that, in addition to being a twist map, $T$ posses
a generating function, $h(x,x')$.  For convenience, 
we will work with a lift of $T$,
call it $\widetilde{T}$, and will use coordinates in $\R \times \R$ 
rather than on the cylinder.  First we will use the generating function
to construct some periodic orbits.

\begin{figure}
	\begin{center}
		\includegraphics[height=3.75cm]{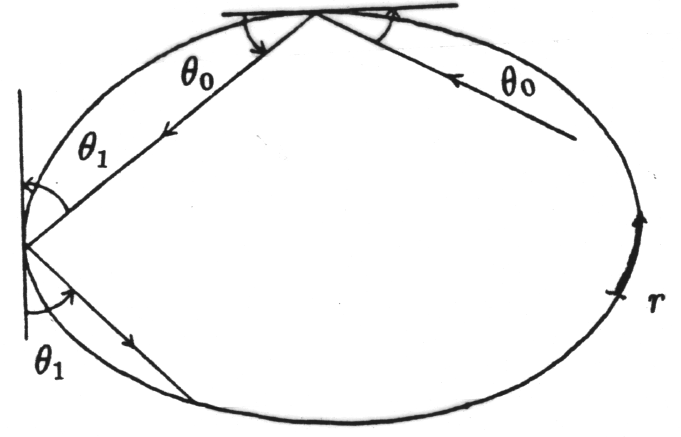}
	\end{center}
	\vspace*{-2\belowdisplayskip}

	\caption{\em The billiard ball dynamical system.
	  {\rm \protect\cite{Birk:orb}} 
	}
\label{fig:billiards}
\end{figure}

     A periodic orbit is characterized by its period and by the number
of times it winds around the cylinder before closing. Suppose we 
want an orbit which, in $q$ steps, makes $p$ turns.  Such an 
orbit would appear on the universal cover as a sequence of points 
$\{\cdots (x_0,p_0), (x_1,p_1), \cdots (x_{q-1}, p_{q-1}), (x_p,q_p), 
\cdots \}$ with $x_{j+q} = x_j + p$.  We could seek it
by trying to find a sequence of position coordinates, 
\begin{equation}
    X = \{x_0, x_1, \ldots, x_{q-1}, x_q ;\,x_{q} = x_0 + p\},
\label{eqn:state-def}
\end{equation}
such that the function
\begin{equation}
    L_{p,q}(X) = \sum_{j=0}^{q-1} h(x_j, x_{j+1})
\label{eqn:Lpq}
\end{equation}
was minimized.  We will call such a sequence a 
{\em p-q minimizing state}. If we could find one, 
then, automatically, 
we could compute the desired kind of periodic orbit.
To see how, consider the condition that (\ref{eqn:Lpq}) be extremal:
\begin{equation}
    \pderiv{L_{p,q}}{x_j} = \pderiv{h}{x}(x_j, x_{j+1})
			+ \pderiv{h}{x'}(x_{j-1},x_j) = 0
			\qquad \mbox{\em for\ } j=0,1,\cdots,q-1.
\label{eqn:ELg}
\end{equation}
We will call these  the {\em Euler-Lagrange equations}.
Now, if $X$ were the  projection of some periodic orbit,
we would be able to recover the missing momentum coordinates
in two ways; we could use either
\begin{displaymath}
    p_j = \pderiv{h}{x'}(x_{j-1},x_j)
    \qquad {\rm or } \qquad
    p_j = - \pderiv{h}{x}(x_j,x_{j+1}).
\end{displaymath}
The condition (\ref{eqn:ELg}) is that these two be equal, so that
if we can find a sequence like (\ref{eqn:state-def}) we have
found the desired periodic orbit.  Arguments like this were first
made by Birkhoff, who used them to construct periodic orbits 
for the map given by the motion of a point particle in a convex,
rigid walled box.  This system can be reduced to an area preserving
twist map by considering the particle's collisions  the  wall
and using coordinates
given by a length, $r$ measured along the perimeter of the domain, 
and the variable $\sigma = -\cos(\theta)$ where $\theta$ is the
angle the particle's path makes with the tangent to the wall,
see figure (\ref{fig:billiards}).
In this system the generating function is just the negative of the
length of the path traced by the ball, and so the minimizing
periodic orbit with $p=2, q=5$ is just the orbit which corresponds
to the longest inscribed star.  Besides the minimizing periodic 
orbit, there
is another, a {\em minimax} orbit.  To see how this orbit
arises take one point of the
minimizing orbit and slide it along the boundary, allowing the other
points to shift so as to keep the total length 
of the star as large as possible.
At first the length must decrease; we have assumed that 
the initial, undistorted star was the longest possible.
Eventually, though, the length of the distorted star will have
to stop decreasing and begin to increase because eventually the 
verticies will reach a configuration which is a cyclic permutation 
of the original star. The configuration
for which the length again begins to increase must also be a 
sationary point of $L_{p,q}$; it satisfies
the Euler-Lagrange equations and so it too corresponds to a 
genuine periodic orbit.

     The action-minimizing periodic orbits, which are called
{\em Birkhoff orbits}, are distinguished by the 
numbers $p$ and $q$
used in their construction.  The rational number $\frac{p}{q}$,
which is the orbit's average angular speed, is called the
{\em rotation number} of the orbit.  More generally, an orbit
$(x_0,p_0), (x_1,p_1),\ldots$ on the
universal cover is said to have rotation number $\alpha$ if 
\begin{equation}
	\alpha = \lim_{n \rightarrow \infty} 
	\frac{x_n - x_0}{n}.
\label{eqn:rot-num}
\end{equation}
This limit does not always exist.  Most of the points in
the stochastic regions of the standard map do not have 
well-defined rotation numbers, though all of the orbits
lying on invariant circles do; orbits on non-resonant
circles have irrational $\alpha$.

This observation prompted Mather, in \cite{Ma:top}, to try 
to find orbits that had irrational rotation numbers, but
were not part of invariant tori.  He succeeeded dramatically,
discovering whole, complicated sets of such orbits and
revealing an unexpected, rich structure in the phase space.

\begin{figure}
	\begin{center}
		\includegraphics[height=6.0cm]{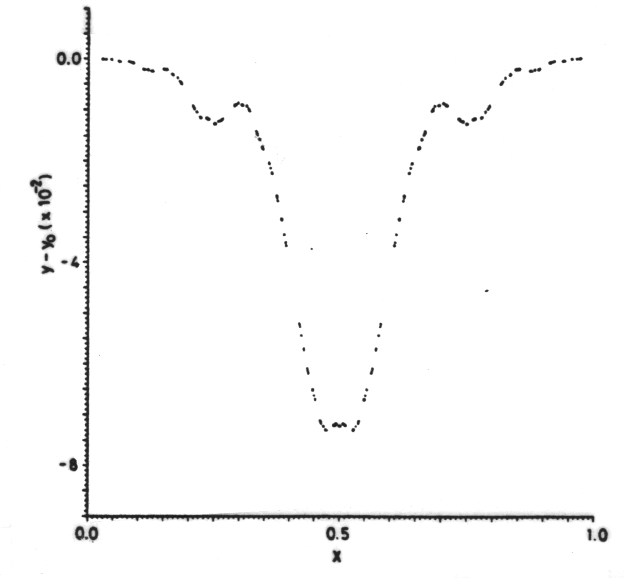}
	\end{center}
	\vspace*{-2\belowdisplayskip}
    \caption{\em A cantorus for the standard map.  The vertical
		axis is measured in units of 
		$y=p - \frac{k}{4\pi} \sin(2\pi x) $, where 
		$k={\rm 1.001635}$ is the size of the perturbation 
		and the rotation number is 
		$\approx \frac{1}{\gamma^2}$
		where $\gamma = \frac{1 + \protect \sqrt{5} }{2}$ 
		is the golden mean.
		\rm \protect \cite{MMP:tran}
		}
\label{fig:cantorus}
\end{figure}

We can construct one of Mather's sets 
by taking a limit of minimizing, Birkhoff
periodic orbits.  That is, we take a sequence of rational
numbers $\{p_0/q_0, p_1/q_1 \cdots\}$ which has an irrational
$\omega$ as a limit, construct the 
corresponding Birkhoff minimizing orbits, and see 
whether they accumulated to any interesting limit set.
Katok, \cite{Kat:rem}, has shown that they do.  If 
there is an invariant circle with rotation number $\omega$,
then the Birkhoff orbits accumulate on it.  If there is
no invariant circle, then the orbits accumulate on a
cantorus, a set which looks like an invariant circle
with a countable set of holes cut out of it, see figure
(\ref{fig:cantorus}).  

The cantori have many properties reminiscent of irrational invariant 
circles; orbits lying in the cantorus are dense and the motion
on the cantorus, is, by a continuous change of coordinate,
equivalent to rotation by the angle $\omega$.
Also, the cantorus has the same kind of smoothness\thesisfootnote{
   A theorem of Birkhoff states that the invariant circles 
   are Lipschitz graphs.
}
as an invariant circle.  If $(\theta_0,p_0)$
and $\theta_1,p_1)$ are any two points from the cantorus then 
there is a constant $L$, independent of the $\theta$'s,
such that
\begin{displaymath}
    |p_0 - p_1| \leq L |\theta_0 - \theta_1|,
\end{displaymath}
that is, the momenta are Lipschitz functions of 
the positions. 

Katok's scheme for
approximating the cantorus by a of periodic
orbits is different from the approach first used by Mather,
but it is much better suited to  numerical experiment;
all computational investigations of cantroi depend on 
approximnation by periodic orbits e.g. \cite{MMP:tran,MP:Newt,Grn:meth}.

%% file: numres/intro2b.tex
\section{Higher dimensional analogs}

In this section we formulate the numerical investigations reported in
the rest of the chapter.  Our studies are based on the 
Katok and Bernstien's paper, \cite{KB:Birk} in which they
study certain $n$-dimensional symplectic maps generated by a function
$H_\epsilon(\bfx, \bfx')$ and
prove the existence of action-minimizing periodic orbits.
For these orbits, which are defined by analogy with the 
Birkhoff orbits on the cylinder, the role of the rational rotation
number $\frac{p}{q}$ is played by a {\em rotation vector},
$\frac{\bfp}{q}$ where $q$ is the length of the orbit and
$\bfp \in \Zn, \; \bfp = (p_0, p_1, \ldots, p_n)$ gives the number
of times the orbit winds around each of the coordinate directions.
As above, each rational vector has a corresponding type of 
$p,q\/$-minimizing state,
\begin{displaymath}
    X = {\bfx_0, \bfx_1, \ldots, 
	\bfx_{q-1}, \bfx_q; \, \bfx_q = \bfx_0 + \bfp }
\end{displaymath}
an action functional, $L_{p,q}$, some Euler-Lagrange equations,
\begin{eqnarray}
    L_{p,q}(X) & = & \sum_{j=0}^{q-1} H_\epsilon( \bfx_j, \bfx_{j+1}) 
    \label{eqn:KB-Lpq} \\
    \pderiv{L_{p,q}}{\bfx_j} & = & 
	\pderiv{H_\epsilon}{\bfx'}(\bfx_{j-1}, \bfx_j ) +
	\pderiv{H_\epsilon}{\bfx}(\bfx_j, \bfx_{j+1}),
    \label{eqn:KB-ELg}
\end{eqnarray}
and at least one minimizing periodic orbit.
Katok and Bernstien's maps are small perturbations of 
some completely integrable
system whose unperturbed generating function, $H_0(\bfx,\bfx')$,
satisfies $H_0(\bfx,\bfx') = h(\bfx' - \bfx)$ where $h(\bfu)$
is strictly convex, i.e., the Hessian matrix of $h$, 
\begin{equation}
    \frac{\partial^2 h}{\partial \bfu^2} = 
    \left[ \begin{array}{cccc}
	\frac{\partial^2 h}{\partial u_0^2} & 
	\frac{\partial^2 h}{\partial u_0 \partial u_1} & 
	\cdots & 
	\frac{\partial^2 h}{\partial u_0 \partial u_{n-1}} \\
	\frac{\partial^2 h}{\partial u_1 \partial u_0} & 
	\frac{\partial^2 h}{\partial u_1^2} & \cdots &  \vdots \\
	\vdots & & \ddots &  \\
	\frac{\partial^2 h}{\partial u_{n-1} \partial u_0 } &
	\cdots & & 
	\frac{\partial^2 h}{\partial u_{n-1}^2} 
    \end{array} \right],
\label{eqn:convexity}
\end{equation}
is positive definite.  This condition is a higher dimensional
analog of the twist condition, but is not the only possible
generalization; Herman, in \cite{Herm:non}, gives another.
In the next section we will present some explicit 4-d
symplectic maps and their generating functions and
section~\ref{sec:pix} we show some pictures of minimizing periodic
orbits and discuss how their shapes and stability depend 
on the size of the perturbation.

The real question here is ``Are there cantori in 4-d symplectic
maps?''  On the analytic side, the answer seems to be ``maybe.''
Katok and Bernstien are able to show that if a sequence
of rational rotation vectors
$\{ \frac{\bfp_0}{q_0},\frac{\bfp_1}{q_1},\ldots\,\}, 
\;\; \bfp_i \in \Zn, q \in \Z$, converges to some irrational rotation
vector, $\bfomega = (\omega_1, \omega_2, \cdots \omega_n)$, then the 
corresponding sequence of Birkhoff orbits also has
a limit.  Unfortunately their results on the properties of the limiting
set are not as strong as those available for twist maps.  They cannot
say what the limiting set looks like or  much about the motion on it.
They are able to establish that the momenta should be H\"{o}lder
continuous functions of the positions, but with 
index $\alpha = \half$, that is if, $(\bthet_0, \bfp_0)$ and
$(\bthet_1, \bfp_1)$ are points from this limit set then,
except, perhaps for a single isolated point,
\begin{equation}
    ||\bfp_0 - \bfp_1|| \leq C ||\bthet_0 - \bthet_1||^{\frac{1}{2}},
\label{eqn:holder}
\end{equation}
for some constant $C$, independent of the $\bthet_i$.
We present some ambiguous numerical investigations aimed at 
verifying or improving this smoothness estimate, but are 
unable to report any definite results.

Finally, in section~\ref{sec:Hed} we discuss a 
pathology forseen by Hedlund.  Hedlund's examples  complicate any
discussion of the behaviour of very long orbits and are an 
obstacle to both analytic and numerical 
investigation of higher dimensional cantori.
We report on some qualitative investigations designed to see
whether Hedlund's pathology actually occurs.  

\subsection{the maps and orbits}

We follow \cite{KB:Birk} and study maps which are generated by functions
of the form
\begin{equation}
H_\epsilon(\bfx,\bfx') = h(\bfx'-\bfx) - V_\epsilon(\bfx,\bfx'),
\label{eqn:gen}
\end{equation}
where $h(\bfx'-\bfx): \Rn \rightarrow \R$, the unperturbed part of the
generating function, satisfies (\ref{eqn:convexity}) and the perturbation
$V_\epsilon(\bfx,\bfx') : \Rn \times \Rn \rightarrow \R$, 
is a small, $C^2$ function satisfying 
$V_\epsilon(\bfx+\bfm, \bfx'+\bfm) = V_\epsilon(\bfx,\bfx') \; 
\forall \bfm \in \Zn$.
We will study 4-d sympectic maps generated by (\ref{eqn:gen}) with
\begin{displaymath}
h(\bfx,\bfx') =  \half \norm{\bfx'-\bfx}^2, \qquad 
V_\epsilon(\bfx,\bfx') =  \epsilon V(\bfx). 
\end{displaymath}
Where
\begin{equation}
	   V(\bfx) =  \left\{
				\begin{array}{lll}
				    \mbox{one of} \\
					V_{trig}(\bfx) & = &
					-\frac{1}{M_{trig}}\;
					\{\; \half (\sin{2\pi x_0} + 
					\sin{2\pi x_1}) +  
					\sin{2\pi(x_0+x_1)}\; \},  \\ \\
					V_{poly}(\bfx) & = &
					-\frac{1}{M_{poly}}\;
					\{\;[x_0^2(1-x_0)^2
					(x_0-\tqrt)(\oqrt - x_0)]\;
					[x_1^2(1-x_1)^2]\;\}, \\
				    \mbox{or}  \\
					V_{ff}(\bfx) & = &
					-\half \,
					\{\; \half (c( x_0) + 
					c( x_1) ) +  
					c(x_0+x_1)\; \},   \\ \\
				   \mbox{with} & &
				   c(x) \;=\; \left\{ \begin{array}{cr}
				       1 -24x^2 +32x^3 & {\rm if}\;
				       x \bmod 1 \leq \half, \\ 
				       9 -48x +72x^2 -32x^3 
				       & {\rm if}\;
				       x \bmod 1 > \half.
				\end{array}
				\right.
				\end{array}
				\right.
\label{eqn:zamples}
\end{equation}
Call the
first perturbation the {\em trigonometric} perturbation, the
second the {\em polynomial} perturbation\thesisfootnote{
    The $x_i$ appearing in the definition of $V_{poly}$ are
    all taken $\bmod$ 1. }
and the third the 
{\em fast-Froschl\'{e}}.  
The constants $M_{trig}$
and $M_{poly}$ are chosen so that $\max_{x \in \Tn}|V(\bfx)| = 1 $.
$V_{ff}(\bfx)$ is a polynomial approximation to a map originally
introduced as a model of star motion in elliptical galaxies
\cite{Fro:int}.
The real Froschl\'{e} map has cosines where ours has $c(x)$
and has three independent constants, one for each of the terms.
Since its introduction the map has been popular as a model
for chaotic Hamiltonian dynamics e.g. 
\cite{Fro:num1,Fro:num2,Bagley:Diff,KM:orbs,MMS:Conv}.

All our examples use ``standard-like'' perturbations, ones where
$V_\epsilon(\bfx, \bfx')$  depends 
on $\bfx$ but not on its successor, $\bfx'$.  We made this choice
of perturbation because it simplifies the map. Using~(\ref{eqn:Hgen})
we obtain
\begin{eqnarray}
\bfp'(\bfx,\bfp) & = & \bfp - \epsilon \, \frac{\partial V}
				      {\partial \bfx}(\bfx), 
				       \nonumber  \\
\bfx'(\bfx,\bfp) & = &  \bfx +  \bfp - 
			\epsilon \,  \frac{\partial V}
					{\partial \bfx}(\bfx). 
\label{eqn:map}
\end{eqnarray}
\begin{figure}[hp]
	\begin{center}
		\includegraphics[height=18.0cm]{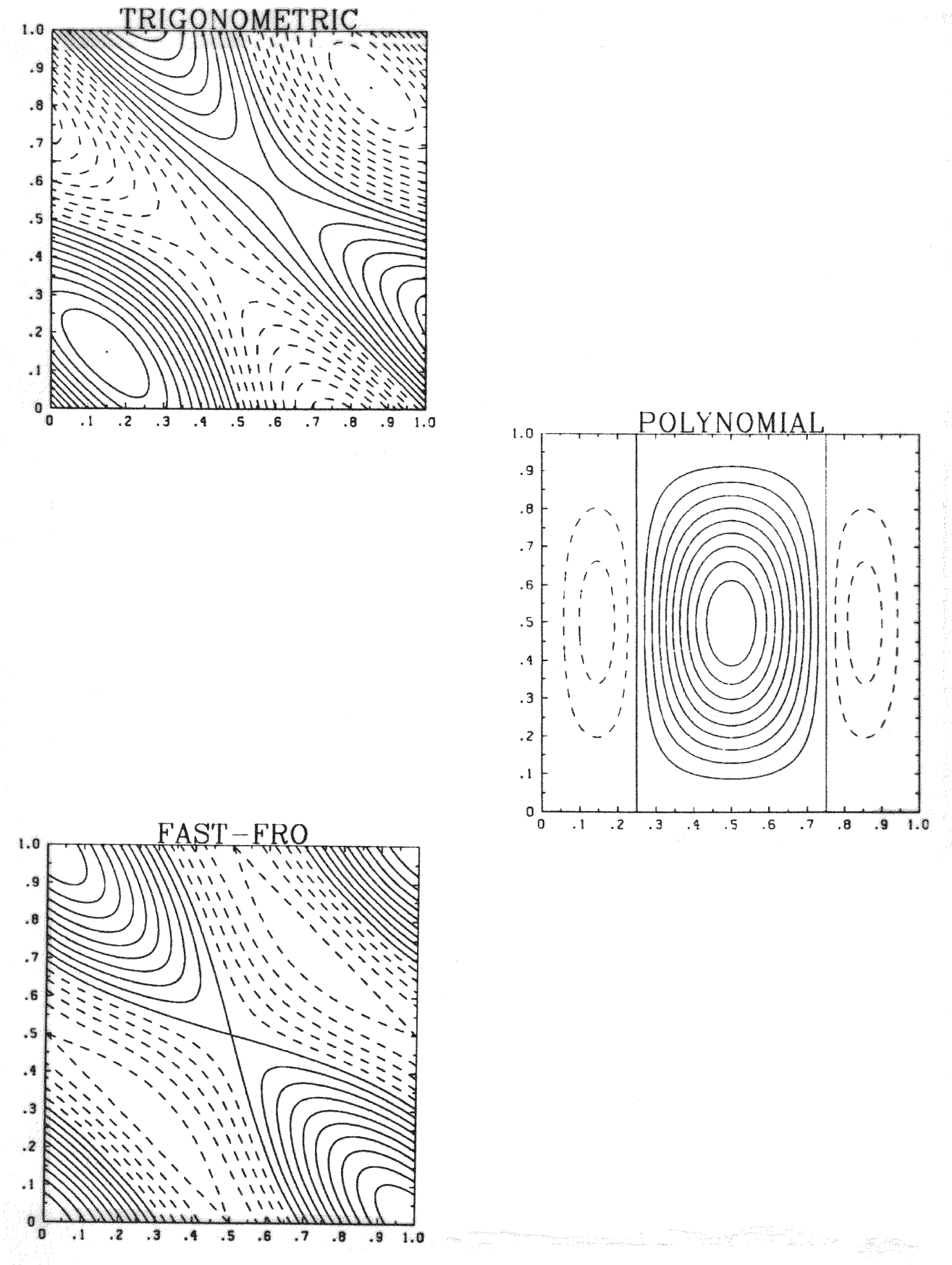}
	\end{center}
	\vspace*{-2\belowdisplayskip}
    \caption{\em Contour maps of $-V_\epsilon(\protect \bfx)$ for the
		{\em (a)} trigonometric, 
		{\em (b)} polynomial, and 
		{\em (c)} fast-Froeschl\'{e} perturbations.
		The conour interval is 0.1 and the contours 
		corresponding to negative values are dashed.
	}
\label{fig:pert-maps}
\end{figure}

%% file: numres/shapes.tex
\subsection{shapes of orbits and Lyapunov exponents}
\label{sec:pix}
\label{subsec:Lyap}

\begin{figure}[hp]
	\begin{center}
		\includegraphics[height=13.0cm]{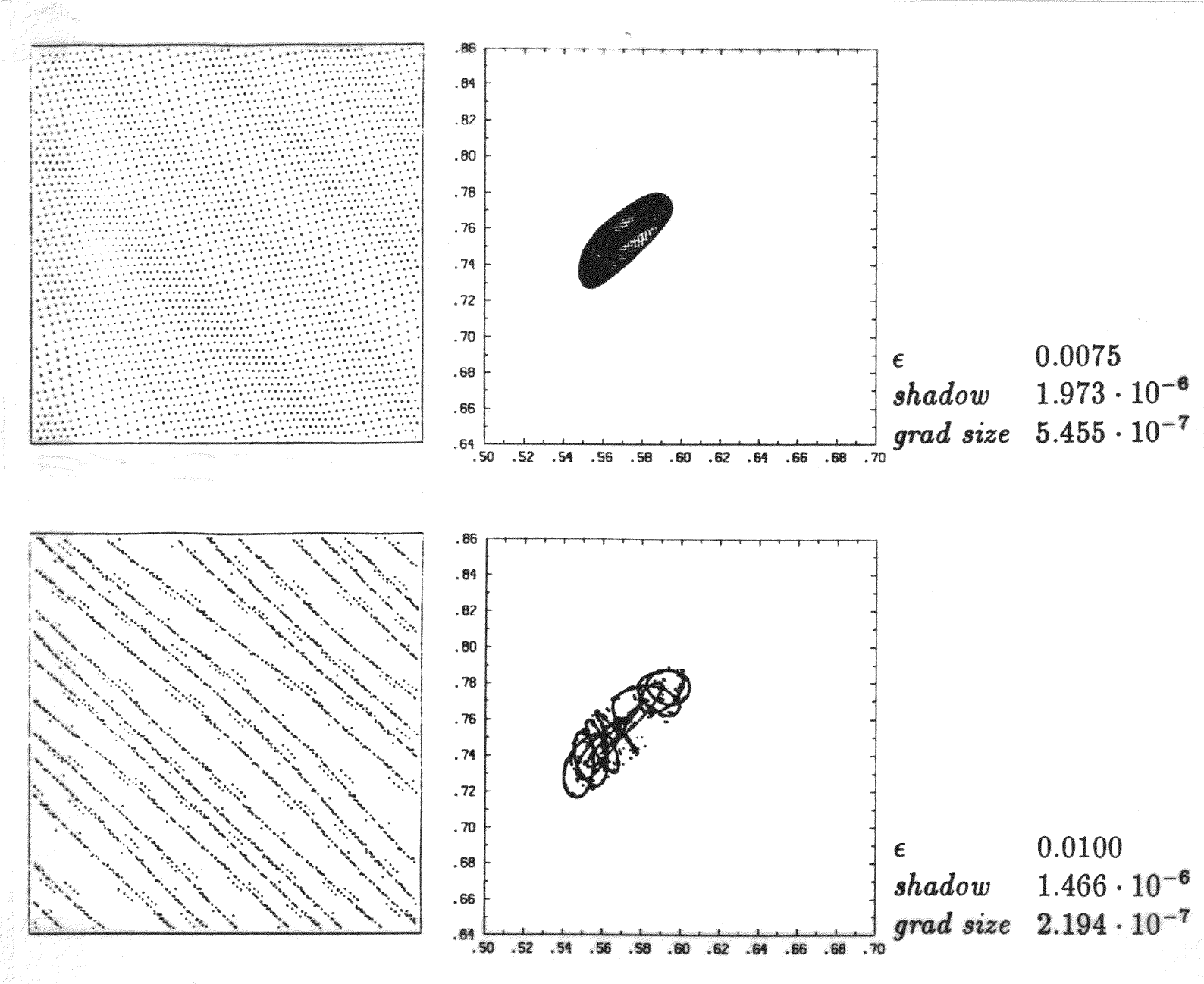}
	\end{center}
	\vspace*{-2\belowdisplayskip}
	\caption{ \em Birkhoff orbits for the trigonometric 
	perturbation and the rotation vector (1432,1897)/2513.
	This panel illustrates the collapse along filaments.
	Notice how the $\epsilon = {\rm 0.0075}$ state has
	momenta seeming to lie on a smooth surface.}
\label{fig:orb-spnf1}
\end{figure}

\begin{figure}[hp]
	\begin{center}
		\includegraphics[height=13.0cm]{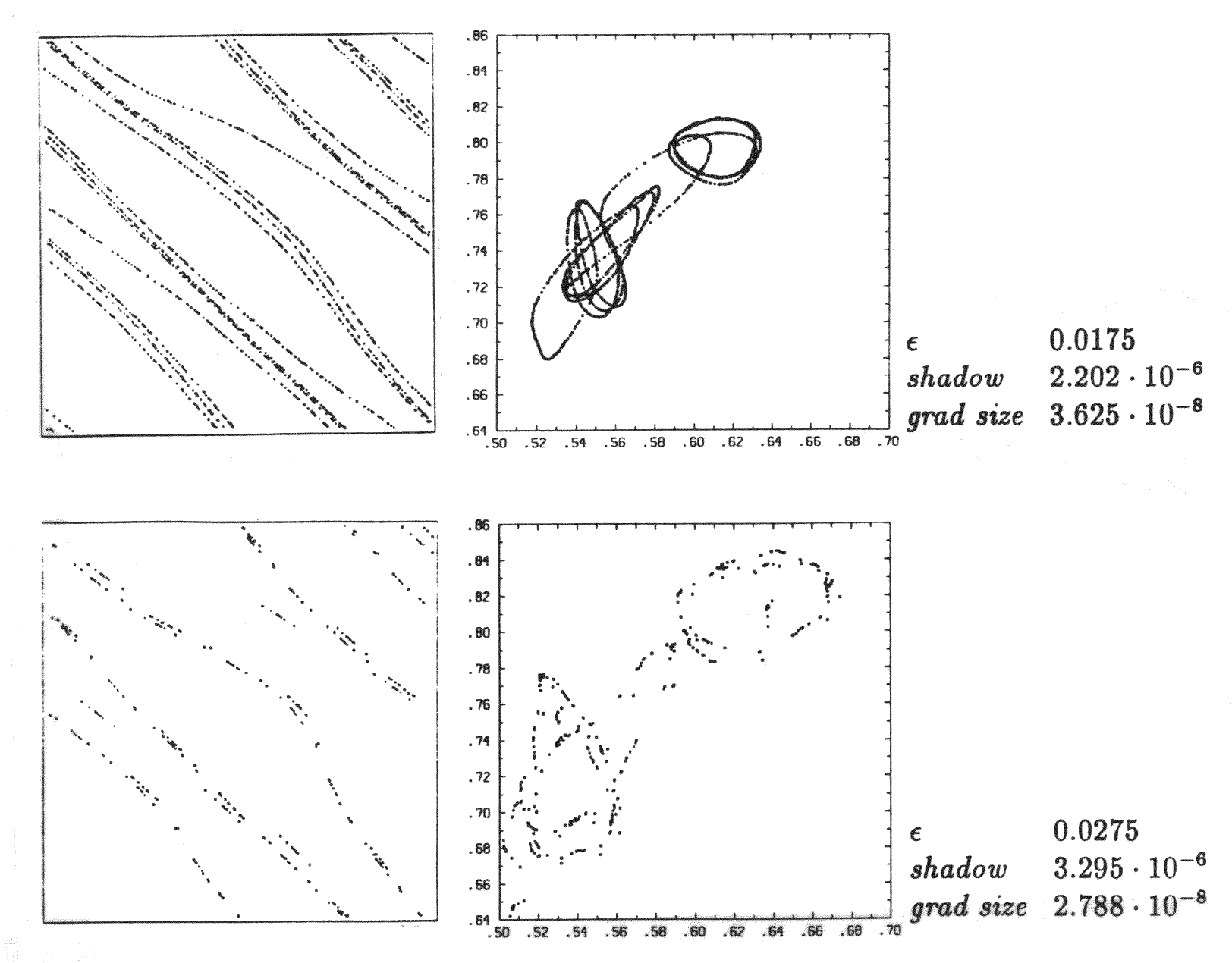}
	\end{center}
	\vspace*{-2\belowdisplayskip}

	\caption{ \em Birkhoff orbits for the trigonometric 
	perturbation and the rotation vector (1432,1897)/2513.
	This pair shows the appearance of Cantor-like clumping
	along the filaments.}
\label{fig:orb-spnf2}
\end{figure}

\begin{figure}[hp]
	\begin{center}
		\includegraphics[height=13.0cm]{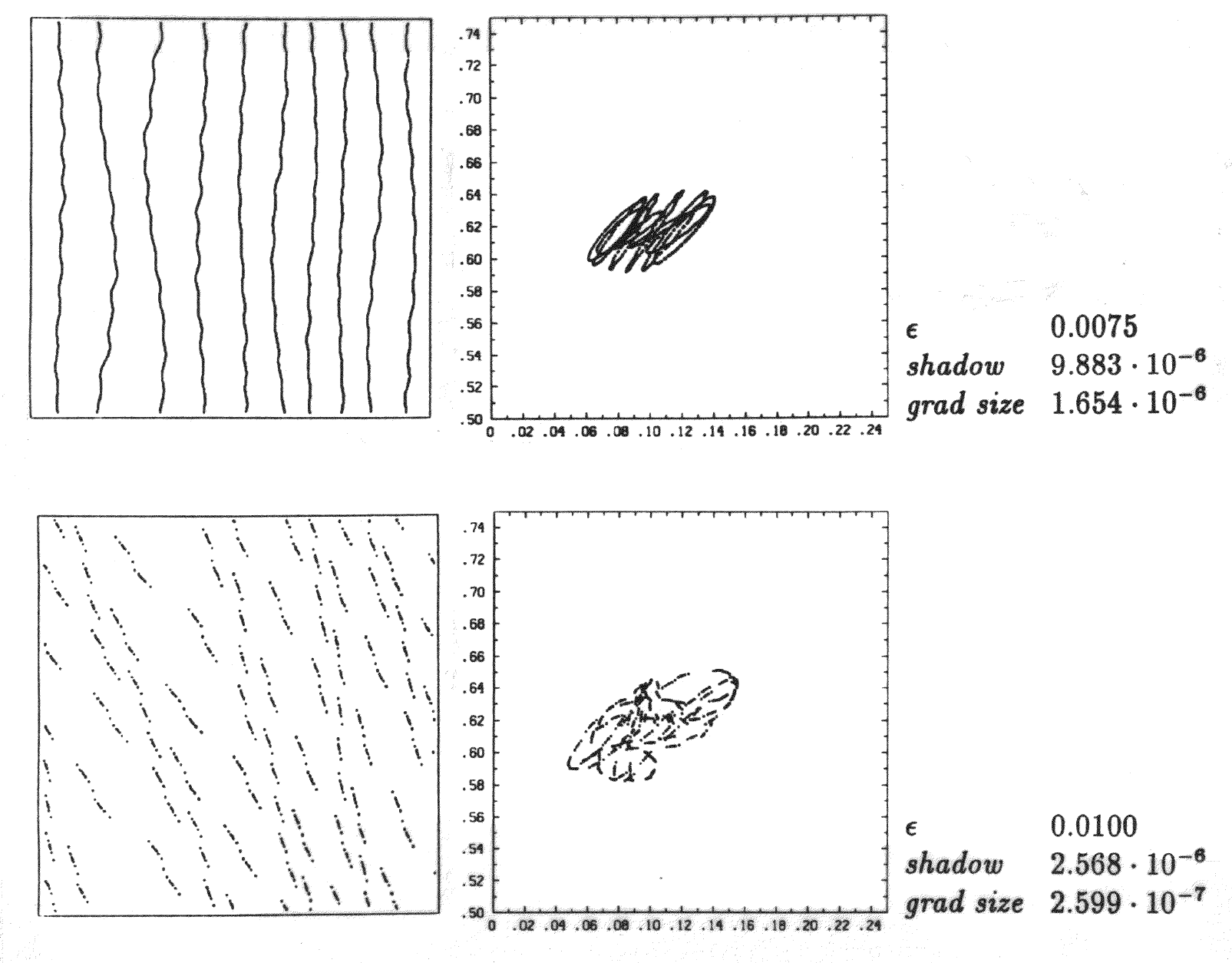}
	\end{center}
	\vspace*{-2\belowdisplayskip}
	\caption{ \em Weakly perturbed Birkhoff orbits for the 
	trigonometric 
	perturbation and the rotation vector (377, 2330)/3770). }
\label{fig:orb-aunf1}
\end{figure}

\begin{figure}[hp]
	\begin{center}
		\includegraphics[height=13.0cm]{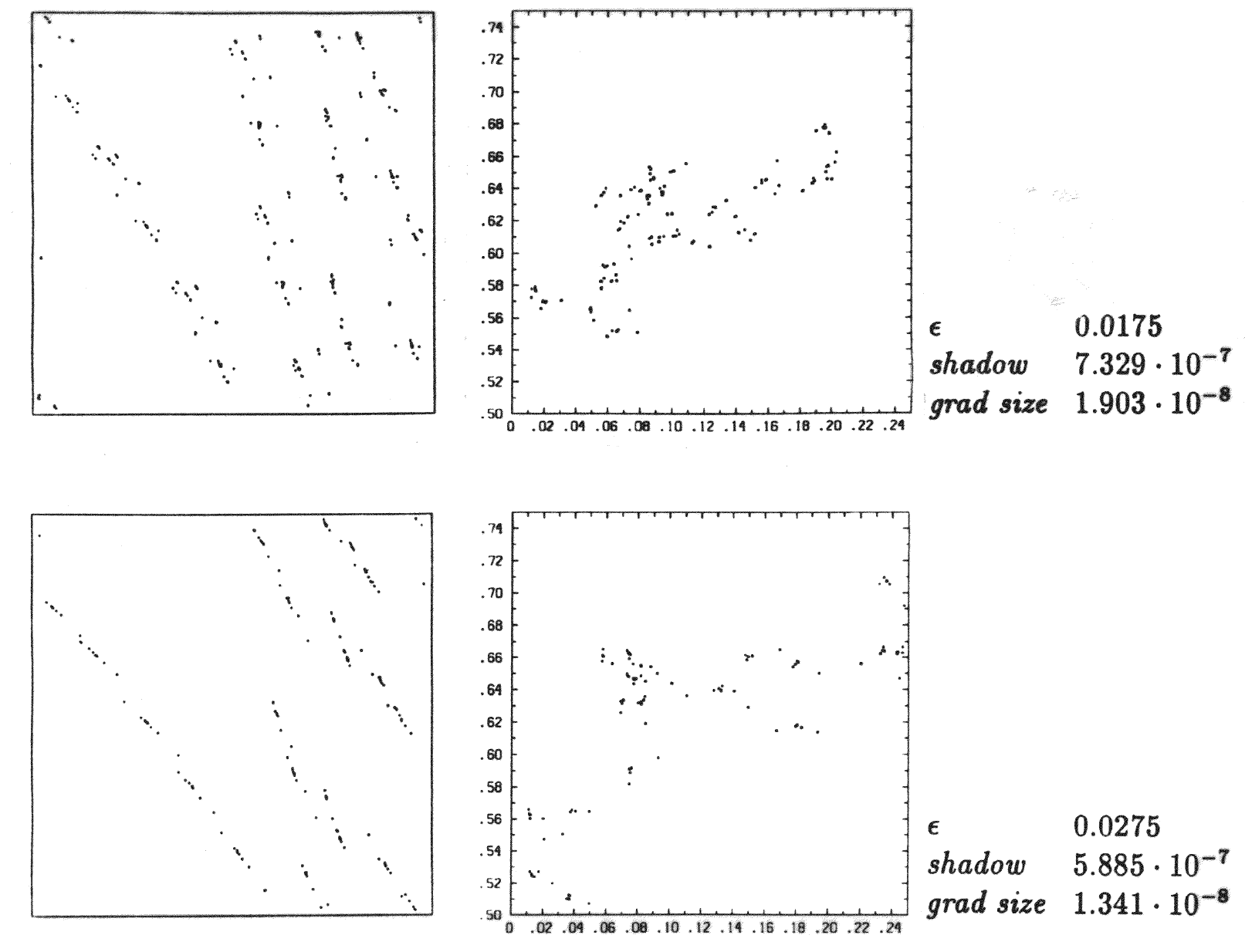}
	\end{center}
	\vspace*{-2\belowdisplayskip}
	\caption{ \em Strongly perturbed Birkhoff orbits for 
	the trigonometric 
	perturbation and the rotation vector (377, 2330)/3770). }
\label{fig:orb-aunf2}
\end{figure}

\begin{figure}[hp]
	\begin{center}
		\includegraphics[height=13.0cm]{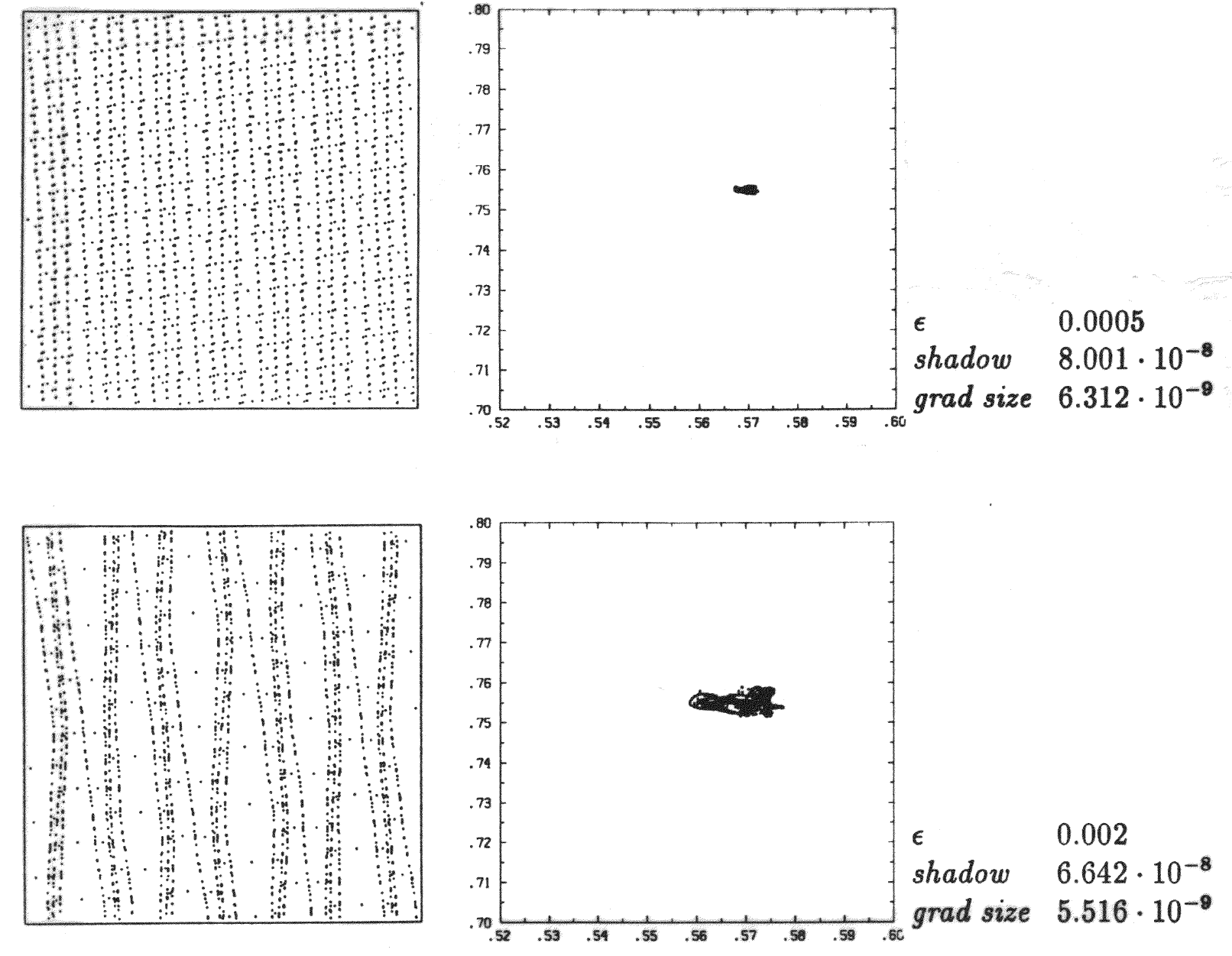}
	\end{center}
	\vspace*{-2\belowdisplayskip}
	\caption{ \em Birkhoff orbits for the polynomial 
	perturbation and the rotation vector (1432,1897)/2513.
	Note that the momenta remain very near their unperturbed
	values.}
\label{fig:orb-sphd1}
\end{figure}

\begin{figure}[hp]
	\begin{center}
		\includegraphics[height=13.0cm]{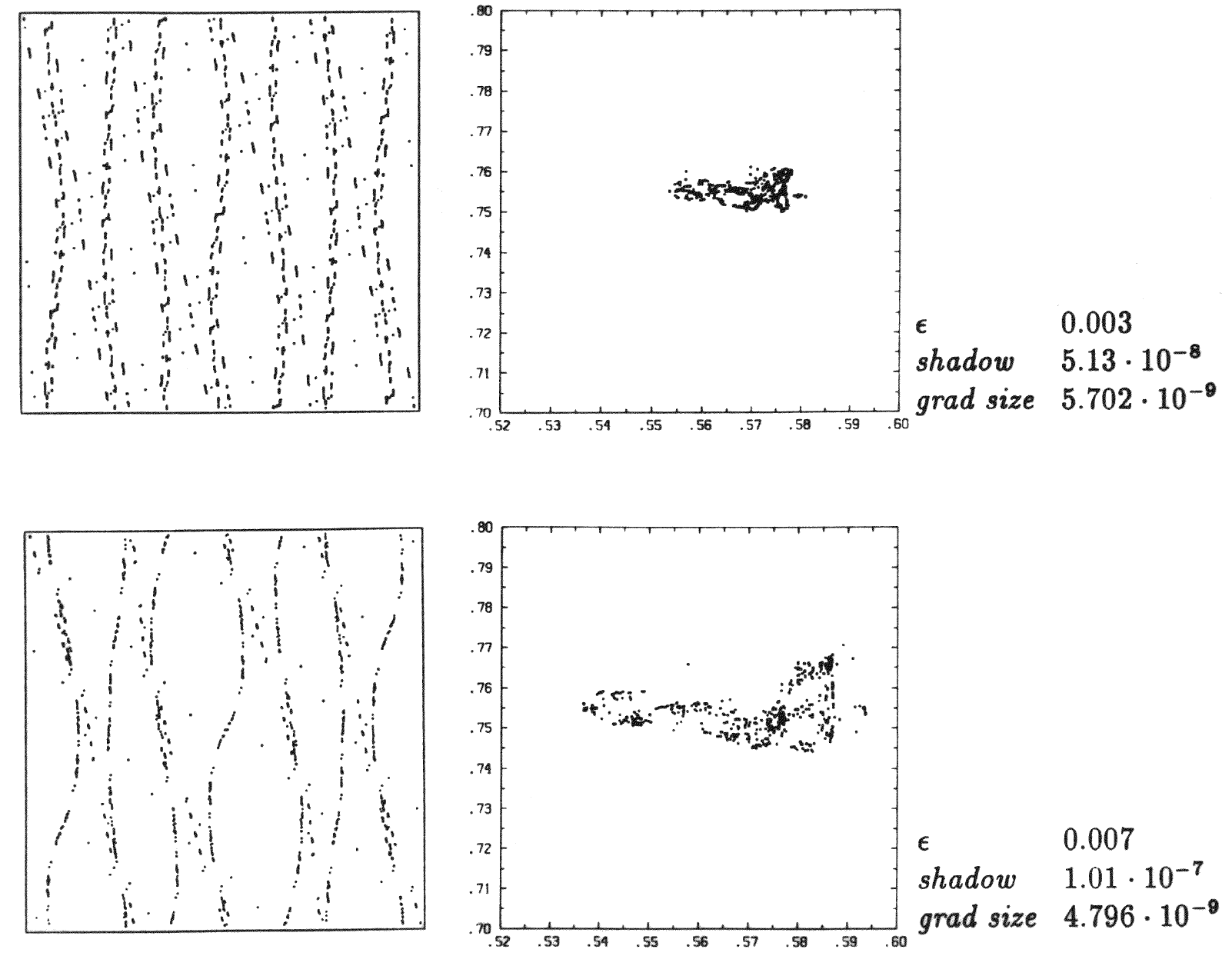}
	\end{center}
	\vspace*{-2\belowdisplayskip}

	\caption{ \em Birkhoff orbits for the polynomial 
	perturbation and the rotation vector (1432,1897)/2513.
	This pair shows the appearance of Cantor-like clumping
	along the filaments.}
\label{fig:orb-sphd2}
\end{figure}

\begin{figure}[hp]
	\begin{center}
		\includegraphics[height=13.0cm]{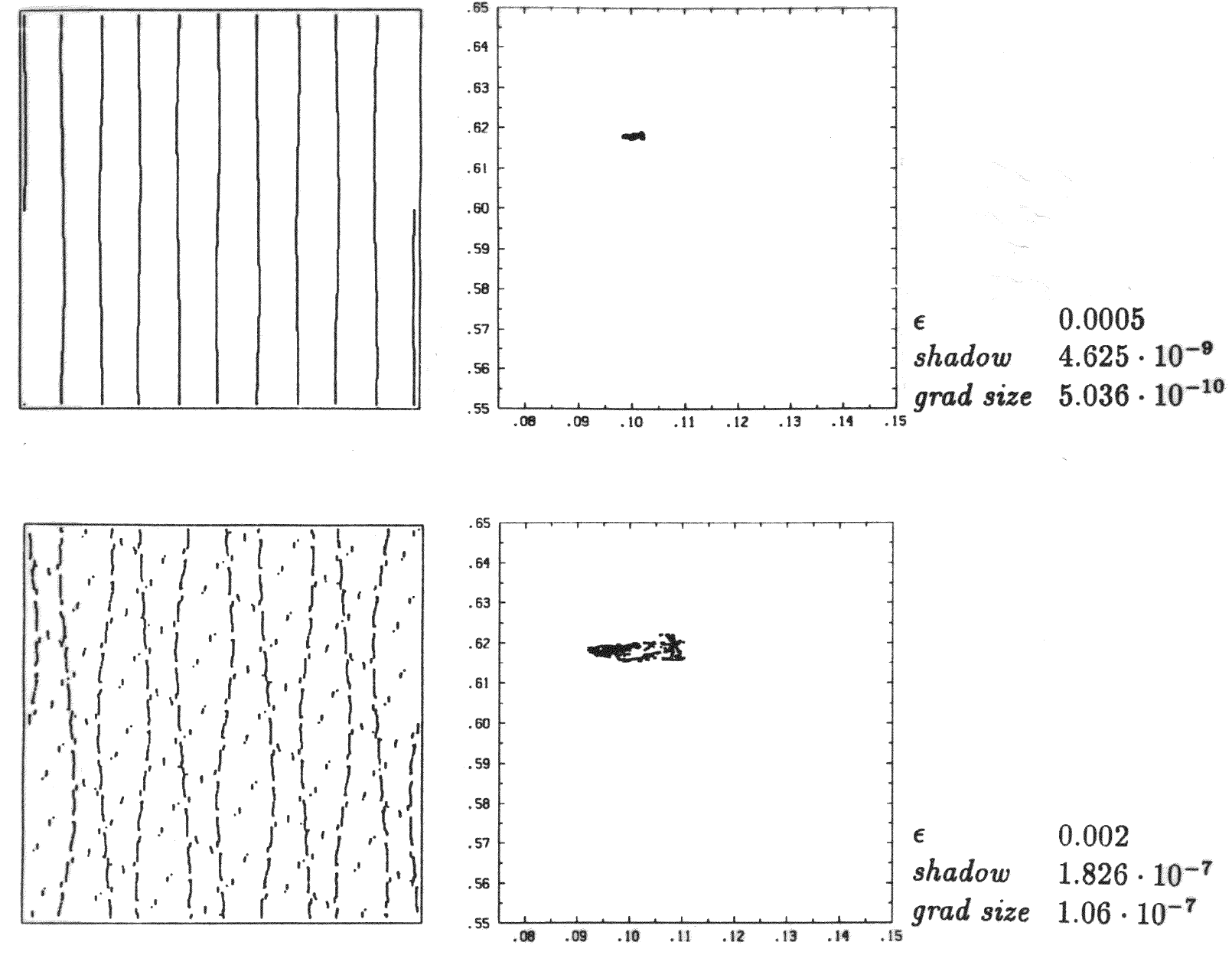}
	\end{center}
	\vspace*{-2\belowdisplayskip}

	\caption{ \em Birkhoff orbits for the 
	polynomial 
	perturbation and the rotation vector (377, 2330)/3770). }
\label{fig:orb-auhd1}
\end{figure}

\begin{figure}[hp]
	\begin{center}
		\includegraphics[height=13.0cm]{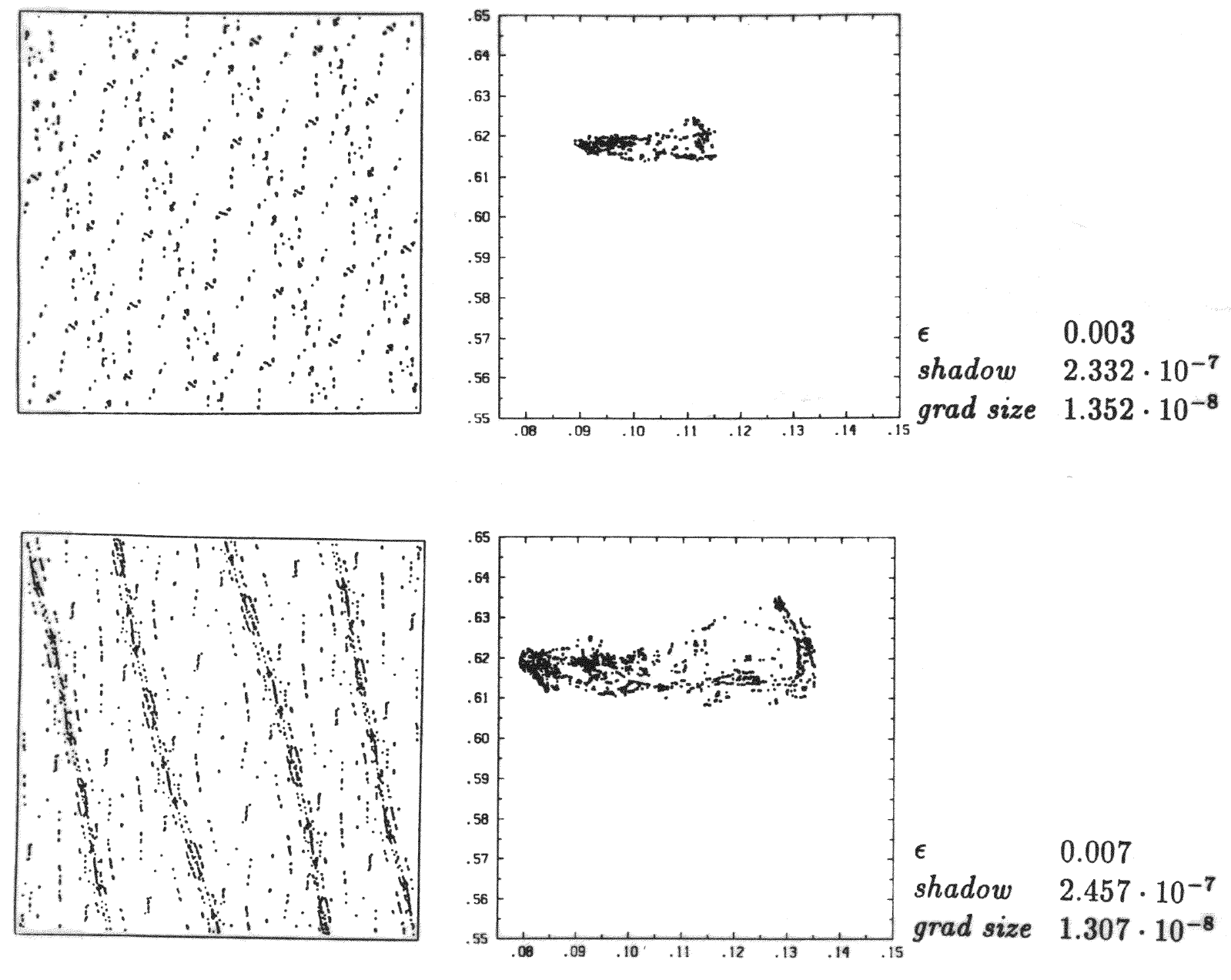}
	\end{center}
	\vspace*{-2\belowdisplayskip}

	\caption{ \em Birkhoff orbits for the polynomial 
	perturbation and the rotation vector (377, 2330)/3770). }
\label{fig:orb-auhd2}
\end{figure}

\begin{figure}[hp]
	\begin{center}
		\includegraphics[height=13.0cm]{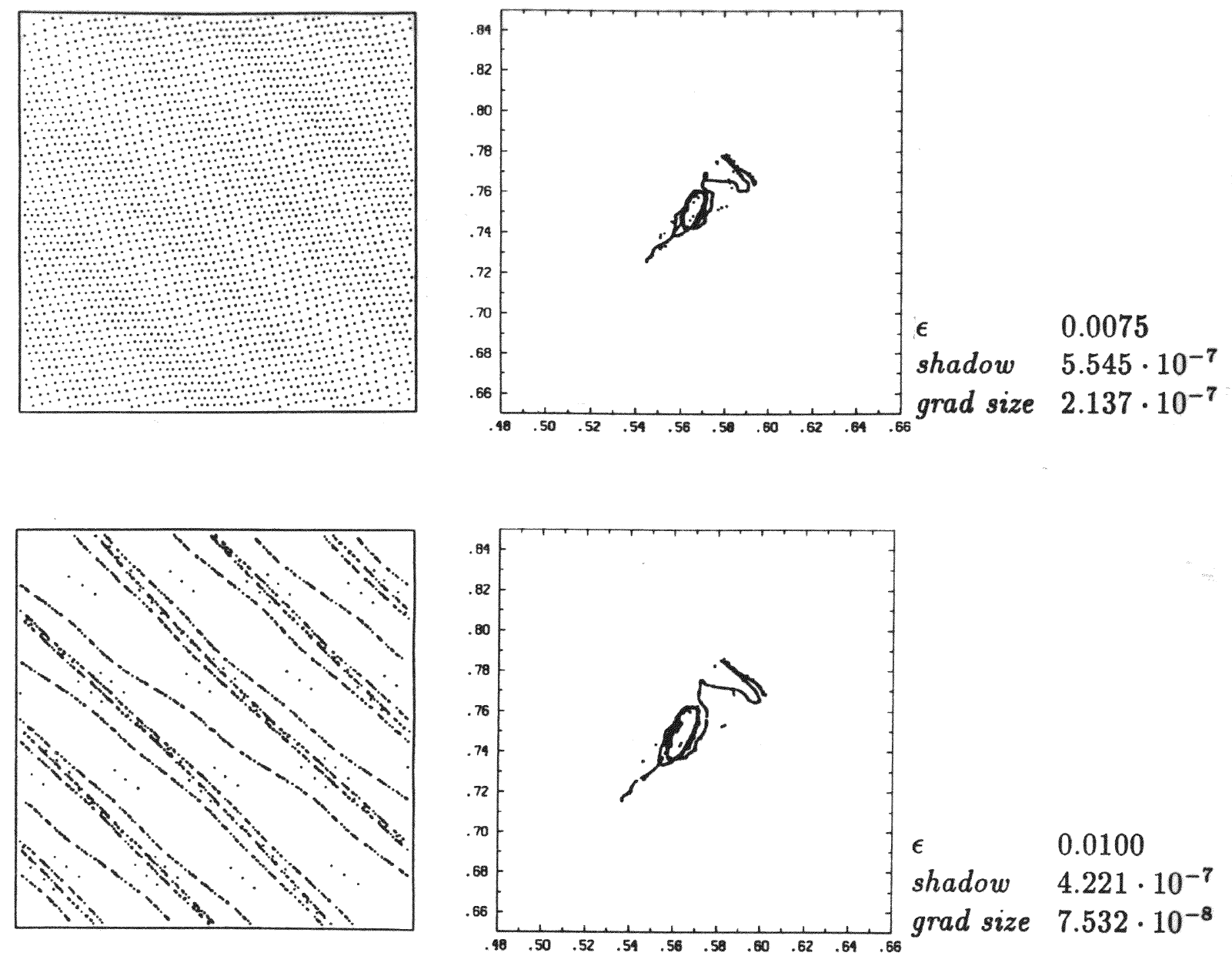}
	\end{center}
	\vspace*{-2\belowdisplayskip}

	\caption{ \em Birkhoff orbits for the fast-Froeschl\'{e} 
	perturbation and the rotation vector (1432,1897)/2513.
	Notice how even the $\epsilon = {\rm 0.0075}$ state 
	seems to have its moment concentrated on a curve. }
\label{fig:orb-spff1}
\end{figure}

\begin{figure}[hp]
	\begin{center}
		\includegraphics[height=13.0cm]{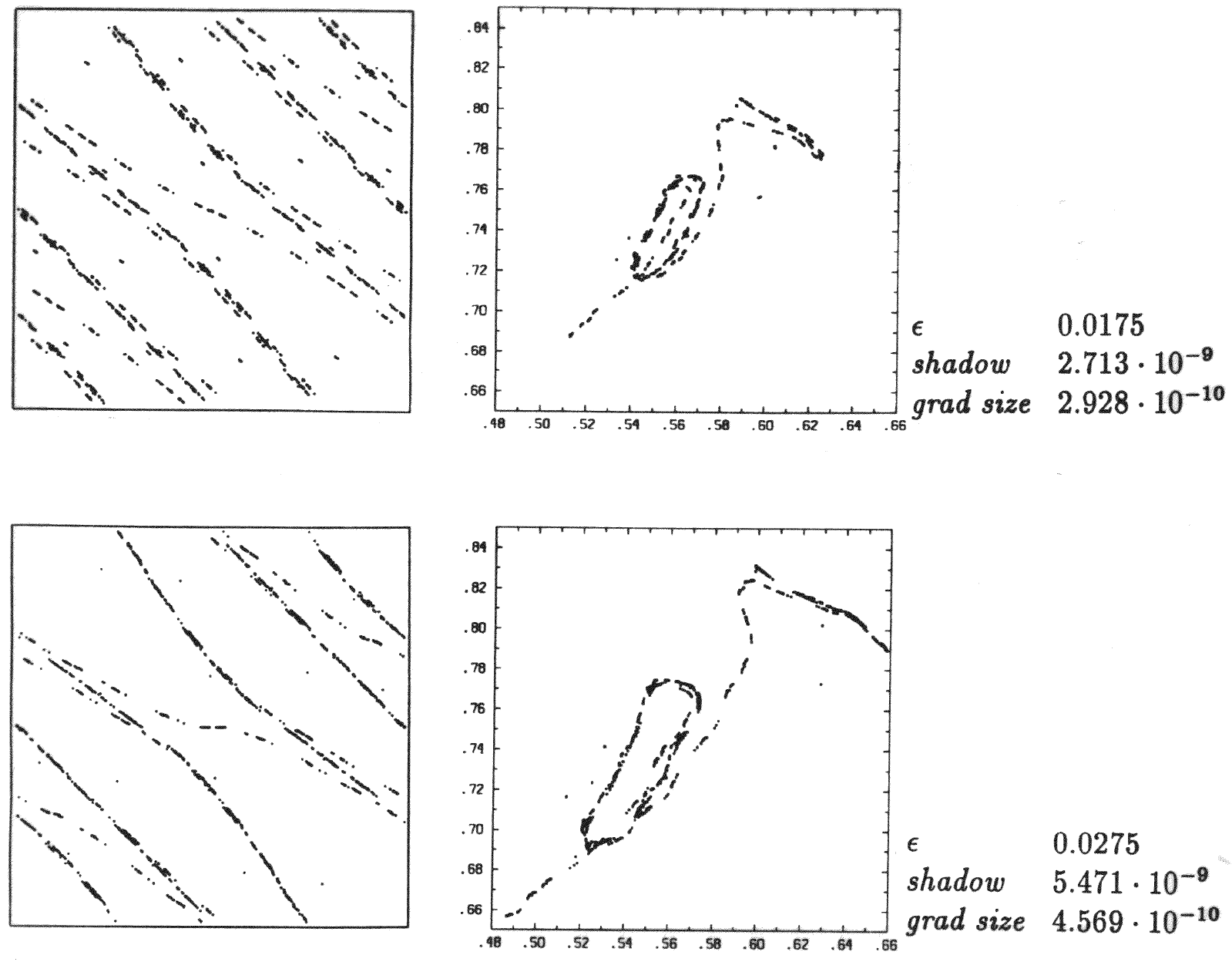}
	\end{center}
	\vspace*{-2\belowdisplayskip}

	\caption{ \em Birkhoff orbits for the fast-Froeschl\'{e} 
	perturbation and the rotation vector (1432,1897)/2513.}
\label{fig:orb-spff2}
\end{figure}

     Figures (\ref{fig:orb-spnf1})--(\ref{fig:orb-spff2}) present several 
families of approximate Birkhoff orbits.  Each orbit is 
displayed as a
pair of projections; one, on the left, is the projection into the 
angular coordinates, the other, on the right, shows the momenta.
Both projections are computed from a p,q-periodic state which is an
approximate solution to the Euler-Lagrange equation (\ref{eqn:KB-ELg}).
The angular projection of a point $\bfx_j$ is an ordered pair 
$(\theta_{j,0}, \theta_{j,1})$, with 
\begin{displaymath}
	\theta_{j,i}  =  x_{j,i} \bmod {\rm 1};
\end{displaymath}
The horizontal is the $\theta_0$ direction and the vertical the
$\theta_1$; both angles lie between 0.0 and 1.0.
The momenta, which are calculated as 
\begin{equation}
	\bfp_j = - \frac{\partial H_\epsilon}
			{\partial \bfx}(\bfx_j, \bfx_{j+1}),
\label{eqn:p-def}
\end{equation}
are arranged similarly; the horizontal is 
the $p_0$ direction and the vertical the $p_1$.

\subsubsection{measures of quality}

Beside each pair 
rotation vector in the form $(p_0, p_1)/q$,
and two measures of the quality of the orbit, {\em shadow} and 
{\em grad size}. 
The first of these measures how closely our orbit, which has
its momenta given by (\ref{eqn:p-def}), approaches the ideal
\begin{eqnarray*}
    (\bfx_{j+1}, \bfp_{j+1}) & = & F_\epsilon(\bfx_j, \bfp_j),\\
    & = & (\bfp'(\bfx_j, \bfp_j),\,\bfp'(\bfx_j, \bfp_j));
\end{eqnarray*}
the value {\em shadow} is 
\begin{eqnarray*}
   & \max_{0 \leq j \leq q-1} & 
    \norm{ (\bfx_{j+1}, \bfp_{j+1}) - F_\epsilon(\bfx_j, \bfp_j)} \\
 =  & \max_{0 \leq j \leq q-1} & 
	\sqrt{
	    \norm{\bfx_{j+1} - \bfx'(\bfx_j, \bfp_j)}^2 +
	    \norm{\bfp_{j+1} - \bfp'(\bfx_j, \bfp_j)}^2 
	} \\
 =  & \max_{0 \leq j \leq q-1} & \bigStrut
	\sqrt{
	    \sum_{k=0}^1 (x_{j+1,k} - x'(\bfx_j, \bfp_j)_k)^2 +
	    		  (p_{j+1,k} - p'(\bfx_j, \bfp_j)_k)^2 
	}.
\end{eqnarray*}
    Most of the states displayed here have {\em shadow} $\approx
10^{-6} $.  The other measure, {\em grad size}, is 
\begin{displaymath}
\left[ \;
	\frac{1}{q} \: \sum_{i=0}^{q-1} \: \left\|  
						\frac{ \partial L_{p,q}}
						{ \partial \bfx_i}
						\right\|^2 \;
\right]^{\frac{\scriptscriptstyle 1}{\scriptscriptstyle 2}};
\end{displaymath}    
it is essentially the norm of the gradient of the action functional,
normalized by the length of the state. 

\subsubsection{shapes}
We display 
orbits for all three perturbations and  for
two rotation vectors, (1432,1897) /2513 and (2330,377) /3770.
The first is an approximation to a irrational vector called 
the {\em spiral mean}, the second approximates $(\frac{1}{10}, \gamma)$,
where $\gamma$ is the golden mean.  Both approximations come from the
Farey triangle scheme of Kim and Ostlund, \cite{KO}, see appendix
\ref{app:num-methods} for details.

For small $\epsilon$ the orbit is well distributed over the angular 
variables and the
momenta look as though they lie on a torus.  With increasing perturbation
the orbits abruptly contract and concentrate along one dimensional filaments.
The system of filaments depends on both the perturbation and the rotation
vector; in figure (\ref{fig:orb-spnf1}b) the (1432,1897)/2513 orbit 
has contracted onto a system of three curves, each of which winds around
the torus once in each angular direction; we will call these curves of
type (1,1).  In figure (\ref{fig:orb-sphd2}b) the same rotation vector
and the polynomial perturbation lead
to a union of seven curves, each of type (0,1). On the other hand,
this same perturbation forces the (2330,377)/3770 state to 
concentrate along a curve of type  (4,1).

\subsubsection{Lyapunov exponents}

\begin{figure}
	\begin{center}
		\includegraphics[height=13.5cm]{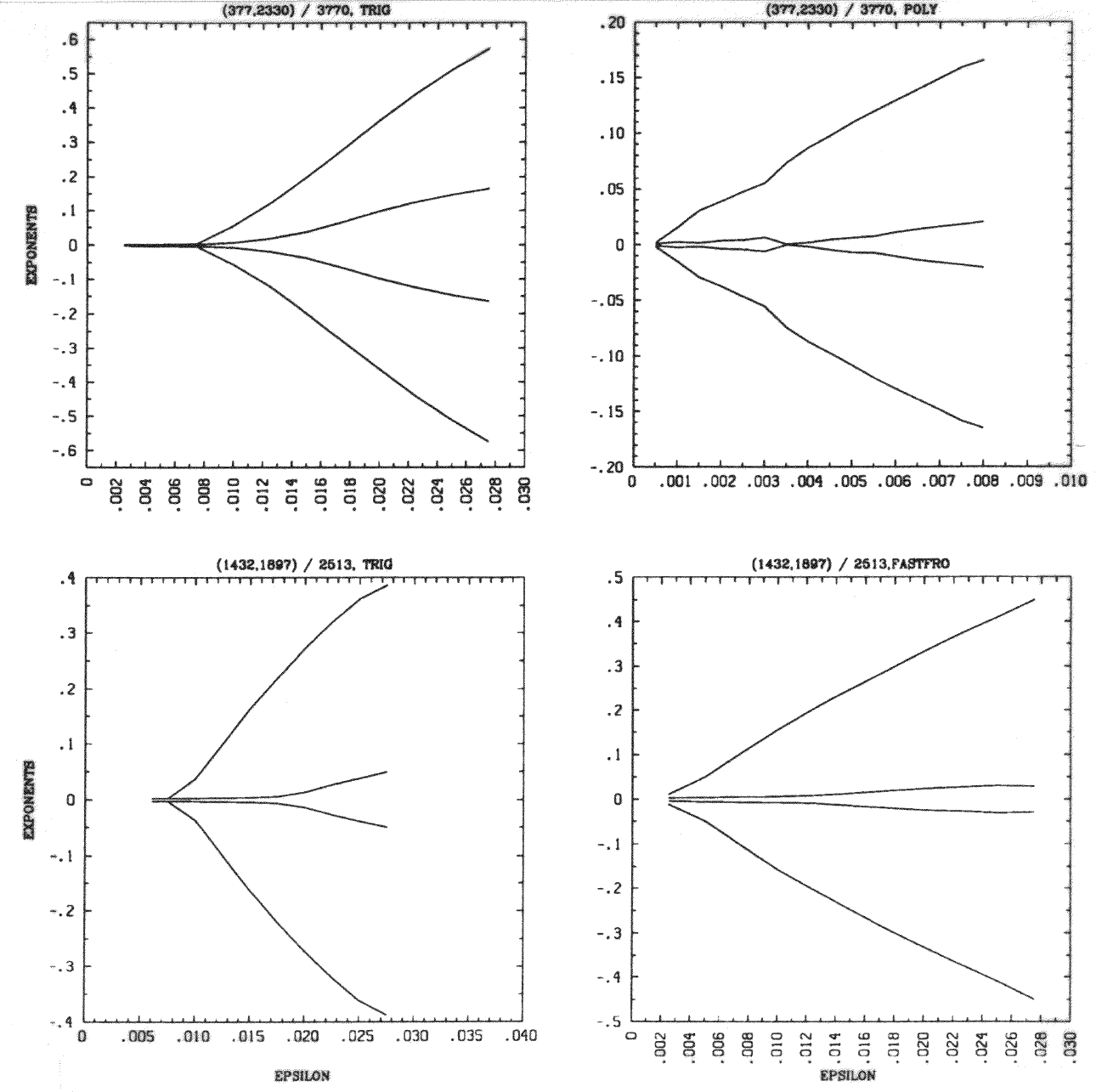}
	\end{center}
	\vspace*{-2\belowdisplayskip}
	\caption{ \em   The Lyapunov exponents for the rotation 
			vector (377,2330)/3770 and the
			trigonometric and
			polynomial perturbations.
			Also those for the vector
			(1432,1897)/2513 with the 
			trigonometric and
			fast-Froeschl\'{e}  
			perturbations. }
\label{fig:Lyaps}
\end{figure}

The qualitative behaviour of the orbits is correlated with their 
stability properties.  The Lyapunov exponents measure the 
exponential rate of divergence of nearby trajectories (see, e.g.,
\cite{Osc}) and, for a periodic
orbit, are just the eigenvalues\thesisfootnote{
     The accurate, direct calculation of the matrix product in 
     (\ref{eqn:DF}) is usually not possible; see appendix 
     \ref{app:num-methods} for a discussion.
} of
\begin{equation}
     DF^q_{\epsilon, (x_0, p_0)} \;\; = \;\;
     DF_{\epsilon, (x_{q-1}, p_{q-1})} \circ  
     DF_{\epsilon, (x_{q-2}, p_{q-2})} \circ \cdots \circ 
     DF_{\epsilon, (x_0, p_0)} 
\label{eqn:DF}
\end{equation}
where $DF_{\epsilon,(x,p)}$ is the Jacobian of the map.
From \ref{eqn:map} we can calculate
\begin{displaymath}
DF_{\epsilon,(x,p)} \; = \; \left[ 
		    \begin{array}{cc}
			\frac{\displaystyle \partial \bfx'}
			     {\displaystyle \partial \bfx}  &
			\frac{\displaystyle \partial \bfx'}
			     {\displaystyle \partial \bfp}  \\
			\frac{\displaystyle \partial \bfp'}
			     {\displaystyle \partial \bfx} \bigStrut &
			\frac{\displaystyle \partial \bfp'}
			     {\displaystyle \partial \bfp} \bigStrut 
		    \end{array}
		\right]
	\; = \;
		\left[
		   \begin{array}{cc}
		       {\bf I} -
		       \frac{\displaystyle \partial^2 V_\epsilon}
			    {\displaystyle \partial \bfx^2} &
		       -{\bf I} \\
		       -\frac{\displaystyle \partial^2 V_\epsilon}
			     {\displaystyle \partial \bfx^2} \bigStrut &
		       {\bf I}  \bigStrut 
		    \end{array}
		\right]
\end{displaymath}
where {\bf I} is the {\em d\/}-dimensional identity matrix and
$\partial^2 V_\epsilon / \partial x^2$ is the Hessian of the
perturbation.
Each of the $ DF_{\epsilon, (x_i, p_i)} $ is
a real symplectic matrix and so the entire product is 
real and sympectic too.
The eigenvalues of $DF^q_{\epsilon, (x_0,p_0)}$ 
thus occur in reciprocal pairs 
\mbox{ $(\lambda_0, 1/\lambda_0)$ } and 
\mbox{ $(\lambda_1, 1/\lambda_1)$ }, \cite{Arnold:Math-Meth};
for the unperturbed map, all four are equal to 
one.  As the perturbation increases first one pair, then the other,
depart from the unit circle.  At about the same parameter value for
which the first pair leaves the circle we see the minimizing state
contract along the filaments.   For large enough perturbation both
pairs are non-zero and the distribution along the direction of the
filaments is also Cantor-like.  See figure~(\ref{fig:Lyaps}) for
the exponents of most of the orbits presented here. 

     At about the same value of the perturbation for which the states
begin to concentrate along filaments, the first pair of Lyapunov
exponents departs form the unit circle.  The eigenvector corresponding
to the largest exponent projects to a vector transverse to 
the filaments.
As we increase the perturbation further the states begin to form into
clumps along the direction of the filaments until, in the last panels
of each series of orbits, the orbits are concentrated near points.

\subsection{non-existence of tori: a prelude}

     Notice that the very perturbed orbits look as though they are full
of holes, as though there are some parts of the torus they cannot visit.
One might imagine that this is just a consequence of the finite lengths
of the our orbits, that if we had orbits with ten times as many 
points some of them would be bound to land in the holes.  We can show
that, for sufficiently large perturbations, the holes are genuine; there
are neighborhoods which all minimizing Birkhoff orbits must avoid.  

  Suppose $V_\epsilon(\bfx)$ is a $C^2$, standard-like perturbation 
to the generating function 
$H_0(\bfx,\bfx') = \half \, \norm{\bfx'-\bfx}$.
Suppose further that $V_\epsilon(\bfx)$ has a
minimum at \mbox{$\bfx=\bfx_{min}$}.  
Then there is an $\epsilon_c$, such that for
$\epsilon > \epsilon_c$, all minimizing states must avoid a region 
containing $\bfx_{min}$.
\newline
{\bf Proof} \hskip 1\parindent   A globally mininimizing state, 
$X$,
must be an extremum of $L_{p,q}$ such that every small, local,
variation, $x_i \rightarrow x_i + \delta$ increases the action.  
That means  that $X$ must satisfy the Euler-Lagrange equations 
(\ref{eqn:KB-ELg}) and also that
\begin{eqnarray}
\frac{\partial^2  L_{p,q}}{\partial \bfx_i^2} & = &
    \left[ 
	\begin{array}{cc}
	    2 - \epsilon \, \frac{{\displaystyle \partial^2 V}}
			 {{\displaystyle \partial x_0^2}}(\bfx_i) &
	    -\epsilon \, \frac{{\displaystyle \partial^2 V}}
			     {{\displaystyle \partial x_0
				     \partial x_1}}(\bfx_i) \\
	    -\epsilon \, \frac{{\displaystyle \partial^2 V}}
			     {{\displaystyle \partial x_0
				 \partial x_1}}(\bfx_i) \bigStrut  &
	    2 - \epsilon \, \frac{{\displaystyle \partial^2 V}}
		  {{\displaystyle \partial x_1^2}}(\bfx_i) 
	\end{array}
    \right],
\label{eqn:maxHess}
\end{eqnarray}
is positive definite.  Because $\bfx_{min}$ 
is a minimum, the eigenvalues,
\mbox{$\mu_0(\epsilon) \leq \mu_1(\epsilon)$},
of the Hessian of $-V_\epsilon(x_{max})$
are negative.  If one of them is less than $-2$ then (\ref{eqn:maxHess})
cannot be satisfied.  Since the $\mu_i$ are decreasing functions of 
$\epsilon$ we need only find that value, $\epsilon_c$, for which
$\mu_0(\epsilon_c) = -2$.

     For the trigonometric perturbation 
\mbox{$\epsilon_c \approx$ 0.03856}; for the polynomial perturbation 
\mbox{$\epsilon_c \approx$ 0.04167}.  The appearance of the states
suggests that neither of these is a very good estimate; the region
near the maximum is completely devoid of points long before 
$\epsilon = \epsilon_c$.   The real interest of an argument like
the one above is that it can provide an estimate of the size of
perturbation needed to destroy all the original invariant tori;
since the whole next chapter is devoted to such estimates, 
we leave the subject for now.

%% file: numres/regularity.tex
\newcommand{\normdx}{{\norm{{\scriptstyle \Delta} \bfx}}}
\newcommand{\normdp}{{\norm{{\scriptstyle \Delta} \bfp}}}
\subsection{smoothness }

    We would like to be able to say that very long periodic orbits 
approximate a Cantor set which we could view as the tattered
remnant of an invariant torus.   Such a remnant would have a kind
of smoothness; two points which lay lie very close to each other in 
the angular variables should not have wildly different
momenta. What we need is a result like the theorem of Birkhoff, 
generalized by Katok \cite{Kat:rem}, which 
says that for points in a Mather set,
the momenta are Lipschitz functions of the coordinates, i.e.
\mbox{ $ \norm{p_i - p_j} \leq C \, \norm{x_i - x_j} $ }
where $C$ is a constant.
Katok and Bernstien \cite{KB:Birk} looked for such a result and,
as mentioned above, were able
to show that, except perhaps at one point, the momenta are
H\"{o}lder continuous with index $ 1 / 2$, that is,
\begin{displaymath}
	\norm{\bfp_i - \bfp_j} \leq 
	\norm{\bfx_i - x_j}^{\alpha}  \qquad \alpha = \half.
\end{displaymath}
for some constant $C$ independent of the $\bfx_i$.

    Hoping to verify or improve their estimate, we computed pairs
\mbox{ $(L, \normdx ) $ }, where 
\mbox{ $ L = \normdp / \normdx $ }, and displayed  them
on logarithmic axes.  If some kind of H\"{o}lder continuity applies,
then
\begin{displaymath}
	L = \frac{ \normdp }{\normdx } \leq
		C \, \normdx^{\alpha - 1},
\end{displaymath}
so
\begin{eqnarray*}
	\log L & \leq & \log C + (\alpha - 1) \log \normdx.
\end{eqnarray*}

\begin{figure}
	\begin{center}
		\includegraphics[height=19.0cm]{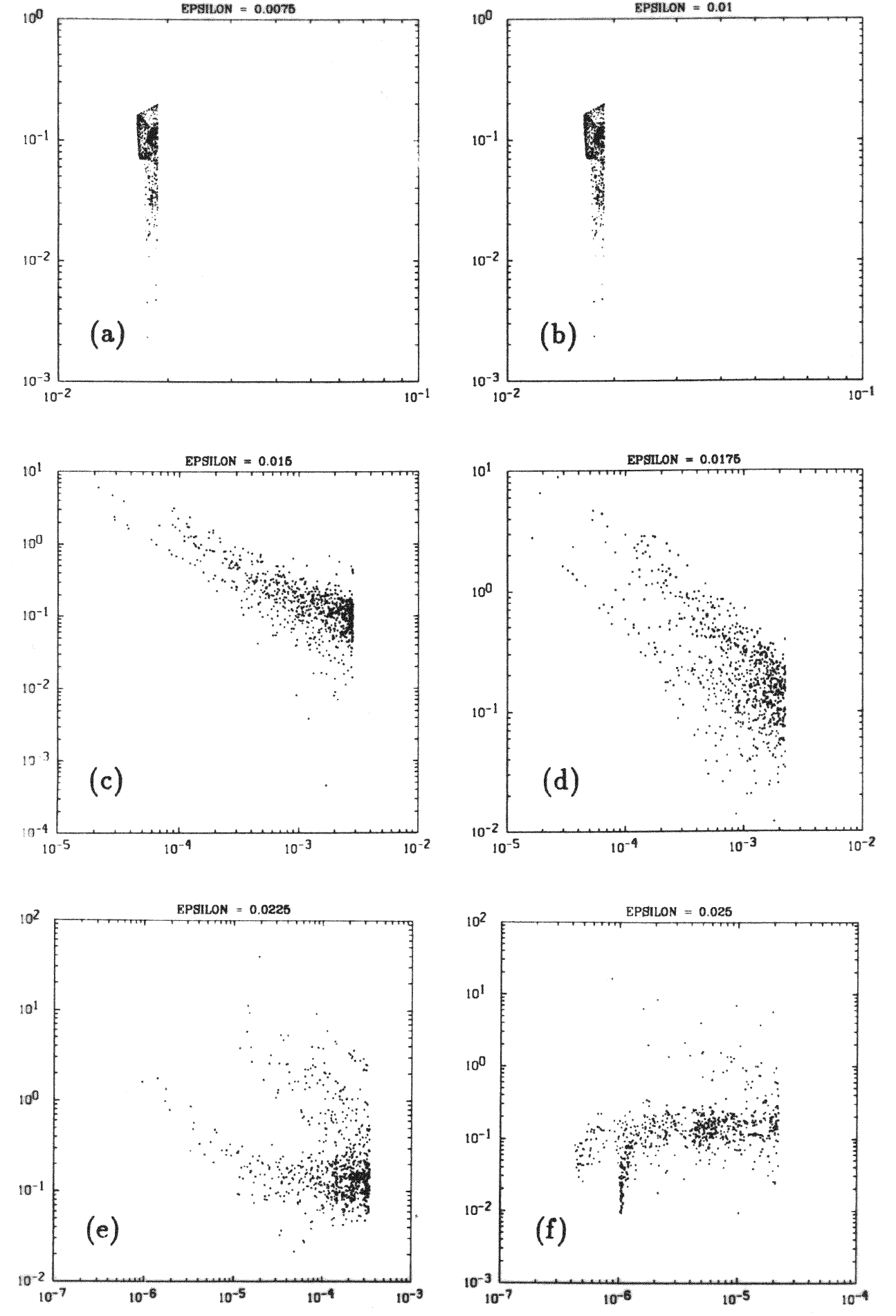}
	\end{center}
	\vspace*{-2\belowdisplayskip}
    \caption{\em  Pairs \protect {$( L, \normdx )$} calculated for the 800
		most closely spaced pairs of points in states
		of the rotation vector (1432,1897)/2513 with the 
		trigonometric perturbation.
	    }
\label{fig:smooth}
\end{figure}

     We can tell whether our orbits are compatible with Lipschitz continuity
by looking at the upper envelope of $(L, \normdx ) $.  If the
envelope is a decreasing function of $\normdx $ then the  H\"{o}lder
index is less than one and the momenta are not 
Lipschitz functions.  If the
envelope is flat or sloping upward then the continuity is 
Lipschitz or better.
Figure (\ref{fig:smooth}) shows some collections of 
\mbox{$(L, \normdx )$} pairs.  
The results are ambiguous at best.  At very small 
perturbation the upper envelope has a positive slope, see figure
(\ref{fig:smooth} parts a and b).  For intermediate values of $\epsilon$, 
those for which the orbit has contracted  
into filaments but has not yet begun to concentrate in points, the 
situation
looks worse; the largest values of $L$ occur for the 
smallest values of $\normdx$, see figures (\ref{fig:smooth}parts c and d).
This would seem to doom any hope that $p$ is a Lipschitz function
of $x$.  Note, however, that the upper envelope has a slope of
$-1$.  This suggests that $\normdp \approx \mbox{\em const}$.
On the other hand, we have, from Katok and Bernstien, that $\bfp$ is 
H\"{o}lder $\half$ .  It is thus possible that the lack of smoothness
may come from not having enough points.
    At very large $\epsilon$, those for which the orbit has 
    contracted into
a few small clumps, $(L, \normdx )$ begins to have an increasing 
envelope again.  Unfortunately, it is just at these very short distances
that we must begin to doubt the quality of our orbits.  Typically we
have {\em shadow} = $10^{-6}$ and so must expect the 
$\bfx$'s, $\bfp$'s and their differences to be uncertain at about 
that level too.

    Finally, we note that the uncertainty in the $\bfp$'s 
    could expalin the 
behavior at intermediate $\epsilon$.  If the components of 
$\bfp$'s are uncertain beyond $\sigma_p$, their differences 
are uncertain to $\sqrt{2}\,\sigma_p$.
Then, no matter what the continuity properties of $p$, for small enough
$\normdx$ we should expect to see $\normdp \approx const.$
This explaination is not vompletely satisfactory in that it 
fails to explain
why some of the graphs in figure~\ref{fig:smooth} seem to have two
different populations of constant $\normdp\/$'s.

%% file: numres/hedlund.tex
\section{Hedlund's examples}
\label{sec:Hed}

    In this section we will worry 
about whether the shapes of our states 
have anything to say about the shapes of much longer states 
with similar rotation vectors.  A central premise of our program 
of rational approximation
is that they do; unfortunately, except for the two  
dimensional case (twist maps on the cylinder), 
we cannot prove this.  We cannot even show
that states with the {\em same} rotation vector 
must have the same shape.
Consider the family of minimizing states with rotation vectors, 
\begin{displaymath}
	\frac{\bfp_0}{q_0}, \, 
	\frac{2\bfp_0}{2q_0}, \, 
	\ldots , \, 
	\frac{n\bfp_0}{nq_0}, \ldots 
	\qquad n \in {\bf Z}^+,
\end{displaymath}
where $\bfp_0 / q_0$ is in lowest terms.
For each of these states there is certainly one solution to the 
Euler-Lagrange equation which
is just a concatenation of $n$ copies of the $\bfp_0 / q_0$ minimizing
state, but there may also be other solutions, some of which may have
lesser total action.

    To see how this can happen,  we consider the problem of finding
{\em minimal geodesics}, curves of smallest possible length, on either
the two (or three) dimensional torus.  This problem arises, 
for example, in the motion of a free
particle in a system with periodic boundary conditions and
could be reduced to a symplectic map via a surface of section,
but in the discussion below it will be simpler to think about
continuous time and smooth trajectories.   
We will work with two different representations of the problem,
one on the two (or three) dimensional torus and another made by 
periodically extending the torus to get the plane (or ${\rm R^3}$). 
In either representation, we will allow the metric to be other than 
the usual Euclidean one.

In the $\Rn$ version of the problem, a minimal geodesic is a curve,
\mbox{$ \gamma : \R \rightarrow \Rn$}, parameterized in terms of, say, 
arc length and for which every 
finite segment is the
shortest possible curve connecting its endpoints.  
Our special interest will be the {\em periodic geodesics};  
on the torus these are curves which wind around and eventually begin to
retrace themselves.  In $\Rn$ they appear
as curves for which  $\exists \tau \in \R$ such that
\begin{equation}
     \gamma(t + \tau) = \gamma(t) + \bfm,\;\;\bfm \in \Zn
\label{eqn:Hed}
\end{equation}
and 
we may classify them according to $m$, which gives the number of times
$\gamma$ winds around each of the coordinate directions on the torus
before repeating itself.   Hedlund studied these curves on
the two dimensional torus and, in \cite{Hed:ex}, showed that
for every pair $(m_0, m_1) \in {\rm \Z^2}$, there is a minimal periodic 
geodesic which winds $m_0$ times around the $\theta_0$ direction 
and $m_1$ times in the $\theta_1$ direction before closing.

He also made an observation which connects the geodesic problem
to the problem of finding Birkhoff periodic orbits.  He asked 
whether, for example, the minimizing periodic geodesic for the pair
(10,20) could be different from the which traces 10 times
over the (1,2) geodesic.  He found that it could not.
The corresponding statement for Birkhoff orbits is that the pathology
outlined at the beginning of the section does not occur for 
two dimensional twist maps of the annulus.

In the last section of his paper, Hedlund demonstrated that 
one cannot expect the analagous result in higher dimension. 
He presented an explicit example of a metric on ${\bf T^3}$ for which
the shortest geodesic of type $(ni,\,nj,\,nk)$ is not {\em n} copies
of the shortest $(i,\,j,\,k)$ geodesic.  Victor Bangert \cite{Bang:min}
has proved that a metric on $\Tn$ has at least $n+1$ minimal geodesics
and has given some principles for the design of Hedlund-type examples.

\begin{figure}[bp]
	\begin{center}
		\includegraphics[height=6.0cm]{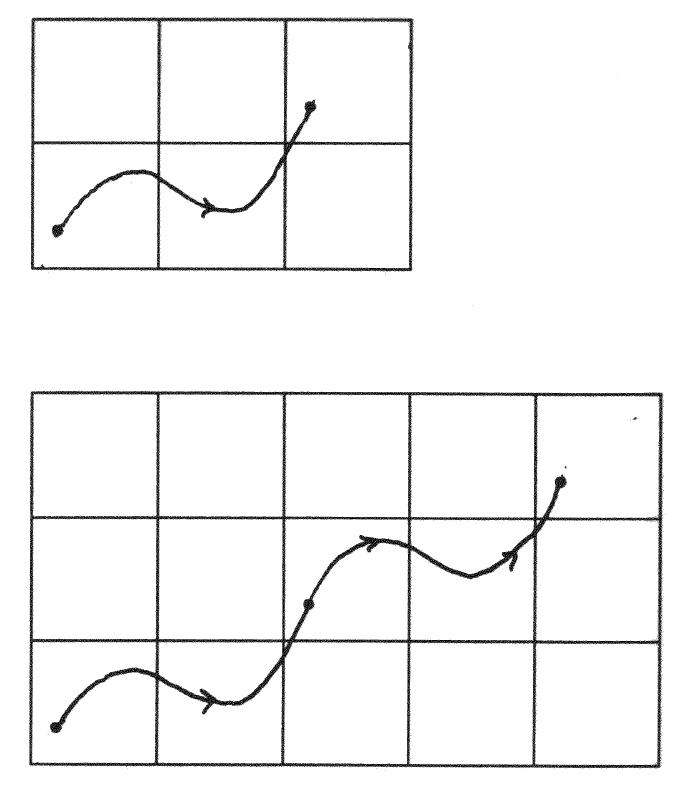}
	\end{center}
	\vspace*{-2\belowdisplayskip}
    \caption{\em Some minimizing periodic geodesics for
    the two dimensional torus; the shortest curve of type
    (2,4) is just 2 copies of the shortest one of type
    (1,2). }
\label{fig:hed1}
\end{figure}

\begin{figure}[bp]
	\begin{center}
		\includegraphics[height=15.0cm]{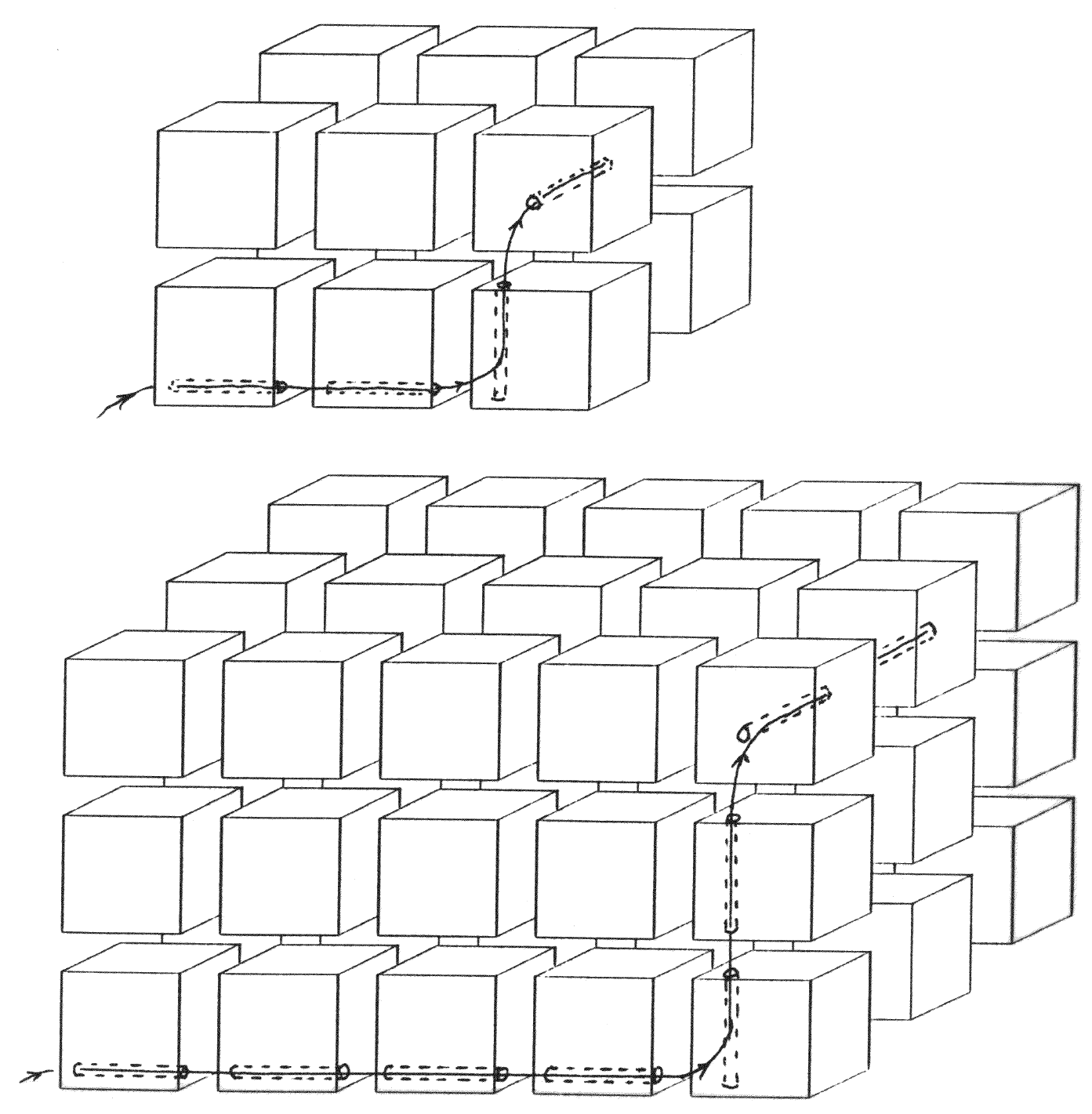}
	\end{center}
	\vspace*{-2\belowdisplayskip}
    \caption{\em Some minimizing periodic geodesics for
    a Hedlund example on
    the three dimensional torus; the shortest curve of type
    (2,4,2) is {\em not} 2 copies of the shortest one of type
    (1,2,1). }
\label{fig:hed2}
\end{figure}

    Figures (\ref{fig:hed1}) and (\ref{fig:hed2}) contain 
the main ideas.  Bangert sets up the metric so it has
certain non-intersecting 
lattices of ``tunnels,'' tubes in the middle of which the 
metric is so small that the length of a segment is,
at most, say, $1 / 100$ of its Euclidean length.  
Outside the tunnels the
metric is such that the length of a segement is a bit longer than its
Euclidean length. In Bangert's examples the tunnels run along the lines
$(0, t, \half),\;(\half, \half, t),$ and $(t, 0, 0),\; t \in \R$ 
and along
all their $\Zn$ translates.  Under these rather severe conditions he is 
able to show that a minimizing geodesic must spend essentially all its
time inside the tunnels, venturing out only to leap from one system of
tunnels to another.


     A minimizing, periodic geodesic then has only three 
short segments lying outside the tunnels, no matter how long it is.  
Note that such a geodesic
strays a long way from the straight line which connects its 
endpoints; the latter is a minimizing periodic geodesic for the 
flat, Euclidean metric.
In the language of Birkhoff orbits, Hedlund's pathology would occur  
if some few {\em p-q} periodic states turned out to have such 
tiny actions that all very long states would be composed of
a few segments, with each segment containing many copies of 
the few economical states.   Although we cannot preclude this 
possibility, we feel it is unlikely.  Hedlund and Bangert's 
examples require that the curves through the
tunnels be much, much shorter than their Euclidean lengths, 
consequently, their metrics are very far from flat.  By contrast, our
generating functions are close to the unperturbed ones.  
We might thus hope that our minimizing states are obliged 
to stay close to the unperturbed states.  Katok has shown, 
in \cite{Kat:min}, that if the perturbed states stay within 
some bounded distance of the unperturbed distance and if the
bound is independent of the length of the state, then Hedlund's
pathology does not occur.

\begin{figure}
	\begin{center}
		\includegraphics[height=13.0cm]{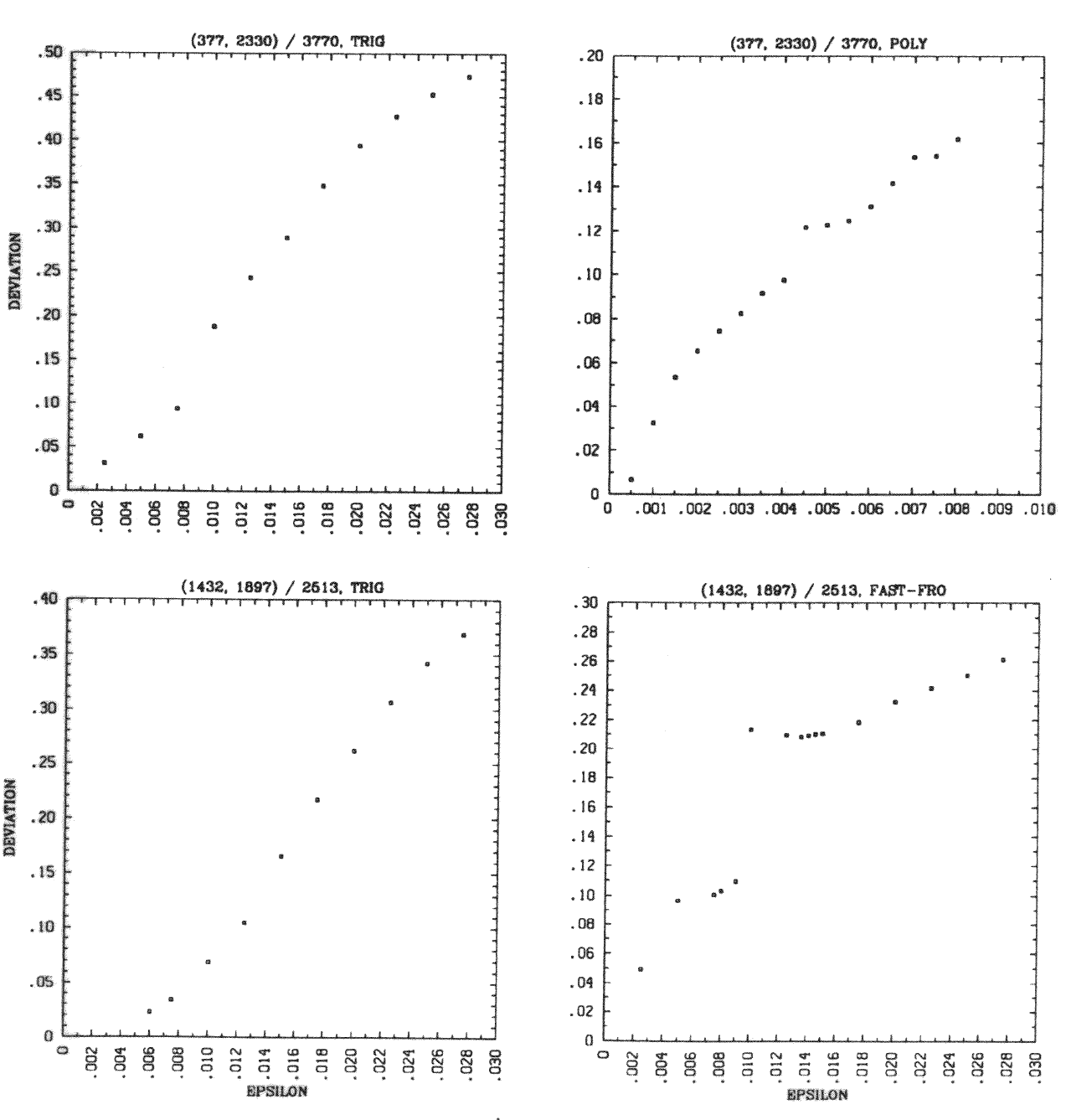}
	\end{center}
	\vspace*{-2\belowdisplayskip}
    \caption{\em The largest displacement between a point in a 
    perturbed minimizing state and the position it would occupy
    in the abscence of the perturbation. Note the abrubt jumps
    in the deviations for the fast-Froeschl\'{e} example.}
\label{fig:deev}
\end{figure}

\begin{figure}
	\begin{center}
		\includegraphics[height=18.0cm]{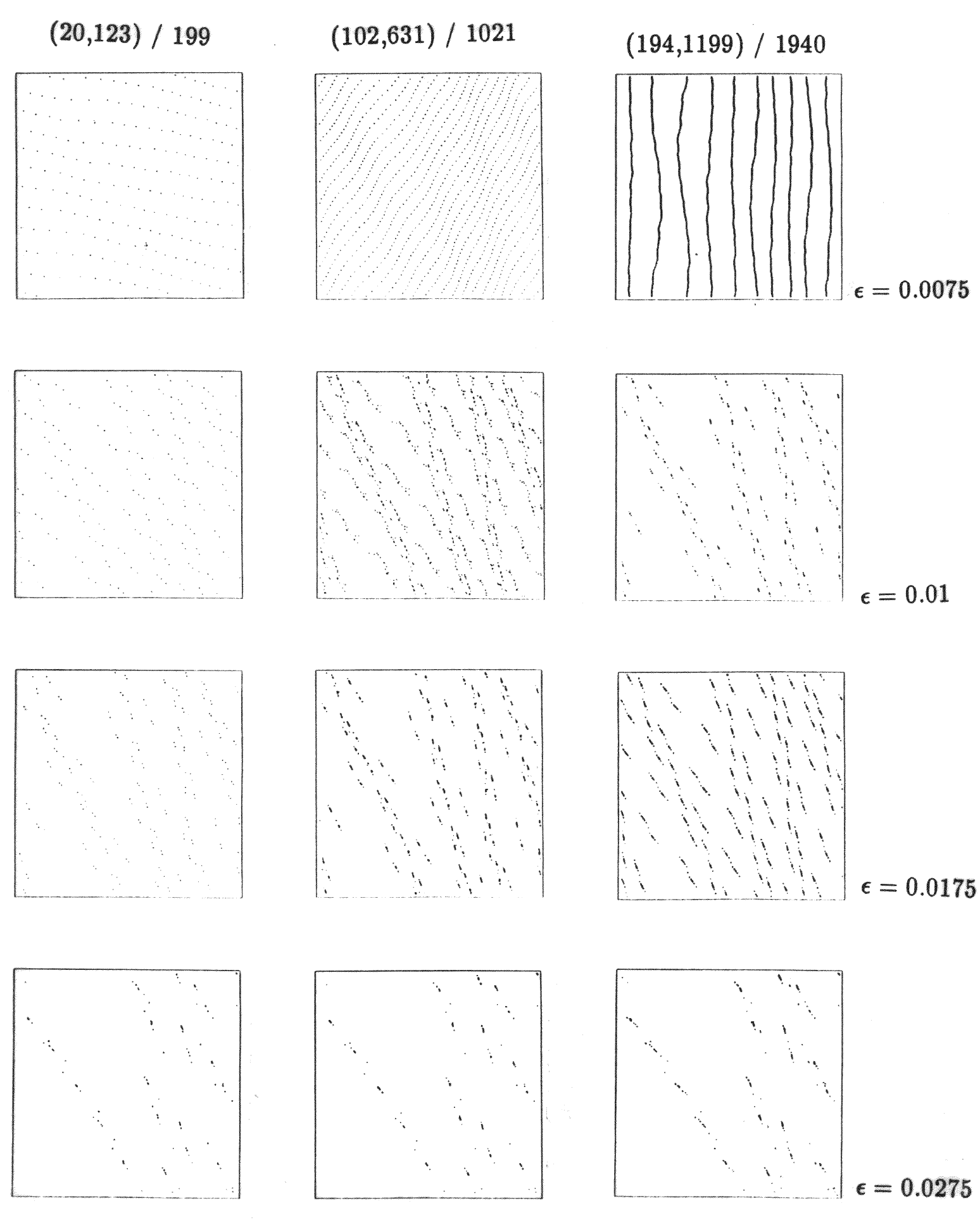}
	\end{center}
	\vspace*{-2\belowdisplayskip}
    \caption{\em A series of orbits whose rotation vectors approximate
		 (377,2330) / 3770.  
	    }
\label{fig:hed_st}
\end{figure}

    We undertook two studies to investigate these issues.  In the first,
figure (\ref{fig:deev}), we measured the deviation of our minimizing
states from the straight line connecting $x_0$ to $x_q$.  The distance
always remains smaller than the diameter of the torus, $1 / \sqrt{2}$.
In the second study we used the Farey triangle algorithm of 
Kim and Ostlund,
(see appendix \ref{app:num-methods}), to get a sequence of 
rotation vectors tending to $(377, 2330) / 3770$.  The states 
for these vectors are displayed in figure \ref{fig:hed_st}.  
The longest orbits look very much like the shortest.
We also did some experiments on families of rotation vectors of the
form\thesisfootnote{
	An unperturbed minimizing state is 
	{\em n} copies of the unperturbed $p_0 /q_0$
	state and our procedures for constructing 
	perurbed minimizing states are 
	such that this shorter, internal periodicity
	would be retained throughout the calculation.
	We tried to circumvent this problem by adding
	a small, random displacement to each of the
	points in the starting guess, see appendix \ref{app:num-methods}.
}
$n\bfp_0 / nq_0$;
The longer states were indistinguishable from the shorter ones.

%% file: Chap3.tex
\input{converse/intro3}
\input{converse/cylinder}
\input{converse/prism}
\input{converse/higher}
\input{converse/results}

%% file: converse/intro3.tex
\chapter{The Frontier of Chaos}

     Our first investigations aimed at the question ``What remains after
invariant tori have been destroyed?''  Our next set asks the more basic
``How could we tell if the tori were there?''  
To answer this question we might  follow
Kolmogorov, Arnold and Moser and seek to find perturbations so small
that some tori would be guaranteed to exist.  
Conversely, we could try to 
find perturbations so large that no invariant tori remain. 
Numerical evidence suggests that the first approach will be hard; 
tori seem to persist well beyond the point where traditional KAM arguments 
break down.\thesisfootnote{
	Several authors have now proved machine-assisted, constructive
	KAM theorems for specific maps; these are in much better 
	agreement with non-rigorous numerical predictions.  See
	e.g. \cite{Alessandra}, \cite{Rana:thesis}, and \cite{LlR:small}. 
}
We will 
adopt the latter strategy; we will try to fill in the blanks in the
following ``converse KAM'' theorem :

\noindent
{\bf Theorem }  
{\em For the n-dimensional symplectic 
	twist map ${\rm F}_{\epsilon} : \An \rightarrow \An$,
	\begin{displaymath}
		{\rm F}_{\epsilon}(x,r) = (x',r') 
					= {\framebox[10em]{\rule{0cm}{1ex}}}
	\end{displaymath}
	depending on the parameters, $\epsilon $, 
	we are guaranteed that no KAM tori exist 
	for any \mbox{$ \epsilon \in {\rm S}_{\rm F} = 
	\{ {\framebox[10em]{\rule{0cm}{1ex}}} \} $. }
}

\noindent
{\bf Proof } \framebox[29em]{\rule{0cm}{1ex}}

Herman, in \cite{Herm:big} first saw that one might get a better 
notion of where invariant tori exist by looking at the edge of 
the region where they do not.  He considered maps, 
$T_{\epsilon}: \T \times \R \rightarrow \T \times \R$,
of the form\thesisfootnote{
	Our examples are not of this form, but, after a 
	change of coordinate,  their inverses are. 
} 
\begin{equation}
	T_{\epsilon}(x,p) = (x',p') = (x + p, p + \epsilon f(x+p)),
\label{eqn:Herm}
\end{equation}
small perturbations to the integrable system, and envisioned a kind of
cartography of non-integrability.  By choosing different $f\/$'s he 
could consider different directions in the space of perturbations.
For each fixed $f$ he could increase the value of $\epsilon$ until it
reached a size, $\epsilon = \epsilon_c(f)$, such that no invariant tori 
remained.  By calculating pairs $(f, \epsilon_c(f))$ he could map out
the edge of non-integrability, the frontier of chaos.
\begin{figure}
	\begin{center}
		\includegraphics[height=3.75cm]{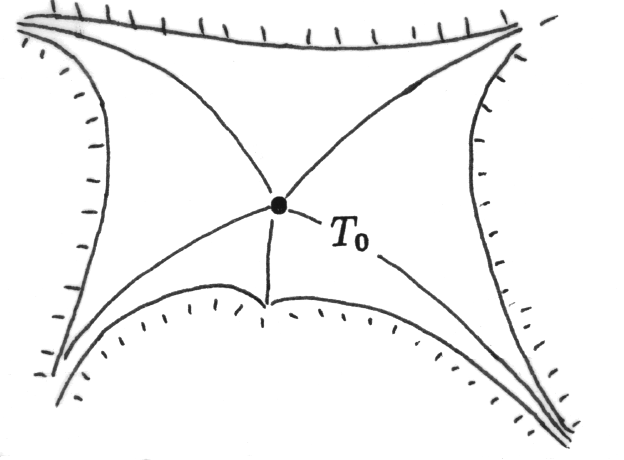}
	\end{center}
	\vspace*{-2\belowdisplayskip}

	\caption{\em 
		The space of near-integrable 
		maps, showing the frontier of 
		non-integrability around $T_0$, 
		an integrable system.}
	\label{fig:frontier}
\end{figure}

    We will concentrate on ways to get rigorous bounds for 
$\epsilon_c(f)$ but will not make a very extensive 
survey\thesisfootnote{
	Jacob Wilbrink, in \cite{Wilb:circ}, used
	a non-rigorous existence criterion to survey 
	a whole one parameter family of maps.  
} of $f$'s.  The rest of the chapter is organized by dimension
of the phase space and sharpness of non-existence criteria.
In the next section we review converse KAM theorems for
area-preserving twist maps on the cylinder, and in section
\ref{sec:machines} we explain
how to prove them with a digital computer.
In \ref{sec:higher-d} we formulate some criteria for
higher dimensional systems and finally, in section
\ref{sec:conv-pix}, apply them to an example.  

%% file: converse/cylinder.tex
\section{Converse KAM results on the cylinder}
\label{sec:cone-families}

Most of the ideas presented here originated with Herman's  paper
\cite{Herm:big}.  Mather picked up these techniques and made applications
to the standard map, \cite{Ma:circ}, and to billiards, \cite{Ma:balls}.
He also introduced a different, more generally applicable
criterion based on the existence of action-minimizing states.
MacKay and Percival augmented Herman's argument with rigorous computation
and discovered a connection between Herman's work and Mather's
action criterion.\thesisfootnote{  
	Recently, Rafael de la Llave (personal communication) 
	has developed an extremely promising criterion based 
	on the construction of hyperbolic orbits. 
	}
The presentation below owes a great deal
to their excellent paper, \cite{MP:Conv}, and to
\cite{Stark:ex}, which came out of Stark's thesis.

\subsection{definitions and a first criterion}

     We will study maps given by (\ref{eqn:Herm})
and try to find criteria which preclude the existence of the kind of 
tori produced by the KAM theory.  We cannot, of course, rule out the 
existence of tori in the broadest sense.  No matter how large the
perturbation, some tori may  remain in the
islands around elliptic periodic points.  In the two dimensional case
we will restrict our attention to the kind of circles which wind once 
around the cylinder; such circles\thesisfootnote {
	These circles are also called {\em rotational} because the
	restriction of the map to such a circle gives a motion conjugate
	to a rotation.
}
can be smoothly deformed into the curve $p=0$.
In higher dimension we will consider those tori which can be smoothly 
deformed into the torus $p = (0, 0, \ldots, 0)$.
\begin{figure}
	\begin{center}
		\includegraphics[height=4.0cm]{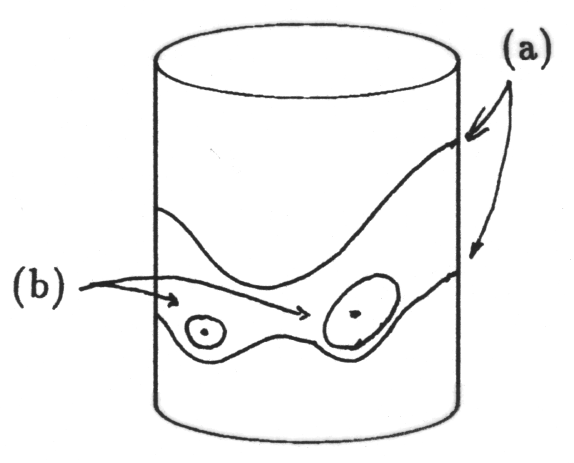}
	\end{center}
	\vspace*{-2\belowdisplayskip}
	\caption{ \em
		The cylinder and several invariant circles, 
		some (a) rotational and some (b) encircling a
		periodic orbit.}
	\label{fig:circles}
\end{figure}

    Maps given by (\ref{eqn:Herm}) are automatically area and orientation 
preserving. We will add the further restrictions that the perturbation,
$f$, be differentiable, periodic, and have average value zero, i.e.
\begin{displaymath}
	f(x) = f(x+1), \qquad \int^1_0 f(x) \; dx = 0.
\end{displaymath}

The restriction on the average value is essential; if it is not met  
$T_{\epsilon}$ has no invariant tori at all.  To see why consider a curve,
$(x,\Gamma_0(x))$, and its image, $(x,\Gamma_1(x))$, where $\Gamma_1$
is given implicitly by 
\begin{displaymath}
	\Gamma_1(x') = p'(x,\Gamma_0(x)),  
\end{displaymath}
or
\begin{equation}
	\Gamma_1(x + \Gamma_0(x)) = \Gamma_0(x) + \epsilon f(x).
\label{eqn:impl}
\end{equation}
Preservation of area and orientation gaurantee that the area
between the two is independent of $\Gamma_0$
since, if we consider another curve, $\Gamma_0'$, and its image,
$\Gamma_1'$, we can write

\begin{displaymath}
    \int_0^1 \Gamma_0' - \Gamma_0 \; = \; 
    \int_0^1 \Gamma_1' - \Gamma_1 
\qquad{\rm so} \qquad
    \int_0^1 \Gamma_0' - \Gamma_1' \; = \;
    \int_0^1 \Gamma_0 - \Gamma_1
\end{displaymath}
and hence we can calculate it for any curve we like. Using 
$\Gamma_0(x) = p_0$ and equation (\ref{eqn:impl}) we get
\begin{displaymath}
    \Gamma_1(x + p_0 ) = p_0 + \epsilon f(x), 
    \qquad{\rm or} \qquad
    \Gamma_1(x) = p_0 + \epsilon f(x - p_0).
\end{displaymath}
Thus we find
\begin{displaymath}
	\Delta \Gamma(x) 
		\; \equiv \;  \Gamma_1(x) - \Gamma_0(x) 
	        \; = \; \epsilon f(x-p_0).
\end{displaymath}
\begin{figure}[t]
	\begin{center}
		\includegraphics[height=4.0cm]{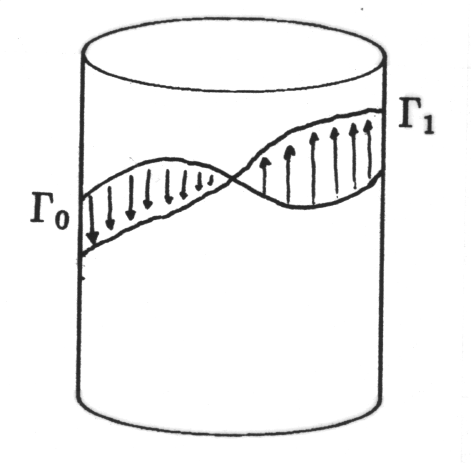}
	\end{center}
	\vspace*{-2\belowdisplayskip}

	\caption{\em A curve and its image. The area between the two is shaded.
	\label{fig:Calabi}}
\end{figure}
The area between the two curves is then 
\begin{eqnarray*}
	\int_0^1 \Delta \Gamma(x) dx & =  & 
	\int_0^1 \epsilon f(x-p_0),
\end{eqnarray*}
the average value of $f$.
Now suppose $\Gamma_0^{inv}$ is an invariant circle. That  means
$\Gamma_1^{inv} = \Gamma_0^{inv}$. Then
\begin{displaymath}
	\int^1_0  \Delta \Gamma(x) dx = 0 
\end{displaymath}
and we have our first and simplest test for the non-existence of
invariant circles.  Unfortunately this is not a very decisive 
criterion; it leaves open the possibility of circles for any value
of $k$ in the Taylor-Chirikov standard map. To do any better we must
more carefully consider the geometry of invariant circles, a
task we turn to next.

\subsection{Lipschitz cone families and their refinement}
\label{sec:cone-fields}
     The first thing to notice is that invariant circles divide the 
cylinder into two disjoint pieces. Orbits which begin below an invariant
circle must always remain below it.  One might hope to turn this
observation into a non-existence criterion, say, by starting an orbit
at some point $(\theta_0, p_0)$ and evolving it forward. If the orbit
eventually attains arbitrarily large momenta then the map has no 
invariant circles.  Chirikov \cite{Chkv:Osc} calls orbits with indefinitely
increasing momentum ``accelerator modes'' and notes that they exist in the 
standard map for $k \geq 2\pi$.

Rigorous implementation of this strategy is hard.
The simple calculation described above does not work because
one can never be sure that a computational error will not carry the orbit
across a genuine invariant circle.  Simply following an orbit cannot establish
the non-existence of circles.  One might instead try to follow an orbit
and say that if it never rises above a certain momentum $p = p_{max}$
then it must be trapped beneath an invariant circle.  That is, one 
might try to prove the {\em existence} of circles.  

\begin{figure}
	\begin{center}
		\includegraphics[height=3.75cm]{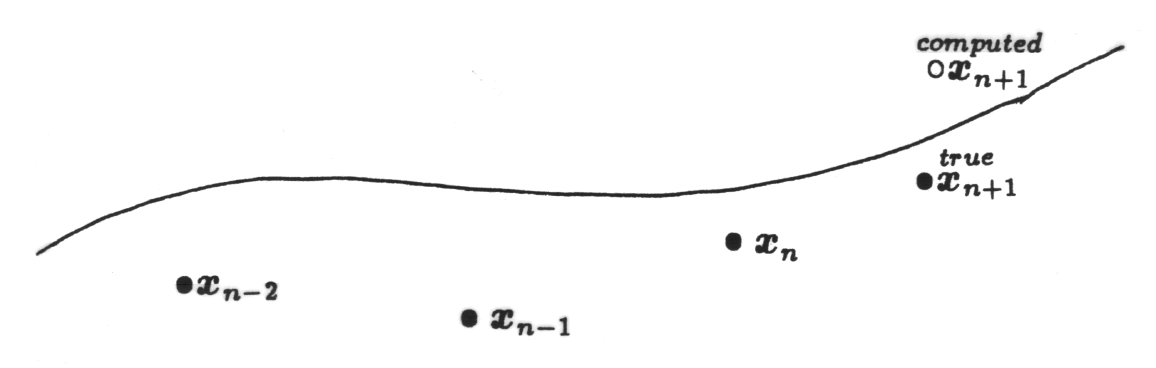}
	\end{center}
	\vspace*{-2\belowdisplayskip}
	\caption{\em 
		Numerical error may carry a point 
		across an invariant circle.}
	\label{fig:error}
\end{figure}

     From an analytic point of view this seems like a good idea. 
A theorem of Birkhoff \cite{Birk:circ} says that if the twist map is
continuously differentiable and if there are two values of the 
momentum, $p_1$ and $p_2, \;\;p_1 < p_2$, such that any orbit which begins 
with momentum less than $p_1$ never attains a momentum greater than $p_2$
then there is an invariant circle
somewhere in the band $p_1 < p < p_2$.
Further, the circle\thesisfootnote{
	\cite{Ma:circ} gives a sketch of the proof of this theorem.
}
is a the graph of some Lipschitz function, 
$\Gamma(\theta)$.  
\begin{figure}[bh]
\parbox[b]{0.33\textwidth}{ \caption{ \em
		If orbits with initial momentum
		less than $p_1$ never rise 
		above $p = p_2$ there is 
		an invariant circle.
		\label{fig:Birk-band}}}
\hfill
\parbox[b]{0.6\textwidth}{
		\centering
		\raisebox{0.2\height}{\includegraphics[height=3.5cm]{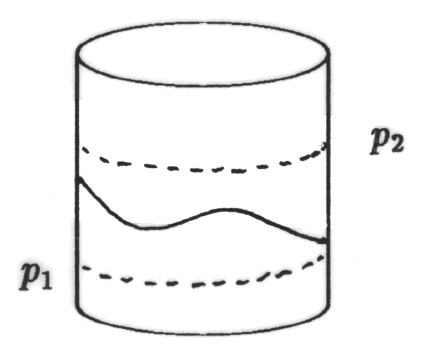}}
		}
\end{figure}

    Despite this analytic support, we cannot get a good existence
criterion either.  Not only is computational error again a problem, but
we must also worry about the cantori.  Although they are not true 
barriers to the diffusion of phase points, they can be formidable
partial barriers\thesisfootnote{
	For the golden cantorus  of the standard map, with 
	$k = {\rm 1.0}$, \cite{MMP:tran} 
	find the mean crossing time to be on the order of 
	${\rm 10}^6$ iterations.
}.
Even if we could calculate an orbit with perfect 
precision we could never be sure that it was permanently trapped below a
particular $p_{max}$.  To get a really useful criterion we must pay 
closer attention to
Birkhoff's theorem, particularly to the part where he tells us that 
rotational invariant circles are the graphs of Lipschitz functions.

\begin{figure}
	\begin{center}
		\includegraphics[height=3.75cm]{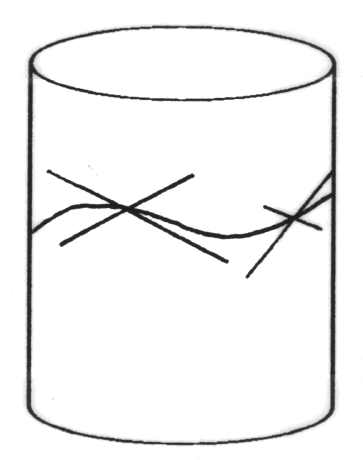}
	\end{center}
	\vspace*{-2\belowdisplayskip}
	\caption{\em An invariant curve and with some Lipschitz cones.}
	\label{fig:cones}
\end{figure}

     Suppose the invariant circle has rotation number $\omega$, then
we will say that it is the graph of $\Gamma_{\omega}(\theta)$.
Since $\Gamma_{\omega}$ is Lipschitz we have
\begin{equation}
	|\Gamma_{\omega}(\theta + {\scriptstyle \Delta}\theta) -
	 \Gamma_{\omega}(\theta)| \leq
	 L \, |{\scriptstyle \Delta}\theta|,
\label{eqn:Lip}
\end{equation}
where $L$ is a constant independent of $\theta$.  On the graph this
means that a vector tangent to the circle
is confined inside a cone, see figure (\ref{fig:cones}).
Since $\Gamma_{\omega}$ is only a Lipschitz function it need
not have a well-defined tangent at every point. That is, although
(\ref{eqn:Lip}) implies that both the right and left limits, 
\begin{eqnarray*}
(\Gamma_{\omega}')_{right} & \equiv & 
    \lim_{\delt{\theta} \searrow 0}
    \frac{ |\Gamma_{\omega}(\theta + {\scriptstyle \Delta}\theta) -
	   \Gamma_{\omega}(\theta)|}{ |{\delt{\theta}}| } \\
(\Gamma_{\omega}')_{left} & \equiv & 
    \lim_{\delt{\theta} \nearrow 0}
    \frac{ |\Gamma_{\omega}(\theta + {\scriptstyle \Delta}\theta) -
	   \Gamma_{\omega}(\theta)|}{ |{\delt{\theta}}| }
\end{eqnarray*}
must exist, they need not be the same.
Nonetheless, both limits must be smaller than $L$, and so 
both the vectors $(1, (\Gamma_{\omega}')_{left})$  and
$(1, (\Gamma_{\omega}')_{right})$ are in the 
cones\thesisfootnote{ 
	Indeed, a	
	Lipschitz function is absolutely continuous and so
	has a derivative defined almost everywhere, 
	see e.g. \cite{Tit:real}. 
	}
pictured in figure (\ref{fig:cones}).

The constant $L$ is a property of 
$\Gamma_{\omega}$ and is defined only along the curve. 
We could, instead, draw a cone at every point, $(\theta, p)$, 
such that if an invariant circle passes through 
$(\theta,p)$ its tangent must lie inside.
We will call such a system of cones a {\em cone family} 
and represent it with two \mbox{$\theta$-periodic} functions, 
$L_+(\theta,p)$ and $L_-(\theta,p)$; a vector 
tangent to a circle through $(\theta, p)$ may only have  
slope, $\ell$, with $L_-(\theta,p) \leq \ell \leq L_+(\theta,p)$. 
The simplest possible cone family is
\begin{equation}
	L_-(\theta,p) = L_{0-},\qquad L_+(\theta,p) = L_{0+}.
\label{eqn:naive}
\end{equation}
We will call this a {\em naive} or {\em uniform} cone family.
We can always get such a family by taking, at the worst, 
$-L_{0-} = L_{0+} = \infty$.  Often, as we shall see, we
can do much better.

\begin{figure}[ht]
	\begin{center}
		\includegraphics[height=4.0cm]{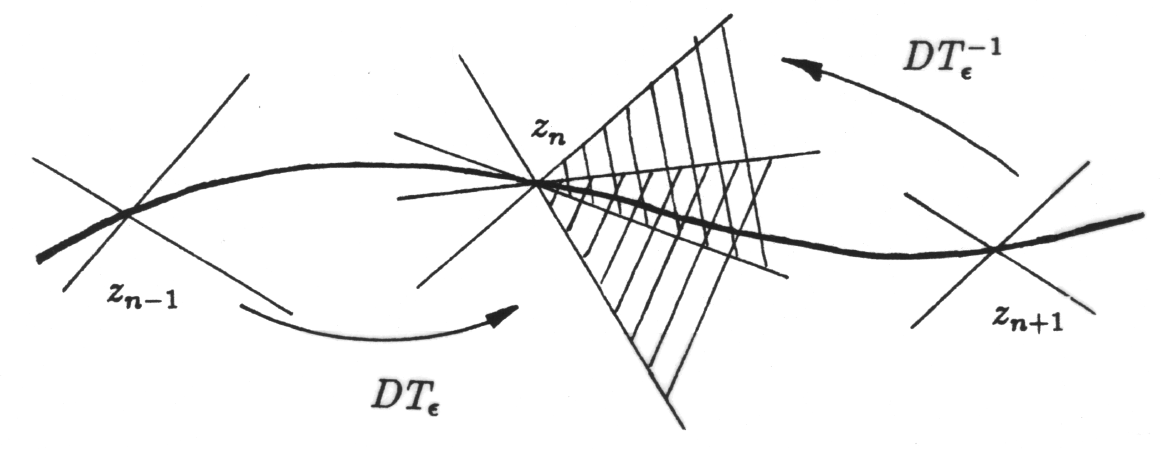}
	\end{center}
	\vspace*{-2\belowdisplayskip}
	\caption{\em Refining the cone family.  The inverse image
			 of the cone at $z_{n+1}$ and the forward image
			 of the cone at $z_{n-1}$ intersect in a new, smaller
			 cone at $z_n$.}
	\label{fig:cone-refinement}
\end{figure}

Each tangent vector lying inside the cone family is ostensibly a
permissible tangent to an invariant curve but the dynamics may
preclude some of the slopes permitted by the naive cone condition.
Consider the action of the map on a tangent vector, say the
vector $\nu$ with footpoint $(\theta,p)$. 
\begin{displaymath}
	\nu' =  DT_{\epsilon, (\theta,p)} \nu
\end{displaymath}
is its image and has footpoint $(\theta',p')$.  We can apply the map
$DT_{\epsilon}$ to all the vectors allowed by the Lipschitz cone at some
point $z_n = (\theta_n, p_n)$ and examine their images at 
$z_{n+1} = (\theta_{n+1}, p_{n+1}) = T_{\epsilon}(z_n) $.
In this way we can use the map on tangent vectors to define a map on
cones.  The image of the cone from $z_n$
will not usually coincide with the cone at $z_{n+1}$.  This means we can 
eliminate part of the cone at $z_n$, for if
there were an invariant graph above $\theta_{n}$ its tangent vector would
have to be one of the ones whose images lie inside the naive cone at 
$z_{n+1}$. We could make a similar argument involving $DT^{-1}_{\epsilon}$
and $z_{n-1}$ and so refine the cone at $z_n$ even further,
see figure (\ref{fig:cone-refinement}).

    More formally, we can use the map to recursively define a sequence 
of cone families, \mbox{  
    $ C_n(\theta,p) \equiv \{L_{n-}(\theta,p), L_{n+}(\theta,p)\} $
} by
\begin{eqnarray}
C_0 & = & \{ L_{0-}, L_{0+} \}  \nonumber \\
	 C_{n+1}(\theta,p) & = & DT^{-1}_{\epsilon} 
		\left\{ C_n(T_{\epsilon}(\theta, p)) \right\}
				\; \cap \;
				C_n(\theta,p)
				\; \cap \; 
				DT^1_{\epsilon} 
				\left\{ C_n(T^{-1}_{\epsilon}
				(\theta,p)) \right\}
\label{eqn:refine}
\end{eqnarray}
where $C_0$ is the naive cone family, (\ref{eqn:naive}).  The 
vectors permitted by the $n$th cone family  have $n$ allowed
images and preimages.  For twist maps this refinement procedure 
produces increasingly restrictive cone families \cite{Stark:ex}.
If it ever happens that $C_n(\theta,p)$ is empty, i.e. that 
the intersection in (\ref{eqn:refine}) contains no vectors, 
then no invariant circle can pass through the point $(\theta,p)$.

\begin{figure}[t]
	\parbox[b]{0.6\textwidth}{
		\includegraphics[height=9.5cm]{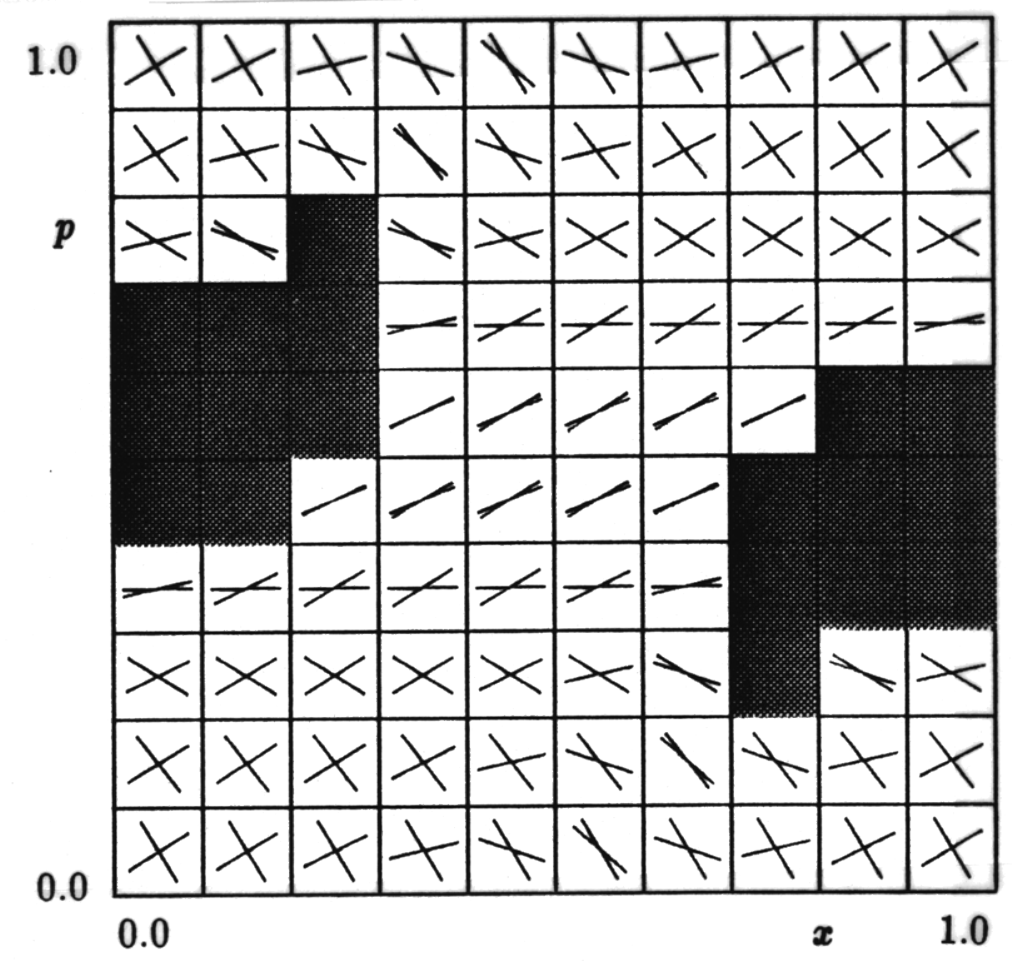}
	}
	\hfill
	\parbox[b]{0.3\textwidth}{ 
		\caption{ \em
			A piecewise constant cone family
			for the standard map with $k=$1.0.
			No invariant circles can pass 
			through the shaded squares.
		}
	}
	\label{fig:cone-boxes}
\end{figure}

     Cone crossing arguments turn out to be quite successful, 
though they need a little more elaboration to be suitable for 
computation.
So far we have seen how to prove that no invariant circle can pass 
through a particular point, now let us use this to prove non-existence 
of circles.  Because a rotational invariant circle must cross
every vertical line, we can establish non-existence by proving 
that no circle can cross a particular
vertical line $\{ (\theta,p)| \theta = \theta_0, p \in [0,1) \} $.
To do that we divide the phase space up into finitely many pieces.
For example, each piece might be a rectangle of the form 
$ R_{ij} =\{(\theta,p)|
\mbox{ $p \in [p_j, p_{j+1}] $ } 
\mbox{ $\theta \in [\theta_j, \theta_{j+1}] $ } \} $ 
We can use this decomposition to construct a sequence of 
piecewise constant cone families, see figure (\ref{fig:cone-boxes}).
\begin{eqnarray}
  C_n(R_{ij}) \equiv  \{ L_{n-}(R_{ij}), L_{n+}(R_{ij}) \} 
    & \qquad & C_0(R_{ij}) = \{ -L, +L \} \nonumber \\
    & \qquad & L_{n -}(R_{ij}) = 
	\lb_{R_{ij}}{L_{n -}(\theta, p)}, \nonumber \\
    & \qquad & L_{n +}(R_{ij}) = 
	\ub_{R_{ij}}{L_{n +}(\theta_0, p)}.
\label{eqn:box-refine}
\end{eqnarray}
where the notations ``u.b.'' and ``l.b.'' mean ``upper bound'' and
``lower bound.''  If the rectangles are small enough,
refinements like (\ref{eqn:box-refine}) can eventually produce a 
whole vertical strip of empty cones.

Finally, we note that the foregoing serves to prove non-existence for 
a single map.  In practice one wants non-existence results for a
whole class of maps, for example, for all the standard maps with
parameters $k_{min} \leq k \leq k_{max}$.  One need only modify
(\ref{eqn:box-refine}) a little, taking the bounds over both 
$R_{ij}$ and $k$.

     Stark has shown that such a program, allied with some
extra observations, can reveal non-existence of circles with only 
a finite amount of work.  He shows, for example, that if one
has a family of maps depending on parameters and one studies a
compact set of the parameters for which no invariant circles exist, 
then the cone-crossing criterion will demonstrate their non-existence
after only a finite amount of computation\thesisfootnote{
    Here ``finite'' means that one could do the calculations 
    to some finite precision and refine the cone families for some
    finite number of steps.  
}.

\subsection{some new coordinates and two more criteria}
\label{subsec:two-crit}

Here we will begin to explain one way to implement the ideas of 
the previous section on a digital computer.  In the process we 
will reformulate the cone-crossing
criterion in a way that obscures its geometric origin\thesisfootnote{
	See \cite{MP:Conv} for a more direct implementation.
}
but reveals a connection to minimizing states.  The first step is to
recast the map in terms of delay coordinates;
we have been considering 
\mbox{$ T_{\epsilon}(\theta,p) = (\theta', p')$}, let us now speak of 
\mbox{$g_{\epsilon} : \T \times \T \mapsto \T \times \T$} so that 
$g_{\epsilon}(\theta_n, \theta_{n+1}) = (\theta_{n+1}, \theta_{n+2})
$ where the $\theta's$
are angular coordinates of successive points in an orbit.  
We will also need a lift of $g$, 
\mbox{$G_{\epsilon} : \R \times \R \rightarrow \R \times \R$, 
$G_{\epsilon}(u,v) = (u',v')$}.
As before, $T_{\epsilon}$ and $G_{\epsilon}$ are related by 
an action generating function, 
$H_{\epsilon}(u,v)$, where
\begin{eqnarray*}
  H_{\epsilon}(x_n, x_{n+1}) = \half (x_{n+1} - x_n)^2 - \epsilon V(x_{n+1}),
    & \qquad & V(x) = -\int_0^x f(y) \; dy, \\
    & \qquad & \partial_1 H_{\epsilon} (x_n, x_{n+1}) = -p_n, \\ 
    & \qquad & \partial_2 H_{\epsilon} (x_n, x_{n+1}) = p_{n+1},
\end{eqnarray*}
and
\begin{eqnarray*}
    G_{\epsilon}( x_{n-1}, x_n ) & \equiv & (x_n, x_{n+1}), \\
    x_{n+1} & = & x'(x_n, p_n),  \\
	    & = & x'(x_n, \; \partial_2 H_{\epsilon}(x_{n-1}, x_n)). 
\end{eqnarray*}

In terms of these coordinates an invariant circle appears as a 
curve $x_{n+1} = \gamma(x_n)$ satisfying 
\begin{eqnarray*}
	\gamma(u + 1) & = & \gamma(u) + 1, \\ 
	G_{\epsilon}(x_n, \gamma(x_n)) & = & ( x_{n+1}, x_{n+2}) 
			      = (\gamma(x_n), \gamma(\gamma(x_n))).
\end{eqnarray*}
\begin{figure}
\parbox[b]{0.6\textwidth}{
	\centering
	\includegraphics[height=6.25cm]{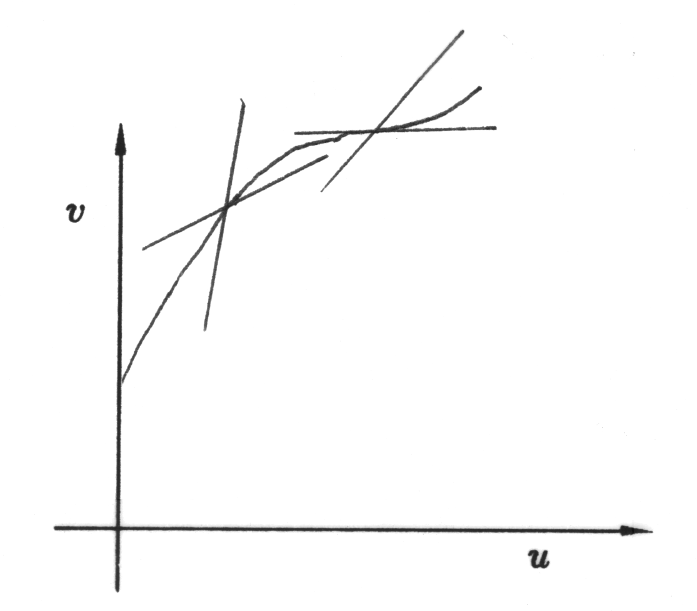}
}
\hfill
\parbox[b]{0.33\textwidth}{
	\caption{\em An invariant curve and 
		some Lipschitz cones in the 
	    delay coordinate system.}
	\label{fig:delay-cone}
}
\end{figure}
The most naive Lipschitz cone, (\ref{eqn:naive}) with 
$L_{0 \pm} = \pm \infty$,
appears here as  \mbox{ $ 0 \leq \ell \leq \infty $ }
where $\ell$ is the slope of $\gamma$.  The lower bound of zero is just the 
requirement that the original map, when restriced to an invariant curve,
be order preserving.  

For examples like (\ref{eqn:Herm}) $u'$ and $v'$ have very simple forms:
\begin{eqnarray}
	u'(u,v) & = & v, \nonumber \\
	v'(u,v) & = & v + (v - u) + \epsilon f(v), \nonumber \\
		& = & 2v -u + \epsilon f(v).
\label{eqn:Gdef}
\end{eqnarray}
$G_{\epsilon}$'s action on tangent vectors is equally simple:
\begin{equation}
	\left[ \begin{array}{c}
			\delta u' \\
			\delta v' 
		\end{array} 
	\right] =  
	\left[ \begin{array}{cc}
			0 & 1 \\
			-1 & 2 - \epsilon \frac{d^2V}{dx^2}
	\end{array} \right]
	\left[ \begin{array}{c}
			\delta u\\
			\delta v
	\end{array} \right].
\label{eqn:DGact}
\end{equation}
For later convenience we will refer to 
\mbox{$2 - \epsilon \frac{d^2 V}{dx^2}(x)$} as $\beta(x)$.

If we take a tangent vector, $[1,\ell ]$, representing a slope of $\ell $
then (\ref{eqn:DGact}) tells us that its image will represent a 
slope $\ell'$ given by:
\begin{eqnarray}
     \ell' & = & \frac{\delta v'}{\delta u'}, 	\nonumber \\
	   & = & \frac{\beta(v) \delta v}{\delta v} -
		 \frac{\delta u}{\delta v},	\nonumber \\
	   & = & \beta(v) - \frac{1}{\ell }.
\label{eqn:slope-evolution}
\end{eqnarray}
Preservation of order   requires both $\ell$ and $\ell'$ be positive.
Combining that with (\ref{eqn:slope-evolution}) we obtain our first real 
criterion.

\noindent
{\bf Criterion 1}  \hskip 1\parindent
{\em If there are any values $v \in [0,1]$ for which $\beta(v) < 0$
then the map $G_{\epsilon}(u,v)$ to which $\beta$ corresponds has no rotational 
invariant circles.  For the standard map this criterion says $k_c \leq 2$. }

We can squeeze one further analytic criterion out of (\ref{eqn:slope-evolution})
by noticing that $\ell'$ will surely be negative if ever $\ell$ is 
very small, and that, always, $\ell' < \max_{v \in [0,1]}{\beta(v)}$. 
Suppose we have $m$ and $M$ such that \mbox{$0 \leq m \leq \beta(v) \leq M$} 
holds everywhere. Then 
\begin{equation}
	\ell' \leq M - \frac{1}{\ell}
\label{eqn:softDG}
\end{equation}
and $\ell' \geq 0$ together imply 
\begin{equation}
	0 \leq M - \frac{1}{\ell} 
	\qquad {\rm or} \qquad 
	\ell \geq \frac{1}{M}.
\label{eqn:global-cone}
\end{equation}

Inequality (\ref{eqn:global-cone}) is a global restriction on slopes,
a new lower bound for the uniform Lipschitz cone family.  We could thus run
through the argument again, this time requiring $\ell' \geq \frac{1}{M}$.
Having done that we would have a better, narrower cone family and could repeat
the argument yet again $\ldots$ better to carry this process straight to its
conclusion and realize that our estimates will stop improving when we find
a slope, $\ell_-$, such that 
\begin{displaymath}
	\ell_- = M - \frac{1}{\ell_-}.
\end{displaymath}
This has two roots. The least of them is just the $\ell_-$ we wanted; the
larger one is a global upper bound on slopes.  It comes from the remark 
above, that $\ell' \leq M$.  Since every vector tangent to an invariant
curve is the image of some other tangent we can conclude $\ell \leq M$.
Once that's done we can argue $\ell' \leq M - \frac{1}{M}$ and so on.
Finally we attain 
\begin{eqnarray}
	\ell_- \leq \ell \leq \ell_+  
	    & \qquad {\rm where} \qquad &
	    \ell_- = \frac{M - \sqrt{M^2 -4}}{2}, \nonumber \\
	&&  \ell_+ = \frac{M + \sqrt{M^2 -4}}{2}.
\label{eqn:globalLip}
\end{eqnarray}
Armed with this best of all possible uniform cones, we are able to make a 
genuine, dynamical cone crossing argument.

\noindent
{\bf Criterion 2} (``Mather $\frac{\scriptstyle 4}{\scriptstyle 3}$'' )
\hskip 1\parindent
{\em  If $m \leq \beta(v) \leq M$ and $\ell_+$ and $\ell_-$ are the 
bounds of the uniform cone family given by 
\mbox{\rm (\ref{eqn:globalLip})}, then
there are no rotational circles if }
\begin{equation}
	\ell_- > m - \frac{1}{\ell_+}.
\label{crit:four-thirds}
\end{equation}

\noindent
{\bf Remark} \hskip 1\parindent 
{\em For the standard map, $m = (2 - k)$ and $M = (2 + k)$
and so {\rm (\ref{crit:four-thirds}) } implies that 
$k_c \leq \frac{\scriptstyle 4}{\scriptstyle 3}$.
}

\noindent
{\bf Proof}  \hskip 1\parindent 
The idea is to concentrate on those states which 
contain the point where $\beta$ attains its minimum,
where $\beta(v) = m$.
Visits to this point are most punishing to the slopes of
tangent vectors; they lead to the smallest possible values of 
$\ell'$ in (\ref{eqn:slope-evolution}).  If $m$ is so small
that even the slope from the upper edge of the uniform family,
$\ell_+$, is diminished to an untenable value, then 
certainly no others can survive.

\subsection{non-existence for minimalists}
\label{sec:min-crit}

We will now reformulate Criterion 2 in the language of minimizing states.  The
new version will prove more fruitful for higher dimensional generalizations.
Here again we follow MacKay and Percival, who demonstrated that their cone 
crossing criterion is equivalent to the action-difference 
criterion put forward by Mather in \cite{Ma:crit}.

We begin by assuming that an invariant circle exists, then we 
deduce some facts about the minimizing orbits lying on it.  Then,
to prove non-existence, we will do a calculation that contradicts
these facts.  Define a {\em minimizing state} to be  
sequence $ \{ \cdots x_{n-1}, x_{n}, x_{n+1}, \cdots \} $
such that every finite segment $x_n, x_{n+1}, \cdots , x_m$ is a 
minimum of the action functional,
\begin{equation}
	W_{m,n}(X) = \sum_{j=m}^{n-1} H_{\epsilon}(x_j, x_{j+1}),
\label{eqn:action-func}
\end{equation}
where $H_{\epsilon}$ is the action generating function 
and we consider variations
which leave $x_n$ and $x_m$ fixed.  Mather's action-difference 
idea is to note that if an irrational invariant circle exists then every 
orbit on it is minimizing and has the same action. That is, if we take
two states arising from orbits on the circle,
$X^a = \{ \cdots, x_0^a, x_1^a, \cdots \}$ and
$X^b = \{ \cdots, x_0^b, x_1^b, \cdots \}$ and take the limit
\begin{equation}
	\lim_{n \rightarrow \infty}
	     {\sum_{j=-n}^{n-1} H_{\epsilon}( x^a_j, x^a_{j+1} ) 
				- H_{\epsilon}( x^b_j, x^b_{j+1})}
\label{eqn:actDiffA}
\end{equation}
it should come out to be zero\thesisfootnote{
	Showing that the action difference (\ref{eqn:actDiffA}) vanishes is
	different,
	and harder, than showing that the {\em average} values
	of the actions are the same.  While the latter follows from the 
	ergodicity of irrational rotation, Mather's
	result requires a more delicate examination of the action
	functional. See \cite{Ma:crit} for details.
}.
He suggests that to test the existence of an invariant circle having 
irrational rotation number $\omega$ one should approximate $\omega$
by a sequence of rational numbers, $\frac{p_n}{q_n}$, and
use the rational numbers to construct the two sequences of Birkhoff periodic 
orbits, the minimax and minimizing orbits.  These sequences
accumulate on two distinct sets on the putative invariant circle. If the
circle is really present, orbits on the two sets should have the same action
and so the limit
\begin{equation}
	\Delta W_{\omega} \equiv
	\lim_{\frac{p_n}{q_n} \rightarrow \omega}
	{ \Delta W_{p_n,q_n}} =
	W_{(p_n,q_n) \; minimax} - W_{(p_n,q_n) \; minimizing}
\label{eqn:actDiffB}
\end{equation}
should tend to zero. If it tends to some other value then no circle
with rotation number $\omega$ exists.

Another way to state this argument is to say that  every orbit
on a rotational invariant circle must have the same action, the action 
corresponding to the limit of the minimizing Birkhoff orbits.  Thus
every  state $X = \{\cdots, x_{-1}, x_0, x_1, \cdots \}$ arising from an
orbit $\{\cdots, (x_{-1},p_{-1}), (x_0,p_0), (x_1, p_1), \cdots \}$ lying
in an invariant circle must be minimizing; every finite segment
snipped out of such a state must be a non-degenerate minimum over all
segments having the same endpoints\thesisfootnote{
	The reader may wonder why the states lying on an invariant
	circle do not belong to a one parameter family, and ask how they
	can lead to non-degenerate minima.  The answer is that we 
	consider only variations which leave the endpoints of finite 
	segments fixed; if we allowed them to move the minima would 
	be degenerate.
}.

     The foregoing suggests a strategy for proving converse KAM theorems.
One chooses an auspicious starting point, $x_0$, for which the 
perturbation
to the generating function is large, and considers every possible state
containing it.  This is not quite so huge a task as it sounds.
Since the map, $G_{\epsilon}(u,v)$, determines the
whole state once, say, $x_0$ and $x_1$ have been given, we need
only consider all possible successors, $x_1$.  For each $x_1$ we 
work out the state, $X$, and the variation of the action
over finite segments, $\{x_{-1}, x_0, \cdots, x_n\}$, 
\begin{displaymath}
    \begin{array}{ccccc}
    \delta W_{-1,n} & = & { \displaystyle \sum_{j = 1}^{n-1} 
			\frac{\partial W_{-1,n}}{\partial x_j } 
			\delta x_j }
		  & + & \half 
			{\bf \delta x}^T{\bf D}^2W_{-1,n} {\bf \delta x} \\
		  & = &  0 & + & \half 
	{ \displaystyle	  \sum_{ j,k = 1}^{ n-1} 
			\frac{ \partial^2 W_{-1,n}}
			     { \partial x_j \partial x_k } 
			     \delta x_j \delta x_k }.
    \end{array}
\end{displaymath}
The term linear in $\delta x_j$ is automatically zero because 
$X$ is a minimizing state. For our examples,
(\ref{eqn:Herm}), the quadratic term can be represented by the
symmetric matrix,
\begin{displaymath}
{\bf D}^2W_{-1,n} = 
\left[ 
  \begin{array}{cccccc}
    2 + \epsilon \frac{\displaystyle df}{\displaystyle dx}(x_0) & 
    -1  &
    \;0  & \cdots & \cdots &  \; 0  \\
    -1 &
    2 + \epsilon \frac{\displaystyle df}{\displaystyle dx}(x_1)&
    -1 &
    \cdots & \cdots &  \;0 \\
    \;\:\vdots & & & \ddots &        & \;\:\vdots \\
    \;\:\vdots & & &        & \ddots & \;\:\vdots \\
    \;0 & \cdots & \cdots & 
    -1 & 2 + \epsilon \frac{\displaystyle df}
			   {\displaystyle dx}(x_{n-2}) & -1 \\
    \; 0 & \cdots & \cdots & \cdots &
    -1 & 2 + \epsilon \frac{\displaystyle df}
			   {\displaystyle dx}(x_{n-1}) 
  \end{array}
\right],
\end{displaymath}
which we shall call $M_n(X)$, or $M_n$ for short.

If $X$ is minimizing then $M_n$ is positive definite.  
Since $M_n$ is so simple it is easily rendered into diagonal form,
a form which makes it simple to calculate the determinant. We can write 
\begin{eqnarray*}
\left[ 
  \begin{array}{ccccc}
    2 + \epsilon \frac{df}{dx}(x_0) & -1  &  0 & 0 & \cdots \\
    -1 & 2 + \epsilon \frac{df}{dx}(x_1) & -1  & 0 & \cdots \\
    0 & -1 & 2 + \epsilon \frac{df}{dx}(x_2) & -1  & \cdots \\
    \vdots & \vdots & \vdots & \vdots & 
  \end{array}
\right]  
&  \rightarrow &
\left[ 
  \begin{array}{ccccc}
    d_0 & 0 & 0 & 0 & \cdots \\
    0 & d_1 & 0 & 0 & \cdots \\
    0 & 0 & d_2 & 0 & \cdots \\
    \vdots & \vdots & \vdots & \vdots & 
  \end{array}
\right]  
\end{eqnarray*}
where the $d_j$ are computed recursively using
\begin{eqnarray}
  d_0 & =  &2 + \epsilon \frac{df}{dx}(x_0),  \nonumber \\
  d_{j+1} & =  & \beta(x_{j+1}) - \frac{1}{d_j}, 
	\;\; {\rm where} \;\;
      \beta(x_{j+1}) = 2 + \epsilon \frac{df}{dx}(x_{j+1}).
\label{eqn:diag-evolution}
\end{eqnarray}

If ever one of the $d_j$ is negative we may conclude that $M_j$
is not positive definite and so does not arise from a minimizing
state.  Notice the similarity between
the evolution equation for the diagonal entries, 
(\ref{eqn:diag-evolution}), and the one for slopes,
(\ref{eqn:slope-evolution}). As we refined the limits on slopes,
so we can refine those on diagonal entries.
We obtain a $d_-$ such that if $d_j < d_-$ then some later
$d_k, \; k > j$ is sure to be negative. We also get $d_+$,
a global upper bound on the $d_j$.  We can thus modify 
(\ref{eqn:diag-evolution}) so that we begin with 
$d_{-1} = d_{+}$, so  $d_0 = \beta(x_0) - \frac{1}{d_{+}}$.
The original prescription corresponds to taking $d_{-1} = \infty$.

%% file: converse/prism.tex
\section{Rigorous Computing}
\label{sec:machines}

     In this section we will see how to implement the action criterion of
the last section on a digital computer.  Since we will eventually want to
treat maps in spaces of arbitrary dimension we will outline some of the
procedures in  greater generality than required for the cylinder. The most 
important part will be a technique for rigorously bounding the image of
a set.

\subsection{two reductions and a plan}
\label{subsec:plan}

As in section (\ref{sec:cone-fields}), we need only show that 
no invariant circle crosses a particular vertical line.  
In the language of the previous section this means our problem is
reduced  to showing that some particular $x_0$ cannot appear as a 
member of any minimizing state.  We can get a further reduction by
noticing that our examples satisfy
\begin{displaymath}
    p'(\theta, p+1) = p'(\theta,p) + 1; 
\end{displaymath}
their dynamical structure is periodic in $p$ as well as in $\theta$.
So, if an invariant circle passes through the point 
$(\theta,p)$, there is also one through $(\theta,p+1)$; if 
no invariant circles pass through some vertical segment
$\mbox{$I_0 \equiv \{(\theta,p) |$} 
	\theta = \theta^{\star}, \; p \in [0,1]\}$, then there
cannot be any at all.  
\begin{figure}[th]
	\begin{center}
		\includegraphics[height=3.75cm]{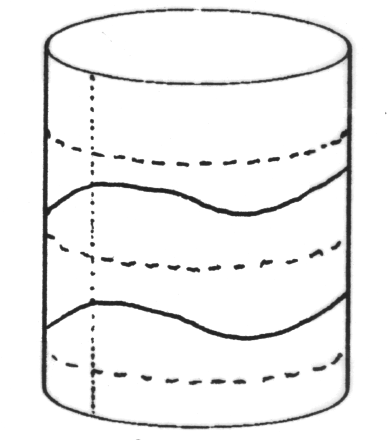}
	\end{center}
	\vspace*{-2\belowdisplayskip}
    \caption{ {\em Rotational invariant circles must cross every vertical
		line, and, for our examples, must be periodic in $p$ as
		well as $\theta$. }
		\label{fig:circle-reductions} }
\end{figure}
Studying a segment like $I_0$ is equivalent to studying a collection of 
states $\{X |\;  x_0 = x^{\star}, \; \mbox{$x_1 \in [0,1]$} \}$, where 
$x^{\star}$ is a lift of $\theta^{\star}$.
     With these reductions in hand, we are ready to plan the main
computation.  Our goal will be to prove:

\noindent
{\bf Theorem}  \newline
{ \em There is an $x^{\star} \in [0,1]$
and an interval of parameter values, 
\mbox{$I_{\epsilon} \equiv [\epsilon_{-},\epsilon_{+}]$},
such that none of the maps, $G_{\epsilon},\; \epsilon \in I_{\epsilon}$,
have a minimizing state with $x_0 = x^{\star}$. }

\noindent
{\bf Plan for the proof:}

\begin{description}
	\item[(i)]{ Formally extend the phase space to include
		    the parameter $\epsilon$ and use the 
		    map $G_{\epsilon}(u,v)$ to define 
		    a new one, $G: \R \times \R \times \R \rightarrow
		    \R \times \R \times \R$, where 
		    \begin{eqnarray}
			 G(\epsilon, u, v) = 
				(\epsilon, G_{\epsilon}(u,v)).
			\label{eqn:G-extension}
		    \end{eqnarray}
		  }
	\item[(ii)]{ Select a starting point $x^{\star}$.  For
		     examples (\ref{eqn:Herm}) we will want 
		     $x^{\star}$ such that $\beta(x^{\star})$
		     is a minimum, a choice which is independent of 
		     $\epsilon$.}
	\item[(iii)]{ Divide the interval [0,1] into a collection 
		      of closed intervals, $I_j$,
		      \mbox{$
				\displaystyle
				\union_{j=1}^{N} I_j = {\rm [0,1]}
			    $}.
		      Using the $I_j$, which 
		      should intersect only at their endpoints,
		      we can construct a collection 
		      of sets in the extended phase space,
		      $S_j \equiv \{(\epsilon,u,v)|
		      \mbox{ $\epsilon \in I_{\epsilon}, \;$}
		      \mbox{$u=x^{\star}$}, 
		      \mbox{$\; v \in I_j\}$}$.
		      In practice this division is done by the 
		      program itself.It begins by trying to 
		      prove the theorem on the whole interval
		      at once, and gets 
		      either, ``Yes, the theorem is true,'' or 
		      ``Maybe it's true.'' 
		      If the answer is ``maybe'' it splits
		      the interval in half and tries the two pieces 
		      separately.  If one of them
		      yields ``maybe'' it gets subdivided  too \ldots.
		      The process of subdivision will go on forever if the
		      theorem is false, but if it is true the work
		      of Stark suggests that the cutting will stop after 
		      finitely many steps. }
	\item[(iv)]{ For each piece $I_j$, try to prove that no minimizing
		     state with $x_0 = x^{\star}$  can have $x_1 \in I_j$.}
\end{description}

    The last step is where the computation comes in; we will use
an argument
like the one at the end of section (\ref{sec:min-crit}), but here we 
calculate upper bounds\thesisfootnote{
	We will often want to evaluate upper bounds, as opposed to
	maxima.  The former are realizeable on computers, the latter
	may not be. }
$\bar{d}_k$ for the {\em k\/}th diagonal entry in 
(\ref{eqn:diag-evolution}).
\begin{eqnarray}
   \bar{d}_0 & = & \ub_{\epsilon \in I_{\epsilon}} \beta(x^{\star})
			- \frac{1}{d_+},       \nonumber \\
   \bar{d}_1 & = & \ub_{(\epsilon, u,v) \in S_j} \beta(v)
			- \frac{1}{\bar{d}_0}, \nonumber  \\
   \bar{d}_2 & = & \ub_{(\epsilon,u,v) \in G(S_j)} \beta(v) 
			- \frac{1}{\bar{d}_1}, \nonumber  \\
	     & \vdots & \nonumber \\
    \bar{d}_{n+1} & = & \ub_{(\epsilon,u,v) \in G^n(S_j)} \beta(v) 
	    		- \frac{1}{\bar{d}_n}.   
\label{eqn:ub-evolution}
\end{eqnarray}
Finding a way to calculate the kind of bound which appears in 
the definition of $\bar{d}_2$, an upper bound over an image of $S_j$,
is the last hurdle in the argument.  What we need is a procedure to
rigorously bound the image of a set. In the next section we will 
explain a quite general scheme due to MacKay and Percival.  

\subsection{bounding images of prisms}
\label{sec:prism-bondage}

\begin{figure}[t]
\raisebox{-0.45\height}{\includegraphics[height=5.0cm]{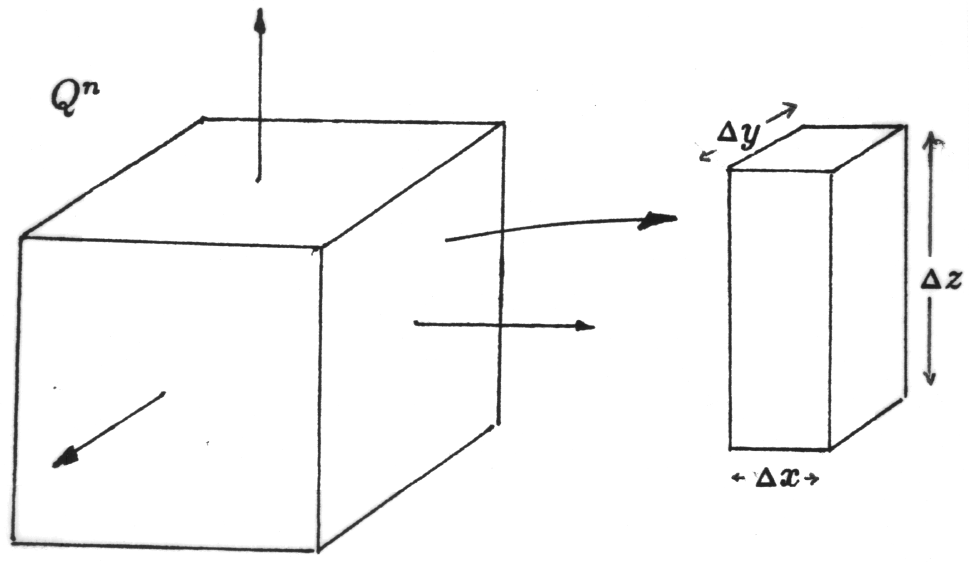}}
\hfill $\displaystyle \equiv 
		(x_c, P) , \; 
		P = \left[ \begin{array}{ccc} 
					\frac{{\scriptstyle \Delta}x}{2}
					& 0 & 0 \\
					0 &
					\frac{{\scriptstyle \Delta}y}{2}
					& 0 \\
					0 & 0 &
					\frac{{\scriptstyle \Delta}z}{2}
			    \end{array}
		     \right]$
\vskip 0.5in
\caption{ \em The {\em n}-dimensional hypercube $Q^n$ is mapped
	  to the prism by the matrix $P$. 
	  \label{fig:prism-pix} }
\end{figure}

For concreteness, and to get an algorithm straightforward
enough to be realized as a computer program, we will concentrate
on sets with a prescribed form,
parallelepipeds, or {\em prisms} for short. 
An {\em n}-dimensional prism is specified by a center point,
$x_c$, and an ${n \times n}$ matrix, $P$.
The prism is the set 
\begin{equation}
	\{ x \in \Rn | x = x_c + P\eta, \; \eta \in Q^n \}, 
\label{eqn:prism-def}
\end{equation}
where $Q^n$ is the {\em n}-dimensional hypercube, 
$\{ \eta \in \Rn | -1 \leq \eta_j \leq 1 \}$, 
see figure (\ref{fig:prism-pix}).  
Our principal technical tool is the following result.

\noindent
{\bf Lemma} (\cite{MP:Conv})\hskip 1\parindent {\em
Suppose $\Phi : \Rn \rightarrow \Rn$ is a ${\rm C}^1$ map.  Then the
$\Phi$ - image of the prism $S \equiv (x_c, P)$ is contained in the prism
$({x_c}', P')$ where ${x_c}'$ is arbitrary, $P' = A \circ  W$ for
an arbitrary invertible matrix $A$, and $W$ the diagonal matrix
\begin{displaymath}
W = \left[ \begin{array}{cccc}
		w_1 & 0 & \cdots & 0 \\
		0 & w_2 & \cdots & 0 \\
		\vdots && \ddots & \vdots \\
		0 & 0 & \cdots & w_n 
	   \end{array}
     \right]
\end{displaymath}
with
\begin{equation}
w_j = \ub \left(
		  |(\Phi(x_c) - {x_c}')_j| + \ub_{x \in S} 
		      \sum_{k=1}^{n} \left|
					[A^{-1} \circ 
					 D\Phi_x \circ
					 P]_{jk}
				     \right|
	  \right).
\label{eqn:w-def}
\end{equation}
}
\begin{figure}
	\begin{center}
		\includegraphics[height=3.75cm]{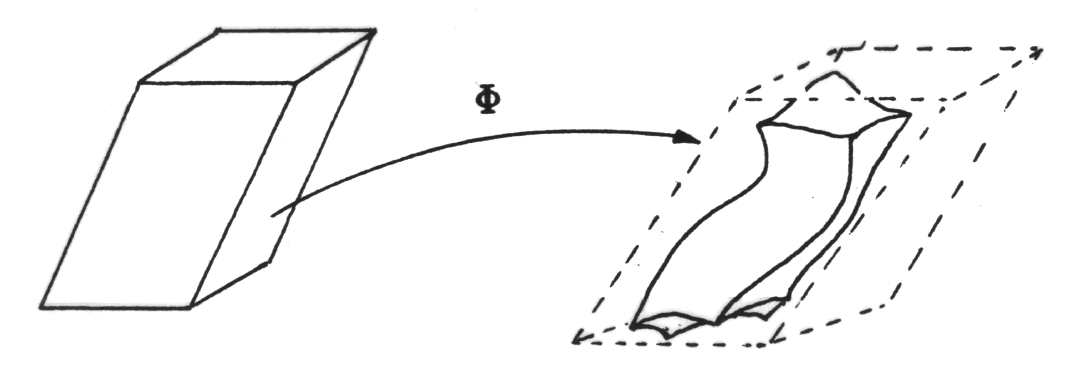}
	\end{center}
	\vspace*{-2\belowdisplayskip}
    \caption{{\em 
	A prism, its image, and a prism bounding the image.}
	\label{fig:prismatic-image}
    }
\end{figure}

\noindent
{\bf Remark} \hskip 1\parindent  {\em The lemma seems unnecessarily 
general; we are left to choose the matrix $A$ and the new center point,
$x_c$ completely arbitrarily.  If we choose them unwisely the new prism 
will surround the image of $S$, but may be much larger than necessary.
Usually we will want
\begin{displaymath}
	{x_c}' \approx \Phi(x_c), 
	\qquad {\rm and} \qquad
	A \approx D\Phi_{x_c} \circ  P.
\end{displaymath}
The freedom allowed by the lemma will make it easy to
handle errors in computing $\Phi(x_c)$ and cases where 
$D\Phi_{x_c}  P$ is singular or nearly singular. }

\begin{figure}
\parbox[b]{0.55\textwidth}{
	\centering
	\includegraphics[height=7cm]{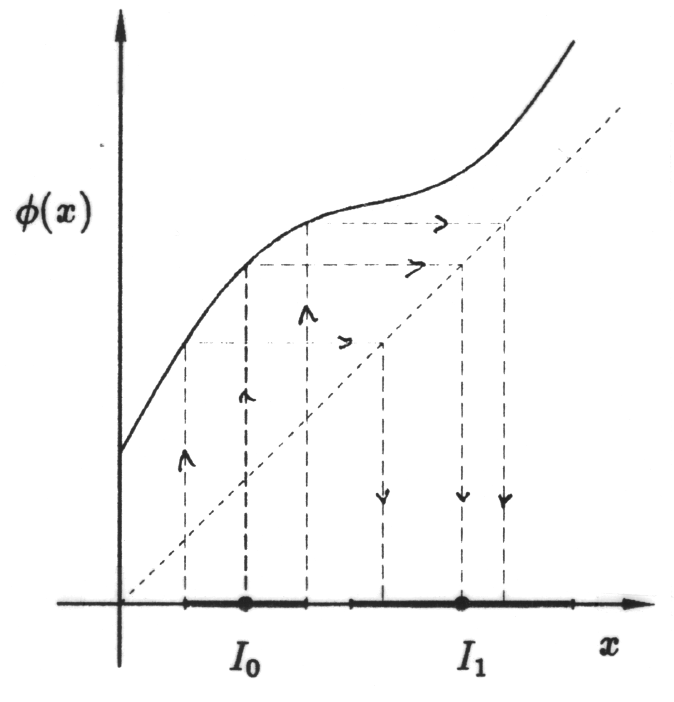}
}
\hfill
\parbox[b]{0.40\textwidth}{
	\caption{ {\em The bounding lemma applied to a lift of 
		the circle map,
		$\phi(x) = x + \Omega + \frac{\epsilon}{2 \pi} 
			\sin{(2 \pi x)}$,
		with \mbox{$\Omega = {\rm 0.3}$},
		\mbox{$\epsilon = {\rm 0.8}$}.  The interval
		$I_1$, at right, is the one given by the lemma;
		it contains the image of $I_0$.
	} \label{fig:bounding} }
}
\end{figure}

\noindent
{\bf Example} (Proof of the lemma for one dimensional maps) \newline
We start in with a one dimensional example, see figure
(\ref{fig:bounding}).  Here the map is some  $C^1$ function, 
$\phi : \R \rightarrow \R$,
and a prism, $S$, is just an interval 
$ x_c - {\scriptstyle \Delta}x 
\leq x \leq x_c + {\scriptstyle \Delta}x $.
We can use the computer to find $\bar{\phi}(x)$, a numerical 
approximation to $\phi(x)$ for which
$|\phi(x) - \bar{\phi}(x)| \leq \delta$. 
Then, choosing  ${x_c}' = \bar{\phi}(x_c)$
and\thesisfootnote{
	The choice of $A$ is meant to suggest the form
	required by the higher dimensional theorem.
	If $\phi'(x_c) = 0$ we will have to make 
	another choice; any constant will do.}
$A = \phi'(x_c) \delt{x}$, we find
\begin{eqnarray*}
    \ub  |x^{\prime}_c - \phi(x_c)| &\leq & \delta, \\
	A^{-1} & = & \frac{1}{\phi'(x_c) \delt{x}}, \\
	W & = & \frac{\delta}{|\phi'(x_c) \delt{x}|} +
		 \ub_{x \in S} 
		    \left|
		       \frac{\phi'(x) \delt{x}}{\phi'(x_c) \delt{x}}
		    \right|, \\
	  & = & \frac{\delta}{|\phi'(x_c) \delt{x}|} +
		\ub_{x \in S} \left| 
				\frac{\phi'(x)}{\phi'(x_c)} 
			      \right|, 
\end{eqnarray*}
and
\begin{equation}
P' \equiv \delt{x'} = A \circ W \geq \delta + 
		\delt{x} (\max_{x \in S}| \phi'(x)|).
\end{equation}
Now let us check some point $x \in S$, and see that its image
is inside the prism $S' = (x'_c, P')$.  Since $x$ is in $S$ 
we can write
\mbox{$x = x_c + \eta \, \delt{x}$} with
\mbox{$\rm -1 \leq \eta \leq 1$}.
If $\phi(x)$ is in $S'$, then, 
\begin{displaymath}
  x'_c - \delt{x'} \leq \phi(x) \leq x'_c  + \delt{x'}
  \qquad {\rm or} \qquad
  |\phi(x) - x'_c| \leq  \delt{x'}.
\end{displaymath}
To see that this is true, consider 
$\gamma(t) = \phi(x_c +t \eta \, \delt{x} )$. $\gamma(t)$
is a $\rm C^1$  function from [0,1] to $\R$ with $\gamma(0) = \phi(x_c),
\; \gamma(1) = \phi(x)$.  By the Mean Value Theorem there is a 
$t_0 \in {\rm [0,1]}$ such that 
\begin{eqnarray*}
    \gamma(1) - \gamma(0)
    & = &
    \frac{d \gamma}{dt}(t_0),  \\
    \phi(x) - \phi(x_c)
    & = &
    \frac{d}{dt}(\phi(x_c + t_0 \eta \, \delt{x})), \\
    & = & 
    \eta \, \delt{x}\,\phi'(x_c + t_0 \eta \, \delt{x}).
\end{eqnarray*}
Rewriting this,
\begin{eqnarray}
    | \phi(x) - x'_c  |
    & =  &
    | \phi(x_c) - x'_c + 
	\eta \, \delt{x} \phi'( x_c + t_0 \eta \, \delt{x}) | ,
						\nonumber \\
    & \leq &
    | \phi(x_c) - x'_c | + 
	|\delt{x} \phi'(x_c + t_0 \eta \, \delt{x})|, \nonumber \\
    & \leq & \delt{x'}, \label{eqn:prism-bd1}  
\end{eqnarray}
even as the lemma claimed. 

\noindent
{\bf Proof}  (The general case) \newline 
The argument is much the same as the 
1-dimensional argument above.  Here the assertion of the theorem is
that every point in the initial prism, $S = (x_c,P)$, has its image
in $S' = (x'_c, P')$.  If one writes a point, $x \in S$, as
$x = x_c + P\eta, \; \eta \in Q^n$ then the theorem says
\begin{equation}
	{P'}^{-1} ( \Phi(x_c + P \eta) - x'_c) = \eta', 
	\qquad  \eta' \in  Q^n.
\label{eqn:assert}
\end{equation}
If we take (\ref{eqn:assert}) one component at a time we find
\begin{equation}
    \left| [{P'}^{-1}
	(\Phi(x_c + P \eta) - x'_c)]_j \right| \leq {\rm 1 }.
\label{eqn:assert-equiv}
\end{equation}

To prove this for the $j$th component we consider a function
$\gamma_j: {\rm [0,1]} \rightarrow \R$, 
\mbox{$\gamma_j(t) = [P'^{-1}\Phi(x_c + t\,P\eta)]_j$}. $\gamma_j(t)$ has
the same smoothness as the map and so
the Mean Value Theorem says \mbox{$\exists t_0 \in {\rm [0,1]}$}
such that 
\begin{eqnarray*}
   \gamma_j(1) - \gamma_j(0) 
   & = &
   \frac{d\gamma_j}{dt}(t_0),  \\
   {\rm or } \qquad 
   [P'^{-1}(\Phi(x_c + P\eta) - \Phi(x_c))]_j
   & = & 
   \left[P'^{-1} \circ D\Phi_{(x_c + t_0\,P\eta)}
			\circ P\eta \right]_j.
\end{eqnarray*}
Arguing as we did in the sequence (\ref{eqn:prism-bd1});
\begin{eqnarray*}
   \left| [P'^{-1}(\Phi(x_c + P\eta) - x'_c)]_j \right|
   & = & 
   \left| \left[ W^{-1} \circ A^{-1} \left\{
				(\Phi(x_c) - x'_c) +  
				D\Phi_{\gamma(t_0)}  \circ
				P \eta \right\}
	\right]_j \right|,\\
   & = & \frac{1}{w_j} 
   	\left| \left[ A^{-1} \left\{ 
				(\Phi(x_c) - x'_c) +  
				D\Phi_{\gamma(t_0)}  \circ
			 	P \eta \right\}
	\right]_j \right|,\\
   & \leq & \frac{1}{w_j}  \left\{
	\begin{array}{l}
   		\left| [A^{-1}(\Phi(x_c) - x'_c)]_j \right|  \\
   		\qquad \qquad + \;\; {\displaystyle 
			\sum_{k=1}^n \left| [A^{-1} 
			\circ D\Phi_{\gamma(t_0)}  
			\circ P]_{jk} \right| 
		    }
	\end{array}  \right\},\\
   & \leq & {\rm 1},
\end{eqnarray*}
which is just the thing required by (\ref{eqn:assert-equiv}).

\subsection{choices for the matrix A}
\label{sec:A-choice}

Although we usually take 
$A \approx D\Phi_{x_c} \circ  P$ we may sometimes need to make a 
different choice to avoid a singular $A$.  Indeed, the very first
prisms  we consider, the ones \pagebreak of the form 
\mbox{$I_{\epsilon} \times x^{\star} \times  I_j$}, have zero
width in the $u$ direction and so have singular matrices, $P$.
In this section we will illustrate two schemes
for fattening up the matrix $D\Phi_{x_c} \circ  P$. 
	The first, the {\em fixed-form} scheme, is borrowed
	directly from \cite{MP:Conv}.  The second, called,
	the {\em column-rotor}, is a slight generalization of 
	theirs.
These techniques 
have not been carefully optimized and are probably not the
best.  They work well enough and, in any case, 
are not the most time consuming part of the algorithm.

\begin{figure}[hb]
	\parbox[b]{0.45\textwidth}{
		\caption{  \em The fixed-form fattener applied to 
			the image of
			a singular, vertical prism.  The map 
			is the delay-embedded version of the 
			standard map with ${\rm k = 0.8}$.
			The new prism, shown in grey, 
			fits snuggly in the $u$ 
			direction but is much more generous
			in the $v$ direction.
			\label{fig:ff-fattener} 
		}
	}
	\hfill
	\parbox[b]{0.5\textwidth}{
		\centering
		\includegraphics[height=7cm]{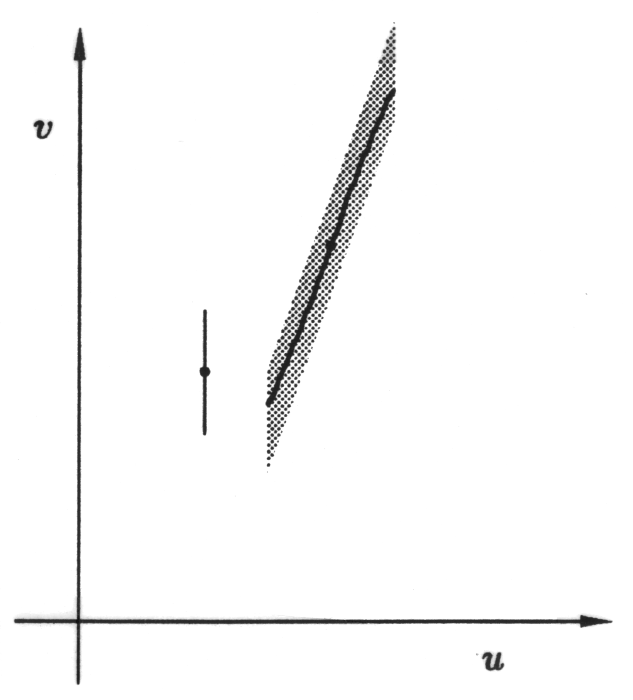}
	}
\end{figure}

\noindent
{\bf Fattener 1} (fixed-form) \hskip 1\parindent
Require the new matrix to have a particular form.  Suppose, for example,
that the initial prism, $P$, and the derivative of the map, 
$D\Phi_{x_c}$, are
\begin{eqnarray*}
P = \left[ 
		\begin{array}{cc} 
			0 & 0 \\
			0 & \frac{{\scriptstyle \Delta} y}{2} 
		\end{array}
    \right],
      \;\;
	D\Phi_{x_c} = \left[ 
			\begin{array}{cc} 
				0 & 1 \\
				-1 & \beta(x_c) 
			\end{array} 
		\right],
	& \; \; {\rm and \; so} \; \;  &
	D\Phi_{x_c} \circ  P =  
	    \left[
		\begin{array}{cc} 
		    0 & \frac{{\scriptstyle \Delta} y}{2} \\
		    0 & \frac{{\scriptstyle \Delta} y}{2} \beta(x_c) 
	    	\end{array}
	    \right],
\end{eqnarray*}
We might then look for a matrix $A$ of the form
\begin{displaymath}
A = 
\left[
	\begin{array}{cc} 
		0 & a_{12} \\
		1 & a_{22}
	\end{array}
\right].
\end{displaymath}
Figure  (\ref{fig:ff-fattener}) shows an application of this 
scheme.
\newpage

\noindent
{\bf Fattener 2} (column-rotor) \hskip 1\parindent  This method 
deals with matrices whose columns, when viewed as vectors,
are all very nearly parallel.  Such matrices will be close to singular,
and must be expected to arise if the dynamics are hyperbolic.  If
we neglect the fattening steps the matrix of the prism
bounding $\Phi^n(S_0)$ looks like 
\begin{equation}
P_n  \approx 
   D\Phi_{\Phi^{n-1}(x_c)} \circ D\Phi_{\Phi^{n-2}(x_c)} 
   \circ \cdots \circ D\Phi_{x_c} \circ P.
\label{eqn:priz-prod}
\end{equation}
If  any of the Lyapunov exponents are positive the columns of the
matrix product (\ref{eqn:priz-prod}) will be parallel to each other and
to the eigenvector corresponding to the largest eigenvalue 
of $D\Phi^n_{x_c}$. 
The idea of this scheme is to rotate the columns with 
respect to one another so as to guarantee a certain minimum angle
between each pair.
\begin{figure}
	\parbox[b]{0.5\textwidth}{
		\centering
		\raisebox{0.05\height}{\includegraphics[height=8.5cm]{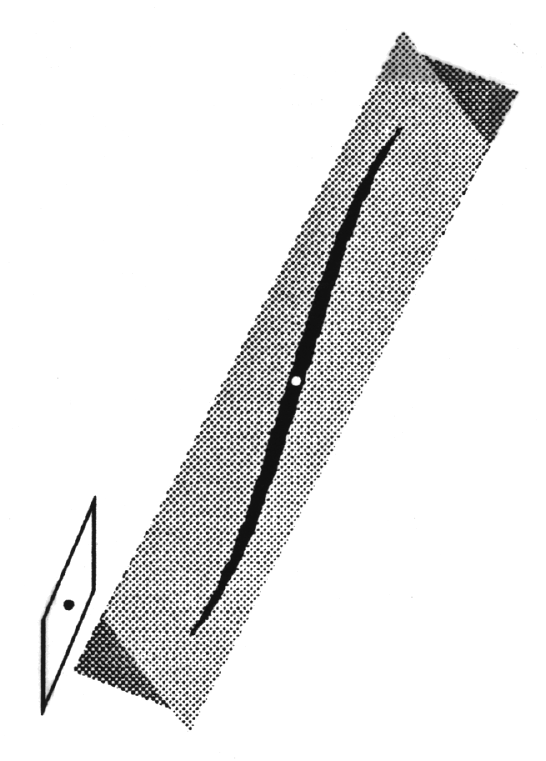}}
	}
\hfill 
	\parbox[b]{0.4\textwidth}{ 
	    \caption{  \em
	    The column-rotor scheme applied to a
	    narrow prism.  The initial prism is at the lower left;
	    it is outlined in black and its center is marked with a dot.
	    The prism's true image is solid black. A bounding prism,
	    produced with the column-rotor scheme using an angle of
	    ${\rm 27}^{\circ}$, is shown in light grey, the darker
	    prism beneath used an angle of ${\rm 90}^{\circ}$.
	    }

	    \label{fig:cr-fattener} }
\end{figure}
In two dimensions, (see figure (\ref{fig:cr-fattener})),
this is an entirely satisfactory program. 
In three and more dimensions it
is possible to find linearly dependent collections 
of column vectors each 
pair of which is separated by a sizeable angle - one could have a triple
of coplanar vectors, for example.  Such collections do not seem to
arise in our calculations, and we have made no special provisions
to avoid them.  The details  of column rotation are described in 
appendix (\ref{app:spec}).

%% file: converse/higher.tex
\section{On to higher dimension}
\label{sec:higher-d}

    Here we develop some new results.
The forms of the arguments will be much the same as in the preceding 
sections, but the maps, tori, and cones will exist in higher dimensional
spaces.  The general results for higher dimensional invariant tori 
are not so strong as for circles on the cylinder,
so we must make a few new 
restrictions and will obtain somewhat weaker results.   We will see 
how to generalize the cone-crossing and action criteria 
and then show an application to the example with the 
trigonometric perturbation,
(\ref{eqn:zamples}).

\subsection{maps and tori} 

     As above, we will consider only small perturbations of integrable
systems.  We will have 2$n\/$-dimensional symplectic maps, 
\mbox{$f_{\epsilon} : \Tn \times \Rn \rightarrow \Tn \times \Rn $},
of the form
\begin{eqnarray}
    f_{\epsilon}( \bthet, \bfp ) 
	& = & 
	(\bthet'(\bthet,\bfp), \bfp'(\bthet,\bfp)) \nonumber \\
    \bthet' 
	& = & 
	\bthet + \bfp - \frac{\partial V_{\epsilon}}{\partial \bthet } 
	\nonumber \\
    \bfp'
	& = &
	\bfp - \frac{\partial V_{\epsilon}}{\partial \bthet} 
\label{eqn:models}
\end{eqnarray}
where $V_{\epsilon}(\bthet) : \Tn \rightarrow \R$  is some periodic function
with at least two continuous derivatives and $\epsilon$ is drawn from
some, perhaps multi-dimensional, parameter space.  We will work mostly
with a lift, 
\mbox{$F_{\epsilon}:\Rn \times \Rn \rightarrow \Rn \times \Rn$}.
As we noted in chapter 2, maps
like (\ref{eqn:models}) are the higher dimensional analogs of standard-type
maps.

    The generating function for a map like (\ref{eqn:models}) is 
\begin{eqnarray}
    H_{\epsilon}(\bfx, \bfx') 
	& = &  
		\frac{1}{2} \norm{\bfx-\bfx'}^2 - 
		V_{\epsilon}(\bfx)  \nonumber \\
	& = & 
		\sum_{j=1}^{n} (x_j' - x_j)^2 - 
		V_{\epsilon}(\bfx).
\label{eqn:high-standard}
\end{eqnarray}
Although $H_{\epsilon}(\bfx, \bfx')$ is formally very similar
to the generating functions used earlier in the chapter
it is not quite the same;
the perturbion, $V_{\epsilon}$,
depends on $\bfx$ rather than $\bfx'$. As we shall 
see, this 
makes no real difference in the formulation of non-existence criteria.
We make this small change because the examples of chapter 
\ref{chap:numres}
have generating functions like (\ref{eqn:high-standard}).

     As on the cylinder, we will not be able to prove the non-existence 
of all possible types of tori, only those which are invariant graphs,
sets of the form 
$\{(\bthet,\bfp)| \bthet \in \Tn, \; \bfp = \bpsi(\bthet)\}$
for some $\bpsi : \Tn \rightarrow \Rn$.
In higher dimension we must add the further requirement that the 
graphs be
{\em Lagrangian}, that is, they must have\thesisfootnote{
	Equivalently, a Lagrangian torus is one on
	whose tangent space the symplectic two-form, 
	\mbox{$\omega = \sum_{j=1}^n dp_j \wedge d\theta_j$}, 
	vanishes. }
\begin{equation}
	\frac{\partial \psi_i}{\partial \theta_j} 
	=
	\frac{\partial \psi_j}{\partial \theta_i}.
\label{eqn:tori-def}
\end{equation}
On the cylinder we have the mighty theorem of Birkhoff
to assure us that any
rotational invariant circle must be a graph. Unfortunately,
for $n > {\rm 1}$ we have no such assurance; there may
be ``accidental'' invariant tori which are graphs, but not
Lagrangian graphs, and there may even be rotational 
invariant tori which are not graphs at all.  
Still, (\ref{eqn:tori-def}) is not a disastrous 
restriction.  Our techniques are fully complementary to traditional 
KAM theory in that
constructive versions of KAM produce just the sort of tori we 
can preclude, invariant, Lagrangian graphs.

Herman, in \cite{Herm:non}, has announced some results 
along the lines of a higher dimensional version of
Birkhoff's theorem, but they are not so comprehensive as 
the original.  He has, however, shown that a
Lagrangian graph, invariant under a map like (\ref{eqn:models}),
is  Lipschitz. This will prove helpful when we try to obtain
global inequalities like (\ref{eqn:globalLip}).

\subsection{Lipschitz cones: old formulae in new guises}

    Both the cone-crossing and action minimizing criteria have 
higher dimensional analogs. We will briefly examine the former 
because of its intuition-pleasing geometric roots, then concentrate
on the latter.  Most of the formulae will bear a strong formal
resemblance to the ones from the first part of the chapter.

    As on the cylinder, we begin by switching to  a map $g$ acting
on the delay coordinates, 
$g_{\epsilon}(\bthet_i, \bthet_{i+1}) = (\bthet_{i+1}, \bthet_{i+2})$, and a
lift, $G_{\epsilon}: \Rn \times \Rn \rightarrow \Rn \times \Rn$ with
$G_{\epsilon}(\bfu,\bfv) = (\bfu',\bfv')$.  In these 
coordinates the derivative of the map is
\begin{equation}
    DG_{\epsilon} \; = \;
       \left[ \begin{array}{cc}
		\displaystyle
			\frac{\partial \bfu'}{\partial \bfu} 
			\bigStrut &
		\displaystyle
			\frac{\partial \bfu'}{\partial \bfv} 
			\bigStrut \\
		\displaystyle
			\frac{\partial \bfv'}{\partial \bfu} 
			\bigStrut &
		\displaystyle
			\frac{\partial \bfv'}{\partial \bfv} 
			\bigStrut
		\end{array} \right] 
	=
		\left[ \begin{array}{cc}
			0 & \id  \\
			-\id & 2\id - 
		\displaystyle
		\frac{\partial^2 V_{\epsilon}}{\partial \bfx^2}(\bfv) 
		\end{array} \right],
\label{eqn:high-Jacob}
\end{equation}
where $\id$ is the $n \times n $ identity matrix and 
$\frac{\partial^2 V_{\epsilon}}{\partial \bfx^2}$ is the matrix of
second partial derivatives of $V_{\epsilon}$.
An invariant graph, $\bfp = \bpsi(\theta)$, appears as
a hypersurface
\begin{eqnarray*}
 \bfv 	& = & \bLam(\bfu),  \\
	& = & \bfu + \psi(\bfu ) - 
	    \frac{\partial V_{\epsilon}}{\partial \bfx} (\bfu). 
\end{eqnarray*}
$\Veps(\bfu)$ and $\bpsi(\bfu)$ and are periodic extensions and 
\mbox{$\bLam(\bfu + \bfm) = \bLam(\bfu) + \bfm$}
$\forall \bfm \in \Zn$.
The geometric object corresponding to a vector tangent to an invariant 
circle is now a hyperplane tangent to the graph.  A vector, 
$(\delta \bfu, \delta \bfv)$,
lying in this hyperplane has
\begin{equation}
    \delta \bfv = \bfL \delta \bfu \qquad {\rm where} \qquad \bfL =
    \left[ \begin{array}{ccc}
	\pderiv{\Lambda_1}{u_1} & \pderiv{\Lambda_1}{u_2} & \cdots \\
	\pderiv{\Lambda_2}{u_1} & \pderiv{\Lambda_2}{u_2} & \cdots \\
	\vdots & \vdots & \ddots 
    \end{array} \right]
\label{eqn:slope-obj}
\end{equation}
so that the tangent plane is the subspace spanned by the $n$ vectors
\begin{eqnarray*}
   & (1, 0, \ldots 0, \pderiv{\Lambda_1}{u_1}, 
	\pderiv{\Lambda_2}{u_1}, \ldots \pderiv{\Lambda_n}{u_1}), & \\
   & (0, 1, \ldots 0, \pderiv{\Lambda_1}{u_2}, 
	\pderiv{\Lambda_2}{u_2}, \ldots \pderiv{\Lambda_n}{u_2}), & \\
   & \vdots & 
\end{eqnarray*}
These are conveniently represented
in block form as $[\id,\bfL]$ where $\id$ is the $n \times n$ identity
matrix and $\bfL$ is as in equation(\ref{eqn:slope-obj}). The action
of the map on the hyperplane is given by
\begin{equation}
    DG_{\epsilon} \circ 
    \left[ \begin{array}{c} \id \\ \bfL \end{array} \right]
	=
    \left[ \begin{array}{cc} 
		0 & \id \\
		-\id & \bbeta
    \end{array} \right]
    \left[ \begin{array}{c} \id \\ \bfL \end{array} \right]
	=
    \left[ \begin{array}{c} \bfL \\ 
			\bbeta \bfL - \id 
    \end{array} \right],
\label{eqn:slopeObj-action}
\end{equation}
where 
$\bbeta = 2\id - \frac{\partial^2 V_{\epsilon}}{\partial \bfx^2}(\bfv)$. 
The new tangent hyperplane must then have 
\begin{equation}
	\bfL' = \bbeta - \bfL^{-1}.
\label{eqn:slopeObj-evolution}
\end{equation}
In the two dimensional slope evolution equation,
(\ref{eqn:slope-evolution}), existence of an invariant circle
meant both the slopes $\ell$ and $\ell'$ had to be positive.
Here the existence of an invariant Lagrangian graph implies 
that the matrices $\bfL$ and $\bfL'$ are positive definite.
On the cylinder we were able to study equation 
(\ref{eqn:slope-evolution}) and obtain a narrower global Lipschitz 
cone; where first we had $0 \leq \ell \leq \infty$ we eventually got
$\ell_{-} \leq \ell \leq \ell_{+}$, with $\ell_{\pm}$ given by equation
(\ref{eqn:globalLip}).  There is a higher dimensional analog of this
best global Lipschitz cone, but we defer it until section
\ref{subsec:narrow}.

\subsection{minimalism revisited}

	We now turn to the higher dimensional generalization of the 
action criterion. The first thing we need is a higher dimensional 
version of the theorem of Mather which told us that invariant circles
are composed entirely of minimizing orbits.  The necessary 
result, which says that every orbit on an 
invariant Lagrangian graph is minimizing, has been proven by Katok,
\cite{Kat:min}, and by MacKay, Meiss and Stark, \cite{MMS:Conv}.
With this result in hand we can proceed as before.  We consider
finite segments, $\bfx_{-1}, \bfx_0, \ldots x_n$ taken out of 
minimizing states.  The action functional is still
\begin{eqnarray*}
    W_{-1,n} & = & \sum_{j=-1}^{n-1} H_{\epsilon}( \bfx_j, \bfx_{j+1}), \\
	    & = & \sum_{j=-1}^{n-1} \frac{1}{2} 
		  \norm{\bfx_{j+1} - \bfx_j}^2 -
		  V_{\epsilon}(\bfx_j).
\end{eqnarray*}
and the second variation of $W_{-1,n}$ is, in block form,
\begin{displaymath}
    \left[ \begin{array}{ccccccc} 
	\bbeta(\bfx_0 ) & -\id &  \emskip 0 & 
		     \emskip 0 & & \cdots & \emskip 0 \\
	-\id & \bbeta(\bfx_1) & -\id &  \emskip 0 &
			      & \cdots & \emskip 0  \\
	\emskip 0 & -\id & \bbeta(\bfx_2) & -\id & 
			 & & \emskip 0 \\
	\emskip \vdots & & \ddots & \ddots & \ddots & & \emskip \vdots \\
	\emskip 0 & & & & -\id & \bbeta(\bfx_{n-2}) & -\id \\
	\emskip 0 & & \cdots & & \emskip 0  
			     & -\id & \bbeta(\bfx_{n-1}) 
    \end{array} \right],
\end{displaymath}
which is readily block-diagonalized to 
\begin{displaymath}
    \left[ \begin{array}{ccc}
	    \bfd_0 & 0 & \cdots \\
	    0 & \bfd_1 & \cdots \\
	    \vdots & \vdots & \ddots
    \end{array} \right].
\end{displaymath}
The diagonal blocks, $\bfd_j$, are given recursively by 
\begin{eqnarray}
    \bfd_0 & = & \bbeta(\bfx_0),  \nonumber \\
	\bfd_{j+1} & = & \bbeta(\bfx_{j+1}) - \bfd_j^{-1}, 
	\qquad \qquad
    \bbeta(\bfx_{j+1}) = 2\id - 
	\frac{\partial^2 V_{\epsilon}}{\partial \bfx^2}(\bfx_{j+1}).
\label{eqn:block-evolution}
\end{eqnarray}
Our concern is that the $\bfd_j$ be positive definite.  It is here 
that blithe, formal, generalization fails us; there are no sensible
formal analogs for results like equations (\ref{eqn:softDG}), 
(\ref{eqn:globalLip}) and (\ref{crit:four-thirds}).  Instead we need
to invent a way to test whether the least eigenvalue of $\bfd_j$ 
is positive. We will develop a whole suite of estimates for this
eigenvalue, then use them and a plan like the one in section 
\ref{subsec:plan} to prove the non-existence of Lagrangian graphs.

    All the matrices we will be discussing are real and symmetric,
hence, Hermitian.  For a particular matrix, $M$, we will need to 
define
$\lambda_-(M)$, the least eigenvalue of $M$, $\lambda_+(M)$,
the largest eigenvalue, and $\Tr{M} = \sum_{j=1}^{dim(M)} M_{jj}$, 
the trace.  The following lemma will be our main tool.

\noindent
{\bf Lemma}
\hskip 1\parindent {\em
For real, symmetric, $n \times n$, positive definite matrices 
$\bbeta$, $\bfd$, and $\bfd'$ with 
\begin{equation}
	\bfd' = \bbeta - \bfd^{-1}
\label{eqn:lemma}
\end{equation}
the following suite of inequalities hold:
\begin{eqnarray}
    \lamM{\bfd'} & \leq & \frac{1}{n} \, \Tr{\bbeta} -
	\frac{n}{\Tr{\bfd}},  \label{eqn:ineq1}  \\
    \lamM{\bfd'} & \leq & \lamP{\bbeta} - 
	\frac{1}{\lamM{ \bfd }},  \label{eqn:ineq2}  \\
    \lamM{\bfd'} & \leq & \lamM{\bbeta} - 
	\frac{1}{\lamP{ \bfd }}.  \label{eqn:ineq3} 
\end{eqnarray}
}

\noindent
{\bf Proof} \hskip 1\parindent
The first inequality, which is due to Herman, comes from the 
observations that for a positive definite, Hermitian matrix,
$M$, $\lamM{M} \leq \frac{1}{n}\Tr{M}$ and 
$\Tr{M^{-1}} \leq \frac{n^2}{\Tr{M}}$. Both these inequalities
are strict except for the degenerate case where all the 
eigenvalues are the same.
The other two inequalities depend on 
\begin{displaymath}
	\lamP{M} =  \max_{\nu \in \Rn,\;\norm{\nu} = 1}
				\iprod{\nu}{M\nu}
\end{displaymath}
and
\begin{displaymath}
	\lamM{M} =  \min_{\nu \in \Rn ,\; \norm{\nu} = 1}
				\iprod{\nu}{M\nu},
\end{displaymath}
where the norm and inner product are
the usual Euclidean norm in $\Rn$ and ordinary dot product, 
$\iprod{\bfu}{\bfv} = \sum_{j=1}^{n} u_j v_j$.
Given these equations we can obtain inequalities about the least
eigenvalue of $\bfd'$ in  (\ref{eqn:lemma}) by evaluating
$\iprod{\nu}{\bfd' \nu}$ on particular vectors.
If, for example, one takes $\nu$ to be the unit eigenvector
corresponding to the  smallest eigenvalue of $\bfd$ one finds
\begin{eqnarray*}
    \lamM{\bfd'} \; \leq \; \iprod{\nu}{\bfd' \nu}
	& = &
	\iprod{\nu}{\bbeta \nu} - \iprod{\nu}{\bfd^{-1} \nu}, \\
	& = & 
	\iprod{\nu}{\bbeta \nu} - \frac{1}{\lamM{\bfd}}, \\
	& \leq &  
	\lamP{ \bbeta } - \frac{1}{\lamM{\bfd}}.
\end{eqnarray*}
This is inequality (\ref{eqn:ineq2}) of the lemma. 
Inequality (\ref{eqn:ineq3})
comes from an identical argument with $\nu$ the unit
eigenvector corresponding to the least eigenvalue of $\bbeta$.

\subsection{global estimates: narrowing the cones}
\label{subsec:narrow}

Here we see how to use our inequalities to reduce the range
of permissible $\lamM{\bfd_j}$.
On the face of it, we must allow $ 0 \leq \lamM{ \bfd } \leq \infty$,
but inequalities (\ref{eqn:ineq1}) and (\ref{eqn:ineq2}) 
have the correct form to allow an iterative refinement like the
one in section \ref{subsec:two-crit}.  Since $\Tr{\bbeta(\bfv)}$,
and $\lamP{\bbeta(\bfv)}$ are continuous, $\Zn$-periodic functions,
they have well defined minima and maxima, say,
\begin{eqnarray*}
	t \;\leq & \Tr{\bbeta} & \;\leq T, \\
	b \;\leq & \lamP{\bbeta} & \;\leq B. 
\end{eqnarray*}
Inequalities (\ref{eqn:ineq1}) and (\ref{eqn:ineq2}) then 
imply that the $\bfd_j$ from a minimizing state must satisfy 
\begin{eqnarray}
	\TrM \leq \Tr{\bfd_j} \leq \TrP,
	    & \; {\rm with } \; &
	    \TrM = \lb \left\{
	   	 	\frac{T - \sqrt{T^2 -4n^2}}{2}
			\right\}, \nonumber \\
	&&  \TrP = \ub \left\{
			\frac{T + \sqrt{T^2 -4n^2}}{2}
			\right\},
\label{eqn:trace-cone}
\end{eqnarray}
and
\begin{eqnarray}
	\lambda_{-min} \leq \lamM{\bfd_j} \leq \lambda_{-max},
	    & \; {\rm with } \; &
	    \lambda_{-min} =  \lb \left\{
				  \frac{ B - \sqrt{B^2 -4}}{2}
				  \right\}, \nonumber \\
	&&  \lambda_{-max} =  \ub \left\{
				  \frac{B + \sqrt{B^2 -4}}{2}
				  \right\}.
\label{eqn:lambda-cone}
\end{eqnarray}
We can also get some analytic use out of inequality (\ref{eqn:ineq3})
by combining it with (\ref{eqn:lambda-cone}).
\begin{eqnarray*}
    \lamP{\bfd} 
	& \leq & \Tr{\bfd} \, - \, (n-1)\lamM{\bfd} \nonumber \\
	& \leq & \Tr{\bfd} \, - \, (n-1)\lambda_{-min}. 
\end{eqnarray*}
Hence,
\begin{eqnarray}
    \lamM{\bfd'} 
	& \leq & \lamM{\bbeta} \; - \; 
		 \frac{1}{\lamP{\bfd}} \nonumber \\
	& \leq & \lamM{\bbeta} \; - \; 
		 \frac{1}{\Tr{\bfd}\,-\,(n-1) \lambda_{-min}}.
\label{eqn:ineq4}
\end{eqnarray}
This profusion of inequalities makes possible a whole host 
of ``Mather $\frac{4}{3}$'' arguments; Herman, in \cite{Herm:non},
gave the one based on (\ref{eqn:ineq1}) and  (\ref{eqn:trace-cone}).
In the next section we show how to apply his criterion, along with 
other, new ones, to a specific example.

%% file: converse/results.tex
\section{A converse KAM theorem}
\label{sec:conv-pix}

Here we use the arguments above on a specific system, the 
trigonometric example from chapter \ref{chap:numres}.
We will use the same example to illustrate some\thesisfootnote{
	Appendix \ref{app:spec} gives a detailed discussion
	of the algorithms used and includes a specification
	of the functions and data structures.  The code 
	itself is in appendix \ref{app:code}.
}
of the issues
in proving a machine-assisted converse KAM theorem and will
show the results of several calculations.

\subsection{analytic preliminaries}

The plan for a converse KAM theorem, section \ref{subsec:plan},
requires a starting point, $\bfxStar$, and the constants 
$t,\,T,\,b$ and $B$ from equations (\ref{eqn:trace-cone}) and 
(\ref{eqn:lambda-cone}).  For the example at hand, 
\begin{eqnarray*}
    \bbeta(\bfv) & = &
	2\id - \epsilon \frac{\partial^2 V_{trig}}{\partial \bfx^2}, \\
    & = &
	2\id -  \frac{\epsilon}{M_{trig}}
			\left[ \begin{array}{cc}
			     \{ \frac{ \sin 2\pi v_0 }{2} + 
				\sin 2\pi(v_0 + v_1) \} &
			     \sin 2\pi(v_0 + v_1)  \\
			     \sin 2\pi(v_0 + v_1)  &
			     \{ \frac{ \sin 2\pi v_1 }{2} + 
				\sin 2\pi(v_0 + v_1) \} 
			\end{array} \right]
\end{eqnarray*}
and so 
\begin{eqnarray}
    \Tr{\bbeta (\bfv)} & = & 4 -  \frac{\epsilon }{M_{trig}} \left\{
			\half \{\sin 2\pi v_0 + \sin 2\pi v_1\}
		 	- 2 \sin 2\pi(v_0 + v_1) \right\} \\
    \lamM{\bbeta(\bfv)} & = & 
		\half \left\{ \begin{array}{l}
		\Tr{\bbeta(\bfv)}  \;-   \\
		\; \; \; \frac{\epsilon}{M_{trig}}
		\sqrt{ \oqrt \left( \sin2\pi v_0 + \sin2\pi v_1 \right)^2 
		+ 4\sin^2 2\pi(v_0+v_1)}
	    \end{array} \right\}
\label{eqn:zamp-expers}
\end{eqnarray}
Both $\Tr{\bbeta}$ and $\lamM{\bbeta}$ achieve their extrema on the 
line $v_0 = v_1$.  The symmetries of $V_{\epsilon}$ also ensure that 
\begin{eqnarray*}
    t - 4 =
	\epsilon \; \min \Tr{\frac{\partial^2 \Vtrg}{\partial \bfx^2}} 
    = -\epsilon \; \max \Tr{\frac{\partial^2 \Vtrg}{\partial \bfx^2}} 
    = 4 - T \\ 
    b - 2 = \epsilon \; \min 
		\lamM{ \frac{\partial^2 \Vtrg}{\partial \bfx^2}} 
    = -\epsilon \; \max 
		\lamM{ \frac{\partial^2 \Vtrg}{\partial \bfx^2}} 
    = 2 - B 
\end{eqnarray*}
We find the approximate positions of the extrema using Newton's
method, then evaluate the bounds $t,\,T$ etc..  From these we can 
calculate the ranges of permissible $\lamM{\bfd_j}$.  

     The choice of the starting point, $\bfxStar$, depends on 
which of the inequalities (\ref{eqn:ineq1}) - (\ref{eqn:ineq3})
we expect to be most fruitful.  Good use of inequality 
(\ref{eqn:ineq1}) would require that $\bfxStar$  be a 
place where $\Tr \bbeta$ attains its minimum; this choice 
immediately gives $\epsilon_c \leq {\rm 0.0435}$.  Best
use of inequalities (\ref{eqn:ineq2}) and (\ref{eqn:ineq3}) 
requires $\bfxStar$ at a place where 
\begin{equation}
	\lamM{\bbeta} = b.
\label{eqn:xstar-def}
\end{equation}
This turns out to be the best choice;
it immediately gives $\epsilon_c \leq 0.0278$.  Note that we need
not be particularly rigorous about finding $\bfxStar$.
Indeed, we are free to choose it anywhere we like; we just get much
better results if (\ref{eqn:xstar-def}) is satisfied. 

\subsection{the computations}

    Once $\bfxStar$ is chosen, we can set up the extended phase 
space, $I_{\epsilon} \times \Rn \times \Rn$, extend $G_{\epsilon}$
to $G$ as in (\ref{eqn:G-extension}), and proceed with 
a proof.  The plan is the same as in section \ref{subsec:plan},
except that here the role of the intervals, $I_j$, is played
by  rectangles in the unit square.  That is, we first ask
``Can any $\bfx \in [{\rm 0,1}] \times [{\rm 0,1}]$ follow 
$\bfxStar$ in a minimizing state?''  If the answer is 
``no'' then we are finished, if not we cut the square in half
and ask the same question for each piece.  Once the rectangle
of potential successors is smaller than the whole square we can
iterate the argument for several steps, bounding image prisms as 
in section \ref{sec:prism-bondage}.  This yields a sequence of prisms 
in the extended phase space, $S_0, \, S_1, \cdots$, with
\begin{eqnarray*}
    S_0 & = & 
    I_{\epsilon} \times \{\bfxStar\} 
		 \times \{ \mbox{\em successor rectangle} \} 
    \equiv (\bfx_{c,0}, P_0) \\
    S_1 & = & (\bfx_{c,1}, P_1)  \supset G(S_0) \\
    & \vdots &
\end{eqnarray*}
Beginning with
\begin{displaymath}
	\ub \lamM{\bfd_{-1}} \equiv \lambda_{-max}
	\qquad {\rm and} \qquad
	\ub \Tr{\bfd_{-1}} \equiv \TrP
\end{displaymath}
we proceed, at each step evaluating the whole suite 
\begin{eqnarray}
    \lamM{\bfd_{j+1}} & \leq & 
	\ub_{(\epsilon,\bfu,\bfv) \in S_{j+1}}
		\left( \frac{1}{n} \, 
		\Tr{\bbeta(\bfv)}\right) -
		\frac{n}{\ub(\Tr{\bfd_j})}  \label{eqn:sineq1} \\ 
    \lamM{\bfd_{j+1}} & \leq & 
	\ub_{(\epsilon,\bfu,\bfv) \in S_{j+1}}
		\left(\lamP{\bbeta(\bfv)}\right) -
		\frac{1}{\ub(\lamM{\bfd_j})}  \label{eqn:sineq2} \\ 
    \lamM{\bfd_{j+1}} & \leq  & 
	\ub_{(\epsilon,\bfu,\bfv) \in S_{j+1}}
		\left(\lamM{\bbeta(\bfv)}\right) -
		\frac{1}{\ub(\Tr{\bfd_j}) - 
			 \lambda_{-min}} \label{eqn:sineq3} 
\end{eqnarray}
and choosing the best upper bound.  Computing (\ref{eqn:sineq1}) automatically
gives the bound on $\Tr{\bfd_j}$ used in (\ref{eqn:sineq3}).
These estimates do not, of course, keep improving forever.
Eventually either one of the $\ub \lamM{\bfd_j}$ falls below 
$\lambda_{-min}$ or one of the prisms $S_j$ gets so large
that the inequalities (\ref{eqn:sineq1}) - (\ref{eqn:sineq2}) are vacuous.
At that point one either quits or cuts the initial prism 
in half\thesisfootnote{
	The choice of which cut to make, whether along the $\epsilon$,
	$v_0$ or $v_1$ axis, depends on the shape of the final $S_j$.
}
and starts over.

\subsection{results}
\label{subsec:results}

Table (\ref{tab:results}) summarizes our results.  We were able to show
that the last few of the minimizing states of section 
\ref{sec:pix}
persist beyond the point where no invariant tori remain.  
\begin{table}[hbt]
\centering
    \begin{tabular}{|l|rrrr|} 
   	\hline \hline
    	\multicolumn{1}{|c|}{$\ub \epsilon_c \leq$} &
    	\multicolumn{1}{c|}{\em longest} &
    	\multicolumn{1}{c|}{\em deepest} &
    	\multicolumn{1}{c|}{\em prisms} &
    	\multicolumn{1}{c|}{\em time (sec.)} \\ \hline
    	0.0278  & 3 {\ \ } & 10 {\ \ } & 39 {\ \ } & 500 {\ \ } \\
    	0.0276  & 4 {\ \ } & 11 {\ \ } & 64  {\ \ } & 759 {\ \ } \\
    	0.0274  & 4 {\ \ } & 13 {\ \ } & 156 {\ \ } & 2698 {\ \ } \\
    	0.0272  & 6 {\ \ } & 21 {\ \ } & 933 {\ \ } &  $\sim$ {\ \ } \\
	\hline \hline
    \end{tabular}

    \caption{ \em
    A sequence of bounds on $\epsilon_c$ and some
    details about the computations which verified them.  The table
    includes: {\em longest}, the length of the longest sequence
    of image prisms considered; {\em prisms} the total number
    of prisms on which the algorithm succeeded; {\em deepest},
    the number of refining cuts needed to make the smallest
    successful prism and {\em time} the execution time in seconds.
    All computations were done on a Sun4.
    }
\label{tab:results}
\end{table}

The figures on the 
following pages show some of the systems of prisms used in the proofs.
The dark grey rectangles are sets which cannot contain a successor to 
$x^{\star}$, the light grey regions may be ignored on 
account of symmetry, (see section \ref{subsec:symmetry}).  
As one might expect, those states which go from $\bfxStar$ to 
neighborhoods near the the maximum of $\Vtrg$, (those which correspond
to rectangles in the upper right corner), are harder to prove
non-minimizing.  To succeed on such a rectangle the program must extend
the corresponding state far enough to evaluate several 
$\ub \lamM{ \bfd_j }$.  Since the prism-bounding algorithm always gives
an $S_{j+1}$ bigger than the true image of $S_j$, the initial prisms 
must be small.
\begin{figure}[p]
    \begin{center}
    	\includegraphics[height=16cm]{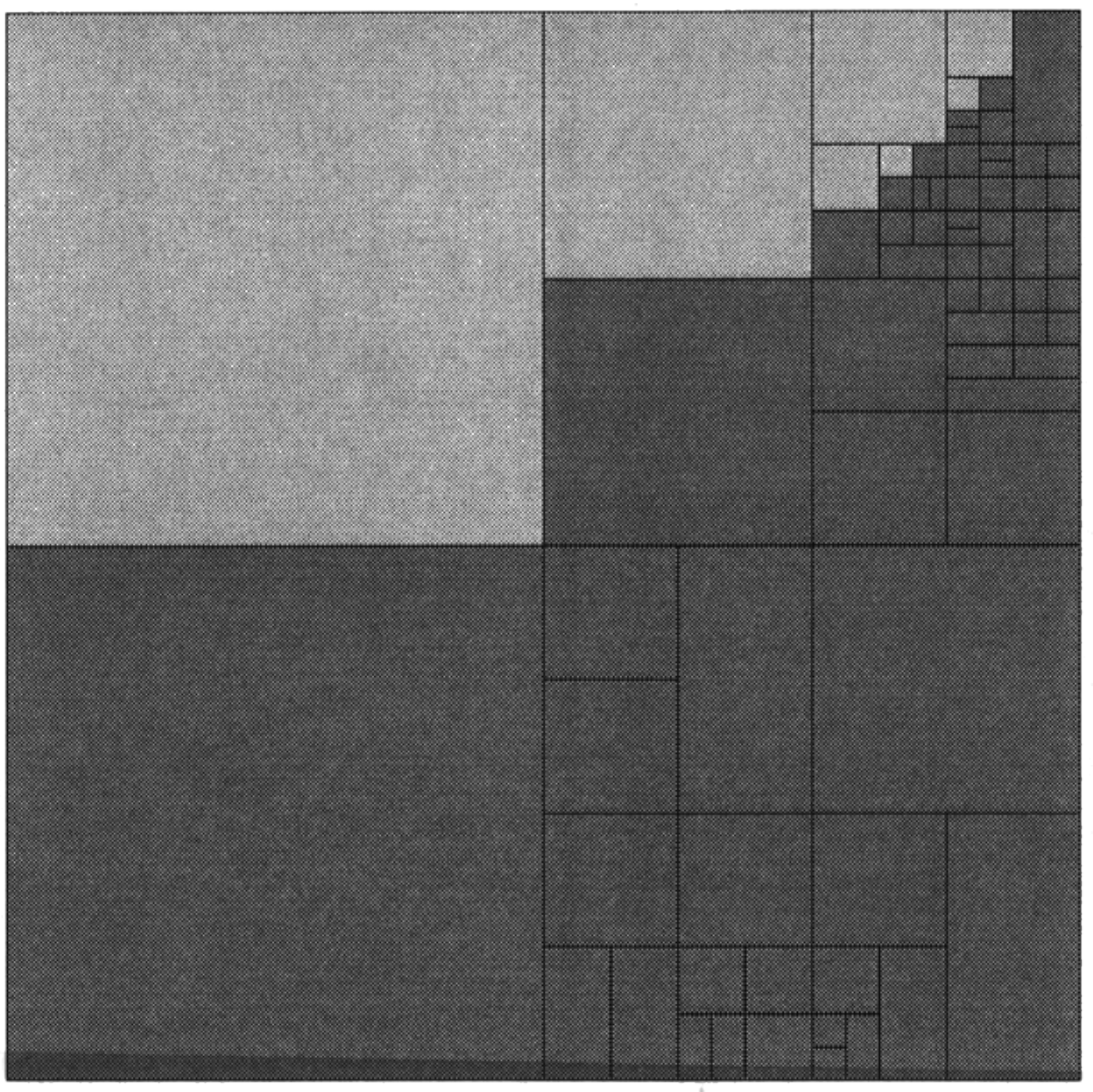}
    \end{center}
    \vspace*{-2\belowdisplayskip}
    \caption{\em
	The system of prisms used to show 
	$\epsilon_c \leq {\rm 0.0276}$.
    }
\label{fig:theorem1}
\end{figure}
\begin{figure}[p]
    \begin{center}
    	\includegraphics[height=16cm]{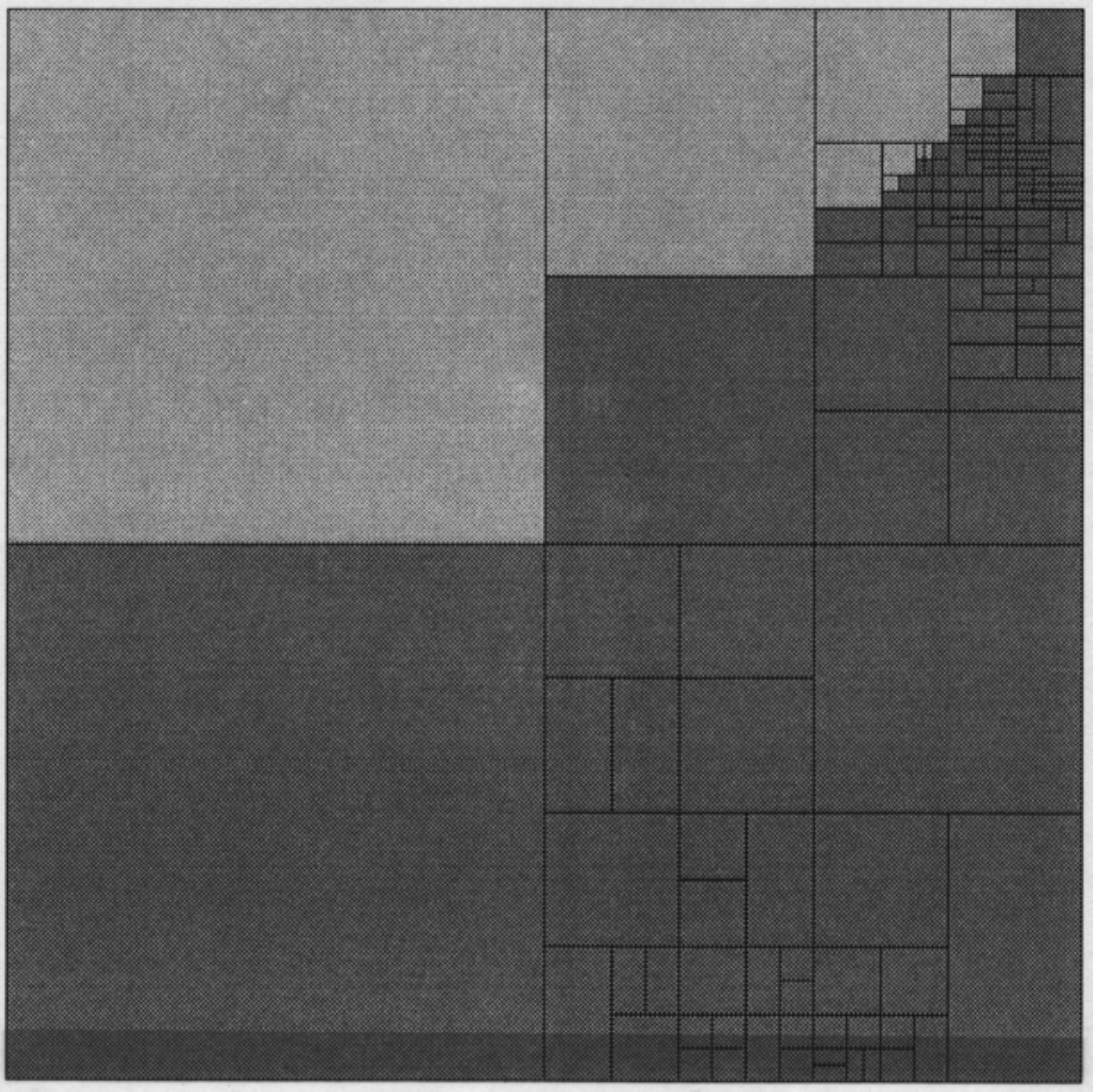}
    \end{center}
    \vspace*{-2\belowdisplayskip}
    \caption{ $\epsilon_c \leq {\rm 0.0274}$ }
\label{fig:theorem2}
\end{figure}
\begin{figure}[p]
    \begin{center}
    	\includegraphics[height=16cm]{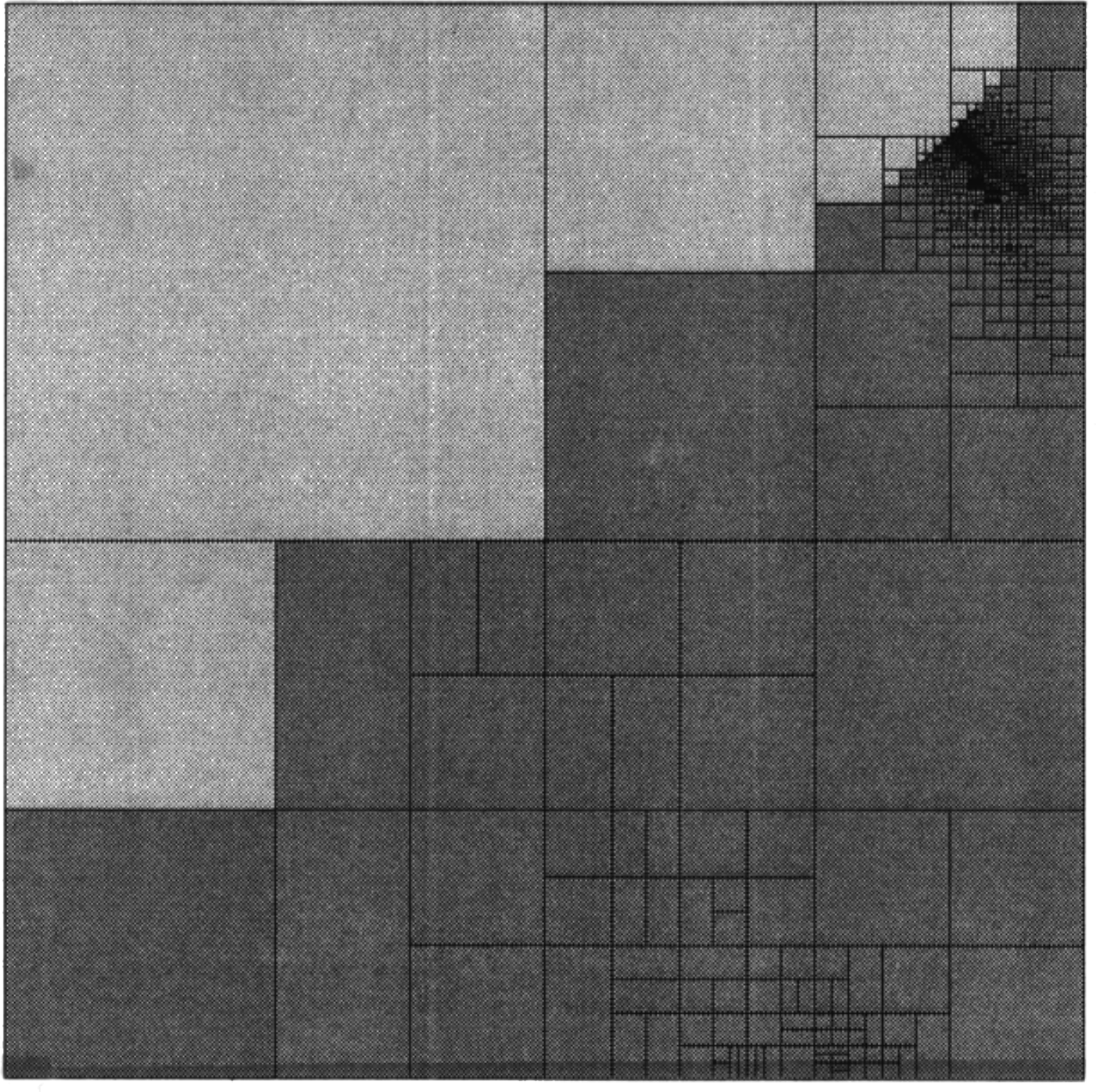}
    \end{center}
    \vspace*{-2\belowdisplayskip}
    \caption{ $\epsilon_c \leq {\rm 0.0272}$ }
\label{fig:theorem3}
\end{figure}

\subsection{using symmetry}
\label{subsec:symmetry}

\begin{figure}[bh]
    \begin{center}
    	\includegraphics[height=6.0cm]{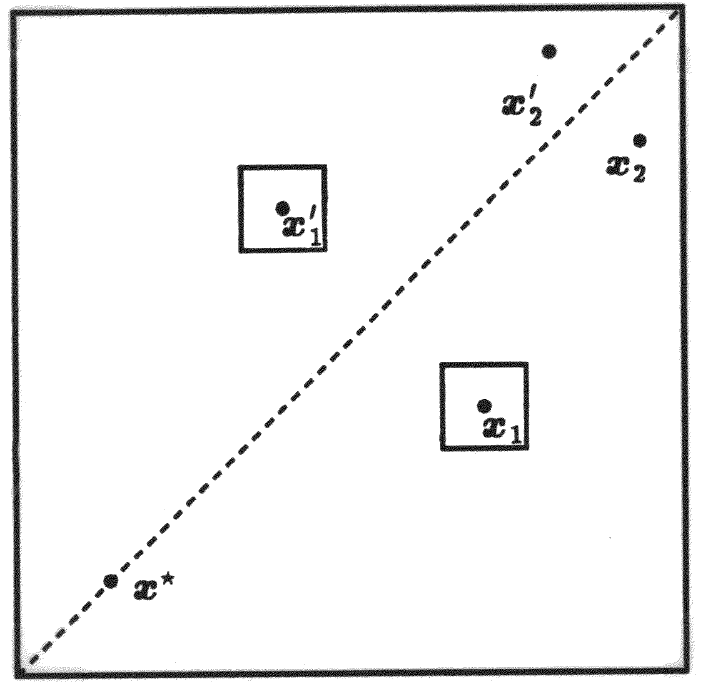}
    \end{center}
    \vspace*{-2\belowdisplayskip}
    \caption{\em Two symmetrically related states have the same action. }
\label{fig:symmetry}
\end{figure}

    In figures (\ref{fig:theorem1}) -- (\ref{fig:theorem3}) 
we were able to ignore around half the possible successors.
To see why, notice that $\Vtrg$ is unchanged by the interchange
of its $v_0$ and $v_1$ arguments.  Two segements, such as 
$\{\cdots,\bfxStar, \bfx_1, \bfx_2, \cdots\}$ and 
$\{\cdots,\bfxStar, \bfx_1', \bfx_2', \cdots\}$ in figure
(\ref{fig:symmetry}), will have the same action because they
are each other's images under the interchange 
$x_{j,0}  \rightleftharpoons x_{j,1}$. Here, 
the interchange is just a reflection about the line\thesisfootnote{
	One must take some care here. The interchange is really 
	a reflection through the diagonal line containing
	$\bfxStar$.  Our program always arranges that
	$\bfxStar$ is in the square 
	$[{\rm 0,1}] \times [{\rm 0,1}]$ and on 
	the line $x_0 = x_1$.
}
$x_0 = x_1$.  
So, refering to figure (\ref{fig:symmetry}), if we prove 
that no minimizing state can pass from $\bfxStar$  through 
the box around  $\bfx_1$, we are automatically assured that
none can go through the box around $\bfx_1'$ either.

%% file: AppendA.tex
\input{appendices/numMeth/introA}
\input{appendices/numMeth/methods1}
\input{appendices/numMeth/methods2}

%% file: appendices/numMeth/introA.tex
\appendix
\chapter{ Approximate Numerical Methods }
\label{app:num-methods}

In this appendix we review the numerical methods used to 
obtain the results of chapter 2.  The first section describes 
the methods used to calculate the minimizing states; the next
section discusses Kim and Ostlund's scheme for approximating 
irrational vectors by rational ones and the last section 
explains how we found the Lyapunov exponents pictured in 
figure (\ref{fig:Lyaps}).

%% file: appendices/numMeth/methods1.tex
\section{Methods of minimization}

All our minimization schemes solve the Euler-Lagrange equations
(\ref{eqn:KB-ELg}).
For each rotation vector, $\bfp / q$ and perturbation 
we produce a sequence of 
states $ \{ X_0, X_1, \ldots
X_k, \ldots \} $ each of which satisfies (\ref{eqn:KB-ELg})
for a particular value of $ \epsilon = \epsilon_j $.
We usually begin with a state whose first point, $\bfx_0$,
lies on the minimum
of the perturbation to the generating function 
(that is, on a maximum of 
$V_\epsilon(\bfx)$) and whose other points are
\mbox{ $ \bfx_j = \bfx_0 + 
\frac{ \scriptstyle j}{ \scriptstyle q} \bfp$ }.  Such a
state is globally minimizing for the unperturbed generating function so
we set $\epsilon_0 = 0 $.
We then increase the size of the perturbation, $\epsilon_j$, in small
steps and use  $X_j$ as a starting point to calculate $X_{j+1}$ using 
either a \mbox{gradient-flow} scheme or Newton's method.

   The former involves integrating the system of differential equations 
\begin{displaymath}
	\frac{d\bfx_i}{d\tau} = 
	\frac{\partial L_{p,q}}{\partial \bfx_i}, \qquad
\end{displaymath}
through a long interval of the formal ``time,'' $\tau$.  This method 
is very slow; it crawls down to the minimum with exponentially 
decreasing speed.  On the other hand it is extremely reliable and seems
very rarely to converge to a state other  than the global minimum.  
Newton's method is much faster, but somewhat prone to converge to 
extrema other than the minimum.  It works by producing a sequence 
of approximate states ${ Y_0, Y_1, \ldots } $ according to the recursive
scheme :
\begin{eqnarray}
Y_0 = \mbox{\em some initial guess,} & & 
		Y_{i+1} = Y_i + D_i  \nonumber \\
	    & & D_i =  -{\bf H}^{-1} d(L_{p,q})
\label{eq:Newt}
\end{eqnarray}
where ${\bf H}^{-1}$ is the inverse of the Hessian of the action 
functional and $d(L_{p,q})$ is the functional's gradient. Since 
{\bf H} has $ (qd)^2 $
entries, solving (\ref{eq:Newt}) could be an O($(qd)^2$) process,
but our Hessian,
\begin{displaymath}
\left[ 
  \begin{array}{cccccc}
    2{\bf I} - \epsilon \, {\bf V}_0 & 
    -{\bf I}  &
    \;\:0  & \cdots & \cdots &  -{\bf I}  \\
    -{\bf I} &
    2{\bf I} - \epsilon \,  {\bf V}_1 &
    -{\bf I} &
    \cdots & \cdots &  \;\:0 \\
    \;\:\vdots & & & \ddots &        & \;\:\vdots \\
    \;\:\vdots & & &        &  \;\:\vdots \\
    \;\:0 & \cdots & \cdots & 
    -{\bf I} & 2{\bf I} - \epsilon \,
			{ \bf V}_{q-2} &
			-{\bf I} \\
    -{\bf I} & \cdots & \cdots & \cdots &
    -{\bf I} & 2{\bf I} - \epsilon \,
			{\bf V}_{q-1} 
  \end{array}
\right] ,
\end{displaymath}
where 
\begin{displaymath}
	{\bf I} = \left[ \begin{array}{cc}
				1 & 0 \\
				0 & 1 
			 \end{array} 
		  \right] , \qquad
	{\bf V}_j \; \equiv \;
		\frac{\partial^2  V}{\partial \bfx^2} (\bfx_j) = 
		\left[ \begin{array}{cc}
		    \frac{\displaystyle \partial^2 V}
			 {\displaystyle \partial x_0^2} &
		    \frac{\displaystyle \partial^2 V}
			 {\displaystyle \partial x_0 \partial x_1} \\
		    \frac{\displaystyle \partial^2 V}
			 {\displaystyle \partial x_0  \partial x_1}  
			 \bigStrut &
		    \frac{\displaystyle \partial^2 V}
			 {\displaystyle \partial x_1^2}  \bigStrut
		\end{array} \right](x_j) ,
\end{displaymath}
has only a few terms off the diagonal.  We implemented two schemes
to solve~(\ref{eq:Newt}), one which does Gauss-Jordan elimination 
\cite{NumRecp}
and another, rather more complicated 
algorithm which generalizes the 1-d
work of Percival and Metsel \cite{MP:Newt}.  We tried the latter
because we hoped it would be more numerically stable; it was not,
and ran a bit more slowly than the Gauss-Jordan program.

%% file: appendices/numMeth/methods2.tex
\newcommand{\dfnu}[1] {DF^q_{(x_0, p_0)}\,\nu_{#1}}

\section{Rational approximation of irrational vectors }

    The problem of approximating a single real number by a sequence
of rationals is completely solved by the simple continued fraction 
algorithm \cite{Khin:num,Rob:num}.  We write
\begin{eqnarray}
\omega & = & 
    a_0 + 
	 \frac{\displaystyle 1}{ a_1 + 
	   \frac{\displaystyle 1 \cfStrut}{\displaystyle a_2 +
	     \frac{\displaystyle 1 \cfStrut}{\displaystyle a_3 +
	       \frac{\displaystyle 1 \cfStrut}{\displaystyle a_4 + 
					\cfStrut	\ddots}
								}
								  }
								   }
\label{eqn:cfrac}
\end{eqnarray}
where the $a_i$, called the {\em partial quotients} of $\omega$,
are positive integers.  We compute them recursively according to 
\begin{eqnarray*}
r_0 = \omega & & a_i = \mbox{Int}[r_i] \\
	     & & r_{i+1} = \frac{1}{r_i - a_i}.
\end{eqnarray*}
If $\omega$ is rational then all but finitely many of the $a_i$ 
are zero, but if $\omega$ is irrational then the sequence 
never terminates.
Truncating the expansion~(\ref{eqn:cfrac}) after finitely many $a_i$
gives a sequence of rational approximations $ \frac{p_0}{q_0},
\frac{p_1}{q_1}, \ldots  $ with many desirable properties.  
Each $\frac{p_i}{q_i}$ is a best approximation in the sense that
the only rationals closer to $\omega$  have larger denominators.
Further, the sequence contains infinitely many $\frac{p_i}{q_i} $
such that $ | \, \omega - {p_i}/{q_i} \, | \leq 1 / \sqrt{5}\,q^2 $.
Indeed, the extremely good convergence 
of this sequence can be a problem.
If one wants many approximations with modest 
denominators one must either
study numbers which, like the golden mean,  have very slowly growing 
$q_i$, or introduce other approximation algorithms which produce more
slowly converging sequences.

    One such algorithm depends on the Farey tree construction of the
rationals.  In a Farey tree one represents the rational number
$\frac{p}{q} $ as an ordered pair \mbox{$(p,q)$}. 
The endpoints of the unit
interval are thus \mbox{$(0,1)$} and \mbox{$(1,1)$}.  The
construction proceeds
by successively splitting intervals with endpoints \mbox{$(p_l,q_l)$} 
and \mbox{$(p_r, q_r)$} into two {\em daughter} 
intervals by inserting an 
interior point at   
\mbox{$((p_l+p_r), (q_l+q_r))$.}
The number \mbox{$((p_l+p_r),(q_l+q_r))$} is called the {\em mediant}
of $(p_l,q_l)$ and $(p_r,q_r)$.
A sequence of Farey subdivisions which begins from the unit interval
will eventually produce all rational numbers, each rational appearing
as a mediant exactly once and in lowest terms.
\begin{figure}
\parbox[b]{0.15\textwidth}{
	{\em level 0} \\
	{\em level 1} \\
	{\em level 2}
}
\raisebox{-0.1\height}{\includegraphics[height=2.5cm]{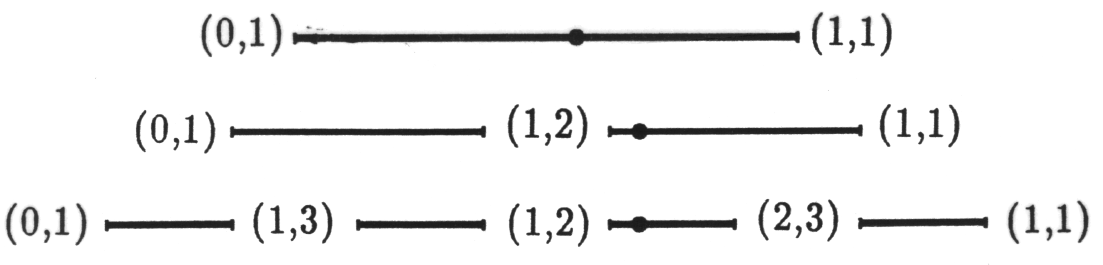}}
\caption{\em Several levels of the Farey tree.  The solid dot shows
	     the position of the golden mean. Its n{\em th} approximation
	     is always the mediant which has the largest sum $p_n+q_n$   
	     of any appearing at at the n{\em th} level.
	     \label{fig:tree}}
\end{figure}
We can use the Farey tree as a tool for rational 
approximation by choosing
$p_n / q_n $ to be the mediant of the 
n{\em th} level interval containing
$\omega$.  
Since an interval in the
n{\em th} level of the tree has length at most $1 / {n+1}$ the 
sequence of Farey approximations must eventually converge.
Since every sequence of Farey approximation
begins with $p_0 / q_0 = \half $ and  each subsequent approximation
requires only a choice of either the left or right daughter interval,
we can represent the sequence of Farey approximations as a binary
address. For example, the address
$llllll \ldots$ would indicate that $\omega$ lies always between
$(0,1)$ and $(1,n)$.

\begin{figure}
     \begin{center}
     	\includegraphics[height=3.75cm]{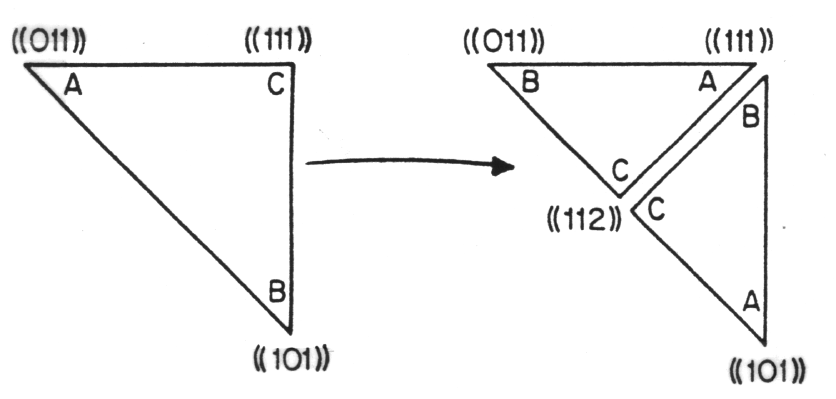}
     \end{center}
     \vspace*{-2\belowdisplayskip}
     \caption{\em  The mediant operation which refines Farey triangles.
		   The parent triangle is represented by an equilateral
		   right triangle.  The algorithm divides this into two 
		   similar, daughter triangles by adding a
		   new point in the middle of the hypotenuse.  The 
		   coordinates of the new point are sums of the  
		   coordinates of the end points of the hypotenuse. 
		   \rm \protect \cite{KO}
		   \label{fig:trimed} }
\end{figure}

\begin{figure}
	\parbox[b]{0.3\textwidth}{\caption{\em
		Five levels of the Farey triangulation, {\em (a)},
		and, {\em (b)}, the corresponding partition of the
		unit square. 
		{\rm \protect \cite{KO} }}}
	\hfill
	\parbox[b]{0.65\textwidth}{
		\raisebox{-0.1\height}{\includegraphics[height=5cm]{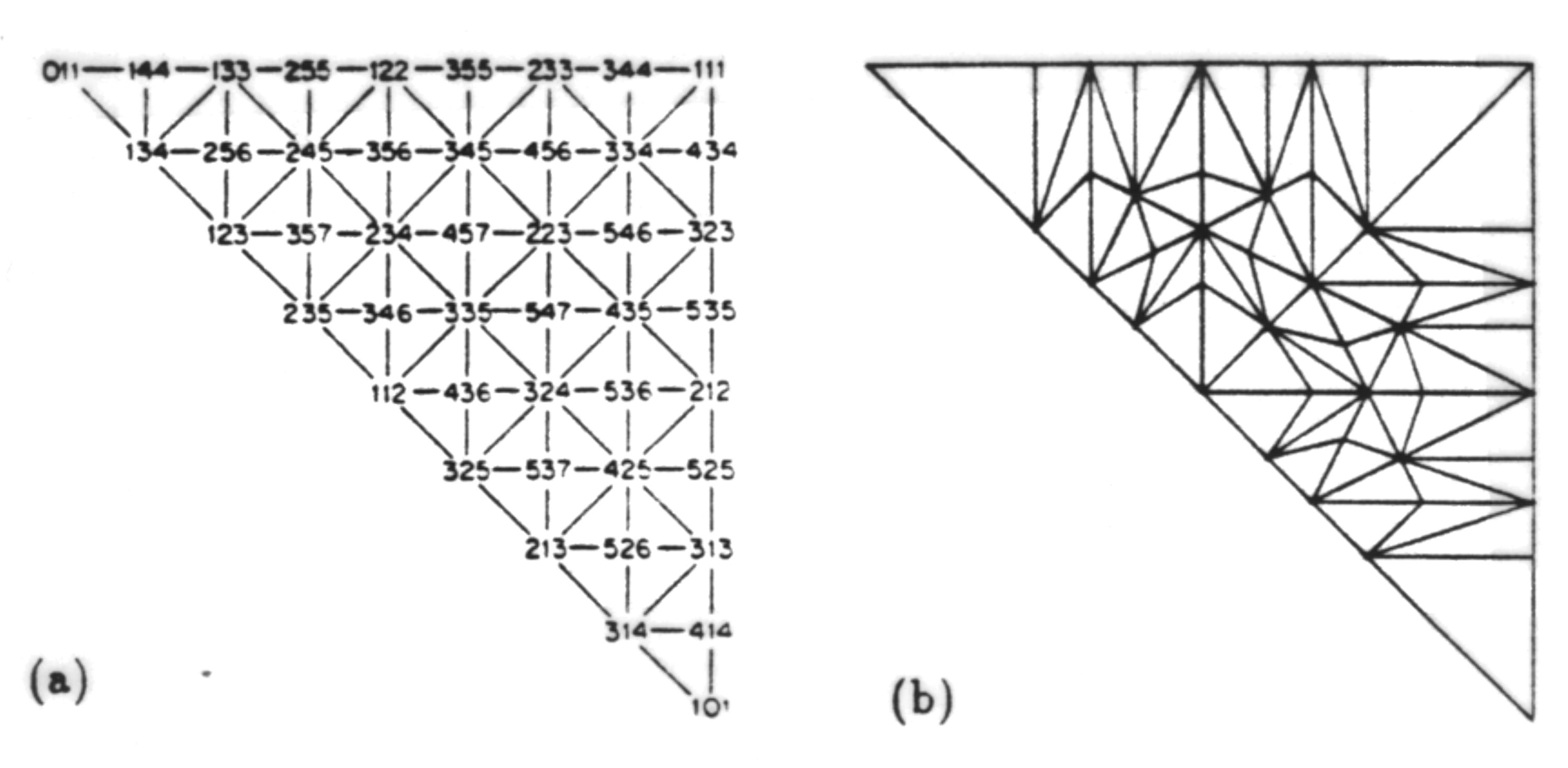}}
	}
\label{fig:tri}
\end{figure}

    Kim and Ostlund \cite{KO} provide a detailed algorithm for 
implementing Farey
approximation on a computer and generalize the idea to solve the problem
of simultaneously approximating two irrationals $(\omega_0, \omega_1)$
by rationals of the form $(p_0 / q, p_1 / q)$\thesisfootnote{
	These are just the sorts of approximations we want; $q$ is 
	the period of our periodic state. },
which they represent as the triple $(p_0, p_1, q)$.  To simplify the 
presentation let us restrict our attention to those vectors for which
$(\omega_0, \omega_1)$ is such that $\omega_0 + \omega_1 \geq 1$;
the other
case is not very different.  The analogs of Farey intervals are 
{\em Farey triangles}, see figure~\ref{fig:trimed}, 
and the act of refinement
again involves adding a point obtained by coordinate-wise addition. 
Even when the vertices of the Farey triangles are viewed as rational
points in ${\bf R}^2$
the 2-d Farey mediant lies on the line connecting its parents so that 
the subdivision into triangles represented in figure~\ref{fig:tri}
reflects a genuine triangular decomposition on the unit square.
Successive subdivisions produce every rational vector, though some
appear twice\thesisfootnote{
	Those vertices in the interior of the triangle 
	\mbox{$(0,1,1)$}, \mbox{$(1,0,1)$}, \mbox{$(1,1,1)$}
	lie on the hypotenuse of two different Farey triangles.
	}.
As in the 1-d Farey approximation scheme, one chooses between 
a right and left daughter at each level of refinement. Irrational
vectors thus have binary addresses.  Kim and Ostlund assert that
the analog of the
golden mean is the vector whose address is $rrrrrrrrr \ldots$;
they call it the {\em spiral mean}.  Its
components are \mbox{$(\tau^{-2}, \tau^{-1})$}, where $\tau$ 
satisfies \mbox{$\tau^3 - \tau - 1 = 0$.} 
    One of the rotation vectors we studied, (1432, 1897) / 2513,
is a an approximation to the spiral mean, and we used the  Farey
triangle algorithm to produce the approximations used in the
sequence of orbits pictured in section~\ref{sec:Hed}.

\section{Lyapunov exponents }

     The Lyapunov exponents displayed in section~{\ref{subsec:Lyap}} 
were found with the algorithm outlined in \cite{BenStr:ex}.
Their method depends on two observations, the first that one can compute
the largest Lyapunov exponent by examining the growth of a vector 
tangent to an orbit, the second that the Lyapunov exponents are constant
on a certain nested family of subspaces of the tangent space. 
To find all the 
exponents one selects a family of linearly independent vectors 
$\nu_0, \nu_1, \ldots, \nu_{2d-1} \in TM_{x_0}$ and  
carries them along the
orbit with the tangent map $DF$.  Unless one makes a fantastically 
improbable choice of initial vectors, each $\nu_i$ will grow with an
exponential rate $\lambda_{max}$,  
\begin{equation}
\lambda_{max} = 
	 \frac{1}{q} \, \log \frac{\left\| DF_{x_0}^q\,\nu_i \right\|}
				 {\left\| \nu_i \right\|},
\label{eqn:lyap}
\end{equation}
equal to the largest Lyapunov exponent.   The 
$\nu_i$ will also become more
and more nearly parallel because
their growth is dominated by that of the eigenvector with the largest 
eigenvalue;
$\dfnu{0}$ will be nearly parallel to this eigenvector.
If we examine those components of $\dfnu{1}$ which are perpendicular to 
$\dfnu{0}$ we should find that they grow 
with a rate given by the next to 
largest
Lyapunov exponent.  Those components of $\dfnu{2}$ 
which are perpendicular
to both $\dfnu{0}$ and $\dfnu{1}$ should grow with a rate given by the
third to largest Lyapunov exponent, and so on.

     In practice the $\dfnu{i}$ are too nearly parallel to 
permit the direct
calculation described above.  Instead one carries out 
the calculation of 
$\dfnu{i}$ in $q$ stages, using the definition of 
$DF^q_{x_0}$, (\ref{eqn:DF}).
Whenever $\dfnu{0}$ gets larger than some modest limit, one performs
a Gram-Schmidt orthogonalization on the vectors,  then normalizes each 
member of the resulting orthogonal collection and keeps a running total
of the logarithms of
the normalization constants.  The Lyapunov exponents are just
\begin{eqnarray*}
	\lambda_i & = & \frac{1}{q} \, \sum_{normalizations} \log n_i,
\end{eqnarray*}
where $n_i$ is a normalization constant for the i{\em th} vector.
We adopted the scheme of \cite{BenStr:ex} only after trying a more
difficult and time consuming  method based on the rate of growth
of the volumes of parallelopipeds. Although this original algorithm
had a pleasing likeness to the definitions of Oseledec's
great paper \cite{Osc}, it gave the same answer as the algorithm
described above, but took quite a bit longer.

%% file: AppendB.tex
\input{appendices/codeSpec/introB}
\input{appendices/codeSpec/overview}
\input{appendices/codeSpec/data}
\input{appendices/codeSpec/algos}

\input{appendices/codeSpec/start}
\input{appendices/codeSpec/traces}
\input{appendices/codeSpec/fatOne}
\input{appendices/codeSpec/fatTwo}

\input{appendices/codeSpec/trunc}

%% file: appendices/codeSpec/introB.tex
\chapter{Converse KAM Methods}
\label{app:spec}

The algorithms used to prove the theorems of section \ref{subsec:results}
have been implemented in the C programming language. 
This appendix descibes the program in some detail.  Section
\ref{sec:overview} gives an overview of a typical computation and
section \ref{sec:types} explains how the basic data: numbers,
intervals, and prisms, are stored in the computer.
Section \ref{sec:detail} carefully describes the crucial algorithms
and serves as an introduction to the parts of the code appearing in 
appendix \ref{app:code}.

%% file: appendices/codeSpec/overview.tex
\section{What the program does}
\label{sec:overview}

   This section expands on the plan for a proof offered in 
section \ref{subsec:plan}.  It first discusses the actual map used,
then gves a more detailed sketch of the computation, ending
with a typical input file and the resulting output.
This section also introduces a convention of typography
and one of nomenclature.
Under the former, bits of text taken directly from computer 
programs will be printed in the {\tt typewriter} typeface.
Under the latter, closely related objects will have similar names.
For efficiency's sake,
I have written two versions of most functions. 
The first, quick and sloppy, is used for exploration.
The second, stately and rigorous, verifies any promising
results suggested by the first.
The quick function usually has some descriptive name, as has 
{\tt bound\_btrace()}, which bounds the trace of the blocks
$\bbeta(\bfx_{i})$.  The rigorous version, {\tt Rbound\_btrace()},
has almost the same name, but for the prefix, {\tt R}, connoting rigor.
A similar convention applies to names of variables; {\tt minLeastLam }
is an approximate value for $\lambda_{-min}$,
the smallest permissible value for the least
eigenvalue of a diagonal block.  The
rigorous estimate of the same number is called {\tt RminLeastLam}.

\subsection{the map}

The program really works with the three-parameter, four-dimensional,
symplectic map,
\begin{eqnarray}
	\bfy' & = & \bfy + \bfJ',  \nonumber \\
	\bfJ'   & = & \bfJ-  \pderiv{V_{abc}}{\bfy}.  \nonumber
\end{eqnarray}
Where
\begin{equation}
	V_{abc}(\bfy)  =  -a \sin (y_0) - b \sin (y_1)
					- c \sin (y_0 + y_1).
\label{eqn:code-map}
\end{equation}
If one takes $a = b = \frac{4\epsilon \pi^2 }{2M_{trig}}, 
\;\; c = \frac{4\epsilon \pi^2}{M_{trig}}$ this map
is conjugate to the trigonometric example via the change of 
coordinates,
\begin{displaymath}
	\bfx =  \frac{\bfy}{2\pi}, 
	\qquad 
	\bfp =  \frac{\bfJ}{2\pi}.
\end{displaymath}
I included the extra parameters because it was easy, and 
left open the possibility of further work.  I used
$\bfy \equiv 2 \pi \bfx$ to avoid having to multiply by $2\pi$
so often.

\subsection{sketch of a computation}

     This section explains what the program does.  First, it 
reads an input file and invoke a host of initialization functions.
These have names like {\tt init$\cdots$()} and do such things as
initialize variables, allocate memory, and copy the input data to 
various output files.  Next, the program chooses the starting point,
$\bfxStar$ and prepares the first, all-encompassing
prism which then becomes the sole member of a linked list of untested prisms.
The rest of the computation is a struggle to get to the end of this list.
It grows shorter whenever the prism-testing algorithm succeeds;
when the program is able to show that none of the points in a 
particular prism could follow $\bfxStar$ in a minimizing state
the successful prism is removed from the list and forgotten.
The list grows longer when the algorithm fails; the offending prism
is divided in two by {\tt refinePrism()} and replaced by the resulting 
pair.

     The program tests a prism in several stages; it begins by 
examining the values of the parameters included in the prism and 
computing $\lambda_{-min}$ and $ \TrM $;
it then invokes a series of prism-testing functions.  The first of these, 
{\tt quick\_try()}, tries to show that the states with
$\bfx_0 = \bfxStar,\; \bfx_{1} = \{\mbox{\em center of the prism}\}$
cannot be minimizing.  If {\tt quick\_try()} fails the prism is
judged hopeless and is immediately halved; if {\tt quick\_try()}
succeeeds the program passes the prism to {\tt try\_Prism()}.
This function does a full, orbit-following, image-bounding test,
but uses only 48-bit, double-precision numbers and does not give
rigorous results.  If {\tt try\_Prism()} succeeds too, then, finally, 
{\tt Rtry\_Prism()} checks the prism rigorously.
Eventually the program either reaches the end of the list, and so 
proves a converse KAM theorem, or founders on a difficult
prism and quits.

\subsection{using the program: a sample}

     The computation which proved $\epsilon_c \leq {\rm 0.0274}$
began when I typed:
\begin{quote}
    {\tt
	converse <trig274.in >\&trig274.out -d30 
    }
\end{quote}
The -d30 sets the maximum {\em depth}; it tells the program 
to quit if it ever fails on a prism which has been subdivided 30 times. 
Other command-line options include:
\begin{description}
    \item[-b {\em filename}] Maintain a backup file.  This
		is esential for long computations; the 
		backup file is updated frequently 
		and contains enough information to
		continue a proof that has been 
		interrupted by some computer disaster.

    \item[-g {\em filename}] Make a graphics file. 
		The program composes a PostScript
		program to draw figures like 
		(\ref{fig:theorem1})-(\ref{fig:theorem3})
		and writes it on {\em filename}.  If
		{\em filename} is the special name,
		{\em off}, then the graphics parts of 
		the program are turned off. 

    \item[-p {\em dp}] Fix the precision used in the rigorous
		parts of the computation to {\em dp} decimal
		places; the example above uses the default, 35.

    \item[-s]  Be stubborn; keep on computing even if some prism
	        cannot be successfully resolved at the maximum depth.
	        This option is good for making pictures and for getting 
	        an idea of how hard a fully successful computation might be.

    \item[-t]   Change the {\em terseness}.  Selecting this option
		makes the program more informative; it prints
		a message whenever it finds a successful prism.
		It also makes the output file much longer, and so 
		I used it only during development of the program.

    \item[-r {\em filename}]  Restore an interrupted computation 
		from a backup file.
\end{description}

The input file, {\tt trig274.in}, looks like:
\singlespace
\begin{tabbing}
\hspace{1\parindent}\= Run on kastor\= with a depth of 30  
\hspace{1\parindent}\= \kill
\>	       \>		\>\em Parameters: \\
\> \tt 0.3085  \> \tt 0.00125  \> $ a_c \;\mbox{\em and}  \; \delt{a} $ \\
\> \tt 0.3085  \> \tt 0.00125  \> $ b_c \;\mbox{\em and}  \; \delt{b} $ \\
\> \tt 0.617   \> \tt 0.0025   \> $ c_c \;\mbox{\em and}  \; \delt{c} $ \\
\\
\> \> \> \em Angles given in units of $2\pi$. \\
\> \tt 1.0	\> \tt 1.0 \> $ \theta_{c,0} \;\mbox{\em and} \; \delt{\theta_0} $\\
\> \tt 1.0	\> \tt 1.0 \> $ \theta_{c,1} \;\mbox{\em and} \; \delt{\theta_1} $\\
\\
\>\tt 0.0274 $<$ epsilon $<$ 0.0276 \\
\>\tt Run on kastor \\
\>\tt May 2nd, 1989  
\end{tabbing}

The parts in the typewriter typeface are copied directly from 
the input file; the parts in italics are additional comments.
The first three lines give the ranges for parameters
$a$, $b$ and $c$. For example, the first line is the pair,
$(a_c, \delt{a})$, which establishes that the initial prism will 
have \mbox{$a_c - \delt{a} \leq a \leq a_c + \delt{a}$.}.
The fifth and sixth lines specify that the prism will have 
$0 \leq \theta_j \leq 2\pi, \; j = 1,2$.  The last few lines 
are comments.  \thesisspace

\ \\
The computation above would yield an output file, 
{\tt trig.out}, looking like:
\singlespace
\begin{quote}
{\tt
apmValidate : null APM value in map.c at line 296. \\
Parameters : \\
a : 3.08500000000000e-01 	 1.25000000000000e-03 \\
b : 3.08500000000000e-01 	 1.25000000000000e-03 \\
c : 6.17000000000000e-01 	 2.50000000000000e-03 \\
\\
Initial region : \\
v[0] : 3.14159265358979e+00 	 3.14159265358979e+00 \\
v[1] : 3.14159265358979e+00 	 3.14159265358979e+00 \\
\\
Comments :\\
\\
0.0274 < epsilon < 0.0276\\
Run on kastor\\
May 2, 1989\\
\\
++++++++++++++++++++++++++++++++++++++++++++++++\\
I find no invariant tori for the range of parameters :\\
0.307250 < a < 0.309750 \\
0.307250 < b < 0.309750 \\
0.614500 < c < 0.619500 \\
\\
Did 322 quick checks, 318 semi-rigorous bounding tries,\\
and 156 rigorous bounding tries. \\
The most deeply refined prism was cut 13 times. \\
The longest semi-rigorous orbit ran for 5 iterations, \\
the longest successful orbit, 4 iterations. \\
Of the 156 successful prisms, 0 fell to the trace criterion, \\
156 to the least eigenvalue test. \\
The best upper bound on the least eigenvalue came from \\
the maxBlam criterion 0.0\% of the time, \\
the minBlam criterion 99.4\% of the time, \\
and from the trace criterion 0.6\% of the time. \\
\\
This investigation took 2697.53 seconds. \\
}
\end{quote}

The first line is an error message from the intialization
phase of the computation, saying that some
variable was not properly allocated;
the program automatically corrects this error. 
The next few lines are
copied directly from the input and the lines after those
give the result:
no tori.  The rest of the file reports details about
the program's performance.
\thesisspace

%% file: appendices/codeSpec/data.tex
\section{Representation of data}
\label{sec:types}

  Here we  explain how data are represented in the program.
This section is fairly technical; it is partly intended as an  
introduction to the program and assumes some knowledge of
C.  Those wishing to avoid
technical details should read only section \ref{sec:apm-intro}, in which
numbers and arbitraty precision arithmetic are discussed.
This leads into a description of {\em intervals} and 
{\em interval arithmetic}, which makes up the next section.
Last, we explain how prisms are represented.

\subsection{numbers and arithmetic}
\label{sec:apm-intro}

     The computations in the rigorous parts of the program use
an arbitrary precision arithmetic library written by Lloyd 
Zussman\thesisfootnote{
	Mr. Zussman's library is licensed under a variant of
	the Free Software Foundation's 
	{\em Gnu EMACS General Public License} and so I am
	obliged to provide a copy of the source code to anyone
	who asks.  Complete source code for my program,
	{\tt converse}, is also available on request.
}.
A desciption of his library and its constituent functions appears
in appendix \ref{app:code}; for now it is enough to know that it
allows one to do arithmetic on numbers represented as finite 
strings of base 10000 ``digits.''  We will call such strings 
{\em APMs}.  Addition, subtraction and multiplication of two 
APMs, say, $x$ and $y$, always yield another number representable
as an APM, but division need not.  The rational number $\frac{x}{y}$
may have an infinite repeating representation in base 10000.
     The division function, {\tt apmDivide()}, deals with this 
problem by allowing the user
to specify the number of decimal places (counting only those
to the right of the decimal point) to which the 
result should be correct.  The special functions,
{\tt apmSin()}, {\tt apmCos()}, and {\tt apmSqrt()}, which I have 
written, use the same stategy.  

     Fixed-precision calculations return a kind of implicit interval.
An answer, $\tilde{a}$, which is accurate to {\em dp} decimal 
places, can be thought of as an interval
guaranteed to contain the true answer, $a$;  
\begin{displaymath}
    \tilde{a} - {\rm 10}^{-dp} \leq a \leq \tilde{a} + {\rm 10}^{-dp}   
\end{displaymath}
The program also uses
functions which do explicit interval arithmetic.   An example is
{\tt Rbd\_sin()} which accepts as its argument an interval,
$[\theta_-, \theta_+] \equiv I_{\theta}$,
and returns an interval, $[s_-, s_+]$, certain to contain $\sin \theta$ for
any $\theta \in I_{\theta}$.
Most of the crucial estimates involve some fixed-precision
calculation and and so the program often uses the variables
\begin{displaymath}
    \mbox{\tt max\_error} =  {\rm 10}^{- dp }, 
\end{displaymath} 
and 
\begin{displaymath}
	\mbox{\tt precision} = dp + \mbox{\tt SAFETY\_DP}.
\end{displaymath}
$dp$ is the number of digits selected with the {\bf -p} option 
and {\tt SAFETY\_DP} is a margin of safety.  All the program's 
intermediate results are calculated to {\tt precision} decimal
places and then, for safety's sake,  regarded as only accurate to
$\pm \, \mbox{\tt max\_error}$.  In the calculations 
summarized in table \ref{tab:results}, $dp = {\rm 35}$ and 
$\mbox{\tt SAFETY\_DP} = {\rm 5}$.

\subsection{intervals and expressions}

The structure representing an interval is 
\begin{typedef}{}
	\member{APM}{ub, lb ; $\}$ Bdd\_apm ;,}
\end{typedef}
\noindent
called a {\em bounded APM}.
The functions {\tt Rbd\_sin()} and 
{\tt Rbd\_cos()} each take one bounded
APM as an argument and return another as the result. 
The only other operations on intervals used by the program are 
addition, subtraction, and multiplication.  This is all  handled 
through two other structures, the {\tt Bapm\_term}, and the 
{\tt Bapm\_expr}.  The former is short for {\em bounded term,} 
the latter
for {\em bounded expression.}  Their full declarations are:
\begin{typedef}{}
    \member{int}{nfactors ;}
    \member{APM}{coef ;}
    \member{Bdd\_apm}{**factors, bound ; $\}$ Bapm\_term ; }
\end{typedef}
\noindent
and
\begin{typedef}{}
    \member{int}{nterms ;}
    \member{APM}{const ;}
    \member{Bdd\_apm}{bound ;}
    \member{Bapm\_term}{*terms ; $\}$ Bapm\_term ; }
\end{typedef}
\noindent
To see the use of these structures, consider trying to 
find a bound on
\begin{displaymath}
	{\rm 2.0} - a \sin(\theta_0) - b \sin( \theta_1 ),
\end{displaymath}
where $a$, $b$, and the $\theta_i$ all belong to intervals.  
One would set up a bounded expression composed of two 
bounded terms:
\begin{displaymath}
    \underbrace{{\rm 2.0}}_{const.}  \;\;
    \underbrace{
    	- \underbrace{a}_{factors} \underbrace{\sin \theta_0}
    }_{Bapm\_term}  \;\;
    \underbrace{
    	- \underbrace{b}_{factors} \underbrace{\sin \theta_1}
    }_{Bapm\_term},
\end{displaymath}
then use {\tt Rbd\_sin()} to set the factors and,
finally, use {\tt Rbd\_expr()} to get bounds on the whole thing.

\subsection{prisms}

     The prisms introduced in section \ref{sec:prism-bondage} are
the fundamental objects of the program; they are stored in
\begin{typedef}{RPrsm}
	\member{int}{in\_torus, n\_cuts ;}
	\member{APM}{*matrix ; }
	\member{char}{*cuts[7] ;}
	\member{Rxtnd\_pt}{*center ;}
	\member{struct Rprsm}{*next ; $\}$ RPrism ; }
\end{typedef}
\noindent
The integer {\tt in\_torus} has one of the values
{\tt NO\_TORI}, {\tt UNTRIED}, {\tt MAYBE}, {\tt ACTIVE},
or {\tt SYMMTRC} according to whether it definitely
does not include points from a minimizing state, has
not yet been tested, has been inconclusively tested, 
is under active consideration or may be disregarded on account
of symmetry.  The integer {\tt n\_cuts} tells how many
subdivisions it took to make this prism and the character
strings {\tt cuts[\ ]} explain how to produce this prism from
the initial, big prism. {\tt center} and {\tt matrix} are
the center point and defining matrix of the prism;
{\tt center} is an example of an {\em extended phase point};
it has seven coordiates in all, three for the parameters 
and two for each of the delay embedded coordinates.
The pointer {\tt next} gives the next Rprism on the
list. 

%% file: appendices/codeSpec/algos.tex
\section{Algorithms}
\label{sec:detail}

     Here we explain and verify the crucial algorithms. 
In the first part of the section we will establish the correctness 
of {\tt apmSin()}, {\tt apmCos()} and {\tt apmSqrt()}, functions 
which we approximate with polynomials gotten by truncating 
Taylor series.  
Next we check the algorithms which set the bounds $\lambda_{-min}$ 
and $\TrM$, then we turn to the computations used to compute 
$\lb \lamM{\bfd_j}$.  In the last part of the section we
examine the prism-bounding algorithms. 

\subsection{special functions}

\subsubsection{sine and cosine}
The real computational work is done by two functions,
{\tt reducedSin()} and {\tt reducedCos()}, which compute the sine 
and cosine of an angle from the interval 
\mbox{$ I_0 \equiv  [0, \frac{\pi}{4}]$}.
These functions and the relations
\begin{eqnarray*}
    \sin(\theta \pm \frac{\pi}{2}) = \pm \cos( \theta ), & &
    \sin (-\theta) = -\sin(\theta ),  \\
    \cos(\theta \pm \frac{\pi}{2}) = \mp \sin( \theta ), & &
    \cos (-\theta) = \cos(\theta ), 
\end{eqnarray*}
allow us to calculate the sine and cosine of any angle.
As mentioned in section \ref{sec:apm-intro}, we must set
{\em dp}, the the number of correct digits we want in the answer.
{\tt setTrigDp({\em dp})} does this; it also chooses the 
order of the Taylor approximation and picks the  number of decimal
places, {\tt trig\_dp}, to which intermediate results are 
calculated. 
To prove that all this works we will estimate the error made 
by {\tt reducedSin()}\thesisfootnote{
	The analysis of {\tt reducedCos()} is much the same.}, 
leaving undetermined {\tt trig\_dp} and the number of terms in 
the polynomials, {\tt trig\_terms}.  We will 
then show how to choose these two and how to reduce an arbitrary
angle $\theta$ to one lying in $[0,\frac{\pi}{4}]$.

     The form of the approximation is
\begin{eqnarray}
    \mbox{\tt reducedSin($\theta$ )}  \approx P_N(\theta)
    & \equiv & 
	\frac{1} {(2N + 1)!} 
	\sum_{j=0}^N 
	\theta^{2j+1} \, (-1)^j \frac{(2N+1)!}{(2j+1)!}
    \nonumber \\
    & \approx & \bigStrut
	\frac{1}{\mbox{\tt sinFactrl}}
		\sum_{j=0}^{N}
		    {\mbox{\tt sinCoef[$j$]}}\, \theta^{2j+1}
\label{eqn:reducedSin}
\end{eqnarray}
where the second line substitutes names used in the code.  
Let us consider an angle, $\theta \in [0, \frac{\pi}{4}]$, 
which is approximately represented by an APM, $\tilde{\theta}$.

\vskip \abovedisplayskip
\noindent
{\bf Proposition} \hspace{1\parindent} {\em If $\Tthet$ is such that 
\mbox{$|\theta - \tilde{\theta}| \leq \epsilon < 1$}, then }
\begin{equation}
    |\sin \theta - P_N( \Tthet )| 
    \leq 
    \epsilon + \frac{\theta^{2N+3}}{(2N+3)!}.
\label{eqn:Taylor}
\end{equation}

\noindent
{\bf Proof} \hspace{1\parindent}  By straightforward computation, 
\begin{eqnarray*}
    |\sin \theta - P_N(\Tthet)|  
    	& \leq &
	|\sin \theta - \sin \Tthet| + |\sin \Tthet - P_N(\Tthet) |, \\
    \bigStrut
    & \leq & |\theta - \Tthet | + \left| 
				\sum_{j=1}^{N} 
				(-1)^j
				\frac{\theta^{2j+1}}{(2j+1)!} 
		  	  \right|, \\
    \bigStrut
    & \leq & \epsilon  +  \frac{\theta^{2N+3}}{(2N+3)!}. 
\end{eqnarray*}

Evaluating long power series like (\ref{eqn:reducedSin}) can take
immense amounts of computer time and memory; if the string of 
digits making up $\Tthet$ has length $\ell$ then the one 
representing $\Tthet^n$ will have length $\approx n\ell$.
So, in the interest of computational speed,
{\tt reducedSin()} truncates some intermediate expressions.  
What it really calculates is a sequence of approximations
to certain polynomials. In the equations below, $[x]_n$ 
is the number given by the truncating $x$ after $n$ 
places to the right of the decimal point, and $tdp$ is
short for {\tt trig\_dp}.
\begin{eqnarray*}
    \bar{S}_0  	& = &  (-1)^N, \\
    \bar{S}_1  	& = & \left[ \Tthet^2 \bar{S}_0 +
			(2N+1)(2N) (-1)^{N-1}\right]_{tdp}, \\
    		& \approx & \Tthet^2(-1)^N + 
			(2N+1)(2N)(-1)^{N-1},  \bigStrut\\
		& \vdots & \\
    \bar{S}_N	& = & \left[ \Tthet^2 \bar{S}_{N-1} + 
				(2N+1)! \right]_{tdp}, \\
		& \approx & \sum_{j=0}^{N}
			    \Tthet^{2j} \, (-1)^j 
			    \frac{(2N+1)!}{(2j+1)!} \bigStrut
\end{eqnarray*}
and, finally, 
\begin{equation}
    \mbox{\tt reducedSin($\Tthet$)} 
	 \equiv  \frac{\Tthet \bar{S}_N}{(2N+1)!} 
		   \label{eqn:reduce-def} 
	 \approx  P_N(\Tthet) 
\end{equation}

Let us consider the additional error introduced by truncation.
Use $S_j$ to denote the exact value of the polynomial 
approximated by $\bar{S}_j$. Then $\bar{S}_0 = S_0$ and
so $S_1$ lies in an interval,
\begin{displaymath}
    \bar{S}_1 - \delta_1 < S_1 < \bar{S}_1 + \delta_1,
\end{displaymath}
with $\delta_1 = {\rm 10}^{-tdp}$. 
Since $S_2 = \Tthet^2 S_1 + C$, where $C$ is a constant, we may
be sure that $S_2$  is in the interval
\begin{eqnarray*}
    [\Tthet^2(\bar{S}_1 - \delta_1) + C, \;
    	      \Tthet^2(\bar{S}_1 + \delta_1) + C] 
    & \subset&  \left[
    		    (\Tthet^2 \bar{S}_1 + C) - \delta_1, \;
    		    (\Tthet^2 \bar{S}_1 + C) + \delta_1
		\right].
\end{eqnarray*}
After truncation we get
\begin{eqnarray*}
    \bar{S}_2 - \delta_2 & < & S_2 < \bar{S}_2 + \delta_2
\end{eqnarray*}
with $\delta_2 = 2\delta_1 $  and after $N$
such steps we are left with an error, $\delta_N = N \,{\rm 10}^{-tdp}$.
Combining this with equations (\ref{eqn:Taylor}) and 
(\ref{eqn:reduce-def}) we get 
\begin{equation}
    |\mbox{\tt reducedSin($\Tthet$)} - \sin{\theta}| \leq
	|\Tthet - \theta| + 
	\frac{N\delta_1}{(2N+1)!} +
	\frac{|\theta|^{2N+3}}{(2N+3)!}
\label{eqn:trig-error}
\end{equation}
     The only unknown quantity here is the
difference between $\theta $ and its APM representation $\Tthet$.
Suppose we can arrange for this to be at least as small as
${\rm 10}^{-tdp}$.
To ensure $dp$ decimal places of accuracy in our answer we
need only choose  $N$ large enough that
$\frac{1}{(2N+3)!} < {\rm 10}^{-(dp+2)}$ and then choose 
{\tt trig\_dp} so large that 
$N\delta_1 \leq {\rm 10}^{-(dp+2)}$ too.

     If we want the sine or cosine of an angle which lies outside 
the interval $I_0$, we must relate it to some calculation
that we can do with the reduced functions.  The program contains a 
very accurate representation\thesisfootnote{
     The current implementation has one good to 
     45 decimal places, but it would be easy to add more.}
of $\pi$, so
it can just subtract the appropriate number of multiples of
$\frac{\pi}{2}$ and, perhaps, reflect about the origin. For very
large angles, the reduction process may lose  so much
precision as to preclude a calculation to the specified accuracy.
In that case the program writes an error message and calculates
the best answer it can.

\subsubsection{square root} 
     The square root function {\tt apmSqrt()} is much simpler.
It takes an argument, $x$, and uses Newton's method to solve 
the equation $y^2 -x = 0$.  Suppose we want $dp$ decimal places of
accuracy in the answer; define $dp+ = dp + 2$.  
{\tt apmSqrt()} recursively calculates a sequence 
$y_j \approx \sqrt{x}$ with 
\begin{eqnarray}
\openup 1.5\jot
    y_0 & = & x  \nonumber \\
    y_{j+1} & =  & \left[ \half( y_j + 
		    	\left[ \frac{x}{y_j}
		    	\right]_{dp+})
		   \right]_{dp+} 
\label{eqn:apmSqrt}
\end{eqnarray}
After the first few steps, the $y_j$ decrease monotonically 
and so we may  write $y_j = \sqrt{x} + r_j$;
the error term, $r_j$, is a small, positive number.
Equation (\ref{eqn:apmSqrt}) then yields the following extremely
conservative estimate:
\begin{eqnarray}
    r_{j+1} & = & y_{j+1} - \sqrt{x}, \nonumber  \\
    \bigStrut
	    & = & \left[ \half(\sqrt{x} + r_j  + 
			\left[ \frac{x}{\sqrt{x} + r_j}
			\right]_{dp+} 
		   \right]_{dp+}  -\sqrt{x}, \nonumber  \\
	     \medStrut
	     & \leq &  (\frac{r_j}{2} + 
			\sqrt{x} + 2\epsilon_{dp+})  - 
			\sqrt{x}, \nonumber \\
	     \medStrut
	     & \leq & \frac{r_j}{2} + 2\epsilon_{dp+} 
\label{eqn:Newt-err}
\end{eqnarray}
where $\epsilon_{dp+} = {\rm 10}^{-dp+}$ is the 
inevitable truncation error.  
If $r_j < \sqrt{x}$, Newton's method actually gives 
$r_{j+1} \sim \frac{r_j^2}{\sqrt{x}}$, but (\ref{eqn:Newt-err}) 
will be good 
enough for us.  It tells us that we must continue computing 
until the difference,
\begin{displaymath}
    y_{j-1} - y_j \; = \; r_{j-1} - r_j 
		  \; > \; \frac{r_j}{2} - 2\epsilon_{dp+},
\end{displaymath}
is less than  ${\rm 10}^{-(dp+1)}$; the last $y_j$ will be 
the answer.

%% file: appendices/codeSpec/start.tex
\subsection{uniform cones and the starting point}

    This section explains how the program evaluates the constants
$\TrM$, $\TrP$, $\lambda_{-min}$ and $\lambda_{-max}$; it also explains
how to get a good value for the starting point $\bfxStar$.
The main technical problem is the correct evaluation of the constants
\begin{displaymath}
    B = \ub \lamP{\bbeta} \qquad \mbox{and} \qquad T = \ub \Tr{\bbeta};
\end{displaymath}
these, together with equations (\ref{eqn:trace-cone}) and 
(\ref{eqn:lambda-cone}), determine everything else.  Finding either
$B$ or $T$ is a matter of maximizing a function on 
$[0,1] \times [0,1] \times \{\mbox{\emph{parameters}}\}$, so it
is enough to explain how to find one of them, say $T$.

    When the program seeks $T$ it sets $a$, $b$ and $c$ to their values  
at the center of the intial prism, then uses Newton's method to 
find a zero of the gradient of $\Tr{\bbeta}$.
For the computations presented
in section \ref{subsec:results}, the search began at 
$(\frac{\pi}{2}, \frac{\pi}{2})$ and continued until it
reached a point $\bfx_T$ such that 
\begin{displaymath}
    \left| \displaystyle 
	\frac{ \partial \Tr{\bbeta(\bfx_T)}}{\partial \bfx} 
    \right|
    \;  < \; (|a_c| + |b_c| + |c_c|) \, \epsilon_{newt},
\end{displaymath}
where $\epsilon_{newt}$ is a small constant. In the code, the 
search is done with ordinary double precision arithmetic and
$\epsilon_{newt}$ is called {\tt NEWT\_TOL} and is equal to ${\rm 10}^{-9}$.
The $\bfx_T$ it finds is very close to the true maximum, and so
a suitable estimate is
\begin{displaymath}
    T \; = \; \Tr{\bbeta(\bfx_T)} + (a_c + b_c + 2c_c) {\rm 10}^{-6} 
		+ (\delt{a} + \delt{b} + 2\delt{c})
\end{displaymath}
where the last term is included to allow for the variation in $a$, $b$
and $c$ over the prism.   The point $\bfx_T$ found by this technique
is the natural starting point for an estimate based on 
Herman's trace condition, so I call it {\em Herman's starting point}.

     The estimate for $B$ works much the same way; a Newton's
method search gives an approximate value for, $\bfx_B$, the 
position where $\max \lamP{\bbeta}$ is attained. $B$ is then
calculated according to  
\begin{eqnarray*}
    B  & = & \lamP{\bbeta(\bfx_B)} + ( a_c + b_c + 2c_c) {\rm 10}^{-6}
		+ ( \delt{a} + \delt{b} + 2 \delt{c} ) 
\end{eqnarray*}
After calculating $B$, the program sets up the starting point, $\bfxStar$,
also called the {\em least-lambda starting point}.  This point is essentially
the same as $\bfx_B$, but is explicitly guaranteed to lie on the
line $x_0 = x_1$ so that the calculation can exploit symmetry, as
explained in section \ref{subsec:symmetry}.

%% file: appendices/codeSpec/traces.tex
\subsection{bounding traces and eigenvalues}

This section explains how the program takes a prism, $P$, and
evaluates the bounds
\begin{eqnarray*}
    \ub_{(\bveps,\bfu,\bfv) \in S} \lamM{\bbeta}, \\
    \ub_{(\bveps,\bfu,\bfv) \in S} \lamP{\bbeta}, \\
    \ub_{(\bveps,\bfu,\bfv) \in S} \Tr{\bbeta}, 
\end{eqnarray*}
where $\bveps \in {\bf R}^3$ stands for the triple of parameters,
$(a,b,c)$.  These are the basic ingredients of the main suite of 
estimates, (\ref{eqn:sineq1}) -- (\ref{eqn:sineq3}).
     Recall that the prism is determined by its center, 
$(\bveps_c, \bfu_c, \bfv_c)$, and by the matrix which maps 
the hypercube, $Q^7$, into the extended phase space. A point 
$\bfeta \in Q^7$ has an image given by
\begin{equation}
    \left[ \begin{array}{l}
		a(\bfeta) \\ b(\bfeta) \\ c(\bfeta) \\
		u_0(\bfeta) \\ u_1(\bfeta) \\ 
		v_0(\bfeta) \\ v_1(\bfeta)
    \end{array} \right]
    \; = \;
    \left[ \begin{array}{l}
		a_c \\ b_c \\ c_c \\ u_{c,0} \\
		u_{c,1} \\ v_{c,0} \\ v_{c,1}
    \end{array} \right]
    \; + \;
    \left[ \begin{array}{ccccccc}
		\delt{a} & 0 & 0 &  & \cdots & & 0 \\
		0 & \delt{b} & 0 &  & \cdots & & 0 \\
		0 & 0 & \delt{c} &  & \cdots & & 0 \\
		\\
		\vdots & \vdots & \vdots & & \ddots & & \\
		\\
		p_{7\,1} & p_{7\,2} & p_{7\,3} 
		& & \cdots & & p_{7\,7} 
    \end{array} \right] 
    \left[ \begin{array}{l}
		\eta_1 \\ \eta_2 \\ \eta_3 \\ 
		\eta_4 \\ \eta_5 \\ \eta_6 \\ \eta_7
    \end{array} \right].
\label{eqn:extendo-prism}
\end{equation}
From this it is easy to show that any $(\bveps, \bfu, \bfv) \in S$
has
\begin{displaymath}
	|v_0 - v_{c,0}| \; \leq \; 
	\sum_{j=1}^7 |p_{6\,j}|
	\qquad {\rm and} \qquad
	|v_1 - v_{c,1}| \; \leq \; 
	\sum_{j=1}^7 |p_{7\,j}|.
\end{displaymath}
Once we have found bounds on the components of $\bfv$, 
we can invoke
{\tt Rbd\_sin()} to get bounds on the functions $\sin(v_0)$,
$\sin(v_1)$ and $\sin(v_0 + v_1)$, then combine those with
$\delt{a}$, $\delt{b}$ and $\delt{c}$ to obtain bounds on 
the expressions appearing in the trace and eigenvalues of
$\bbeta$.

     In the program, all this is done with the  {\tt Bapm\_expr}
machinery described in section \ref{sec:apm-intro}.  The expressions
$a \sin(v_0)$, $b \sin(v_1)$ and $c \sin(v_0 + v_1)$ arise so often
that they are given their own names: {\tt Ra\_sin}, {\tt Rb\_sin}
and {\tt Rc\_sin}; their values are set by 
{\tt Rglobal\_bounds({\em priz})}. In terms of these, 
the estimates we need are:
\begin{eqnarray*}
    \ub_{S} \Tr{\bbeta}
    & = & 4.0 \,+\, {\tt Ra\_sin.bound.ub} \, + \,
		{\tt Rb\_sin.bound.ub} \, + \,
		2\,{\tt Rc\_sin.bound.ub} \\
    \ub_{S} \lamM{\bbeta}
    & = & \half \left\{ 
		\ub \Tr{\bbeta} - \lb \sqrt{ {\tt discrim.lb}}
		\right\}, \\
    \ub_{S} \lamP{\bbeta}
    & = & \half \left\{ 
		\ub \Tr{\bbeta} + \lb \sqrt{ {\tt discrim.ub}}
		\right\}
\end{eqnarray*}
where {\tt discrim} is a bounded APM containing estimates 
over $S$ of the quantity
\begin{equation}
    \left(a\sin(v_0) + b\sin(v_1) \right)^2 + 4 c^2 \sin^2(v_0+v_1).
\label{eqn:discrim-def}
\end{equation}
Note how, in every estimate described above,
we allow each of the terms
$a\sin(v_0) \cdots$ to vary independently; the bounds
we obtain are almost certainly too conservative.

%% file: appendices/codeSpec/fatOne.tex
\subsection{bounding the images of prisms}

     The bulk of the computation is devoted to the kind of prism-bounding 
calculations described in section \ref{sec:prism-bondage}.  In this section
we will see how the program takes a prism in the extended phase space, 
$S = (\bfx_c, P)$, and constructs another, $S' = (\bfx_c', P')$, 
guaranteed to contain $G(S)$.  The computation of $\bfx_c'$ is easy;
$\bfx_c' \approx G(\bfx_c)$ where
\begin{eqnarray}
    G(a,b,b,\bfu,\bfv) \;  \equiv \; (a',b',c',\bfu',\bfv') 
			& = & (a,b,c,\bfu',\bfv'), \nonumber \\	
		\bfu' 	& = & \bfv, \nonumber \\
		\bfv'	& = & 2\bfv - \bfu - \pderiv{V_{abc}(\bfv)}{\bfx}. 
\label{eqn:Gabc-def}
\end{eqnarray}
Although only $\bfv'$ involves any real computation, and so only 
it introduces any error,  we will find it useful
to assign a somewhat larger uncertainty, $\delta_{c}$, to both $\bfu'$ 
and $\bfv'$.  
The computation of $P'$  is much more difficult;
the work falls into two parts: setting up the matrix $A$ and evaluating 
the numbers,
\begin{eqnarray}
w_j  & = & 
	\ub  |[A^{-1}(G(x_c) - {x_c}')]_j| + \ub_{x \in S} 
	      \sum_{k=1}^{7} \left|
				[A^{-1} \circ 
				 DG_x \circ
				 P]_{jk}
	      \right|, \nonumber \\
     & \leq & \rowSum{A^{-1}}{j} \delta_c + \ub_{x \in S} 
	      \sum_{k=1}^{7} \left|
				[A^{-1} \circ 
				 DG_x \circ
				 P]_{jk}
	      \right|,
\label{eqn:Wabc-def}
\end{eqnarray}
The second term, which involves bounds over $\bfx \in S$, 
will be the hard part.
As was mentioned in section \ref{sec:A-choice}, the program uses two
schemes to prepare $A$.  The first, the fixed-form scheme, is specially
suited to prisms with zero volume.  Since all the
prisms on the linked list are of the form
\begin{displaymath}
	\{\mbox{\em parameters}\} 
	\times \{ \bfxStar \} \times \{ \mbox{\em possible successors}\},
\end{displaymath}
all are singular. Accordingly, the fixed-form scheme is always used on the
first step of a round of prism-bounding.  Since the first image is 
non-singular by construction,  the second and subsequent iterates employ 
a different, more
accurate scheme, the column-rotor.  This section describes both schemes
and verifies that they are correctly implemented.

     Most of the work will come in showing that the $w_j$ are calculated
properly, a task simplified by the following definitions and
proposition.

\vskip \abovedisplayskip
\noindent
{\bf Definition} \hskip 1\parindent For any real, $m \times n$,
matrix $A$, define 
\begin{displaymath}
   \rowSum{A}{k} \equiv \sum_{j=1}^{n} |a_{k\,j}|,
\end{displaymath}
the {\em k-th row sum of $A$}, and
\begin{displaymath}
     \matSum{A}
    \equiv \sum_{k=1}^{m} \sum_{j=1}^{n} |a_{k\,j}| 
    = \sum_{k=1}^{m} \rowSum{A}{k}
\end{displaymath}

\noindent
{\bf Proposition} \hskip 1\parindent {\em For any real, 
$m \times n$ matrix $A$ and real, $n \times \l$ matrix $B$, 
the product $C = AB$ satisfies }
\begin{equation}
    \rowSum{C}{k}  \leq  \rowSum{A}{k} \matSum{B} 
    \qquad {\rm and } \qquad
    \matSum{C}    \leq  \matSum{A} \matSum{B} 
\label{eqn:prop}
\end{equation}
{\bf Proof} \hskip 1\parindent By direct calculation:
\begin{eqnarray*}
    \rowSum{C}{k} = \sum_{j=1}^{l} |c_{k\,j}|
    & = &  \sum_{j=1}^{l} \, 
		\left| 
			\sum_{i=1}^{n} a_{k\,i}b_{i\,j}
		\right|, \\
    & \leq & \sum_{j=1}^{l} \; 
		\sum_{i=1}^{n} |a_{k\,i}|\,|b_{i\,j}|, \\
    & \leq & \sum_{i=1}^{n} |a_{k\,i}|\, \rowSum{B}{i}, \\
    & \leq & \sum_{i=1}^{n} |a_{k\,i}| \matSum{B} 
		= \rowSum{A}{k}  \matSum{B}.
\end{eqnarray*}
Then, using the first part of (\ref{eqn:prop}), one finds
\begin{displaymath}
    \matSum{C} = \sum_{k=1}^{m} \rowSum{C}{k} 
    \; \leq \; \sum_{k=1}^{m} \rowSum{A}{k} \matSum{B} = 
    \matSum{A} \matSum{B}.
\end{displaymath}
It also follows from the definitions that
\begin{displaymath}
    \rowSum{( A + B )}{k} \leq \rowSum{A}{k} + \rowSum{B}{k}.
\end{displaymath}

     We will use a block-matrix representation for $DG$, the derivative
of the map;
\begin{equation}
    DG = \left[ \begin{array}{crc}
		    \id & 0 &  \hspace{1em} 0 \\
		    0 & 0 & \hspace{1em} \id \\
		    \bgamma & -\id & \hspace{1em} \bbeta 
	\end{array} \right],
\label{eqn:DG-def}
\end{equation}
where
\begin{displaymath}
    \bbeta( \bfv ) = 
	\left[ \begin{array}{cc}
	    2 - a\sin(v_0) - c\sin(v_0+v_1) & - c\sin(v_0+v_1) \\
	    - c\sin(v_0+v_1) & 2 - b\sin(v_1) - c\sin(v_0+v_1) 
	\end{array} \right]
\end{displaymath}
and
\begin{displaymath}
    \bgamma( \bfv ) =
	\left[ \begin{array}{ccc}
	    \cos(v_0) & 0 & \cos(v_0 + v_1) \\
	    0 & \cos(v_1) & \cos(v_0 + v_1)
	\end{array} \right].
\end{displaymath}
It will also prove convenient to have block forms for the matrix
$P$ and to build a column vector, $\bfw$, out of the $w_j$.
\begin{equation}
    P \equiv 
	\left[ \begin{array}{ccc}
		\Ppp & 0 & 0 \\
		\Pup & \Puu & \Puv \\
		\Pvp & \Pvu & \Pvv 
	\end{array} \right]
    \qquad {\rm and} \qquad
    \bfw \equiv 
	\left[ \begin{array}{c}
		\bfw_p \\ \bfw_u \\ \bfw_v
	\end{array} \right],
\label{eqn:PW-def}
\end{equation}
where $\Ppp$ is 3 \by 3, $\Pup$ and $\Pvp$ are 3 \by 2, and the rest of
the blocks are 2 \by 2.  The elements of $\bfw$ are:
\begin{displaymath}
    w_p = \left[ 
	      \begin{array}{c} w_1 \\ w_2 \\ w_3 \end{array} 
	  \right], \;
    w_u = \left[ \begin{array}{c} w_4 \\ w_5 \end{array} \right]
    \; {\rm and} \qquad 
    w_v = \left[ \begin{array}{c} w_6 \\ w_7 \end{array} \right].
\end{displaymath}

\subsubsection{the fixed-form fattener}

     When using this scheme we force the matrix $A$ to be of the form
\begin{equation}
    A = \left[ \begin{array}{ccc}
		\App & 0 & 0 \\
		\Aup & 0 & \Auv \\
		\Avp & \Avu & \Avv
    \end{array} \right].
\label{eqn:ffA}
\end{equation}
The explicit forms of the blocks will be chosen to simplify 
the calculation of 
the $w_j$. Given (\ref{eqn:ffA}) one can get a formula for $A^{-1}$
in terms of the blocks and their inverses:
\begin{eqnarray}
    A^{-1} & = & \left[ \begin{array}{ccc}
			\AppInv & 0 & 0 \\
			0 & -\AvuInv \Avv \AuvInv & \AvuInv \\
			0 & \AuvInv & 0
		 \end{array} \right]
		 \left[ \begin{array}{ccc}
			\id & 0 & \hspace{1em} 0 \\
			-\Aup \AppInv & \id & \hspace{1em} 0  \\
			-\Avp \AppInv & 0 & \hspace{1em} \id 
		 \end{array} \right] \nonumber \\
	  & = &  \left[ \begin{array}{ccc}
			\AppInv & 0 & 0 \\
			\left\{\begin{array}{l}
			    \AvuInv \Avv \AuvInv \Aup \AppInv \\
			    -\AvuInv \Avp \AppInv 
			\end{array} \right\}
			& -\AvuInv \Avv \AuvInv & \AvuInv \\
			-\AuvInv \Aup \AppInv & \AuvInv & 0
		\end{array} \right] 
\label{eqn:ffA-inv}
\end{eqnarray}
Taking $\App = \Ppp$ and using (\ref{eqn:ffA-inv}), (\ref{eqn:PW-def})
and (\ref{eqn:DG-def}), we get
$ A^{-1}\circ DG \circ P = $
\begin{equation}
    \left[ \begin{array}{ccc}
	\id & 0 & 0 \\
	\left\{
	\begin{array}{l}
	    \AvuInv ( \bgamma \Ppp - \Pup) \\
	    + \AvuInv ( \bbeta \Pvp - \Avp) \\
	    + \AvuInv \Avv \AuvInv( \Avp - \Pup )
	\end{array} \right\} &
	\left\{ \begin{array}{c}
		\AvuInv \bbeta \Pvu - \\
		\AvuInv \Avv \AuvInv \Pvu 
	\end{array} \right\} &
	\left\{
	\begin{array}{c}
	    \AvuInv( \bbeta \Pvv - \Puv ) \\
	    - \AvuInv \Avv \AuvInv \Pvv 
	\end{array} \right\} \\
	\AuvInv( \Pvp - \Aup) & \AuvInv \Pvu & \AuvInv \Pvv
    \end{array} \right].
\label{eqn:monster}
\end{equation}
When computing the $w_j$ we must allow the matrices $\bgamma$
and $\bbeta$, which depend on $a$, $b$, $c$ anf $\bfv$ to vary
over $S$.  All the other blocks, those in $A$ and those in $S$,
are constant.  The form of (\ref{eqn:monster}) suggests the following
choices for the blocks of $A$:
\begin{eqnarray}
    \App & = & \Ppp, \nonumber \\
    \Aup & = & \Pvp, \nonumber \\
    \Avp & = & \bgamma_c \Ppp - \Pup + \bbeta_c \Pvp, \nonumber \\
    \Auv & = & \Pvu + \Pvv, \nonumber \\
    \Avu & = & \bbeta_c ( \Pvu + \Pvv ), \nonumber \\
    \Avv & = & \bbeta_c \Pvv - \Puv, 
\label{eqn:A-blocks}
\end{eqnarray}
where $\bbeta_c$ and $\bgamma_c$ are the values of 
$\bbeta$ and $\bgamma$ at the prism's center.  Note that the
entries in the blocks making up $P$ are exactly represented sas APMs;
so are their sums, products, and differences.  Thus
$\Auv$, $\Aup$ and $\App$ are exact; the other blocks of
$A$, which involve the evaluation of special functions, are
uncertain to the extent that the values of the special functions
are.

The choices (\ref{eqn:A-blocks}) immediately determine most 
of the $w_j$;
    the row sums contributing to $\bfw_p$ are automatically 
equal to one and, unless $\Auv$ is singular,
$\bfw_v =  \scriptstyle
\left[ \begin{array}{c} \scriptstyle 1 \\ 
			\scriptstyle 1 
	\end{array} \right]$.
The program checks the invertibility of $\Auv$ by evaluating its
determinant, an exact calculation.  If ${\rm det}[\Auv]$ were 
to be zero the program would
write an error message and halt; this has never actually happened.
The remaining row sums, those contributing to $\bfw_u$, are
\begin{eqnarray}
\ub \rowSum{A^{-1} \circ DG_x \circ P}{k} & = &
    \ub \left\{ \begin{array}{l}
    	\left[ 
	    \AvuInv( \bgamma - \bgamma_c ) \Ppp
	    + \AvuInv (\bbeta -\bbeta_c) \Pvp 
	\right]_{j\,\star} +  \\
    	\left[ 
	    \AvuInv \bbeta \Pvu 
	    + \AvuInv (\bbeta -\bbeta_c) \Pvv 
	\right]_{j\,\star} 
    \end{array} \right\} \nonumber \\
    & \leq & \rowSum{\AvuInv}{j} \; 
	\ub \left\{ \begin{array}{l}
	     \matSum{{ ( \bgamma - \bgamma_c ) \Ppp
	    + (\bbeta -\bbeta_c) \Pvp }} + \\
	    \matSum{ \bbeta \Pvu 
			+ (\bbeta -\bbeta_c) \Pvv } 
	\end{array} \right\},
    \bigStrut \nonumber \\
    & \leq &
	\rowSum{\AvuInv}{j} 
	    \left\{ \begin{array}{l}
		    \ub(\matSum{\bgamma - \bgamma_c}) \matSum{\Ppp} + \\
		    \ub(\matSum{\bbeta}) \matSum{\Pvu} +  \\
		    \ub(\matSum{\bbeta - \bbeta_c}) 
		    (\matSum{\Pvp} + \matSum{\Pvv} )
	    \end{array} \right\} \hugeStrut
\label{eqn:ff-ests}
\end{eqnarray}
where $k = j+3$, $j = 1,2$ and all upper bounds are taken over $\bfx \in S$.
Out of all the numbers appearing in (\ref{eqn:ff-ests}), only
$\rowSum{\AvuInv}{j}$ and the upper bounds on $\matSum{\beta}$, 
$\matSum{\bbeta -\bbeta_c}$
and $\matSum{\bgamma - \bgamma_c}$ cannot be calculated exactly;
the first can be estimated to any desired precision with the APM
library, the rest are handled with the {\tt Bapm\_term, Bapm\_expr} 
machinery.

%% file: appendices/codeSpec/fatTwo.tex
\subsubsection{the column-rotor scheme}
    This technique fattens matrices $A \approx DG_{x_c} \circ P$, where 
$DG$ and $P$ are as in equations (\ref{eqn:DG-def}) and (\ref{eqn:PW-def}).
Such $A$'s have almost the same form as (\ref{eqn:ffA}), but they have
non-vanishing $\Auu$ blocks.  The method's name comes from the way it
tries to ensure that $A$ is non-singular; it rotates parts of columns
4-7 with respect to each other so as to guarantee that they are not parallel.
For example, the function {\tt Rsubspace\_rot()}, which performs the
rotations, begins by finding the angle between the two, 2-d column
vectors enclosed in braces in the matrix below.
\begin{displaymath}
	\left[ \begin{array}{ccccccc}
	    a_{1\,1} & a_{1\,2} & a_{1\,3} & 0 & \cdots & & \\
	    a_{2\,1} & a_{2\,2} & a_{2\,3} & 0 & \cdots & & \\
	    a_{3\,1} & a_{3\,2} & a_{3\,3} & 0 & \cdots & & \\
	    \vdots & \vdots & \vdots & 
		\left[ \begin{array}{c}
		    a_{4\,4} \\ a_{5\,4}
		\end{array} \right] &
		\left[ \begin{array}{c}
		    a_{4\,5} \\ a_{5\,5}
		\end{array} \right] &
		\begin{array}{c}
		    a_{4\,6} \\ a_{5\,6}
		\end{array}  &
		\begin{array}{c}
		    a_{4\,7} \\ a_{5\,7}
		\end{array}  \\
	    &&&
		\begin{array}{c}
		    a_{6\,4} \\ a_{7\,4}
		\end{array} &
		\begin{array}{c}
		    a_{6\,5} \\ a_{7\,5}
		\end{array} &
		\begin{array}{c}
		    a_{6\,6} \\ a_{7\,6}
		\end{array}  &
		\begin{array}{c}
		    a_{6\,7} \\ a_{7\,7}
		\end{array} 
	\end{array} \right]
\end{displaymath}
If columns 4 and 5 are nearly parallel then so are these two 
vectors; {\tt Rsubspace\_rot()} would rotate the shorter of 
the two through some fixed angle, then go on to check and, 
perhaps rotate, other pairs until the matrix had no
parallel columns.  As we noted in section \ref{sec:A-choice}, 
this technique is not at all optimal. Indeed, it is not even 
certain to produce a non-singular matrix, though, in practice, it
always does.  The column-rotor scheme produces smaller, more snuggly 
fitting bounding prisms than the fixed-form fattener 
and so improves the program's performance.

     The main computational work in this scheme is in inverting the
matrix $A$ and in calculating the $w_j$.  Since, after 
column-rotation, $A$ bears no direct relation to $DG_{x_c} \circ P$, we 
cannot expect any special form for $A^{-1} \circ DG_x \circ P$.  
Instead, we must use the APM library to compute some $\widetilde{A} 
\approx A^{-1}$ directly.  Define\thesisfootnote{ 
    Some of the notation in this section, like $B$ here, is 
    introduced as a guide to the names of variables used in the code.
} a 4 \by 4 matrix $B$ such that
\begin{displaymath}
    \left[ \begin{array}{cc} 
	\Buu & \Buv \\ 
	\Bvu & \Bvv 
    \end{array} \right]
    \left[ \begin{array}{cc} 
	\Auu & \Auv \\ 
	\Avu & \Avv 
    \end{array} \right] = \id. 
\end{displaymath}
Then
\begin{eqnarray}
    A^{-1} & = &  \left[ \begin{array}{cll}		
			\id & 0 & 0 \\
			0 & \Buu & \Buv \\
			0 & \Bvu & \Bvv 
		  \end{array} \right]
    		  \left[ \begin{array}{cll}		
		  	\AppInv & 0 & 0 \\
		  	-\Aup \AppInv & \id & 0 \\
		  	-\Avp \AppInv & 0 & \id 
		  \end{array} \right], \nonumber \\
    & = & \left[ \begin{array}{cll}
		\AppInv & 0 & 0 \\
		\left\{ \begin{array}{l}
			-\Buu \Aup \AppInv \\
			-\Buv \Avp \AppInv
		\end{array} \right\} & \Buu & \Buv  \\
		\left\{ \begin{array}{l}
			-\Bvu \Aup \AppInv \\
			-\Bvv \Avp \AppInv
		\end{array} \right\} & \Bvu & \Bvv  
	  \end{array} \right]
    \; \approx \;  \left[ \begin{array}{ccc}
			\ATpp & 0 & 0 \\
			\ATup & \ATuu & \ATuv \\
			\ATvp & \ATvu & \ATvv
		\end{array} \right].
\label{eqn:crA-inv}
\end{eqnarray}
Note that the lower-left, 4 \by 4 block of $\widetilde{A}$ is just $B$.
Then, again taking $\App = \Ppp $, the product 
$A^{-1} \circ DG_x \circ P$ is 
\begin{equation}
    \left[ \begin{array}{ccc}
	\id & 0 & 0 \\
	\left\{ \begin{array}{l}
		\ATup \Ppp + \ATuu \Pvp + \\
		\ATuv(\bgamma \Ppp - \Pup + \bbeta \Pvp)
	\end{array} \right\} &
	\left\{ \begin{array}{l}
		\ATuu \Pvu + \\
		\ATuv(\bbeta \Pvu - \Puu)
	\end{array} \right\} &
	\left\{ \begin{array}{l}
		\ATuu \Pvv + \\
		\ATuv(\bbeta \Pvv - \Puv)
	\end{array} \right\} \\
	\left\{ \begin{array}{l}
		\ATvp \Ppp + \ATvu \Pvp + \\
		\ATvv(\bgamma \Ppp - \Pup + \bbeta \Pvp)
	\end{array} \right\} &
	\left\{ \begin{array}{l}
		\ATvu \Pvu + \\
		\ATvv(\bbeta \Pvu - \Puu)
	\end{array} \right\} &
	\left\{ \begin{array}{l}
		\ATvu \Pvv + \\
		\ATvv(\bbeta \Pvv - \Puv)
	\end{array} \right\} 
    \end{array} \right]. 
\label{eqn:cr-big}
\end{equation}
Since the fattening scheme does not alter the first three columns,
the blocks $\Aup$ and $\Avp$ have the forms dictated by
$A = DG_{x_c} \circ P$; these are the same as the forms used in 
equation (\ref{eqn:A-blocks}) for the fixed-form scheme.
Equation (\ref{eqn:cr-big}) then simplifies to 
\begin{displaymath}
    \left[ \begin{array}{ccc}
	\id & 0 & 0 \\
	\left\{ \begin{array}{l}
		\ATuv (\bgamma - \bgamma_c)\Ppp + \\ 
		\ATuv (\bbeta  - \bbeta_c)\Pvp
	\end{array} \right\} &
	\left\{ \begin{array}{l}
		\ATuu \Pvu + \\
		\ATuv(\bbeta \Pvu - \Puu)
	\end{array} \right\} &
	\left\{ \begin{array}{l}
		\ATuu \Pvv + \\
		\ATuv(\bbeta \Pvv - \Puv)
	\end{array} \right\} \\
	\left\{ \begin{array}{l}
		\ATuv (\bgamma - \bgamma_c)\Ppp  + \\
		\ATuv (\bbeta  - \bbeta_c)\Pvp
	\end{array} \right\} &
	\left\{ \begin{array}{l}
		\ATvu \Pvu + \\
		\ATvv(\bbeta \Pvu - \Puu)
	\end{array} \right\} &
	\left\{ \begin{array}{l}
		\ATvu \Pvv + \\
		\ATvv(\bbeta \Pvv - \Puv)
	\end{array} \right\} 
    \end{array} \right] 
\end{displaymath}
and the row sums contributing to $\bfw_u$ are 
\begin{eqnarray}
&&
    \ub	\left\{ \begin{array}{l}
	    \rowSum{\ATuv (\bgamma - \bgamma_c) \Ppp +
		    \ATuv (\bbeta - \bbeta_c) \Pvp }{j} + \\
	    \rowSum{\ATuu \Pvu + 
		    \ATuv ( \bbeta \Pvu - \Puu ) }{j} + \\
	    \rowSum{\ATuu \Pvv + 
		    \ATuv ( \bbeta \Pvv - \Puv ) }{j}
	\end{array} \right\}, \nonumber \\
& \leq &
    \ub \rowSum{\ATvu}{j} \,\left\{
	    \ub (\matSum{\bgamma - \bgamma_c}) \matSum{\Ppp} +
    	    \ub (\matSum{\bbeta - \bbeta_c}) \matSum{\Pvp}  
    \right\} + \medStrut \nonumber \\
&&  \ub \matSum{ \ATuu \Pvu + \ATuv (\bbeta \Pvu - \Puu) } + 
					\medStrut  \nonumber \\
&&  \ub \matSum{ \ATuu \Pvv + \ATuv (\bbeta \Pvv - \Puv) }.
\label{eqn:cr-rows}
\end{eqnarray}
All the upper bounds are taken over $\bfx \in S$;
the formulae for $\bfw_v$ are similar. The program calculates
the entries in $\widetilde{A}$ to at least {\tt precision} decimal
places, then treats them as exact in the evaluation of 
$\rowSum{\ATvu}{j}$ and in expressions like 
\begin{equation}
    \ub \matSum{ \ATuu \Pvv + \ATuv (\bbeta \Pvv - \Puv) }.
\label{eqn:cr-bound-rows}
\end{equation}
Upper bounds like (\ref{eqn:cr-bound-rows}) are so important 
that the program includes a special function,
{\tt Rbound\_rows()}, to evaluate them.  
To account for
the small errors ( $\leq {\rm 10}^{-precision}$) in $\widetilde{A}$,  
the program adds {\tt max\_error} to the value of $w_j$ as
computed according to (\ref{eqn:cr-rows}).
Since the entries of $\bbeta$ and  $P$ are all less in 
absolute value than 10, and since {\tt max\_error} is at least five
orders of magnitude bigger than than the largest error in 
$\widetilde{A}$, this is a very conservative estimate.

\subsubsection{matrix inversion}

     Notice that only blocks from the lower-left corner
of $\widetilde{A}$ appear in equation (\ref{eqn:cr-rows}); it
will be enough to calculate just these blocks to 
{\tt precision} decimal places.  The function, {\tt Rgauss()},
which does the calculation, takes a matrix $M$ and uses
the Gauss-Jordan algorithm with full pivoting to produce a
result $\widetilde{M} \approx M^{-1}$ such that 
$M\widetilde{M} = \id + O(\epsilon)$, that is
\begin{displaymath}
	|[M\widetilde{M}]_{i\,j} - \delta_{ij}| \leq \epsilon
\end{displaymath}
where $\delta_{ij}$ is the Kroneker delta function and $\epsilon$ is,
as usual, ${\rm 10}^{-precision}$.

     To apply the Gauss-Jordan algorithm to an {\em n \by n} matrix
$M$ one constructs the \mbox{\em n \by 2n} matrix
\begin{displaymath}
    G   = \left[ \begin{array}{cccc}
		M_{1\,1} & M_{1\,2} & \cdots & M_{1\,n} \\
		M_{2\,1} & M_{2\,2} & \cdots & M_{2\,n} \\
		\vdots & \vdots & \ddots &  \\
		M_{n\,1} & M_{n\,2} & & M_{n\,n}
	  \end{array} \right| \left. \begin{array}{cccr}
	      1 & 0 & \cdots & 0 \\
	      0 & 1 & \cdots & 0 \\
	      \vdots & \vdots & \ddots & \vdots \\
	      0 & 0 & & \hspace{1.5ex} 1
	  \end{array} \right]
\end{displaymath}
made by appending a copy of the identity to the right side
of $M$.
The algorithm transforms the left side of G into the 
identity through a sequence of row operations which simultaneously
transform the right side into $A^{-1}$. The first step is to 
multiply the top row by a constant so that the (1,1) entry is
equal to one, then subtract suitably scaled multiples of the first 
row from each of the others in such a way as to eliminate the
entries in the first column.  After this step the system looks like
\begin{equation}
    G'  = \left[ \begin{array}{cccc}
		1 & \frac{M_{1\,2}}{M_{1\,1}} 
		    & \cdots  & \frac{M_{1\,n}}{M_{1\,1}} \\
		0 & M_{2\,2} - \frac{M_{2\,1}M_{1\,2}}{M_{1\,1}} 
		    &  & \\
		\vdots & & \ddots &  \\
		0 & M_{n\,2} - \frac{M_{n\,1}M_{1\,2}}{M_{1\,1}} 
		    &  & \\
	  \end{array} \right| \left. \begin{array}{cccr}
	      \frac{1}{M_{1\,1}} & 0 & \cdots &  \hspace{1.5ex} 0 \\
	      -\frac{M_{2\,1}}{M_{1\,1}} & 1 & & \\
	      \vdots & \vdots & \ddots &  \\
	      -\frac{M_{n\,1}}{M_{1\,1}} & & & 1
	  \end{array} \right].
\label{eqn:G1-def}
\end{equation}
In the second step one uses multiples of the second row to
eliminate all but the (2,2) entry form the second column
$\ldots$ and so on. The true Gauss-Jordan algorithm with full 
pivoting may rearrange some of the rows and columns so as to 
place large entries on the diagonal of the left-hand block;
also, real implementations use only a single {\em n \by n}
array, gradually replacing the matrix $M$ by its approximate
inverse, $\widetilde{M}$.  The reader interested in the details
of the algorithm should consult either the code, which is
in appendix \ref{app:code}, or the excellent book \cite{NumRecp}.
Here, we will mostly ignore the rearrangenents, because they do 
not affect the error estimates we need.

     The divisions needed to calculate intermediate results like
(\ref{eqn:G1-def}) can only be done approximately so we must calculate
bounds on the errors they introduce.
Suppose all the calculations are done to some fixed precision,
{\tt inv\_dp} and define $\eInv = {\rm 10}^{inv\_dp}$. 
We will need a new symbol, $\GT'$,  to denote the approximate 
value of the matrix $G'$ and will also need to
define $\delta_1$, the largest error made in calculating an entry
of $\GT'$;
\begin{displaymath}
	\delta_1 = \ub_{j,k}|[\GT' - G']_{jk}|.
\end{displaymath}
The second step produces
\begin{equation}
    G'' = \left[ \begin{array}{cccc}
		1 & 0 & \star & \cdots \\
		0 & 1 & \star & \cdots \\
		0 & 0 & \star & \cdots \\
		\vdots & \vdots & \vdots & 
	  \end{array} \right| \left. \begin{array}{cccc}
	      \frac{1}{M_{1\,1}} & 0 & 0 & \cdots \\
	      \star & \frac{M_{1\,1}}{M_{1\,1}M_{2\,2} - 
				      M_{2\,1}M_{1\,2}}
				      & 0 & \cdots \\
	      \star & \star & 1 & \cdots  \\
	      \vdots & \vdots & \vdots &
	  \end{array} \right].
\label{eqn:G2-def}
\end{equation}
Ideally, we would use $G'$ to calculate $G''$ according to 
\begin{displaymath}
    G''_{i\,j} \; = \; 
    \left\{ \begin{array}{cl} \displaystyle
	\frac{G'_{i\,j}}{G'_{2\,2}} & {\rm if} \; i = 2 \\
 	{G'_{i\,j}} -  \displaystyle
	    \frac{G'_{i\,2}G'_{2\,j}}{G'_{2\,2}} 
	    & {\rm if} \; i \neq 2.  \bigStrut
    \end{array} \right. ,
\end{displaymath}
but instead, {\tt Rgauss()} actually calculates 
\begin{equation}
    \GT''_{i\,j} \; = \; 
    \left\{ \begin{array}{cl}
	\left[  \displaystyle
	    \frac{\GT'_{i\,j}}{\GT'_{2\,2}}
	\right]_{inv\_dp} & {\rm if} \; i = 2 \\
	\left[ \GT'_{i\,j} -  \displaystyle
	    \left[	
	        \frac{\GT'_{i\,2}\GT'_{2\,j}}{\GT'_{2\,2}} 
	    \right]_{inv\_dp}
	\right]_{inv\_dp}
	    & {\rm if} \; i \neq 2  \bigStrut
    \end{array} \right.
\label{eqn:GT2-calc}
\end{equation}
From this we must estimate  $\delta_2$, an upper bound on the difference 
between $\GT''$ and $G''$.  {\tt Rgauss()} finds $\delta_2$ in stages, 
as follows:
\begin{description}
   \item[(i)]{ Compute
		\begin{eqnarray*}
		    \delta_{piv} & =  &
			\frac{\delta_1}{|\GT'_{2\,2}| - \delta_1}
		    	+ \eInv \\
		    & \leq & \left[
			\frac{\delta_1}{|\GT'_{2\,2}| - \delta_1}
		    \right]_{inv\_dp} + 2 \eInv. \bigStrut
		\end{eqnarray*}
		This is a bound on the error made by taking
		\begin{displaymath}
		    \frac{1}{G'_{2\,2}} \, = \,
		    \left[ 
			\frac{1}{\GT'_{2\,2}}
		    \right]_{inv\_dp} \, \equiv \,
		    {\tt piv\_inv};
		\end{displaymath}
		{\tt piv\_inv} is the name used in the code.
    }
    \item[(ii)] { \begin{displaymath}
			\delta_r  = \delta_1 |\mbox{\tt piv\_inv}|
			+ \delta_p ( \ub_{k\neq2} |\GT'_{2\,k}|)
			+ \delta_1 \delta_p.
		  \end{displaymath}
		  This is a bound on the error introduced by normalizing 
		  the second row  so that its (2,2) entry is equal to one.
		}
    \item[(iii)] { \begin{eqnarray*}
			\delta_m  & = & 2 \delta_1 + 
			\delta_r \ub_{l \neq 2}|\GT'_{l\,2}| +
			\delta_1 \delta_r, \\
			& \geq & 
			\delta_1 + 
			\delta_1 \, \ub_{k \neq 2}
			    | \mbox{\tt piv\_inv} \, \GT'_{2\,k}| +
			\delta_r \ub_{l \neq 2}|\GT'_{l\,2}| +
			\delta_1 \delta_r. \\
		  \end{eqnarray*}
		  This is a matrix-wide bound on the errors made in 
		  computations like those in (\ref{eqn:GT2-calc}).
		  The inequality is a consequence of the pivoting
		  part of the algorithm, which ensures that
		  $|\mbox{\tt piv\_inv}\,\GT'_{2\,k}| \, \leq \, 1$.
		}
    \item[(iv)] { Finally,
		  \begin{displaymath}
		     \delta_2 = [\delta_m]_{inv\_dp} + \eInv.
		  \end{displaymath}
		}
\end{description}
Similar estimates eventaully give $\delta_n$, a matrix-wide
estimate on the difference between entries of $\widetilde{M}$ 
and the true inverse, $M^{-1}$.  From this we can conclude
\begin{equation}
     \left|  [M\widetilde{M}]_{i\,j} - \delta_{ij} \right|
     \; \leq \; n\delta_n \, \ub_{l,m}|M_{l\,m}| .
\label{eqn:inverse-error}
\end{equation}
Unless $M$ is singular, we can choose {\tt inv\_dp} so as
to make the error (\ref{eqn:inverse-error}) as small as we like.
{\tt Rgauss()} guarantees both $\delta_n$ and the error given by  
(\ref{eqn:inverse-error}) to be less than ${\rm 10}^{-precision}$.

%% file: appendices/codeSpec/trunc.tex
\subsubsection{about truncation}
Both the schemes described above produce matrices, $P'$, whose entries 
are long strings of digits, longer than those of the original matrix, $P$.
To avoid the computational cost of storing and manipulating long strings, the
program truncates the entries in $P'$ to {\tt precision} decimal places;
this introduces a small, readily manageable error.

    Call the truncated prism $P'_{trunc}$; its entries differ from those of
$P'$ by, at most, $\epsilon = {\rm 10}^{-precision}$, so that $\bfx \in S'$
\begin{displaymath}
	\bfx = \bfx'_c + P' \bfeta \qquad \mbox{\em for some} \; \bfeta \in Q^7
\end{displaymath}
differs from 
\begin{displaymath}
	\tilde{\bfx} = \bfx_c' + P'_{trunc} \bfeta
\end{displaymath}
by, at most, $7\epsilon$ in each coordinate.  The simplest way to 
handle this error is to incorporate it into $\delta_c$, the upper 
bound on the difference $|(G_{abc}(\bfx_c) - \bfx_c)_j'|$.
The coordinates of $G_{abc}(\bfx_c)$ are calculated 
out to {\tt precision} decimal places, so we must have 
\begin{displaymath}
 \delta_c \geq 8\epsilon.
\end{displaymath}
Since the program uses 
$\delta_c = \mbox{\tt max\_error} = {\rm 10}^{safety\_dp} \epsilon 
= {\rm 10}^5 \epsilon$,
this condition is abundantly satisfied.

%% file: AppendC.tex
\input{./appendices/code/introC}

\input{./appendices/code/Zussman}
\input{./appendices/code/source}

\thesisspace
\normalsize

%% file: appendices/code/introC.tex
\chapter{Computer Programs}
\label{app:code}

     This appendix contains the most important parts of the C programs used to
prove the results described in chapter 3.  In the interest of economy, we have
deleted most of the non-rigorous and semi-rigorous parts of the code, leaving only
those parts which bear on the correctness of our converse KAM results.  The first
section contains Lloyd Zussman's own description of his arbitrary precision library,
the rest of the appendix has been copied directly from the source files used to compile
the program.

%% file: appendices/code/Zussman.tex
\subsection{Arbitrary precision library}

\singlespace
\scriptsize
\begin{verbatim} 
APM
apmInit(init, scale_factor, base)
long init;
int scale_factor;
short base;
{}
  This routine initializes a new APM value.  The 'init' parameter is a long
  integer that represents its initial value, the 'scale_factor' variable
  indicates how this initial value should be scaled, and 'base' is the base of
  the initial value.  Note that the APM value returned by this routine is
  normally a reclaimed APM value that has been previously disposed of via
  apmDispose(); only if there are no previous values to be reclaimed will this
  routine allocate a fresh APM value (see also the apmGarbageCollect()
  routine).

  Bases can be 2 - 36, 10000, or 0, where 0 defaults to base 10000.

  If the call fails, it will return (APM)NULL and 'apm_errno' will contain a
  meaningful result.  Otherwise, a new APM value will be initialized.

  For example, assume that we want to initialize two APM values in base 10000,
  the first to 1.23456 and the second to 1 E20 ("one times 10 to the 20th
  power"):

      APM apm_1 = apmInit(123456L, -5, 0);
      APM apm_2 = apmInit(1L, 20, 0);

  As a convenience, the following macro is defined in apm.h:

      #define apmNew(BASE)    apmInit(0L, 0, (BASE))

int
apmDispose(apm)
APM apm;
{}
  This routine disposes of a APM value 'apm' by returning it to the list of
  unused APM values (see also the apmGarbageCollect() routine).  It returns
  an appropriate status which is also put into 'apm_errno'.

int
apmGarbageCollect()
{}
  When APM values are disposed of, they remain allocated.  Subsequent calls to
  apmInit() may then return a previously allocated but disposed APM value.
  This is done for speed considerations, but after a while there may be lots of
  these unused APM values lying around.  This routine reclaims the space taken
  up by these unused APM values (it frees them).  It returns an appropriate
  status which is also put into 'apm_errno'.

int
apmAdd(result, apm1, apm2)
APM result;
APM apm1;
APM apm2;
{}
  This routine adds 'apm1' and 'apm2', putting the sum into 'result', whose
  previous value is destroyed.  Note that all three parameters must have been
  previously initialized via apmInit().

  The 'result' parameter cannot be one of the other APM parameters.

  The return code and the 'apm_error' variable reflect the status of this
  function.

int
apmSubtract(result, apm1, apm2)
APM result;
APM apm1;
APM apm2;
{}
  This routine subtracts 'apm2' from 'apm1', putting the difference into
  'result', whose previous value is destroyed.  Note that all three parameters
  must have been previously initialized via apmInit().

  The 'result' parameter cannot be one of the other APM parameters.

  The return code and the 'apm_errno' variable reflect the status of this
  function.

int
apmMultiply(result, apm1, apm2)
APM result;
APM apm1;
APM apm2;
{}
  This routine multiplies 'apm1' and 'apm2', putting the product into 'result',
  whose previous value is destroyed.  Note that all three parameters must have
  been previously initialized via apmInit().

  The 'result' parameter cannot be one of the other APM parameters.

  The return code and the 'apm_errno' variable reflect the status of this
  function.

int
apmDivide(quotient, radix_places, remainder, apm1, apm2)
APM quotient;
int radix_places;
APM remainder;
APM apm1;
APM apm2;
{}
  This routine divides 'apm1' by 'apm2', producing the 'quotient' and
  'remainder' variables.  Unlike the other three basic operations,
  division cannot be counted on to produce non-repeating decimals, so
  the 'radix_places' variable exists to tell this routine how many
  digits to the right of the radix point are to be calculated before
  stopping.  If the 'remainder' variable is set to (APM)NULL, no
  remainder is calculated ... this saves quite a bit of computation time
  and hence is recommended whenever possible.

  All APM values must have been previously initialized via apmInit() (except,
  of course the 'remainder' value if it is to be set to NULL).

  Division by zero creates a zero result and a warning.

  The 'quotient' and 'remainder' variables can't be one of the other APM
  parameters.

  The return code and the 'apm_errno' variable reflect the status of this
  function.

int
apmCompare(apm1, apm2)
APM apm1;
APM apm2;
{}
  This routine compares 'apm1' and 'apm2', returning -1 if 'apm1' is less than
  'apm2', 1 if 'apm1' is greater than 'apm2', and 0 if they are equal.

  It is not an error if 'apm1' and 'apm2' are identical, and in this case the
  return value is 0.

  The 'apm_errno' variable contains the error code.  You must check this value:
  if it is set to an error indication, the comparison failed and the return
  value is therefore meaningless.

int
apmCompareLong(apm, longval, scale_factor, base)
APM apm;
long longval;
int scale_factor;
short base;
{}
  This routine works just like apmCompare(), but it compares the 'apm' value to
  'longval', scaled by 'scale_factor' in 'base'.  The 'apm_errno' variable
  contains the error code.

int
apmSign(apm)
APM apm;
{}
  This routine returns the sign of the 'apm' value: -1 for negative, 1 for
  positive.  The 'apm_errno' variable contains the error code.  You must check
  'apm_errno': if it's non-zero, the function return value is meaningless.

int
apmAbsoluteValue(result, apm)
APM result;
APM apm;
{}
  This routine puts the absolute value of 'apm' into 'result', whose previous
  value is destroyed.  Note that the two parameters must have been previously
  initialized via apmInit().

  The 'result' parameter cannot be the other APM parameter.

  The return code and the 'apm_errno' variable reflect the status of this
  function.

int
apmNegate(result, apm)
APM result;
APM num;
{}
  This routine puts the additive inverse of 'apm' into 'result', whose previous
  value is destroyed.  Note that the two parameters must have been previously
  initialized via apmInit().

  The 'result' parameter cannot be the other APM parameter.

  The return code and the 'apm_errno' variable reflect the status of this
  function.

int
apmReciprocal(result, radix_places, apm)
APM result;
int radix_places;
APM num;
{}
  This routine puts the multiplicative inverse of 'apm' into 'result', whose
  previous value is destroyed.  Note that the two APM parameters must have been
  previously initialized via apmInit().  Since taking the reciprocal involves
  doing a division, the 'radix_places' parameter is needed here for the same
  reason it's needed in the apmDivide() routine.

  Taking the reciprocal of zero yields zero with a warning status.

  The 'result' parameter cannot be the other APM parameter.

  The return code and the 'apm_errno' variable reflect the status of this
  function.

int
apmScale(result, apm, scale_factor)
APM result;
APM apm;
int scale_factor;
{}
  This routine assigns to 'result' the value of 'apm' with its radix point
  shifted by 'scale_factor' (positive 'scale_factor' means shift left).  The
  'scale_factor' represents how many places the radix is shifted in the base of
  'apm' unless 'apm' is in base 10000 ...  in this special case, 'scale_factor'
  is treated as if the base were 10.

  This is a very quick and accurate way to multiply or divide by a power of 10
  (or the number's base).

  The 'result' parameter cannot be the other APM parameter.

  The return code and the 'apm_errno' variable reflect the status of this
  function.

int
apmValidate(apm)
APM apm;
{}
  This routine sets 'apm_errno' and its return status to some non-zero value if
  'apm' is not a valid APM value.

int
apmAssign(result, apm)
APM result;
APM num;
{}
  This routine assigns the value of 'apm' to 'result', whose previous value is
  destroyed.  Note that the two parameters must have been previously
  initialized via apmInit().

  It is not considered an error if 'result' and 'apm' are identical; this case
  is a virtual no-op.

  The return code and the 'apm_errno' variable reflect the status of this
  function.

int
apmAssignLong(result, long_value, scale_factor, base)
APM result;
long long_value;
int scale_factor;
short base;
{}
  This routine assigns a long int to 'result'.  Its second through fourth
  parameters correspond exactly to the parameters of apmInit().  The only
  difference between the two routines is that this one requires that its result
  be previously initialized.  The 'long_value' parameter is a long that
  represents the value to assign to 'result', the 'scale_factor' variable
  indicates how this value should be scaled, and 'base' is the base of the
  value.

  Bases can be 2 - 36, 10000, or 0, where 0 defaults to base 10000.

  For example, assume that we want to assign values to two previously
  initialized APM entities, apm_1 and apm_2.  The base will be base 10000, the
  first value will be set to 1.23456 and the second will be set to 1 E20 ("one
  times 10 to the 20th power"):

      int ercode;

      ercode = apmAssignLong(apm_1, 123456L, -5, 0);
      ...

      ercode = apmAssignLong(apm_2, 1L, 20, 0);
      ...

  The return code and the 'apm_errno' variable reflect the status of this
  function.

int
apmAssignString(apm, string, base)
APM apm;
char *string;
short base;
{}
  This routine takes a character string containing the ASCII representation of
  a numeric value and converts it into a APM value in the base specified.  The
  'apm' parameter must have been previously initialized, 'string' must be
  non-NULL and valid in the specified base, and 'base' must be a valid base.

  The return code and the 'apm_errno' variable reflect the status of this
  function.

int
apmConvert(string, length, decimals, round, leftjustify, apm)
char *string;
int length;
int decimals;
int round;
int leftjustify;
APM apm;
{}
  This routine converts a APM value 'apm' into its ASCII representation
  'string'. The 'length' parameter is the maximum size of the string (including
  the trailing null), the 'decimals' parameter is the number of decimal places
  to display, the 'round' parameter is a true-false value which determines
  whether rounding is to take place (0 = false = no rounding), the
  'leftjustify' parameter is a true-false value which determines whether the
  result is to be left justified (0 = false = right justify; non-zero = true =
  left justify), and the 'apm' paramter is the APM value to be converted.

  The 'string' parameter must point to an area that can hold at least 'length'
  bytes.

  If the 'decimals' parameter is < 0, the string will contain the number of
  decimal places that are inherent in the APM value passed in.

  The return code and the 'apm_errno' variable reflect the status of this
  function.

int
(*apmErrorFunc(newfunc))()
int (*newfunc)();
{}
  This routine registers an error handler for errors and warnings.  Before any
  of the other APM routines return to the caller, an optional error handler
  specified in 'newfunc' can be called to intercept the result of the
  operation.  With a registered error handler, the caller can dispense with the
  repetitious code for checking 'apm_errno' or the function return status after
  each call to a APM routine.

  If no error handler is registered or if 'newfunc' is set to NULL, no action
  will be taken on errors and warnings except to set the 'apm_errno' variable.
  If there is an error handler, it is called as follows when there is an error
  or a warning:

      retcode = (*newfunc)(ercode, message, file, line, function)

  where ...

      int retcode;          /* returned by 'newfunc': should be 'ercode' */
      int ercode;           /* error code */
      char *message;        /* a short string describing the error */
      char *file;           /* the file in which the error occurred */
      int line;                    /* the line on which the error occurred */
      char *function;        /* the name of the function in error */

  Note that your error handler should normally return 'ercode' unless it does a
  longjmp, calls exit(), or in some other way interrupts the normal processing
  flow.  The value returned from your error handler is the value that the apm
  routine in error will return to its caller.

  The error handler is called after 'apm_errno' is set.

  This routine returns a pointer to the previously registered error handler or
  NULL if one isn't registered.

int
apmCalc(result, operand, ..., NULL)
APM result;
APM operand, ...;
{}
  This routine performs a series of calculations in an RPN ("Reverse
  Polish Notation") fashion, returning the final result in the 'result'
  variable.  It takes a variable number of arguments and hence the
  rightmost argument must be a NULL.

  Each 'operand' is either a APM value or a special constant indicating
  the operation that is to be performed (see below).  This routine makes
  use of a stack (16 levels deep) similar to that in many pocket
  calculators.  It also is able to access a set of 16 auxiliary
  registers (numbered 0 through 15) for holding intermediate values.

  The stack gets reinitialized at the start of this routine, so values
  that have been left on the stack from a previous call will disappear.
  However, the auxiliary registers are static and values remain in these
  registers for the duration of your program.  They may also be
  retrieved outside of this routine (see the apmGetRegister() and
  apmSetRegister() routines, below).

  An operand that is an APM value is automatically pushed onto the stack
  simply by naming it in the function call.  If the stack is full when a
  value is being pushed onto it, the bottommost value drops off the
  stack and the push succeeds; this is similar to how many pocket
  calculators work.  Also, if the stack is empty, a pop will succeed,
  yielding a zero value and keeping the stack empty.  The topmost value
  on the stack is automatically popped into the 'result' parameter after
  all the operations have been performed.

  An operand that is one of the following special values will cause
  an operation to be performed.  These operations are described in the
  following list.  Note that the values "V", "V1", and "V2" are used
  in the following list to stand for temporary values:

        APM_ABS                pop V, push absolute value of V
        APM_NEG                pop V, push -V
        APM_CLEAR        empty the stack
        APM_DUP                pop V, push V, push V
        APM_SWAP        pop V1, pop V2, push V1, push V2
        APM_SCALE(N)        pop V, push V scaled by N [ as in apmScale() ]
        APM_PUSH(N)        V = value in register N, push V
        APM_POP(N)        pop V, store it in register N
        APM_ADD                pop V1, pop V2, push (V2 + V1)
        APM_SUB                pop V1, pop V2, push (V2 - V1)
        APM_MUL                pop V1, pop V2, push (V2 * V1)
        APM_DIV(N)        pop V1, pop V2, push (V2 / V1) with N radix places
                        [ as in apmDivide() ], remainder goes into register 0
        APM_RECIP(N)        pop V, push 1/V with N radix places
                        [ as in apmReciprocal() ]

  Since register 0 is used to hold the remainder in a division, it is
  recommended that this register not be used to hold other values.

  As an example, assume that APM values "foo", "bar", and "baz" have
  been initialized via apmInit() and that "foo" and "bar" are to be used
  to calculate "baz" as follows (assume that divisions stop after 16
  decimal places have been calcluated):

        baz = 1 / ((((foo * bar) + foo) / bar) - foo)

  The function call will be:

        bcdCalc(baz, foo, APM_DUP, APM_POP(1), bar, APM_DUP, APM_POP(2),
                APM_MUL, APM_PUSH(1), APM_ADD, APM_PUSH(2), APM_DIV(16),
                APM_PUSH(1), APM_SUB, APM_RECIP(16), NULL);

  Note that the value of "foo" is stored in register 1 and the value of
  "bar" is stored in register 2.  After this call, these registers will
  still contain those values.

int
apmGetRegister(regvalue, regnumber)
APM regvalue;
int regnumber;
{}
  The value in auxiliary register number 'regnumber' is assigned to APM
  value 'regvalue'.  The 'regnumber' parameter must be between 0 and 15,
  inclusive.  The 'regvalue' parameter must have been previously
  initialized via apmInit().

int
apmSetRegister(regvalue, regnumber, newvalue)
APM regvalue;
int regnumber;
APM newvalue;
{}
  The value in auxiliary register number 'regnumber' is assigned to APM
  value 'regvalue', and then the APM value 'newvalue' is stored in that
  same register.  The 'regnumber' parameter must be between 0 and 15,
  inclusive.  The 'regvalue' and 'newvalue' parameters must have been
  previously initialized via apmInit().
\end{verbatim}

%% file: appendices/code/source.tex
\normalsize
\thesisspace
\section{Source code}

	The listings below contain only those functions crucial to the
correct execution of a converse KAM calculation. Some references to 
inessential or semi-rigorous  parts of the code have been left in place 
because we wished to present the important functions exactly as they 
appear in the original source files.

\subsection{special functions}
\subsubsection{the header file apmSpecial.h}
\singlespace
\scriptsize
\normalsize
\subsubsection{apmCos(), etc.}
\scriptsize
\input{./appendices/code/special}

\normalsize
\subsection{interval arithmetic}
\subsubsection{the header file bounding.h}
\scriptsize
\input{./appendices/code/boundingh}
\normalsize
\subsubsection{expressions}
\scriptsize
\input{./appendices/code/bounding}
\normalsize
\subsubsection{bounding trig. functions}
\scriptsize
\input{./appendices/code/bd_trig}

\normalsize
\subsection{starting points and global bounds}
\scriptsize
\input{./appendices/code/start}
\input{./appendices/code/lipschitz}
\input{./appendices/code/lambda}
\input{./appendices/code/trace}

\normalsize
\subsection{control of the computation}
\subsubsection{the header file converse.h}
\scriptsize
\input{./appendices/code/converseh}
\normalsize
\subsubsection{main()}
\scriptsize
\input{./appendices/code/converse}
\normalsize
\subsubsection{Rtry\_prism()}
\scriptsize
\input{./appendices/code/follow}

\normalsize
\subsection{the map}
\subsubsection{the header file map.h}
\scriptsize
\input{./appendices/code/maph}
\normalsize
\subsubsection{mapping functions}
\scriptsize
\input{./appendices/code/map}
\input{./appendices/code/prismatic}
\normalsize

\subsection{images of prisms}
\subsubsection{the header file rows.h}
\scriptsize
\input{./appendices/code/rowsh}
\normalsize
\subsubsection{Rglobal\_bounds()}
\scriptsize
\input{./appendices/code/rowSums}
\normalsize
\subsubsection{column-rotor}
\scriptsize
\input{./appendices/code/cr_fatten}
\input{./appendices/code/cr_rows}
\normalsize
\subsubsection{fixed-form}
\scriptsize
\input{./appendices/code/ff_fatten}
\input{./appendices/code/ff_rows}
\normalsize
\subsubsection{matrix inverter}
\scriptsize
\input{./appendices/code/gauss}

%% file: appendices/code/special.tex
\begin{verbatim} 
# include <stdio.h>
# include <math.h>
# include "apm.h"
# include "apmPrint.h"
# include "apmSpecial.h"

APM    *sinCoef, *cosCoef ;
APM    zero, one, two ;
APM    pi, two_pi, half_pi, threeHalf_pi, eighths_2pi[8] ;
APM    Theta, scratch, xMod2pi, Theta_sq, Answer ;
APM    sinFactrl, cosFactrl, apmOrder ;
APM    approx[2], diff, ub_diff ;
int    trig_dp, specialsInit = NO  ;
int    trig_terms, dp_lost ;
char   pi_str[] = "3.14159265358979323846243383279502884197169399375" ;
char   log_buf[BUF_SZ] ;

/* ++++++++++++++++++++++++ */

initApmSpecials()
{
    int     k ;

        /* Initialize a bunch of APMs.  Theta will be the reduced argument 
           of a trig function; it will be between zero and pi / 4.          */
  
    pi = apmNew( 0 ) ;
    one = apmInit( 1L, 0, 0 ) ;
    two = apmInit( 2L, 0, 0 ) ;
    zero = apmInit( 0L, 0, 0 ) ;
    diff = apmNew( 0 ) ;
    Theta = apmNew( 0 ) ;
    Answer = apmNew( 0 ) ;
    two_pi = apmNew( 0 ) ;
    half_pi = apmNew( 0 ) ;
    scratch = apmNew( 0 ) ;
    ub_diff = apmNew( 0 ) ;
    xMod2pi = apmNew( 0 ) ;
    apmOrder = apmNew( 0 ) ;
    Theta_sq = apmNew( 0 ) ;
    sinFactrl = apmNew( 0 ) ;
    cosFactrl = apmNew( 0 ) ;
    approx[0] = apmNew( 0 ) ;
    approx[1] = apmNew( 0 ) ;
    threeHalf_pi = apmNew( 0 ) ;
    for( k=0 ; k < 8 ; k++ )
        eighths_2pi[k] = apmNew( 0 ) ;

        /* Obtain some rational mutiples of pi.  These will be helpful
           when we go to restrict the domain of the trig functions to 
           between zero and pi / 4 .                                        */

    apmAssignString( pi, pi_str, 0 ) ;

    apmMultiply( scratch, two, two ) ;
    apmDivide( eighths_2pi[0], (PI_DP+2), (APM)NULL, pi, scratch) ;

    for( k=1 ; k < 8 ; k++ ) 
        apmAdd( eighths_2pi[k], eighths_2pi[0], eighths_2pi[k-1] ) ;

    apmMultiply( two_pi, pi, two ) ;
    apmAssign( half_pi, eighths_2pi[1] ) ;
    apmAssign( threeHalf_pi, eighths_2pi[5] ) ;

    setTrigDp( DFLT_TRIG_DP ) ;

    dp_lost = 0 ;
    specialsInit = YES ;

    return( 1 ) ;
}
/* ++++++++++++++++++++++++++ */

setTrigDp( dp ) 

int dp ;
{
    double  j, j_fact, ten_to_dp ;

        /* Check to see that the desired accuracy is compatible
           with our knowledge of pi.                                */

    if( (dp+2) > PI_DP ) {
        fprintf( stderr, 
        "We don't know pi well enough to achieve the desired accuracy. \n" ) ;
        return( 0 ) ;
    }
    else
        trig_dp = dp+2 ;

        /* Assume the argument is between zero and pi / 4.  How many
           terms from the Taylor series do we need to include ?     */

    trig_terms = 1 ;
    ten_to_dp = pow( 10.0, (double)dp ) ;
    for( j = 1.0, j_fact = 1.0 ; j_fact < ten_to_dp ; j += 2.0 )  {
        j_fact *= j * (j + 1) ;
        trig_terms++ ;
        if( trig_terms > MAX_TRIG_TERMS ) {
            fprintf( stderr, "Too many terms required. \n" ) ;
            return(0) ;
        }
    }

    trig_dp += (int)( ceil( log10((double) trig_terms) ) ) ;
    setTrigCoef() ;
    return( dp ) ;
}
/* +++++++++++++++++++++++++++++++++++++ */

reduceArg( x )
/*
    Takes x, chops off enough multiples of two_pi to get it 
    into the interval between zero and two_pi.  Checks that we
    haven't lost an unacceptable amount of precision in doing
    this stage of the reduction.  Then  chops off multiples
    of pi/4 to get the argument into the interval between zero and
    pi/4.  Sets Theta equal to the reduced argument and returns
    an integer indicating in which of eight equally spaced intervals
    x (mod two_pi) lay.  If any precision is lost, dp_lost is set 
    to the number of decimal places lost.
                                                                        */
APM  x ;
{
    int   octant ;
    char  qtnt_str[BUF_SZ] ;

        /* Note that we haven't lost any decimal places yet. */
    dp_lost = 0 ;

        /* Whack out many multiples of two_pi. */
    apmDivide( scratch, 3, (APM)NULL, x, two_pi ) ;
    apmFloorString( qtnt_str, BUF_SZ, scratch ) ;
    apmAssignString( scratch, qtnt_str, 0 ) ;
    apmMultiply( Answer, scratch, two_pi ) ;
    apmSubtract( xMod2pi, x, Answer ) ;
    if( apmSign( xMod2pi ) == -1 ) 
         apmCalc( xMod2pi, xMod2pi, two_pi, APM_ADD, NULL ) ;

    for( octant=0 ;  (octant < 8)  ; octant++ ) {
         if( apmCompare(xMod2pi, eighths_2pi[octant]) < 0 )
             break ;
    }

    switch( octant ) {
        case 0 :
            apmAssign( Theta, xMod2pi ) ;
            break ;

        case 1 :
            apmSubtract( Theta, half_pi, xMod2pi ) ;
            break ;

        case 2 :
            apmSubtract( Theta, xMod2pi, half_pi ) ;
            break ;

        case 3 :
            apmSubtract( Theta, pi, xMod2pi ) ;
            break ;

        case 4 :
            apmSubtract( Theta, xMod2pi, pi ) ;
            break ;

        case 5 :
            apmSubtract( Theta, threeHalf_pi, xMod2pi ) ;
            break ;

        case 6 :
            apmSubtract( Theta, xMod2pi, threeHalf_pi ) ;
            break ;

        case 7 :
            apmSubtract( Theta, two_pi, xMod2pi ) ;
            break ;

        default :
            break ;
    }

        /* Check for loss of precision */
    if( (PI_DP - strlen(qtnt_str)) < trig_dp )
        dp_lost = trig_dp - PI_DP + strlen(qtnt_str) ;
    else
        dp_lost = 0 ;

    return( octant ) ;
}
/* +++++++++++++++++++++++ */

reducedSin()
/*
    Takes the sine of Theta, puts the result in Answer.
                                                                */
{
    int  order, dp_to_find, term_num ; 

    apmAssign( Answer, zero ) ;
    apmMultiply( Theta_sq, Theta, Theta ) ;

    term_num = trig_terms - 1 ;
    for( order = ( 2 * trig_terms - 1 ) ; order > 0 ; order -= 2 ) {

            /* Multiply the old partial sum by Theta squared
               and add in a new coefficient                        */
        apmMultiply( scratch, Answer, Theta_sq ) ;
        apmAdd( Answer, sinCoef[term_num--], scratch ) ;
        apmTruncate( Answer, trig_dp ) ;
    }

        /*  Multiply by the final factor of Theta,
            divide by the factorial, and return  */

    if( dp_lost > 0 ) 
        dp_to_find = trig_dp + 1 - dp_lost ;
    else
        dp_to_find = trig_dp + 1 ;

    apmMultiply( scratch, Answer, Theta ) ;
    apmDivide( Answer, dp_to_find, (APM)NULL, scratch, sinFactrl ) ;
    return ;
}
/* ++++++++++++++++++++++++++++++++++++ */

reducedCos()
/*
    Takes the cosine of Theta, puts the result in Answer.
                                                                */
{
    int  order, dp_to_find, term_num ; 

    apmAssign( Answer, zero ) ;
    apmMultiply( Theta_sq, Theta, Theta ) ;

    term_num = trig_terms - 1 ;
    for( order = ( 2 * trig_terms - 2 ) ; order >= 0 ; order -= 2 ) {

            /* Multiply the old partial sum by Theta squared
               and add in a new coefficient                        */
        apmMultiply( scratch, Answer, Theta_sq ) ;
        apmAdd( Answer, cosCoef[term_num--], scratch ) ;

        apmTruncate( Answer, trig_dp ) ;
    }

        /*  Divide by the factorial, 
            Put the result into Answer, and return  */

    if( dp_lost > 0 ) 
        dp_to_find = trig_dp + 1 - dp_lost ;
    else
        dp_to_find = trig_dp + 1 ;

    apmDivide( scratch, dp_to_find, (APM)NULL, Answer, cosFactrl ) ;
    apmAssign( Answer, scratch ) ;
    return ;
}
/* ++++++++++++++++++++++++++++++++++++ */

apmSin( result, x )

APM   result, x ;
{
    int   octant ;

    if( specialsInit == NO ) {
        fprintf( stderr, 
                 "apmSin() : Please call initApmSpecials(). \n" ) ;
         apmAssignLong( result, 0L, 0, 0 ) ;
        apm_errno = APM_EPARM ;
        return ;
    }
    else 
        apm_errno = APM_OK ;

        /* Reduce the argument, report any loss of precision, and
           note in which octant x (mod two_pi) lay.                  */

    octant = reduceArg( x ) ;
    if( dp_lost > 0 )  {
        fprintf( stderr, 
        "apmSin : Big argument, lost %d decimal places from the answer. \n", 
                                                                  dp_lost ) ;
        apm_errno = APM_WTRUNC ; 
    }
    else
        apm_errno = APM_OK ;

        /* Evaluate the sine.  Which of the two reduced functions
           one uses depends on the octant.                        */

    switch( octant ) {
        case 0 :
            reducedSin() ;
            break ;

        case 1 :
            reducedCos() ;
            break ;

        case 2 :
            reducedCos() ;
            break ;

        case 3 :
            reducedSin() ;
            break ;

        case 4 :
            reducedSin() ;
            apmNegate( scratch, Answer ) ;
            apmAssign( Answer, scratch ) ;
            break ;

        case 5 :
            reducedCos() ;
            apmNegate( scratch, Answer ) ;
            apmAssign( Answer, scratch ) ;
            break ;

        case 6 :
            reducedCos() ;
            apmNegate( scratch, Answer ) ;
            apmAssign( Answer, scratch ) ;
            break ;

        case 7 :
            reducedSin() ;
            apmNegate( scratch, Answer ) ;
            apmAssign( Answer, scratch ) ;
            break ;

        default :
            break ;
    }

    apmAssign( result, Answer ) ;

    return ;
}
/* +++++++++++++++++++++++++ */

apmCos( result, x )

APM   result, x ;
{
    int   octant ;

    if( specialsInit == NO ) {
        fprintf( stderr, 
                 "apmCos() : Please call initApmSpecials() first. \n" ) ;
         apmAssignLong( result, 0L, 0, 0 ) ;
        apm_errno = APM_EPARM ;
        return ;
    }
    else 
        apm_errno = APM_OK ;

        /* Reduce the argument, report any loss of precision, and
           note in which octant x (mod two_pi) lay.                  */

    octant = reduceArg( x ) ;
    if( dp_lost > 0 )  {
        fprintf( stderr, 
        "apmCos : Big argument, lost %d decimal places from the answer. \n", 
                                                                  dp_lost ) ;
        apm_errno = APM_WTRUNC ; 
    }
    else
        apm_errno = APM_OK ;

        /* Evaluate the cosine.  Which of the two reduced functions
           one uses depends on the octant.                          */

    switch( octant ) {
        case 0 :
            reducedCos() ;
            break ;

        case 1 :
            reducedSin() ;
            break ;

        case 2 :
            reducedSin() ;
            apmNegate( scratch, Answer ) ;
            apmAssign( Answer, scratch ) ;
            break ;

        case 3 :
            reducedCos() ;
            apmNegate( scratch, Answer ) ;
            apmAssign( Answer, scratch ) ;
            break ;

        case 4 :
            reducedCos() ;
            apmNegate( scratch, Answer ) ;
            apmAssign( Answer, scratch ) ;
            break ;

        case 5 :
            reducedSin() ;
            apmNegate( scratch, Answer ) ;
            apmAssign( Answer, scratch ) ;
            break ;

        case 6 :
            reducedSin() ;
            break ;

        case 7 :
            reducedCos() ;
            break ;

        default :
            break ;
    }

    apmAssign( result, Answer ) ;
    return ;
}
/* +++++++++++++++++++++++++ */

apmSqrt( result, dp, x ) 
/*
    Find square roots using Newton's method.
                                                */
int    dp ;
APM    x, result ;
{
    int     comp, dp_plus ;
    APM     *this_approx, *next_approx, *temp ;

/*
        Check that all the scratch variables are ready.
*/

    if( specialsInit == NO ) {
        fprintf( stderr, 
                 "apmSqrt() : Please call initApmSpecials() first. \n" ) ;
         apmAssignLong( result, 0L, 0, 0 ) ;
        apm_errno = APM_EPARM ;
        return ;
    }
    else 
        apm_errno = APM_OK ;

/*
        If the argument is zero, just return zero.
        If the argument is negative, whine.
*/
    if( (comp = apmCompare( x, zero )) == 0 ) {
        apmAssign( result, zero ) ;
        return ;
    }
    else if( comp == -1 ) {
        fprintf( stderr, "apmSqrt() : Can't handle negative arguments.\n" ) ;
        apm_errno = APM_EPARM ;
        return ;
    }
    else
        apm_errno = APM_OK ;

/*
        Do up Newton.  The rule is

                y[n+1] = (y[n] + x/y[n]) / 2.0
*/
    dp_plus = dp + 2 ;
    apmAssignLong( ub_diff, 1L, -dp_plus, 0 ) ;

    this_approx = &approx[0] ;
    next_approx = &approx[1] ;

    apmAssign( *this_approx, x ) ;
    apmAssign( *next_approx, zero ) ;
    apmSubtract( diff, *this_approx, *next_approx ) ;
    while( apmCompare( diff, ub_diff) > 0 ) {
        apmDivide( scratch, dp_plus, (APM) NULL, x, *this_approx ) ;
        apmCalc( scratch, scratch, *this_approx, APM_ADD, NULL ) ;
        apmDivide( *next_approx, dp_plus, (APM) NULL, scratch, two ) ;
        apmTruncate( *next_approx, dp_plus ) ; 

        apmCalc( diff, *this_approx, *next_approx, APM_SUB, APM_ABS, NULL ) ;
        m_swap( this_approx, next_approx, temp ) ;
    }

    apmAssign( result, *this_approx ) ;
    return ;
}
/* +++++++++++++++++++++++++++++++++++++++ */

apmFloor( result, arg, base )

int   base ;
APM   result, arg ;
{
    char  buf[BUF_SZ], *cpt ;

    apmConvert( buf, BUF_SZ, 2, NO_ROUND, LEFT_JUST, arg ) ;
    for( cpt = buf ; *cpt != '\0' ; cpt++ )
        if( *cpt == '.' )
            *cpt = '\0' ;

    apmAssignString( result, buf, base ) ;
}
/* ++++++++++++++++++++++++++++++++ */

setTrigCoef() 
{
   int   j, order, coef_num ;
   char  *malloc() ;

   sinCoef = (APM *) malloc( trig_terms * sizeof( APM ) ) ;
   cosCoef = (APM *) malloc( trig_terms * sizeof( APM ) ) ;
   if( (sinCoef == NULL) || (cosCoef == NULL) ) {
       fprintf( stderr, "Trouble allocating %d APMs for coefficients.\n" ) ;
       exit(0) ;
    }

   for( j=0 ; j < trig_terms ; j++ ) {
       sinCoef[j] = apmNew( 0 ) ;
       cosCoef[j] = apmNew( 0 ) ;
    }

    if( (trig_terms % 2) != 0 ) {
        apmAssignLong( sinCoef[trig_terms-1], -1L, 0, 0 ) ;
        apmAssignLong( cosCoef[trig_terms-1], -1L, 0, 0 ) ;
    }
    else {
        apmAssignLong( sinCoef[trig_terms-1], 1L, 0, 0 ) ;
        apmAssignLong( cosCoef[trig_terms-1], 1L, 0, 0 ) ;
    }

    coef_num = trig_terms - 2 ;
    for( order = (2 * trig_terms - 1) ; order > 1 ; order -= 2 ) {
        /* coefficients for the sine */

        apmAssignLong( apmOrder, -((long) order), 0, 0 ) ;
        apmMultiply( scratch, sinCoef[coef_num+1], apmOrder ) ;
        apmAssignLong( apmOrder, (long)(order-1), 0, 0 ) ;
        apmMultiply( sinCoef[coef_num], scratch, apmOrder ) ;

        /* coefficients for the cosine */
        apmMultiply( scratch, cosCoef[coef_num+1], apmOrder ) ;
        apmAssignLong( apmOrder, -(long)(order-2), 0, 0 ) ;
        apmMultiply( cosCoef[coef_num], scratch, apmOrder ) ;

        coef_num-- ;
    }

    apmAssign( sinFactrl, sinCoef[0] ) ;
    apmAssign( cosFactrl, cosCoef[0] ) ;
}
/* +++++++++++++++++++++++++++++++++++++++++++++++ */

apmFloorString( s, n, x )

APM        x ;
int        n ;
char    *s ;
{
    apmConvert( s, n, 1, NO_ROUND, LEFT_JUST, x ) ;
    strip_frac( s ) ;
}
/* +++++++++++++++++++++ */

strip_frac( str ) 

char  *str ;
{
    char  *cpt ;

    for( cpt = str ; cpt != '\0' ; cpt++ )
        if( *cpt == '.' ) {
            *cpt = '\0' ;
            break ;
        }
}
/* +++++++++++++++++++++++ */

apmLogBd( x )

APM   x ;
/*
        Returns an upper bound on the base-10 log of an apm.
*/
{
    int   order ;
    char  *bpt ;

    if( apmCompare( one, x ) <= 0 ) {
        apmFloorString( log_buf, BUF_SZ, x ) ;
        return( strlen( log_buf ) ) ;
    }
    else {
        apmConvert( log_buf, BUF_SZ, (BUF_SZ-4), NO_ROUND, LEFT_JUST, x ) ;
/*
        Skip to the digits beyond the decimal point
*/
        for( bpt=log_buf ; *bpt != '.' ; bpt++ ) ;
        bpt++ ;
/*
        Count the number of zeroes to the right of the decimal point.
*/
        for( order=0 ; (*bpt == '0') ; bpt++, order-- ) ;
        return( order ) ;
    }
}

\end{verbatim}

%% file: appendices/code/boundingh.tex
\begin{verbatim} 
/*
        Data structures for calculating semi-rigorous bounds 
        on expressions. 
                                                                */

typedef   struct {  double    ub, lb ; } Bdd_dbl ;

typedef   struct {  int                nfactors ;
                    double        coef ;
                    Bdd_dbl        **factors, bound ; } Bdd_term ;

typedef  struct {   int                nterms ;
                    double      const  ;
                    Bdd_dbl     bound  ;
                    Bdd_term    *terms ; }  Bdd_expr ;
/*  
        APM partners to the structures above
                                                */

typedef   struct {  APM   ub, lb ; } Bdd_apm ;

typedef   struct {  int                nfactors ;
                    APM                coef ;
                    Bdd_apm        **factors, bound ; } Bapm_term ;

typedef  struct {  int                nterms ;
                   APM                const  ;
                   Bdd_apm        bound  ;
                   Bapm_term        *terms ; }  Bapm_expr ;
/* +++++++++++++++++++++++ */

                                         apmAssign( empty->lb, full->lb) )
                                         empty->lb = full->lb  )
                                         new.lb = apmNew( base ) )

extern  int  RmaxAbs() ;
\end{verbatim}

%% file: appendices/code/bounding.tex
\begin{verbatim} 
# include  <stdio.h>
# include  <math.h>
# include  "apm.h"
# include  "converse.h"
# include  "bounding.h"

APM        Rextrema, Rextremb, Rub, Rlb ;
APM        Rprod[4], *Rlastp = (Rprod + 4) ;
double  prod[4], *lastp = (prod + 4) ;
/* ++++++++++++++++++++++++++++ */
initBounding()
{
    int   j ;

    Rub = apmNew( BASE ) ;
    Rlb = apmNew( BASE ) ;

    Rextrema = apmNew( BASE ) ;
    Rextremb = apmNew( BASE ) ;

    for( j=0 ; j < 4 ; j++ )
        Rprod[j] = apmNew( BASE ) ;
}
/* +++++++++++++++++++++++++++ */

Rbound_term( tpt )
/*
    Take a list of bounded factors and obtain a bound on their
    product.
*/

Bapm_term   *tpt ;
{
    APM             *ppt ;
    Bdd_apm  *facptr, **lastf, **fpt ;

/*
        If there is only one factor, deal with it directly.
*/
    if( tpt->nfactors == 1 ) {
        apmAssign( Rextrema, tpt->factors[0]->ub ) ;
        apmAssign( Rextremb, tpt->factors[0]->lb ) ;
    }
/*
        Handle expressions with more than one factor.
        Since some of the factors may be negative we 
        can't just multiply to gether all the upper
        and lower bounds.
*/
    else {
        apmAssign( Rextrema, tpt->factors[0]->ub ) ;
        apmAssign( Rextremb, tpt->factors[0]->lb ) ;

        fpt = &tpt->factors[1] ;
        for( lastf = tpt->factors + tpt->nfactors ; fpt < lastf ; fpt++ ) {
            facptr = *fpt ;

            apmMultiply( Rprod[0], facptr->ub, Rextrema ) ;
            apmMultiply( Rprod[1], facptr->ub, Rextremb ) ;
            apmMultiply( Rprod[2], facptr->lb, Rextrema ) ;
            apmMultiply( Rprod[3], facptr->lb, Rextremb ) ;

            apmAssign( Rextrema, Rprod[0] ) ;
            apmAssign( Rextremb, Rprod[0] ) ;
            for( ppt = (Rprod+1) ; ppt < Rlastp ; ppt++ ) {
                if( apmCompare( *ppt, Rextrema ) == 1 )
                    apmAssign( Rextrema, *ppt ) ;
                else if( apmCompare( *ppt, Rextremb ) == -1 )
                    apmAssign( Rextremb,  *ppt ) ;
            }
        }
    }

    apmCalc( Rextrema, Rextrema, tpt->coef, APM_MUL, NULL ) ;
    apmCalc( Rextremb, Rextremb, tpt->coef, APM_MUL, NULL ) ;
    if( apmCompare( Rextrema, Rextremb ) == -1 ) {
        apmAssign( tpt->bound.ub, Rextremb ) ;
        apmAssign( tpt->bound.lb, Rextrema ) ;
    }
    else {
        apmAssign( tpt->bound.ub, Rextrema ) ;
        apmAssign( tpt->bound.lb, Rextremb ) ;
    }
}
/* ++++++++++++++++++++++++++++++++++++ */

Rbound_expr( ept ) 
/*
    Obtain bounds on the terms in a bounded expression, add them up,
    and so obtain a bound on the whole.
*/

Bapm_expr    *ept ;
{
    Bapm_term  *tpt, *last_term ;

    apmAssign( Rub, ept->const ) ;
    apmAssign( Rlb, ept->const ) ;

    tpt = ept->terms ;
    for( last_term = tpt + ept->nterms ; tpt < last_term ; tpt++ ) {
        Rbound_term( tpt ) ;
        apmCalc( Rub, Rub, tpt->bound.ub, APM_ADD, NULL ) ;
        apmCalc( Rlb, Rlb, tpt->bound.lb, APM_ADD, NULL ) ;
    }

    apmAssign( ept->bound.ub, Rub ) ;
    apmAssign( ept->bound.lb, Rlb ) ;
}
/* +++++++++++++++++++++++++++++++++ */

RmaxAbs( result, x, y )

APM   result, x, y ;
{
    apmAbsoluteValue( Rub, x ) ;
    apmAbsoluteValue( Rlb, y ) ;

    if( apmCompare( Rub, Rlb ) == 1 ) 
        apmAssign( result, Rub ) ;
    else
        apmAssign( result, Rlb ) ;
}
\end{verbatim}

%% file: appendices/code/bd_trig.tex
\begin{verbatim} 
# include  <stdio.h>
# include  <math.h>
# include  "apm.h"
# include  "apmSpecial.h"
# include  "converse.h"
# include  "bounding.h"
# include  "pi.h"

APM        half, three_halfs ;
APM        Rdelta, Rmax_cos, Rmin_cos  ;
APM        Rmax_x, Rmin_x, Rfloor_x, Rlft_val, Rrght_val ;

Bdd_apm        Rnew_theta ;
/* -------------------------------- */

initTrigBd()
/*
        Set up the APM's defined above.
*/
{
    Rdelta = apmNew( BASE ) ;
    Rmin_x = apmNew( BASE ) ;
    Rmax_x = apmNew( BASE ) ;
    Rfloor_x = apmNew( BASE ) ;
    Rmax_cos = apmNew( BASE ) ;
    Rmin_cos = apmNew( BASE ) ;
    Rlft_val = apmNew( BASE ) ;
    Rrght_val = apmNew( BASE ) ;

    Rnew_theta.ub = apmNew( BASE ) ;
    Rnew_theta.lb = apmNew( BASE ) ;

    half = apmInit( 2L, 0, BASE ) ;
    three_halfs = apmInit( 3L, 0, BASE ) ;
    apmCalc( half, half, APM_RECIP( precision ), NULL ) ; 
    apmCalc( three_halfs, half, three_halfs, APM_MUL, NULL ) ;
}
/* ++++++++++++++++++++++++++++++++++++++ */

Rbd_cos( bound, theta ) 
/*
        Obtain bounds for the cosine function over
        a certain given range of angles.
*/

Bdd_apm   *theta, *bound ;
{
/*
        An APM partner to the function above.  The variables 
        used here are static, and are defined at the top 
        of the file.
*/
/*
        Get some variables equal to theta / TWO_PI.  These will
        help decide whether the interval under consideration
        contains any extrema.
*/
    apmDivide( Rmin_x, precision, (APM)NULL, theta->lb, two_pi ) ;  
    apmDivide( Rmax_x, precision, (APM)NULL, theta->ub, two_pi ) ;  

    apmFloor( Rfloor_x, Rmin_x, BASE ) ;
    apmCalc( Rmin_x, Rmin_x, Rfloor_x, APM_SUB, NULL ) ;
    apmCalc( Rmax_x, Rmax_x, Rfloor_x, APM_SUB, NULL ) ;
    apmSubtract( Rdelta, Rmax_x, Rmin_x ) ;
    if( apmCompare( Rdelta, one ) == 1 ) {
        apmAssign( bound->ub, one ) ;
        apmNegate( bound->lb, one ) ; 
    }

    else {
        apmCos( Rlft_val, theta->lb ) ;
        apmCos( Rrght_val, theta->ub ) ;
        if( apmCompare( Rlft_val, Rrght_val ) == 1 ) {
            apmAssign( Rmax_cos, Rlft_val ) ;
            apmAssign( Rmin_cos, Rrght_val ) ;
        }
        else {
            apmAssign( Rmax_cos, Rrght_val ) ;
            apmAssign( Rmin_cos, Rlft_val ) ;
        }
/*
        Check for extrema.
*/
        if( apmCompare( Rmax_x, one) == 1 ) 
            apmAssign( Rmax_cos, one ) ;

        if( (apmCompare( Rmax_x, three_halfs) == 1) || 
           ((apmCompare( Rmin_x, half) == -1) && 
            (apmCompare( Rmax_x, half) ==  1))  ) apmNegate( Rmin_cos, one ) ;

        apmAdd( bound->ub, Rmax_cos, max_error ) ;
        apmSubtract( bound->lb, Rmin_cos, max_error ) ;
    }

    return ;
}
/* +++++++++++++++++++++++++++++ */

Rbd_sin( bound, theta ) 
/*
    Use the relation sin( x - HALF_PI ) = cos( x ) 
    and the function bd_cos() to obtain a bound on 
    the sines of angles lying in a given range.
*/

Bdd_apm  *theta, *bound ;
{
/*
        Rnew_theta is used here but is declared at the top of
        the file
*/
    apmSubtract( Rnew_theta.ub, theta->ub, half_pi ) ;
    apmSubtract( Rnew_theta.lb, theta->lb, half_pi ) ;

    Rbd_cos( bound, &Rnew_theta ) ;
    return ;
}
\end{verbatim}

%% file: appendices/code/start.tex
\begin{verbatim} 
# include <stdio.h>
# include <math.h>
# include "apm.h"
# include "converse.h"
# include "pi.h"

APM  Rstart_size ;
/* +++++++++++++++++++++ */
setHermStart( priz )

RPrism  *priz ;
{
    double  a, b, c, two_c, x, y ;
    double  jump_sz, jump_scl, dx, dy ;
    double  gx, gy, hxx, hxy, hyy, hdet, tolerance ;

    a = apmtodbl( priz->center->p[0] ) ;
    b = apmtodbl( priz->center->p[1] ) ;
    c = apmtodbl( priz->center->p[2] ) ;
    two_c = 2.0 * c ;

    tolerance = NEWT_TOL * (fabs(a) + fabs(b) + fabs(c)) ;

/*   
        Use Newton's method to try to find a minimum for the
        trace of the matrix beta.
*/

    x = HALF_PI ;
    y = HALF_PI ;

    do {
        /*        components of the gradient.        */

        gx =  -a * cos( x )  -  two_c * cos( x + y ) ;
        gy =  -b * cos( y )  -  two_c * cos( x + y ) ;

        /*         components of the Hessian         */

        hxx =  a * sin( x )  +  two_c * sin( x + y ) ;
        hxy =  two_c * sin( x + y ) ;
        hyy =  b * sin( y )  +  two_c * sin( x + y ) ;
        hdet = hxx * hyy  -  hxy * hxy ;

        /*          A Newton's method step                */
        if( hdet != 0.0 ) {
            dx = (  gx * hyy  -  gy * hxy ) / hdet ;
            dy = ( -gx * hxy  +  gy * hxx ) / hdet ;

            if( (jump_sz = fabs(dx) + fabs(dy)) > MAX_JUMP ) {
                jump_scl = MAX_JUMP / jump_sz ;
                dx *= jump_scl ;
                dy *= jump_scl ;
            }

            x -= dx ;
            y -= dy ;
                
        }
        else {
            fprintf( stderr, "Death during Newton's method. \n" ) ;
            cease() ;
        }
    } while( (fabs(gx) + fabs(gy))  > tolerance ) ;        

/*
        Force the starting point to lie on the line x=y.
*/
    dbltoapm( priz->center->z.u[0], BASE, x ) ;
    dbltoapm( priz->center->z.u[1], BASE, x ) ;

#if DEBUG
    printf( "Herman's starting point : x = %.6e, y= %.6e \n", x, x ) ;
    fflush( stdout ) ;
# endif
}
/* ++++++++++++++++++++++ */

setLLStart( priz )

RPrism  *priz ;
{
/*
        Beware : this function expects to be called AFTER
        setHermStart(), no matter which criterion is in force.
*/
    double  discrim, sqrt_disc, sqrt() ;
    double  a_sin, a_cos, b_sin, b_cos, c_sin, c_cos ;
    double  a, b, c, two_c, x, y ;
    double  jump_sz, jump_scl, dx, dy ;
    double  gx, gy, hxx, hxy, hyy, hdet, tolerance ;
    double  dDisc_dx, dDisc_dy ;

    a = apmtodbl( priz->center->p[0] ) ;
    b = apmtodbl( priz->center->p[1] ) ;
    c = apmtodbl( priz->center->p[2] ) ;
    two_c = 2.0 * c ;

    x = apmtodbl( priz->center->z.u[0] ) ;
    y = apmtodbl( priz->center->z.u[1] ) ;

    tolerance = NEWT_TOL * (a + b + c) ;

    do {
        /*        preliminaries                        */

        a_sin = a * sin( x ) ;
        b_sin = b * sin( y ) ;
        c_sin = two_c * sin( x + y ) ;

        a_cos = a * cos( x ) ;
        b_cos = b * cos( y ) ;
        c_cos = two_c * cos( x + y ) ;

        discrim = ( a_sin - b_sin ) * ( a_sin - b_sin ) +
                                          c_sin * c_sin ;
        sqrt_disc = sqrt( discrim ) ;
        dDisc_dx = a_cos * (a_sin - b_sin) + c_cos * c_sin ;
        dDisc_dy = b_cos * (b_sin - a_sin) + c_cos * c_sin ;

        /*        components of the gradient.        */

        gx = -a_cos - c_cos - dDisc_dx / sqrt_disc ;
        gy = -b_cos - c_cos - dDisc_dy / sqrt_disc ;

        /*         components of the Hessian         */

        hxx = a_sin + c_sin +
              ( a_sin * (a_sin - b_sin) - 
                a_cos * a_cos - c_cos * c_cos + 
                c_sin * c_sin ) / sqrt_disc
              + dDisc_dx * dDisc_dx / (discrim * sqrt_disc) ;

        hxy = c_sin +
              ( a_cos * b_cos + c_sin * c_sin -
                                c_cos * c_cos ) / sqrt_disc
              + dDisc_dx * dDisc_dy / (discrim * sqrt_disc) ;

        hyy = b_sin + c_sin +
              ( b_sin * (b_sin - a_sin) - 
                b_cos * b_cos - c_cos * c_cos + 
                c_sin * c_sin ) / sqrt_disc
              + dDisc_dy * dDisc_dy / (discrim * sqrt_disc) ;

        hdet = hxx * hyy  -  hxy * hxy ;

        /*          A Newton's method step                */
        if( hdet != 0.0 ) {
            dx = (  gx * hyy  -  gy * hxy ) / hdet ;
            dy = ( -gx * hxy  +  gy * hxx ) / hdet ;

            if( (jump_sz = fabs(dx) + fabs(dy)) > MAX_JUMP ) {
                jump_scl = MAX_JUMP / jump_sz ;
                dx *= jump_scl ;
                dy *= jump_scl ;
            }

            x -= dx ;
            y -= dy ;
                
        }
        else {
            fprintf( stderr, "Death during Newton's method. \n" ) ;
            cease() ;
        }
    } while( (fabs(gx) + fabs(gy))  > tolerance ) ;        

/*
        Force the starting point to lie on the line x=y.
*/
    dbltoapm( priz->center->z.u[0], BASE, x ) ;
    dbltoapm( priz->center->z.u[1], BASE, x ) ;

#if DEBUG
    printf( "Least eigenvalue starting point : x = %.6e, y= %.6e \n", x, x ) ;
    fflush( stdout ) ;
# endif
}
/* +++++++++++++++++++++ */

shiftStart( priz )
/* 
    Shift the starting point off the main diagonal. 
*/

RPrism    *priz ;
{
    double    x, y, a, b, amin, bmin ;

    a = apmtodbl( priz->center->p[0] ) ;
    b = apmtodbl( priz->center->p[1] ) ;

    amin = a - apmtodbl( priz->matrix[0] ) ;
    bmin = b - apmtodbl( priz->matrix[MAT_DIM+1] ) ;

    x = apmtodbl( priz->center->z.u[0] ) ;
    y = apmtodbl( priz->center->z.u[1] ) ;

    if( fabs(x - y) < DELTA ) {
        if( amin < bmin ) {
            x += DELTA ;
            y -= DELTA ;
        }
            else {
            x -= DELTA ;
            y += DELTA ;
        }
    }

    dbltoapm( priz->center->z.u[0], BASE, x ) ;
    dbltoapm( priz->center->z.u[1], BASE, y ) ;
}
\end{verbatim}

%% file: appendices/code/lipschitz.tex
\begin{verbatim} 
# include <stdio.h>
# include <math.h>
# include "apm.h"
# include "apmSpecial.h"
# include "converse.h"
# include "bounding.h"
# include "pi.h"

APM        Rdf_sq, Rdf ;
APM        lip_scratch ;
APM        sixteen, eight, four ;
APM        Rdscrm, Rsqrt_disc ;
APM        Rmax_slope, Rmin_slope, Rfirst_slope ;
double  max_slope, min_slope, first_slope ;
RPrism  *earliest ;
Bdd_apm Rmax_btrace, Rmin_btrace, Rfirst_btrace ;
/* +++++++++++++++++++++ */

initLip()
{
/*
    This function depends in detail on the choice of map.
*/

/*
        APM stuff
*/
    four = apmInit( 4L, 0, BASE ) ;
    eight = apmInit( 8L, 0, BASE ) ;
    sixteen = apmInit( 16L, 0, BASE ) ;

    Rmin_slope = apmNew( BASE ) ;        /*  The external APMs */
    Rmax_slope = apmNew( BASE ) ;
    Rfirst_slope = apmNew( BASE ) ;
    Rdf = apmInit( (long)(DEG_FREE), 0, BASE ) ;
    Rdf_sq = apmInit( (long)(DF_SQ), 0, BASE ) ;
    Rstart_size = apmInit( 1L, -START_SZ, BASE ) ;

    Rdscrm = apmNew( BASE ) ;
    Rsqrt_disc  = apmNew( BASE ) ;
    lip_scratch = apmNew( BASE ) ;

    newBapm( Rmax_btrace, BASE ) ;
    newBapm( Rmin_btrace, BASE ) ;
    newBapm( Rfirst_btrace, BASE ) ;

    earliest = conjureRPrism() ;
}
/* +++++++++++++++++++++++ */

setCone( priz ) 

RPrism  *priz ;
/*
        Get the minimum and maximum values for the
        trace of the slope object.  Note that we 
        exploit the symmetry of the potential; the minimum
        and maximum values of the trace of (beta - 2I) have
        the same absolute value.
*/
{
    int  j ;
    APM  *mat_pos ;

    for( j=0 ; j < N_PARMS ; j++ ) 
        apmAssign( earliest->center->p[j], priz->center->p[j] ) ;

    for( j=0 ; j < DEG_FREE ; j++ ) {
        apmAssign( earliest->center->z.v[j],  priz->center->z.u[j] ) ;

    }

    Rglobal_bounds( earliest ) ;
    Rbound_btrace( &Rmin_btrace, earliest ) ;
/*
        Account for the imprecision of the starting point 
        and the variation of the parameters.
*/
    apmAssignLong( lip_scratch, 0L, 0, BASE ) ;
    mat_pos = priz->matrix ;

    for( j=0 ; j < N_PARMS ; j++ ) {
        apmCalc( lip_scratch, lip_scratch, 
                              priz->center->p[j] , Rstart_size,
                              APM_MUL, APM_ADD,
                              *mat_pos,
                              APM_ABS, APM_ADD, NULL ) ;
        mat_pos += 1 + MAT_DIM ;
    }

    apmCalc( Rmin_btrace.lb, Rmin_btrace.lb, lip_scratch, 
                             APM_SUB, NULL ) ;
    apmCalc( Rmin_btrace.ub, Rmin_btrace.ub, lip_scratch, 
                             APM_ADD, NULL ) ;

        /* exploit the symmetry */
    apmSubtract( Rmax_btrace.ub, eight, Rmin_btrace.lb ) ;
    apmSubtract( Rmax_btrace.lb, eight, Rmin_btrace.ub ) ;

    apmCalc( Rdscrm,  Rmax_btrace.lb, APM_DUP, APM_MUL,
                        four, Rdf_sq, APM_MUL, APM_SUB, NULL ) ;

    apmSqrt( Rsqrt_disc, precision, Rdscrm ) ;
    apmAdd( lip_scratch, Rmax_btrace.lb, Rsqrt_disc ) ;
    apmDivide( Rmax_slope, precision, (APM)NULL, lip_scratch, two ) ;

    apmSubtract( lip_scratch, Rmax_btrace.lb, Rsqrt_disc ) ;
    apmDivide( Rmin_slope, precision, (APM)NULL, lip_scratch, two ) ;

    min_slope = apmtodbl( Rmin_slope ) ;
    max_slope = apmtodbl( Rmax_slope ) ;

}
/* +++++++++++++++++++++++++++++++++++++++++++ */

setSlopes( priz )

RPrism  *priz ;
/*
        Recall that our orbit will, at the beginning of
        a round of orbit-following, have just passed through a
        point on the torus whose beta will diminish the 
        slope. This implies that the slope is already smaller 
        than the value of max_slope found above.  Calculate
        a better upper bound on what the slope could be and
        store it in first_slope and Rfirst_slope.
*/
{
    int  j, mat_pos ;

    for( j=0 ; j < N_PARMS ; j++ )  {
        apmAssign( earliest->center->p[j], priz->center->p[j] ) ;

        mat_pos = j * (MAT_DIM + 1) ;
        apmAssign( earliest->matrix[mat_pos], priz->matrix[mat_pos] ) ;
    }

    for( j=0 ; j < DEG_FREE ; j++ ) {
        apmAssign( earliest->center->z.v[j], priz->center->z.u[j] ) ;

        /* 
                Account for imprecision in the starting point.
         */
        mat_pos = STAID_LEN + TWO_DF*MAT_DIM + 
                  N_PARMS + DEG_FREE + j * (MAT_DIM + 1) ;
        apmAssign( earliest->matrix[mat_pos], Rstart_size ) ;
    }
    
    Rglobal_bounds( earliest ) ;
    Rbound_btrace( &Rfirst_btrace, earliest ) ;
    
    apmDivide( lip_scratch, precision, (APM)NULL, Rdf_sq, Rmax_slope ) ;
    apmCalc( Rfirst_slope, Rfirst_btrace.ub, lip_scratch, APM_SUB,
                                            max_error, APM_ADD, NULL ) ;

    first_slope = apmtodbl( Rfirst_slope ) + DBL_ERR ;

}
\end{verbatim}

%% file: appendices/code/lambda.tex
\begin{verbatim} 
# include <stdio.h>
# include <math.h>
# include "apm.h"
# include "apmSpecial.h"
# include "converse.h"
# include "bounding.h"
# include "rows.h"

APM            Rsqrt_disc ;
APM            Ra_term, Rb_term, Rc_term ;
APM            Rtrace_ll, RminBlam_ll, RmaxBlam_ll, Rdenom ;
Bdd_apm            RBtrace, RminLam, RmaxLam ;
RPrism            *earliest ;

Bdd_dbl            discrim ;
Bdd_dbl     a_sq, b_sq, c_sq ;
Bdd_dbl            *lamFacts[2] ;
Bdd_term    ab_term ;

APM            four, lam_scratch ;
Bdd_apm            Rdiscrim ;
Bdd_apm     Ra_sq, Rb_sq, Rc_sq ;
Bdd_apm            *RlamFacts[2] ;
Bapm_term   Rab_term ;

APM            RfirstLeastLam, RminLeastLam, RmaxLeastLam, RsumTinyLams ;
double            firstLeastLam, minLeastLam, maxLeastLam, sumTinyLams ;
/* ++++++++++++++++++++++++++++++ */

initLambda()
{
/*
        Do up the APMs
*/
    Ra_term = apmNew( BASE ) ;
    Rb_term = apmNew( BASE ) ;
    Rc_term = apmNew( BASE ) ;

    Rdenom = apmNew( BASE ) ;
    Rtrace_ll = apmNew( BASE ) ;
    Rsqrt_disc = apmNew( BASE ) ;
    RminBlam_ll = apmNew( BASE ) ;
    RmaxBlam_ll = apmNew( BASE ) ;

    RminLeastLam = apmNew( BASE ) ;
    RmaxLeastLam = apmNew( BASE ) ;
    RsumTinyLams = apmNew( BASE ) ;
    RfirstLeastLam = apmNew( BASE ) ;

    newBapm( Ra_sq, BASE ) ;
    newBapm( Rb_sq, BASE ) ;
    newBapm( Rc_sq, BASE ) ;
    newBapm( RmaxLam, BASE ) ;
    newBapm( RminLam, BASE ) ;
    newBapm( RBtrace, BASE ) ;
    newBapm( Rdiscrim, BASE ) ;

    four = apmInit( 4L, 0, BASE ) ;
    lam_scratch = apmNew( BASE ) ;

    earliest = conjureRPrism() ;
/*
        Set up the terms.
*/

    ab_term.nfactors = Rab_term.nfactors = 2 ;
    ab_term.factors = lamFacts ;
    Rab_term.factors = RlamFacts ;
    ab_term.coef = -2.0 ;
    Rab_term.coef = apmInit( -2L, 0, BASE ) ;
    newBapm( Rab_term.bound, BASE ) ;
        ab_term.factors[0] = &a_sin.bound ;
        ab_term.factors[1] = &b_sin.bound ;
        Rab_term.factors[0] = &Ra_sin.bound ;
        Rab_term.factors[1] = &Rb_sin.bound ;
}
/* ++++++++++++++++++++++++ */

Rbd_Blams( leastBlam, bigBlam, trace )

Bdd_apm         *leastBlam, *trace, *bigBlam ;
/*
        An APM partner to bd_Blams ;
*/
{
/* Bound the terms for the discriminant. */
    RsetSq( &Ra_sq, &Ra_sin.bound ) ;
    RsetSq( &Rb_sq, &Rb_sin.bound ) ;
    RsetSq( &Rc_sq, &Rc_sin.bound ) ;
    Rbound_term( &Rab_term ) ;

/* Bound the discriminant itself. */

        /* lower bound */
    apmCalc( Rdiscrim.lb, Ra_sq.lb, Rb_sq.lb, APM_ADD,
                          four, Rc_sq.lb, APM_MUL, APM_ADD,
                          Rab_term.bound.lb, APM_ADD, NULL ) ;

    if( apmCompare( Rdiscrim.lb, zero ) < 1 )
        apmAssign( Rdiscrim.lb, zero ) ;

        /* upper bound */
    apmCalc( Rdiscrim.ub, Ra_sq.ub, Rb_sq.ub, APM_ADD,
                          four, Rc_sq.ub, APM_MUL, APM_ADD,
                          Rab_term.bound.ub, APM_ADD, NULL ) ;

    if( apmCompare( Rdiscrim.ub, zero ) < 1 )
        apmAssign( Rdiscrim.ub, zero ) ;

/* Do up the final bounds on the eigenvalues. 
                First do those requiring 
                sqrt( discrim.lb ). 
*/
    apmSqrt( Rsqrt_disc, precision, Rdiscrim.lb ) ;
    apmCalc( lam_scratch, trace->ub, Rsqrt_disc, APM_SUB, 
                                        max_error, APM_ADD, NULL ) ;
    apmDivide( leastBlam->ub, precision, (APM)NULL, lam_scratch, two ) ;

    apmCalc( lam_scratch, trace->lb, Rsqrt_disc, APM_ADD, 
                                        max_error, APM_SUB, NULL ) ;
    apmDivide( bigBlam->lb, precision, (APM)NULL, lam_scratch, two ) ;

/*
                Next those requiring
                sqrt( discrim.lb )
*/
    apmSqrt( Rsqrt_disc, precision, Rdiscrim.ub ) ;
    apmCalc( lam_scratch, trace->lb, Rsqrt_disc, APM_SUB, 
                                        max_error, APM_SUB, NULL ) ;
    apmDivide( leastBlam->lb, precision, (APM)NULL, lam_scratch, two ) ;

    apmCalc( lam_scratch, trace->ub, Rsqrt_disc, APM_ADD, 
                                        max_error, APM_ADD, NULL ) ;
    apmDivide( bigBlam->ub, precision, (APM)NULL, lam_scratch, two ) ;

}
/* ++++++++++++++++++++++++++++ */

setLLbounds( priz ) 
/*
    Get bounds on the least eigenvalue of the variation of the action
    functional.  This is equivalent to the summer's estimate of the 
    value of size of the perturbation for which no  minimizing state
    can include the maximum of the perturbation.
*/

RPrism  *priz ;
{
    int            j, mat_pos ;
    APM            *pmat_pos ;

    for( j=0 ; j < N_PARMS ; j++ ) 
        apmAssign( earliest->center->p[j], priz->center->p[j] ) ;

        mat_pos =  j * (MAT_DIM + 1) ;
        apmAssign( earliest->matrix[mat_pos],  priz->matrix[mat_pos] ) ;

    for( j=0 ; j < DEG_FREE ; j++ )  
        apmAssign( earliest->center->z.v[j], priz->center->z.u[j] ) ;

/*

    Rglobal_bounds( earliest ) ;
    Rbound_btrace( &RBtrace, earliest ) ;
    Rbd_Blams( &RminLam, &RmaxLam, &RBtrace ) ;

/*
        Account for the imprecision of the starting point 
        and the variation of the parameters.
*/
    apmAssignLong( lam_scratch, 0L, 0, BASE ) ;
    pmat_pos = priz->matrix ;

    for( j=0 ; j < N_PARMS ; j++ ) {
        apmCalc( lam_scratch, lam_scratch, 
                              priz->center->p[j] , Rstart_size,
                              APM_MUL, APM_ADD,
                              *pmat_pos,
                              APM_ABS, APM_ADD, NULL ) ;
        pmat_pos += 1 + MAT_DIM ;
    }

    apmCalc( RminLam.lb, RminLam.lb, lam_scratch, APM_SUB, NULL ) ;
    apmCalc( RminLam.ub, RminLam.ub, lam_scratch, APM_ADD, NULL ) ;

/*
        Exploit the symmetry of the example.  The 
        largest value for an eigenvalue is 
        4.0 - (leastLam.lb).

        The calculation above assumes that the
        u part of the prism's center contains a
        starting point suitable for a least-eigenvalue
        kind of test, i.e. the point where the least ev
        attains its minimum.  The bdd_apm RmaxLam will
        contain information about the largest ev of beta 
        at the spot where leastLam is small.  To get the
        thing we really want for the calculations
        below we must exploit the symmetry described
        above.
*/
    apmSubtract( RmaxLam.ub, four, RminLam.lb ) ;
    apmCalc( Rdiscrim.ub, RmaxLam.ub, APM_DUP, APM_MUL,
                        four, APM_SUB, NULL ) ;
    apmSqrt( Rsqrt_disc, precision, Rdiscrim.ub ) ;

/* 
        A global lower bound - if the least eigenvalue of
        one of the diagonal blocks (see notes, Jan 10 )
        slips below this value then the next block is 
        sure to have a negative eigenvalue.
*/
    apmSubtract( lam_scratch, RmaxLam.ub, Rsqrt_disc ) ;
    apmDivide( RminLeastLam, precision, (APM) NULL, lam_scratch, two ) ;
    apmCalc( RminLeastLam, RminLeastLam, max_error, APM_SUB, NULL ) ;
    minLeastLam = apmtodbl( RminLeastLam ) ;

        /* 
            A lower bound on the sum of the non-maximal eigenvalues
            of a diagonal block.
        */
    sumTinyLams = minLeastLam ;
    apmAssign( RsumTinyLams, RminLeastLam ) ;

/*
        A global upper bound.
*/
    apmAdd( lam_scratch, RmaxLam.ub, Rsqrt_disc ) ;
    apmDivide( RmaxLeastLam, precision, (APM) NULL, lam_scratch, two ) ;
    apmCalc( RmaxLeastLam, RmaxLeastLam, max_error, APM_ADD, NULL ) ;
    maxLeastLam = apmtodbl( RmaxLeastLam ) ;

}
/* ++++++++++++++++++++++++++++++ */

RsetSq( xsq, x ) 

Bdd_apm           *x, *xsq ;
{
    if( apmCompare( x->ub, zero ) > 0 ) {
        if( apmCompare( x->lb, zero ) > 0 ) {
            apmMultiply( xsq->ub, x->ub, x->ub ) ;
            apmMultiply( xsq->lb, x->lb, x->lb ) ;
        }
        else {
            apmAbsoluteValue( lam_scratch, x->lb ) ;
            if( apmCompare( x->ub, lam_scratch ) > 0 ) {
                apmMultiply( xsq->ub, x->ub, x->ub ) ;
                apmAssign( xsq->lb, zero ) ;
            }
            else {
                apmMultiply( xsq->ub, x->lb, x->lb ) ;
                apmAssign( xsq->lb, zero ) ;
            }
        }
    }
    else {
        apmMultiply( xsq->ub, x->lb, x->lb ) ;
        apmMultiply( xsq->lb, x->ub, x->ub ) ;
    }
}
/* ++++++++++++++++++++++++++++++++ */

setLeastLam( priz ) 

RPrism  *priz ;
/* 
        Calculate an upper bound on the largest eigenvalue of beta
        at the initial point, then use it and the global bound,
        maxLeastLam to set firstLeastLam.
*/
{
    int            j, mat_pos ;

    for( j=0 ; j < N_PARMS ; j++ )  {
        earliest->center->p[j] = priz->center->p[j] ;

        mat_pos = j * (MAT_DIM + 1) ;
        earliest->matrix[mat_pos] = priz->matrix[mat_pos] ;
    }

    for( j=0 ; j < DEG_FREE ; j++ ) 
        earliest->center->z.v[j] = priz->center->z.u[j] ;

    Rglobal_bounds( earliest ) ;
    Rbound_btrace( &RBtrace, earliest ) ;
    Rbd_Blams( &RminLam, &RmaxLam, &RBtrace ) ;

/* 
        Obtain an upper bound on the least
        eigenvalue of the block of the Hessian of
        the action functional corresponding to the
        starting point. As in the functions in follow.c,
        compute a whole suite of estimates and choose
        the best one.

*/

        /*
                Rdenom is a global upper bound
                on the size of the largest eigevalue
                of a diagonal block.   
                    Rdenom = maximum trace - (n-1) * minimum ev.

                It's used together with the least eigenvalue
                of beta (evaluated at the starting point) :

                    LeastLam <= RminBlam.ub - 1.0 / Rdenom
        */
    apmCalc( Rdenom, Rdf, one, APM_SUB, 
                     RminLeastLam, APM_MUL, APM_NEG, 
                     Rmax_slope, APM_ADD, NULL) ;
    apmDivide( lam_scratch, precision, (APM) NULL, one, Rdenom ) ;
    apmSubtract( RminBlam_ll, RminLam.ub, lam_scratch ) ;

        /*
                Here we try to attain a small estimate by 
                saying :
                    LeastLam <= RmaxBlam.ub - 1.0 / maxLeastLam.
        */
    apmDivide( lam_scratch, precision, (APM) NULL, one, RmaxLeastLam ) ;
    apmSubtract( RmaxBlam_ll, RmaxLam.ub, lam_scratch ) ;

        /*
                Finally we make the estimate
                    LeastLam <= first_slope / DEG_FREE
        */
    apmDivide( Rtrace_ll, precision, (APM)NULL, Rfirst_slope, Rdf ) ;

        /* 
                Choose the best (smallest) lower bound.
        */
    apmAssign( RfirstLeastLam, RmaxBlam_ll ) ;
    if( apmCompare( RfirstLeastLam, RminBlam_ll ) == 1 ) 
        apmAssign( RfirstLeastLam, RminBlam_ll ) ;
    if( apmCompare( RfirstLeastLam, Rtrace_ll ) == 1 ) 
        apmAssign( RfirstLeastLam, Rtrace_ll ) ;

    firstLeastLam = apmtodbl( RfirstLeastLam ) ;

}
\end{verbatim}

%% file: appendices/code/trace.tex
\begin{verbatim} 
# include <stdio.h>
# include <math.h>
# include "apm.h"
# include "converse.h"
# include "map.h"
# include "bounding.h"
# include "rows.h"

Bdd_apm            *Rfact_buf[NUM_FACTS] ;
Bapm_expr   Rb_trc ;
Bapm_term   Rtrace_terms[NUM_TERMS] ;
/* ++++++++++++++++++++++++++++++ */

initTrace()
{
    int                j ;                
    Bdd_apm        **Rfpt ;
/*
        Set up the expressions.
*/
    Rb_trc.nterms = NUM_TERMS ;

    Rb_trc.const = apmInit( 4L, 0, BASE ) ;

    newBapm( Rb_trc.bound, BASE ) ;

    Rb_trc.terms = Rtrace_terms ;
 
 /*
        Set up their terms.
*/
    Rfpt = Rfact_buf ;
    for( j=0 ; j < NUM_TERMS ; j++ ) {
        Rtrace_terms[j].nfactors = 1 ;
        Rtrace_terms[j].coef = apmInit( -1L, 0, BASE ) ;
        Rtrace_terms[j].factors = Rfpt ;

        newBapm( Rtrace_terms[j].bound, BASE ) ;

        Rfpt++ ;
    }
/*
        Fix up the constant in the third term . . . it should be
        -2.0.
*/
    apmAssignLong( Rtrace_terms[2].coef, -2L, 0, BASE ) ;
/*
        Associate the factors - which are only pointers 
        to bounded objects - to genuine, properly initialized objects.
*/
    /* first term */
    Rb_trc.terms[0].factors[0] = &Ra_sin.bound ;

    /* second term */
    Rb_trc.terms[1].factors[0] = &Rb_sin.bound ;

    /* third term */
    Rb_trc.terms[2].factors[0] = &Rc_sin.bound ;
}
/* ++++++++++++++++++++++++ */

Rbound_btrace( result, priz )

RPrism   *priz ;
Bdd_apm  *result ;
/*
        An APM partner to bound_btrace.  Some of the variables 
        used here are defined above.
*/
{
        /* Bound the expression */
    Rbound_expr( &Rb_trc ) ;
    apmCalc( Rb_trc.bound.ub, Rb_trc.bound.ub, max_error, APM_ADD, NULL ) ;
    apmCalc( Rb_trc.bound.lb, Rb_trc.bound.lb, max_error, APM_SUB, NULL ) ;

    apmAssign( result->ub, Rb_trc.bound.ub ) ;
    apmAssign( result->lb, Rb_trc.bound.lb ) ;
}
\end{verbatim}

%% file: appendices/code/converseh.tex
\begin{verbatim} 
# ifndef YES
# endif

# ifndef WORKED
# endif

                                   been considered, is too hard to
                                   decide, is under active 
                                   consideration, or is  equivalent
                                   to some symmetrically related,
                                   other prism.                            */

/*
        Data types for non-rigorous, rough calculations
                                                                */

typedef  double   *Tor_pt, *Parm_pt ;

typedef  struct { Tor_pt        u, v ; } Embed_pt ;

typedef  struct { Embed_pt        z ;
                  Parm_pt        p ; } Xtnd_pt ;

typedef  struct  prsm {  int                in_torus, n_cuts ;
                         char                *cuts[N_PARMS+TWO_DF] ;
                         double                *matrix ;
                         Xtnd_pt        *center ;
                         struct prsm        *next ;     } Prism ;
/*
        Data types for rigorous, arbitrary precision, calculations 
                                                                */

typedef  APM   *RTor_pt, *RParm_pt ;

typedef  struct { RTor_pt        u, v ; } REmbed_pt ;

typedef  struct { REmbed_pt        z ;
                  RParm_pt        p ; } RXtnd_pt ;

typedef  struct  Rprsm {  int                in_torus, n_cuts ;
                          APM                *matrix ;
                          char                *cuts[MAT_DIM] ;
                          RXtnd_pt        *center ;
                          struct Rprsm        *next ;     } RPrism ;
/* +++++++++++++++++++++++++++++++++++++++++++++++++++++++++++++++ */

extern  Prism    *conjurePrism() ;
extern  RPrism   *conjureRPrism() ;

/*
     Some variables used throughout the converse KAM calculations 
                                                                    */
extern  int         do_graph, do_backup, restoration ;
extern int        precision, depth, furthest, terse, stubborn ;
extern  int        quick_tries, tries, Rtries, max_steps, max_NTsteps ;
extern  int        HermSuccess, LLSuccess, ll_used[3], most_cuts ;
extern  int        (* fatten)(), (* row_sums)() ; 
extern  int        fxed_form(), Rfxed_form(), col_rotor(), Rcol_rotor() ;
extern  int        ff_rows(), Rff_rows(), cr_rows(), Rcr_rows() ;
extern  APM        Rfirst_slope, Rmin_slope, Rmax_slope, Rdf, Rdf_sq ;
extern  APM        RminLeastLam, RmaxLeastLam, RfirstLeastLam, RsumTinyLams ;
extern  APM     half, max_error, Rstart_size, RSmBlock_err, RBgBlock_err ;
extern  char        *graf_file, *back_name, *rest_name, *parm_names[] ;
extern  double  firstLeastLam, minLeastLam, maxLeastLam, sumTinyLams ;
extern  double  first_slope, min_slope, max_slope ;
extern  double  apmtodbl(), parm_roof[], parm_floor[] ;
extern  double  SmBlock_err, BgBlock_err ;
\end{verbatim}

%% file: appendices/code/converse.tex
\begin{verbatim} 
# include  <stdio.h>
# include  <math.h>
# include  "apm.h"
# include  "converse.h"
# include  "tree.h"

int        do_graph, do_backup, restoration ;
int        precision, depth, err_hndlr, furthest ;
int        stubborn, terse ;
APM        max_error, RSmBlock_err, RBgBlock_err ;
double  SmBlock_err = DF_SQ * DBL_ERR ;
double  BgBlock_err = DEG_FREE * N_PARMS * DBL_ERR ;
/* ++++++++++++++++++++++++++++ */

main (argc, argv)

int   argc ;
char  *argv[] ;
{
    int             verdict, Rverdict, tree_verdict, nsteps ;
    Prism    *image_prism ;
    RPrism   *active_prism, *old_prism ;

    handle_opts( argc, argv ) ;
    active_prism = conjureRPrism() ;
    image_prism = conjurePrism() ;
    commence( active_prism ) ;

        /* Study the current prism, cutting it up if need be */
    while( active_prism != NULL ) {
    /*
        Try a preliminary, non-rigorous calculation to see if 
        prospects are good.  If they are, do a rigorous check.
        If they aren't, try to refine the prism.  If it has already
        been refined enough, just give up.
                                                                    */
        if( do_graph == YES ) 
            graphPrism( active_prism, ACTIVE ) ;
/*
        Check the tree to see if an equivalent prism 
        is already finished.  If so, record the result 
        and press on.  If not, do a detailed analysis.
*/
        tree_verdict = consultTree( active_prism ) ;

# if FANCY_TREE
        if( (tree_verdict == MAYBE) || (tree_verdict == NO_TORI) ) {
            if( do_graph == YES ) 
                graphPrism( active_prism, tree_verdict ) ;
            if( do_backup == YES )
                make_backup( active_prism ) ;

            old_prism = active_prism ;
            active_prism = old_prism->next ;

            old_prism->in_torus = tree_verdict ;
            if( terse == NO ) 
                printRPrism( old_prism, 0 ) ;
            releaseRPrism( old_prism ) ;
        }
# else
        if( tree_verdict == MAY_SKIP ) {
            if( do_graph == YES ) 
                graphPrism( active_prism, SYMMTRC ) ;
            if( do_backup == YES )
                make_backup( active_prism ) ;

            old_prism = active_prism ;
            active_prism = old_prism->next ;

            releaseRPrism( old_prism ) ;
        }
# endif

        else {
            prepare_trial( active_prism ) ;
            verdict = try_prism( active_prism, image_prism, &nsteps ) ;

            Rverdict = UNTRIED ;
            if( verdict == NO_TORI ) {
                Rverdict = Rtry_prism( active_prism, image_prism, &nsteps ) ;
                if( Rverdict == NO_TORI ) {
                    active_prism->in_torus = NO_TORI ;
# if FANCY_TREE
                    colorLeaf( active_prism ) ;
# endif
                    if( terse == NO ) 
                        printRPrism( active_prism, nsteps );
                    if( do_graph == YES )
                        graphPrism( active_prism, NO_TORI ) ;
                    if( do_backup == YES )
                        make_backup( active_prism ) ;

                    old_prism = active_prism ;
                    active_prism = old_prism->next ;
                    releaseRPrism( old_prism ) ;
                }
# if TATTLE
                else {
                    printf( 
                        "Disagreement between try() and Rtry(). \n" ) ;
                    printf( "disputed prism : \n\t" ) ;
                    printRPrism( active_prism, nsteps ) ;
                    fflush( stdout ) ;
                }
# endif

            }

            if( (Rverdict == MAYBE) || (verdict == MAYBE) ) {

                    /* Either refine the prism . . . */
                if( may_refine(active_prism) == YES ) {
                    refinePrism( active_prism, image_prism ) ;
                    if( do_graph == YES ) {
                        graphPrism( active_prism->next, UNTRIED ) ;
                        graphPrism( active_prism, ACTIVE ) ;
                    }
                }

                    /* . . . or give up and move on. */
                else {  
                    if( do_graph == YES ) 
                        graphPrism( active_prism, MAYBE ) ;
                    if( do_backup == YES )
                        make_backup( active_prism ) ;

                    active_prism->in_torus = MAYBE ;
                    moveEdge_o_Chaos( active_prism, nsteps ) ;
                    if( terse == NO ) 
                        printRPrism( active_prism, nsteps ) ;

                    old_prism = active_prism ;
                    active_prism = old_prism->next ;
# if FANCY_TREE
                    colorLeaf( old_prism ) ;
# endif
                    releaseRPrism( old_prism ) ;
                }
            }
        }
    }

    cease() ;
}
\end{verbatim}

%% file: appendices/code/follow.tex
\begin{verbatim} 
# include <stdio.h>
# include <math.h>
# include "apm.h"
# include "apmSpecial.h"
# include "converse.h"
# include "bounding.h"
# include "rows.h"
# include "pi.h"

                           to be used in determining how long 
                           quick_try should go.
                        */

                                ? furthest : ((n/QS_TO_RS)+3) )

/*
        Declarations for some external variables 
        mentioned in converse.h.  The APMs are initialized by
        initFollowing().
*/
/*
    The functions in this file manipulate copies of the data 
    passed to them. The copies are kept in Prisms and RPrisms
    gotten with the conjuring functions by initFollowing().
*/
Prism          *workPriz[2] ;

double     b_buf[DF_SQ], *b_ptrs[DF_SQ] ;
double           parmbuf[2*N_PARMS], coordbuf[2*TWO_DF] ;

Xtnd_pt   xpt_a, xpt_b ;
/*
        Some APM variables needed for orbit
        following and  slope watching.
*/

RPrism     *Rwork[2] ;

APM           f_scratch, Rdenom ;
APM           Rsum, Rmax_sum ;
APM           Rtrace_ll, RmaxBlam_ll, RminBlam_ll ;
double     trace_ll, maxBlam_ll, minBlam_ll ;

/*
        The variables declared below don't really need to 
        be bounded objects (they did in an earlier version of the code),
        but the .ub in their uses makes the code easier to understand.
*/
Bdd_dbl           b_trace, minBlam, maxBlam, leastLam, slope ;
Bdd_apm           Rb_trace, RminBlam, RmaxBlam, RleastLam, Rslope ;

int     is_first_trial = YES ;
int        local_furth, ll_used[3] ;
int        HermSuccess, LLSuccess ;
int        max_steps, max_NTsteps, tries, Rtries, quick_tries, most_cuts ;
/* +++++++++++++++++++++++++++++++++ */

prepare_trial( priz ) 

RPrism  *priz  ;
{
    int   j ; 

    if( areNewParms( priz ) == YES ) {

/*
        Unless this is the very first prism, 
        record the center point - it will be moved by 
        setHermStart() and setLLStart() and will neeed to be
        restored to its correct value.
*/
        if( is_first_trial == NO ) { 
            for ( j=0 ; j < DEG_FREE ; j++ ) {
                    apmAssign( xpt_a.z.u[j], priz->center->z.u[j] ) ;
                    apmAssign( xpt_a.z.v[j], priz->center->z.v[j] ) ;
            }
        }

        setHermStart( priz ) ;
        setCone( priz ) ;
# if USE_LL
        setLLStart( priz ) ;
        setLLbounds( priz ) ;
# endif
# if USE_SHIFT
        shiftStart( priz ) ;
# endif
/*
             Unless this is the very first trial, restore the 
             correct value of the centerpoint before evaluating
             the initial estimates for the slope and least eigenvalue.
*/
        if( is_first_trial == YES )
             is_first_trial = NO ;
        else {
            for ( j=0 ; j < DEG_FREE ; j++ ) {
                    apmAssign( priz->center->z.u[j], xpt_a.z.u[j]  ) ;
                    apmAssign( priz->center->z.v[j], xpt_a.z.v[j]  ) ;
            }
        }

        setSlopes( priz ) ;

# if USE_LL
        setLeastLam( priz ) ; 
# else
        firstLeastLam = 1.0 ;
        minLeastLam = 0.5 ;
        dbltoapm( RfirstLeastLam, BASE, firstLeastLam ) ;
        dbltoapm( RminLeastLam, BASE, minLeastLam ) ;
# endif
    }
}
/* +++++++++++++++++++++++++++++++++ */

initFollowing()
{
/*
        Set up the correct connections between the various
        static variables in this file.
*/
    int              j, all_well ;

    all_well = YES ;
/*
        Set up the working prisms.
*/
    workPriz[0] = conjurePrism() ;
    workPriz[1] = conjurePrism() ;
    if( (workPriz[0] == NULL) || (workPriz[1] == NULL) ) 
        all_well = NO ;

/*
        Set up the APM stuff
*/
    f_scratch = apmNew( BASE ) ;
    Rdenom = apmNew( BASE ) ;

    Rtrace_ll = apmNew( BASE ) ;
    RminBlam_ll = apmNew( BASE ) ;
    RmaxBlam_ll = apmNew( BASE ) ;

    newBapm( Rslope, BASE ) ;
    newBapm( Rb_trace, BASE ) ;
    newBapm( RminBlam, BASE ) ;
    newBapm( RmaxBlam, BASE ) ;  
    newBapm( RleastLam, BASE ) ;

# if (USE_LL == NO)
    apmAssignLong( RleastLam.ub, 1L, 0, BASE ) ;
    apmAssignLong( RleastLam.lb, 1L, 0, BASE ) ;
# endif

    Rsum = apmNew( BASE ) ;
    Rmax_sum = apmNew( BASE ) ;
    dbltoapm( Rmax_sum, BASE, MAX_SUM ) ;
    
    Rwork[0] = conjureRPrism() ;
    Rwork[1] = conjureRPrism() ;
    if( (Rwork[0] == NULL) || (Rwork[1] == NULL) )
        all_well = NO ;

/*
        Set up the extended points - they're used by 
        quick_test(), and are pointed to by the
        "center" attributes of the working prisms. 
*/
    xpt_a.z.u = coordbuf ;
    xpt_a.z.v = coordbuf + DEG_FREE ;
    xpt_a.p   = parmbuf ;

    xpt_b.z.u = coordbuf + TWO_DF ;
    xpt_b.z.v = coordbuf + TWO_DF + DEG_FREE ;
    xpt_b.p   = parmbuf + N_PARMS ;
/*
        Set up pointers to the matrix which receives the 
        changeable parts of the jacobian; the one called 
        "beta" in most of my notes.  
*/
    for( j=0 ; j < (sizeof( b_buf ) / sizeof( double )) ; j++ )
        b_ptrs[j] = &b_buf[j] ;
}
/* ++++++++++++++++++++++++ */

Rtry_prism( initial_priz, final_priz, nsteps ) 

int         *nsteps ;
Prism    *final_priz ;
RPrism   *initial_priz ;
/*
     Rigorously decides whether a prism of initial data may 
     contain any invariant Lagrangian tori, an APM version of
     the routine tryPrism() above.
*/
{
    int                count ;
    RPrism      *priz, *priz_prime, *temp_priz ;

    Rtries++ ;
    priz = Rwork[0] ;
    priz_prime = Rwork[1] ;

/*
        Note that Rtry_prism() does not call setSlopes,setStart or
        setCone.  All that should have been done with a call to 
        prepare_trial().
*/

    isNewPrism = YES ;
    RcopyRPrism( priz, initial_priz ) ;

    fatten = Rfxed_form ;
    row_sums = Rff_rows ;

    *nsteps = count = 1 ;
    apmAssign( Rslope.ub, Rfirst_slope ) ;
    apmAssign( RleastLam.ub, RfirstLeastLam ) ;
    if( apmCompare(Rslope.ub, Rmin_slope) == -1 ) {
        HermSuccess++ ;
        copyRPrism( final_priz, priz ) ;
        return( NO_TORI ) ;
    }
    if( apmCompare(RleastLam.ub, RminLeastLam) == -1 ) {
        LLSuccess++ ;
        copyRPrism( final_priz, priz ) ;
        return( NO_TORI ) ;
    }

# if (USE_RIGOR == NO)
    copyRPrism( final_priz, priz ) ;
    return( NO_TORI ) ;
# endif

    while( big_RPrism( priz ) == NO ) {
/*
                Check the slope.
*/
        count++ ;
/*
                Calculate some bounds useful for both criteria.
*/
        Rglobal_bounds( priz ) ;
        Rbound_btrace( &Rb_trace, priz ) ;

# if USE_LL
                /* mrm's condition  */
        Rbd_Blams( &RminBlam, &RmaxBlam, &Rb_trace ) ;

        apmDivide( f_scratch, precision, (APM)NULL, one,
                                         RleastLam.ub ) ;
        apmSubtract( RmaxBlam_ll, RmaxBlam.ub, f_scratch ) ;

        apmSubtract( Rdenom, Rslope.ub, RsumTinyLams ) ;
        if( apmCompare( Rdenom, zero) > 0 ) {
            apmDivide( f_scratch, precision, (APM) NULL, one, Rdenom ) ;
            apmSubtract( RminBlam_ll, RminBlam.ub, f_scratch ) ;
        }
        else 
            apmAssign( RminBlam_ll, zero ) ;
# endif

                /* Herman's condition */
        apmDivide( f_scratch, precision, (APM) NULL, Rdf_sq, Rslope.ub ) ;
        apmSubtract( Rslope.ub, Rb_trace.ub, f_scratch ) ;

# if USE_LL
        apmDivide( Rtrace_ll, precision, (APM)NULL, Rslope.ub, Rdf ) ;

        Rbest_ll( RleastLam.ub, RmaxBlam_ll,
                    RminBlam_ll, Rtrace_ll ) ;

# endif

/*                 
                Do some truncations to speed things up
*/
# if USE_LL
        apmTruncate( RleastLam.ub, precision ) ;
# endif
        apmTruncate( Rslope.ub, precision ) ;

        if( apmCompare(Rslope.ub, Rmin_slope) == -1 ) {
            *nsteps = count ;
            if( count > max_NTsteps ) 
                max_NTsteps = count ;

            HermSuccess++ ;
            copyRPrism( final_priz, priz ) ;
            return( NO_TORI ) ;
            }
        else if( apmCompare(RleastLam.ub, RminLeastLam) == -1 ) {
            *nsteps = count ;
            if( count > max_NTsteps ) 
                max_NTsteps = count ;

            LLSuccess++ ;
            copyRPrism( final_priz, priz ) ;
            return( NO_TORI ) ;
        }
        else {
            if( count == furthest )
                break ;

            Rprismatic_image( priz_prime, priz ) ;

            m_swap( priz, priz_prime, temp_priz ) ;
        }
# if USE_CR
        if( count > FF_CYCLS )  {
            fatten = Rcol_rotor ;
            row_sums = Rcr_rows ;
        }
# endif
    }
  
    *nsteps = count ;
    copyRPrism( final_priz, priz ) ;
    return( MAYBE ) ;
}
/* +++++++++++++++++++++++++++++++++++++++ */

big_RPrism( Priz )

RPrism  *Priz ;
{
    APM   *Rrpt, *Rend_mat, *Rend_row ;

    Rend_mat = Priz->matrix + MAT_SZ ;
    for( Rrpt = Priz->matrix ; Rrpt < Rend_mat ; ) {
        apmAssignLong( Rsum, 0L, 0, BASE ) ;
        for( Rend_row = Rrpt + MAT_DIM ; Rrpt < Rend_row ; Rrpt++ )
            apmCalc( Rsum, Rsum, *Rrpt, APM_ABS, APM_ADD, NULL ) ;

        if( apmCompare( Rsum, Rmax_sum) == 1 ) 
            return( YES ) ;
    }

    return( NO ) ;
}
/* ++++++++++++++++++++++ */

Rbest_ll( best, minBlam_ll, maxBlam_ll, trace_ll ) 

APM   best, minBlam_ll, maxBlam_ll, trace_ll ;
{
    apmAssign( best, maxBlam_ll ) ;
    if( apmCompare( best, minBlam_ll ) == 1 ) 
        apmAssign( best, minBlam_ll ) ;

    if( apmCompare( best, trace_ll ) == 1 ) 
        apmAssign( best, trace_ll ) ;
}
\end{verbatim}

%% file: appendices/code/maph.tex
\begin{verbatim} 
extern APM        RDeriv[], *Rbeta_ptrs[], *Rgamma_ptrs[] ;
extern double   Deriv[], *beta_ptrs[], *gamma_ptrs[] ;
\end{verbatim}

%% file: appendices/code/map.tex
\begin{verbatim} 
/*
        Functions to perform the extended Froeschle map and to
        calculate its jacobian.  Each function has a rigorous
        and a non-rigorous form; the former always has a name
        beginning with a "R".

        The functions in this file are quite specific - 
        they pertain to maps  of the form

        (p,u,v) -> (p',u',v')

        p' = p
        u' = v
        v' = 2v - u -dV(v)

        where u, v, u' anf v' are all in 2d Euclidean space,
        p is an element of a space of parameters and
            V(v) = -a * sin( v[0] ) + -b * sin( v[1] ) +
                   -c * sin( v[0] + v[1] )

        The parameters a, b, and c are always passed through
        an array called "parms" with 

        a = parms[0], b = parms[1], c = parms[2].
*/
# include <stdio.h>
# include <math.h>
# include "apm.h"
# include "apmSpecial.h"
# include "converse.h"
# include "map.h"

APM        Rmixing_term, Rv_sum, map_scratch ;
APM        *Rbeta_ptrs[DF_SQ]  ;
APM        *Rgamma_ptrs[DF_SQ], RDeriv[MAT_SZ] ;
double        *beta_ptrs[DF_SQ]   ;
double        *gamma_ptrs[DF_SQ], Deriv[MAT_SZ] ;
/* +++++++++++++++++++++++++++++++++++++++++ 
        Rimage() 
+++++++++++++++++++++++++++++++++++++++++ */

Rimage( x_prime, x ) 

RXtnd_pt   *x, *x_prime ;
/*
        Finds the image, x_prime,  of a delay-embedded point, x.
        The parameters of the map are in the parameter-space point
        called "parms".
*/
{
    APM                *x_pt, *xp_pt, *last_x ;
    RParm_pt        parms ;

    parms = x->p ;
    x_pt = x->p ;
    xp_pt = x_prime->p ;
    for( last_x = x_pt + N_PARMS ; x_pt < last_x ; x_pt++ )
        apmAssign( *xp_pt++, *x_pt ) ;

        /* Because of the way delay embedding works, 
           the first member of x_prime is the same as 
           the second member of x .        
        */

    x_pt = x->z.v ;
    xp_pt = x_prime->z.u ;
    for( last_x = x_pt + DEG_FREE ; x_pt < last_x ; x_pt++ )
        apmAssign( *xp_pt++, *x_pt ) ;

        /*  Do up the actual map.  One could 
            write a version of image() which worked for
            any standard-type symplectic map; it would 
            rely on another function, perturb(), to 
            completely define the map.  Instead we
            incorporate the perturbation to the 
            generating function right into our map -
            we hope to save a little time.
        */
    apmAdd( Rv_sum, x->z.v[0], x->z.v[1] ) ;
    apmCos( map_scratch, Rv_sum ) ;
    apmMultiply( Rmixing_term, map_scratch, parms[2] ) ;

    apmCos( map_scratch, x->z.v[0] ) ;
    apmCalc( x_prime->z.v[0], two, x->z.v[0], APM_MUL, 
                                 x->z.u[0], APM_SUB, 
                                 parms[0], map_scratch, APM_MUL,
                                 Rmixing_term, APM_ADD,
                                 APM_ADD, NULL                        ) ;

    apmCos( map_scratch, x->z.v[1] ) ;
    apmCalc( x_prime->z.v[1], two, x->z.v[1], APM_MUL, 
                                 x->z.u[1], APM_SUB, 
                                 parms[1], map_scratch, APM_MUL,
                                 Rmixing_term, APM_ADD,
                                 APM_ADD, NULL                        ) ;
}
/* +++++++++++++++++++++++++++++++++++++++++++
        find_Rbeta()

        In the interest of speed, we provide functions which only
        calculate those parts of the Jacobian that actually
        depend on parms and (u,v).  The other parts are 
        assumed to have been correctly set by a call to
        initJacobian() or initRjacobian(), both of which
        may be found below.
+++++++++++++++++++++++++++++++++++++++++++ */

find_Rbeta( b_block, x )

APM       *b_block[] ;
RXtnd_pt  *x ;
{
    apmAdd( Rv_sum, x->z.v[0], x->z.v[1] ) ;
    apmSin( map_scratch, Rv_sum ) ;
    apmMultiply( Rmixing_term, x->p[2], map_scratch ) ;

    apmSin( map_scratch, x->z.v[0] ) ;
    apmCalc( *b_block[0], x->p[0], map_scratch, APM_MUL,
                          two, APM_SWAP, APM_SUB,
                          Rmixing_term, APM_SUB, NULL ) ;

    apmNegate( *b_block[1], Rmixing_term ) ;
    apmNegate( *b_block[2], Rmixing_term ) ;

    apmSin( map_scratch, x->z.v[1] ) ;
    apmCalc( *b_block[3], x->p[1], map_scratch, APM_MUL,
                          two, APM_SWAP, APM_SUB,
                          Rmixing_term, APM_SUB, NULL ) ;
}
/* ++++++++++++++++++++++++++++++++++++++++++++++++++++ 

    Rgamma() : calculate the dependence of 
    v' on the parameters.  Even as the functions
    above, gamma() and Rgamma() change only those components
    pointed to by elements of a block of pointers.  
+++++++++++++++++++++++++++++++++++++++++++++++++++++++ */

find_Rgamma( g_block, x ) 

APM          *g_block[] ;
RXtnd_pt  *x ;
{
    apmAdd( Rv_sum, x->z.v[0], x->z.v[1] ) ;
    apmCos( Rmixing_term, Rv_sum ) ;

    apmCos( *g_block[0], x->z.v[0] ) ;
    apmAssign( *g_block[1], Rmixing_term ) ;
    apmCos( *g_block[2], x->z.v[1] ) ;
    apmAssign( *g_block[3], Rmixing_term ) ;
}
/* ++++++++++++++++++++++++++++++++++++++++ */

initRjacobian( jac ) 
/*
        Set the constant parts of a jacobian matrix
*/

APM  *jac ;
{
    int            j ;
    APM     *end_jac, *jpt ;

/*        
        If the array of APM's called jac has not yet been 
        initialized, do that first.
*/
    if( apmValidate(jac[0]) != APM_OK ) {
        end_jac = jac + MAT_SZ ;
        for( jpt=jac ; jpt < end_jac ; jpt++ )
            *jpt = apmNew( BASE ) ;
    }

    end_jac = jac + MAT_SZ ;                        /* Set all the entries */
    for( jpt=jac ; jpt < end_jac ; jpt++ )        /* to zero.                */
        apmAssignLong( *jpt, 0L, 0, BASE ) ;

/*        Put the identity in the (p,p) position. */
    jpt = jac ;
    for( j=0 ; j < N_PARMS ; j++ ) {
        apmAssignLong( *jpt, 1L, 0, BASE ) ;
        jpt += MAT_DIM + 1 ;
    }

/*        Put the identity in the (u,v) position. */
    jpt = jac + STAID_LEN + N_PARMS + DEG_FREE ;
    for( j=0 ; j < DEG_FREE ; j++ ) {
        apmAssignLong( *jpt, 1L, 0, BASE ) ;
        jpt += MAT_DIM + 1 ;
    }

/*        Put -1 times the identity in the (v,u) position. */
    jpt = jac + STAID_LEN + (DEG_FREE * MAT_DIM) + N_PARMS ;
    for( j=0 ; j < DEG_FREE ; j++ ) {
        apmAssignLong( *jpt, -1L, 0, BASE ) ;
        jpt += MAT_DIM + 1 ;
    }
}
/* +++++++++++++++++++++ */

initMap()
{
/*
    This function depends in detail on the choice of map.
*/
    beta_ptrs[0] = Deriv + STAID_LEN + (DEG_FREE * MAT_DIM) +
                                        N_PARMS  + DEG_FREE ;
    beta_ptrs[1] = beta_ptrs[0] + 1 ;
    beta_ptrs[2] = beta_ptrs[0] + MAT_DIM ;
    beta_ptrs[3] = beta_ptrs[2] + 1 ;

    gamma_ptrs[0] = Deriv + STAID_LEN + (DEG_FREE * MAT_DIM) ;
    gamma_ptrs[1] = gamma_ptrs[0] + 2 ;
    gamma_ptrs[2] = gamma_ptrs[0] + MAT_DIM + 1 ;
    gamma_ptrs[3] = gamma_ptrs[1] + MAT_DIM ;

/*
        APM stuff
*/
    Rbeta_ptrs[0] = RDeriv + STAID_LEN + (DEG_FREE * MAT_DIM) +
                                                N_PARMS + DEG_FREE ;
    Rbeta_ptrs[1] = Rbeta_ptrs[0] + 1 ;
    Rbeta_ptrs[2] = Rbeta_ptrs[0] + MAT_DIM ;
    Rbeta_ptrs[3] = Rbeta_ptrs[2] + 1 ;

    Rgamma_ptrs[0] = RDeriv + STAID_LEN + (DEG_FREE * MAT_DIM) ;
    Rgamma_ptrs[1] = Rgamma_ptrs[0] + 2 ;
    Rgamma_ptrs[2] = Rgamma_ptrs[0] + MAT_DIM + 1 ;
    Rgamma_ptrs[3] = Rgamma_ptrs[1] + MAT_DIM ;

    initJacobian( Deriv ) ;
    initRjacobian( RDeriv ) ;
/*
        Further APM stuff - constants and scratch variables.
*/
    Rv_sum = apmNew( BASE ) ;
    map_scratch = apmNew( BASE ) ;
    Rmixing_term = apmNew( BASE ) ;
}
/* +++++++++++++++++++++++ */

Rjacobian( xpt )

RXtnd_pt  *xpt ;
{
    find_Rbeta( Rbeta_ptrs, xpt ) ;
    find_Rgamma( Rgamma_ptrs, xpt ) ;
}
\end{verbatim}

%% file: appendices/code/prismatic.tex
\begin{verbatim} 
# include <stdio.h>
# include <math.h>
# include "apm.h"
# include "apmSpecial.h"
# include "converse.h"
# include "bounding.h"
# include "map.h"

int                (* fatten)(), (* row_sums)() ;
APM                Rw[MAT_DIM] ;
double                w[MAT_DIM]  ;
/* ++++++++++++++++++++++++++++++++ */

Rprismatic_image( pz_prime, pz ) 

RPrism   *pz_prime, *pz ;
{
    int           j ;
    APM    *mpt, *end_mat, *wpt, *end_w ;

/*
        Find the image of the center of the prism.
*/
    Rimage( pz_prime->center, pz->center ) ;

    Rjacobian( pz->center ) ;        /* Calculate the derivative 
                                   of the map.                    */

/*
        Fatten the matrix    Deriv * pz->matrix  so that it isn't too
        singular.
*/
    (* fatten) ( pz_prime->matrix, RDeriv, pz->matrix ) ;

/*
        Get upper bounds on the rows of the fattened matrix,
        and use them to get the matrix of a prism gauranteed 
        to enclose the image of pz.
*/
    (* row_sums)( Rw, pz_prime->matrix, RDeriv, pz ) ;

    end_w = Rw + MAT_DIM ;
    end_mat = pz_prime->matrix + MAT_SZ ;
    for( mpt = pz_prime->matrix ; mpt < end_mat ; ) {
        for( wpt = Rw ; wpt < end_w ; wpt++, mpt++ )
            apmCalc( *mpt, *mpt, *wpt, max_error, 
                           APM_ADD, APM_MUL, NULL ) ;
    }

    truncateRPrism( pz_prime, precision ) ;

}
/* +++++++++++++++++++++ */

initPrismatic()
{
    int   j ;

    for( j=0 ; j < N_PARMS ; j++ ) {
        Rw[j] = apmNew( BASE ) ;
        apmAssign( Rw[j], one ) ;
        w[j]  = 1.0 ;
    }

    for( j=N_PARMS ; j < (N_PARMS + DEG_FREE) ; j++ ) 
        Rw[j] = apmNew( BASE ) ;

    for( j=(N_PARMS + DEG_FREE) ; j < MAT_DIM ; j++ ) {
        w[j] = 1.0 + DBL_ERR ;
        Rw[j] = apmNew( BASE ) ;
        apmAdd( Rw[j], one, max_error ) ;
    }
}
\end{verbatim}

%% file: appendices/code/rowsh.tex
\begin{verbatim} 
extern int        isNewPrism ;

extern int        global_bounds(), Rglobal_bounds() ;
extern int        Rbeta_dif_star(), Rgamdif_star() ;
extern double   beta_dif_star(), gamdif_star() ;

extern Bdd_dbl        cos_zero, cos_one, cos_sum ;
extern Bdd_expr a_sin, b_sin, c_sin ;

extern Bdd_apm   Rcos_zero, Rcos_one, Rcos_sum ;
extern Bapm_expr Ra_sin, Rb_sin, Rc_sin ;

extern APM        neg_one, neg_two, Rrow_abs[] ; 
/* +++++++++++++++++++++++++++++++++++++++++++++++++++++++ */
\end{verbatim}

%% file: appendices/code/rowSums.tex
\begin{verbatim} 
# include <stdio.h>
# include <math.h>
# include "apm.h"
# include "apmSpecial.h"
# include "converse.h"
# include "bounding.h"
# include "rows.h"

APM          neg_one, neg_two ;
APM          Rrows[DEG_FREE], Rrow_abs[DEG_FREE] ;

Bdd_dbl          a, b, c, cos_zero, cos_one, cos_sum ;
Bdd_dbl   sin_zero, sin_one, sin_sum, theta ;
Bdd_dbl          *row_factors[NUM_FACTS] ;
Bdd_term  row_terms[NUM_TERMS] ;
Bdd_expr  beta_dif[3], gamma_dif[3] ;
Bdd_expr  a_sin, b_sin, c_sin ;

Bdd_apm          Ra, Rb, Rc, Rcos_zero, Rcos_one, Rcos_sum ;
Bdd_apm   Rsin_zero, Rsin_one, Rsin_sum, Rtheta ;
Bdd_apm          *Rrow_factors[NUM_FACTS] ;
Bapm_term Rrow_terms[NUM_TERMS] ;
Bapm_expr Rbeta_dif[3], Rgamma_dif[3] ;
Bapm_expr Ra_sin, Rb_sin, Rc_sin ;
/* ++++++++++++++++++++++++++++++ */

initRowSums()
/*
        Set up the expressions and terms as described in my notes
        from 11/14.
*/
{
    int              j, k ;
    Bdd_dbl   **dpt ;
    Bdd_apm   **apt ;
    Bdd_term  *tpt ;
    Bapm_term *Rtpt ;

/*
    Set up some APM's to be used to hold intermediate 
    results.
*/
    newBapm( Ra, BASE ) ;
    newBapm( Rb, BASE ) ;
    newBapm( Rc, BASE ) ;
    newBapm( Rtheta, BASE ) ;

    newBapm( Rcos_zero, BASE ) ;
    newBapm( Rcos_one, BASE ) ;
    newBapm( Rcos_sum, BASE ) ;
    newBapm( Rsin_zero, BASE ) ;
    newBapm( Rsin_one, BASE ) ;
    newBapm( Rsin_sum, BASE ) ;

    neg_one = apmInit( -1L, 0, BASE ) ;
    neg_two = apmInit( -2L, 0, BASE ) ;

    for( j=0 ; j <DEG_FREE ; j++ )  {
        Rrows[j] = apmNew( BASE ) ;
        Rrow_abs[j] = apmNew( BASE ) ;
    }

/*
        Set the number of terms in the bounded expressions
*/
    a_sin.nterms = Ra_sin.nterms = 1 ;
    b_sin.nterms = Rb_sin.nterms = 1 ;
    c_sin.nterms = Rc_sin.nterms = 1 ;

    beta_dif[0].nterms = Rbeta_dif[0].nterms = 2 ;
    beta_dif[1].nterms = Rbeta_dif[1].nterms = 1 ;
    beta_dif[2].nterms = Rbeta_dif[2].nterms = 2 ;

    gamma_dif[0].nterms = Rgamma_dif[0].nterms = 1 ;
    gamma_dif[1].nterms = Rgamma_dif[1].nterms = 1 ;
    gamma_dif[2].nterms = Rgamma_dif[2].nterms = 1 ;

/*        
        Assign terms
*/

    tpt = row_terms ;
    Rtpt = Rrow_terms ;
    for( j=0 ; j < 3 ; j++ ) {
        beta_dif[j].terms = tpt ;
        Rbeta_dif[j].terms = Rtpt ;
        tpt += beta_dif[j].nterms ;
        Rtpt += Rbeta_dif[j].nterms ;

        gamma_dif[j].terms = tpt ;
        Rgamma_dif[j].terms = Rtpt ;
        tpt += gamma_dif[j].nterms ;
        Rtpt += Rgamma_dif[j].nterms ;
    }

    a_sin.terms = tpt++ ;
    Ra_sin.terms = Rtpt++ ;

    b_sin.terms = tpt++ ;
    Rb_sin.terms = Rtpt++ ;

    c_sin.terms = tpt++ ;
    Rc_sin.terms = Rtpt++ ;

/*
        Set nfactors.
*/
    
    Rbeta_dif[0].terms[0].nfactors = beta_dif[0].terms[0].nfactors = 1 ;
    Rbeta_dif[0].terms[1].nfactors = beta_dif[0].terms[1].nfactors = 1 ;
    Rbeta_dif[1].terms[0].nfactors = beta_dif[1].terms[0].nfactors = 1 ;
    Rbeta_dif[2].terms[0].nfactors = beta_dif[2].terms[0].nfactors = 1 ;
    Rbeta_dif[2].terms[1].nfactors = beta_dif[2].terms[1].nfactors = 1 ;

    Rgamma_dif[0].terms->nfactors = gamma_dif[0].terms->nfactors = 1 ;
    Rgamma_dif[1].terms->nfactors = gamma_dif[1].terms->nfactors = 1 ;
    Rgamma_dif[2].terms->nfactors = gamma_dif[2].terms->nfactors = 1 ;

    a_sin.terms->nfactors = Ra_sin.terms->nfactors = 2 ;
    b_sin.terms->nfactors = Rb_sin.terms->nfactors = 2 ;
    c_sin.terms->nfactors = Rc_sin.terms->nfactors = 2 ;
/*
        Assign factors. 
*/

    dpt = row_factors ;
    apt = Rrow_factors ;
    for( j=0 ; j < 3 ; j++ ) {
    /*
        beta_dif
    */
        for( k=0 ; k < beta_dif[j].nterms ; k++ ) {
            beta_dif[j].terms[k].factors = dpt ;
            Rbeta_dif[j].terms[k].factors = apt ;

            dpt += beta_dif[j].terms[k].nfactors ;
            apt += Rbeta_dif[j].terms[k].nfactors ;
        }

    /*
        gamma_dif
    */
        for( k=0 ; k < gamma_dif[j].nterms ; k++ ) {
            gamma_dif[j].terms[k].factors = dpt ;
            Rgamma_dif[j].terms[k].factors = apt ;

            dpt += gamma_dif[j].terms[k].nfactors ;
            apt += Rgamma_dif[j].terms[k].nfactors ;
        }
    }

    a_sin.terms->factors = dpt ;
    Ra_sin.terms->factors = apt ;
    dpt += 2 ;
    apt += 2 ;

    b_sin.terms->factors = dpt ;
    Rb_sin.terms->factors = apt ;
    dpt += 2 ;
    apt += 2 ;

    c_sin.terms->factors = dpt ;
    Rc_sin.terms->factors = apt ;

/*
        Set up those of the "bound" attributes which are 
        bounded APM's.
*/

    for( j=0 ; j < NUM_TERMS ; j++ ) {
        newBapm( Rrow_terms[j].bound, BASE ) ;
    }

    for( j=0 ; j < 3 ; j++ ) {
        newBapm( Rbeta_dif[j].bound, BASE ) ;
        newBapm( Rgamma_dif[j].bound, BASE ) ;
    }

    newBapm( Ra_sin.bound, BASE ) ;
    newBapm( Rb_sin.bound, BASE ) ;
    newBapm( Rc_sin.bound, BASE ) ;

/*
        Set up the terms and expressions. 
*/

/*  a_sin */

    a_sin.const = 0.0 ;
    Ra_sin.const = apmNew( BASE ) ;
        a_sin.terms->coef = 1.0 ;
        Ra_sin.terms->coef = apmInit( 1L, 0, BASE ) ;
            
            a_sin.terms->factors[0] = &a ;
            a_sin.terms->factors[1] = &sin_zero ;
            Ra_sin.terms->factors[0] = &Ra ;
            Ra_sin.terms->factors[1] = &Rsin_zero ;

/*  b_sin */

    b_sin.const = 0.0 ;
    Rb_sin.const = apmNew( BASE ) ;
        b_sin.terms->coef = 1.0 ;
        Rb_sin.terms->coef = apmInit( 1L, 0, BASE ) ;
            
            b_sin.terms->factors[0] = &b ;
            b_sin.terms->factors[1] = &sin_one ;
            Rb_sin.terms->factors[0] = &Rb ;
            Rb_sin.terms->factors[1] = &Rsin_one ;

/*  c_sin */

    c_sin.const = 0.0 ;
    Rc_sin.const = apmNew( BASE ) ;
        c_sin.terms->coef = 1.0 ;
        Rc_sin.terms->coef = apmInit( 1L, 0, BASE ) ;
            
            c_sin.terms->factors[0] = &c ;
            c_sin.terms->factors[1] = &sin_sum ;
            Rc_sin.terms->factors[0] = &Rc ;
            Rc_sin.terms->factors[1] = &Rsin_sum ;

/*  beta_dif */
 
        /* beta_dif[0] =  (2.0 - a * sin(v[0]) - c * sin(v[0] + v[1]) ) 
                   -{ 2.0 - ac * sin(vc[0]) - cc * sin(vc[0] + vc[1])  
            Where ac, cc, vc[0], and vc[1] are the values of these
            numbers at the center of the prism.  The whole second 
            term ( enclosed in braces ) is an entry in the jacobian 
            of the map                                                        
                                                                */
    Rbeta_dif[0].const = apmNew( BASE ) ;
        beta_dif[0].terms[0].coef = -1.0 ;
        Rbeta_dif[0].terms[0].coef = neg_one ;

            beta_dif[0].terms[0].factors[0] = &a_sin.bound ;
            Rbeta_dif[0].terms[0].factors[0] = &Ra_sin.bound ;

        beta_dif[0].terms[1].coef = -1.0 ;
        Rbeta_dif[0].terms[1].coef = neg_one ;

            beta_dif[0].terms[1].factors[0] = &c_sin.bound ;
            Rbeta_dif[0].terms[1].factors[0] = &Rc_sin.bound ;

        /* beta_dif[1] =     -2.0 * c * sin.bound( v[0] + v[1] )
                         - { -2.0 * cc * sin.bound( vc[0] + vc[1] ) } 
         */
    Rbeta_dif[1].const = apmNew( BASE ) ;
        beta_dif[1].terms[0].coef = -2.0 ;
        Rbeta_dif[1].terms[0].coef = neg_two ;

            beta_dif[1].terms[0].factors[0] = &c_sin.bound ;
            Rbeta_dif[1].terms[0].factors[0] = &Rc_sin.bound ;

        /* beta_dif[2] = 2.0 - b * sin.bound(v[1]) - c * sin(v[1] + v[0])  
                -{ 2.0 - bc * sin.bound(vc[1]) - cc * sin(vc[1] + vc[0]) }
        */
    Rbeta_dif[2].const = apmNew( BASE ) ;
        beta_dif[2].terms[0].coef = -1.0 ;
        Rbeta_dif[2].terms[0].coef = neg_one ;

            beta_dif[2].terms[0].factors[0] = &b_sin.bound ;
            Rbeta_dif[2].terms[0].factors[0] = &Rb_sin.bound ;

        beta_dif[2].terms[1].coef = -1.0 ;
        Rbeta_dif[2].terms[1].coef = neg_one ;

            beta_dif[2].terms[1].factors[0] = &c_sin.bound ;
            Rbeta_dif[2].terms[1].factors[0] = &Rc_sin.bound ;

/* gamma_dif */

        /* gamma_dif[0] = da * ( cos(v[0]) - cos(vc[0]) ) 
           Where da is half the prism's width as measured
           along the a-axis and vc is as above.                  */

    Rgamma_dif[0].const = apmNew( BASE ) ;

        Rgamma_dif[0].terms[0].coef = apmNew( BASE ) ;
            
            gamma_dif[0].terms[0].factors[0] = &cos_zero ;
            Rgamma_dif[0].terms[0].factors[0] = &Rcos_zero ;

        /* gamma_dif[1] = db * ( cos(v[1]) - cos(vc[1]) )  */

    Rgamma_dif[1].const = apmNew( BASE ) ;

        Rgamma_dif[1].terms[0].coef = apmNew( BASE ) ;
            
            gamma_dif[1].terms[0].factors[0] = &cos_one ;
            Rgamma_dif[1].terms[0].factors[0] = &Rcos_one ;

        /* gamma_dif[2] = dc * ( cos(v[0]  + v[1]) - 
                                 cos(vc[0] + vc[1]) )  */

    Rgamma_dif[2].const = apmNew( BASE ) ;

        Rgamma_dif[2].terms[0].coef = apmNew( BASE ) ;
            
            gamma_dif[2].terms[0].factors[0] = &cos_sum ;
            Rgamma_dif[2].terms[0].factors[0] = &Rcos_sum ;
}
/* +++++++++++++++++++++++++++++++++ */

Rglobal_bounds( pz )

RPrism        *pz ;
{
    int             j ;
    APM             *apt, *end_row ;

    apmAdd( Ra.ub, pz->center->p[0], pz->matrix[0] ) ;
    apmSubtract( Ra.lb, pz->center->p[0], pz->matrix[0] ) ;

    apmAdd( Rb.ub, pz->center->p[1], pz->matrix[MAT_DIM+1] ) ;
    apmSubtract( Rb.lb, pz->center->p[1], pz->matrix[MAT_DIM+1] ) ;

    apmAdd( Rc.ub, pz->center->p[2], pz->matrix[2*MAT_DIM+2] ) ;
    apmSubtract( Rc.lb, pz->center->p[2], pz->matrix[2*MAT_DIM+2] ) ;

    apt = pz->matrix + STAID_LEN + (DEG_FREE * MAT_DIM) ;
    for( j=0 ; j < DEG_FREE ; j++ ) {
        apmAssign( Rrows[j], zero ) ;
        for( end_row=apt + MAT_DIM ; apt < end_row ; apt++ ) {
            apmCalc( Rrows[j], Rrows[j], *apt, 
                                        APM_ABS, APM_ADD, NULL ) ;
        }
    }

    apmAdd( Rtheta.ub, pz->center->z.v[0], Rrows[0] ) ;
    apmSubtract( Rtheta.lb, pz->center->z.v[0], Rrows[0] ) ;
    Rbd_sin( &Rsin_zero, &Rtheta ) ;
    Rbd_cos( &Rcos_zero, &Rtheta ) ;
        
    apmAdd( Rtheta.ub, pz->center->z.v[1], Rrows[1] ) ;
    apmSubtract( Rtheta.lb, pz->center->z.v[1], Rrows[1] ) ;
    Rbd_sin( &Rsin_one, &Rtheta ) ;
    Rbd_cos( &Rcos_one, &Rtheta ) ;
        
    apmCalc( Rtheta.ub, Rtheta.ub, pz->center->z.v[0], Rrows[0],
                                        APM_ADD, APM_ADD, NULL ) ;
    apmCalc( Rtheta.lb, Rtheta.lb, pz->center->z.v[0], Rrows[0],
                                        APM_SUB, APM_ADD, NULL ) ;
    Rbd_sin( &Rsin_sum, &Rtheta ) ;
    Rbd_cos( &Rcos_sum, &Rtheta ) ;

    Rbound_expr( &Ra_sin ) ;
    Rbound_expr( &Rb_sin ) ;
    Rbound_expr( &Rc_sin ) ;

}
/* +++++++++++++++++++++++++++++++++ */

Rbeta_dif_star( answer, deriv ) 

APM   answer, *deriv ;
{
     APM    *dpt ;

     dpt = deriv + STAID_LEN + (MAT_DIM*DEG_FREE) + N_PARMS + DEG_FREE ;
     apmSubtract( Rbeta_dif[0].const, two, *dpt++ ) ; 
     apmMultiply( Rbeta_dif[1].const, neg_two, *dpt ) ;
     dpt += MAT_DIM ;
     apmSubtract( Rbeta_dif[2].const, two, *dpt ) ;

    Rbound_expr( &Rbeta_dif[0] ) ;
    Rbound_expr( &Rbeta_dif[1] ) ;
    Rbound_expr( &Rbeta_dif[2] ) ;

    RmaxAbs( answer, Rbeta_dif[0].bound.ub, Rbeta_dif[0].bound.lb )  ;
    RmaxAbs( Rrow_abs[0], Rbeta_dif[1].bound.ub, Rbeta_dif[1].bound.lb ) ;
    RmaxAbs( Rrow_abs[1], Rbeta_dif[2].bound.ub, Rbeta_dif[2].bound.lb ) ;

/*
        Add max_error to the answer to account for the uncertainties
        in beta**(center).
*/
    apmCalc( answer, answer, Rrow_abs[0], Rrow_abs[1], max_error, 
                                APM_ADD, APM_ADD, APM_ADD, NULL ) ; 
}
/* +++++++++++++++++++++ */

Rgamdif_star( answer, deriv, pmat )

APM    answer, *deriv, *pmat ;
{
    APM    *apt, *Rda, *Rdb, *Rdc ;

    Rda = pmat ;
    Rdb = pmat + MAT_DIM + 1 ;
    Rdc = pmat+ (2 * MAT_DIM) + 2 ;

    apmAssign( Rgamma_dif[0].terms[0].coef, *Rda ) ;
    apmAssign( Rgamma_dif[1].terms[0].coef, *Rdb ) ;
    apmMultiply( Rgamma_dif[2].terms[0].coef, two, *Rdc ) ;

    apt = deriv + STAID_LEN + (DEG_FREE * MAT_DIM) ;
    apmCalc( Rgamma_dif[0].const, *Rda, APM_NEG, *apt, APM_MUL, NULL ) ;
    apt += MAT_DIM + 1 ; 
    apmCalc( Rgamma_dif[1].const, *Rdb, APM_NEG, *apt, APM_MUL, NULL ) ;
    apt++ ;
    apmCalc( Rgamma_dif[2].const, two, APM_NEG, *Rdc, *apt, 
                                       APM_MUL, APM_MUL, NULL ) ;

    Rbound_expr( &Rgamma_dif[0] ) ;
    Rbound_expr( &Rgamma_dif[1] ) ;
    Rbound_expr( &Rgamma_dif[2] ) ;

    RmaxAbs( answer, Rgamma_dif[0].bound.ub, Rgamma_dif[0].bound.lb )  ;
    RmaxAbs( Rrow_abs[0], Rgamma_dif[1].bound.ub, Rgamma_dif[1].bound.lb ) ;
    RmaxAbs( Rrow_abs[1], Rgamma_dif[2].bound.ub, Rgamma_dif[2].bound.lb ) ;
/*
        Add max_error to the answer to account for the uncertainties
        in beta**(center).
*/
    apmCalc( answer, answer, Rrow_abs[0], Rrow_abs[1], max_error,
                                APM_ADD, APM_ADD, APM_ADD, NULL ) ; 
}
\end{verbatim}

%% file: appendices/code/cr_fatten.tex
\begin{verbatim} 
# include  <stdio.h>
# include  <math.h>
# include  "apm.h"
# include  "apmSpecial.h"
# include  "converse.h"
# include  "bounding.h"
# include  "rows.h"
# include  "pi.h"

                                   recorded here in units of pi.    */

APM         Rcthet, Rsthet, Rsmall_sinsq ;
APM         Rarea, Rsin_sq, Rnorm_one, Rnorm_two, Rsign ;
APM         Rnorm_prod, Rsign, Rx, Ry ;
double   cthet, sthet, small_sinsq ;
/* ++++++++++++++++++++++++++++++++++++++++++ */

initRotor()
{
    Rcthet = apmNew( BASE ) ;
    Rsthet = apmNew( BASE ) ;

    Rx = apmNew( BASE ) ;
    Ry = apmNew( BASE ) ;
    Rarea = apmNew( BASE ) ;
    Rsign = apmNew( BASE ) ;
    Rsin_sq = apmNew( BASE ) ;
    Rnorm_one = apmNew( BASE ) ;
    Rnorm_two = apmNew( BASE ) ;
    Rnorm_prod = apmNew( BASE ) ;
    Rsmall_sinsq = apmNew( BASE ) ;

    cthet = cos( PI * THETA_ROT ) ;
    sthet = sin( PI * THETA_ROT ) ;
    small_sinsq = sthet * sthet ;

    dbltoapm( Rx, BASE, THETA_ROT ) ;
    apmMultiply( Ry, pi, Rx ) ;
    apmCos( Rcthet, Ry ) ;
    apmSin( Rsthet, Ry ) ;
    apmMultiply( Rsmall_sinsq, Rsthet, Rsthet ) ;
}
/* ++++++++++++++++++++++++++++++++++++++++++ */

Rcol_rotor( Amat, Deriv, Prizmat )

APM  *Amat, *Deriv, *Prizmat ;
/*
     Prepares the matrix called "A" in my notes.  Mostly we want to
     have A = DF*Priz, but  we want to ensure that A is not singular.
     In the interest of speed we have coded the calculations below with 
     pointers.  Our hope is that the resulting function will scream along 
     at ultrasonic speed.  Unfortunately it is quite unreadable.
*/
{
    int                     j, k ;
    APM                     *Aend, *Dend, *Pend ;
    register APM     *Apt, *Dpt, *Ppt ;

/*
        Copy the few terms which appear in the top rows of Amat.
*/
    Aend = Amat + N_PARMS * (MAT_DIM + 1) ;
    for( Apt = Amat, Ppt = Prizmat ; Apt < Aend ; Apt += (MAT_DIM + 1 ),
                                                  Ppt += (MAT_DIM + 1 ) ) 
        apmAssign( *Apt, *Ppt ) ;

/*
        Clear out those parts of Amat which change from iteration to 
        iteration. 
*/

    Aend = Amat + MAT_SZ ;
    for( Apt = Amat + STAID_LEN ; Apt < Aend ; Apt++ )
        apmAssignLong( *Apt, 0L, 0, BASE ) ;

/*
        Set the (u,p) part of A
        It's equal to the (v,p) part of Prizmat.
*/
    
    Aend = Amat + STAID_LEN + (DEG_FREE * MAT_DIM) ;
    Ppt = Prizmat + STAID_LEN + (DEG_FREE * MAT_DIM) ;
    for( Apt = Amat + STAID_LEN ; Apt < Aend ; Apt += TWO_DF ) {
        for( Pend = Ppt + N_PARMS ; Ppt < Pend ; Ppt++ )
            apmAssign( *Apt++, *Ppt ) ;

        Ppt += TWO_DF ;
    }

/*
        Set the (v,p) part - three terms.
*/
        /*  First term - equal to Deriv(v,p) * Prizmat(p,p)  */

    Dpt = Deriv + STAID_LEN + (DEG_FREE * MAT_DIM) ;
    Apt = Amat + STAID_LEN + (DEG_FREE * MAT_DIM) ;
    for( Aend = Apt + (DEG_FREE*MAT_DIM) ; Apt < Aend ; Apt += TWO_DF ) {
        Ppt = Prizmat ;
        for( Dend = Dpt + N_PARMS ; Dpt < Dend ; Dpt++ ) {
            apmCalc( *Apt, *Apt, *Dpt, *Ppt, APM_MUL, APM_ADD, NULL ) ;
            Apt++ ;
            Ppt += MAT_DIM + 1 ;
        }

        Dpt += TWO_DF ;
    }

        /*  Second term - equal to negative Prizmat(u,p) */

    Ppt = Prizmat + STAID_LEN ;
    Apt = Amat + STAID_LEN + (DEG_FREE * MAT_DIM) ;
    for( Pend = Ppt + (DEG_FREE * MAT_DIM) ;  Ppt < Pend ; Ppt += TWO_DF ) {
        for( Aend = Apt + N_PARMS ; Apt < Aend ; Apt++ )
            apmCalc( *Apt, *Apt, *Ppt++, APM_SUB, NULL ) ;

        Apt += TWO_DF ;
    }

        /*  Third term - equal to Deriv(v,v) * Prizmat(v,p)  */

    Dpt = Deriv + STAID_LEN + (DEG_FREE * (MAT_DIM + 1)) + N_PARMS ;
    Dend = Deriv + MAT_SZ ;
    Apt = Amat + STAID_LEN + (DEG_FREE * MAT_DIM) ;
    while( Dpt < Dend ) {
        Ppt = Prizmat + STAID_LEN + (DEG_FREE * MAT_DIM) ;
        Pend = Prizmat + MAT_SZ ;
        while( Ppt < Pend ) {
            Aend = Apt + N_PARMS ;
            while( Apt < Aend ) {
                apmCalc( *Apt, *Apt, *Dpt, *Ppt, APM_MUL, APM_ADD, NULL ) ;
                Apt++ ;
                Ppt++ ;
            }

            Dpt++ ;
            Ppt += TWO_DF ;
            Apt -= N_PARMS ;
        }

        Dpt += N_PARMS + DEG_FREE ;
        Apt += MAT_DIM ;
    }

/*
        (u,u) part
        equals Priz(v,u)
*/

    Apt = Amat + STAID_LEN + N_PARMS ;
    Aend = Amat + STAID_LEN + (DEG_FREE * MAT_DIM) ;
    Ppt = Prizmat + STAID_LEN + (DEG_FREE * MAT_DIM) + N_PARMS ;
    while( Apt < Aend ) {
        Pend = Ppt + DEG_FREE ;
        while( Ppt < Pend )  {
            apmAssign( *Apt++, *Ppt ++) ;
        }

        Apt  += N_PARMS + DEG_FREE ;
        Ppt  += N_PARMS + DEG_FREE ;
    }

/*
        (u,v) part
        equals Priz(v,v)
*/

    Apt = Amat + STAID_LEN + N_PARMS + DEG_FREE ;
    Aend = Amat + STAID_LEN + (DEG_FREE * MAT_DIM) ;
    Ppt = Prizmat + STAID_LEN + (DEG_FREE*MAT_DIM) + N_PARMS + DEG_FREE ;
    while( Apt < Aend ) {
        Pend = Ppt + DEG_FREE ;
        while( Ppt < Pend ) 
            apmAssign( *Apt++, *Ppt++ ) ;

        Apt  += N_PARMS + DEG_FREE ;
        Ppt  += N_PARMS + DEG_FREE ;
    }

/*
        The (v,u) part - equal to Deriv(v,v) * Priz(v,u) - Priz(u,u) ,
*/

        /* First term */
    Apt = Amat + STAID_LEN + (DEG_FREE * MAT_DIM) + N_PARMS ;
    Aend = Apt + (DEG_FREE * MAT_DIM) ;
    Dpt = Deriv + STAID_LEN + (DEG_FREE*MAT_DIM) + N_PARMS + DEG_FREE ;
    while( Apt < Aend ) {
        Ppt = Prizmat + STAID_LEN + (DEG_FREE * MAT_DIM) + N_PARMS ;
        Pend = Ppt + DEG_FREE ;
        while( Ppt < Pend ) {
            Dend = Dpt + DEG_FREE ;
            while( Dpt < Dend )  {
                apmCalc( *Apt, *Apt, *Dpt++, *Ppt, APM_MUL, 
                                             APM_ADD, NULL ) ;

                Ppt += MAT_DIM ;
            }
            Apt++ ;
            Dpt -= DEG_FREE ;
            Ppt -= (DEG_FREE * MAT_DIM) - 1 ;
        }
        Dpt += MAT_DIM ;
        Apt += N_PARMS + DEG_FREE ;
    }
        
        /* Second term */
    Apt = Amat + STAID_LEN + (DEG_FREE * MAT_DIM) + N_PARMS + DEG_FREE ;
    Ppt = Prizmat + STAID_LEN + N_PARMS ;
    Pend = Ppt + (MAT_DIM * DEG_FREE) ;
    while( Ppt < Pend ) {
        Aend = Apt + DEG_FREE ;
        while( Apt < Aend )  {
            apmCalc( *Apt, *Apt, *Ppt, APM_SUB, NULL ) ;
            Apt++ ;
            Ppt++ ;
        }

        Ppt += N_PARMS + DEG_FREE ;
        Apt += N_PARMS + DEG_FREE ;
    }

/*
        (v,v) part - equals Deriv(v,v) * Priz(v,v) - Priz(u,v) 
*/

        /* First term */
    Apt = Amat + STAID_LEN + (DEG_FREE * MAT_DIM) + N_PARMS + DEG_FREE ;
    Aend = Apt + (DEG_FREE * MAT_DIM) ;
    Dpt = Deriv + STAID_LEN + (DEG_FREE*MAT_DIM) + N_PARMS + DEG_FREE ;
    while( Apt < Aend ) {
        Ppt = Prizmat + STAID_LEN + (DEG_FREE*MAT_DIM) + 
                                    N_PARMS + DEG_FREE ;
        Pend = Ppt + DEG_FREE ;
        while( Ppt < Pend ) {
            Dend = Dpt + DEG_FREE ;
            while( Dpt < Dend )  {
                apmCalc( *Apt, *Apt, *Dpt++, *Ppt, APM_MUL, 
                                             APM_ADD, NULL ) ;

                Ppt += MAT_DIM ;
            }
            Apt++ ;
            Dpt -= DEG_FREE ;
            Ppt -= (DEG_FREE * MAT_DIM) - 1 ;
        }
        Dpt += MAT_DIM ;
        Apt += N_PARMS + DEG_FREE ;
    }
        
        /* Second term */
    Apt = Amat + STAID_LEN + (DEG_FREE * MAT_DIM) + N_PARMS + DEG_FREE ;
    Ppt = Prizmat + STAID_LEN + N_PARMS + DEG_FREE ;
    Pend = Ppt + (MAT_DIM * DEG_FREE) ;
    while( Ppt < Pend ) {
        Aend = Apt + DEG_FREE ;
        while( Apt < Aend ) {
            apmCalc( *Apt, *Apt, *Ppt, APM_SUB, NULL ) ;
            Apt++ ;
            Ppt++ ;
        }

        Ppt += N_PARMS + DEG_FREE ;
        Apt += N_PARMS + DEG_FREE ;
    }

# if USE_ROT
/*
        Do up the rotations.
*/
    for( j=0 ; j < TWO_DF ; j++ )
        for( k=(j+1) ; k < TWO_DF ; k++ )
        Rsubspace_rot( Amat, j, k ) ;
# endif
}
/* +++++++++++++++++++++++++++++ */

Rsubspace_rot( Amat, col_one, col_two )

int        col_one, col_two ;
APM        *Amat ;
{
    APM   *Apt, *Apt2 ;

    Apt = Amat + STAID_LEN + N_PARMS + 
                 (col_two - col_one - 1) * MAT_DIM  +
                 col_one ;
    Apt2 = Apt + col_two - col_one ;

    apmCalc( Rarea, *Apt, Apt2[MAT_DIM], APM_MUL,
                    Apt[MAT_DIM], *Apt2, APM_MUL, 
                                  APM_SUB, NULL ) ;
    apmCalc( Rnorm_one, *Apt, APM_DUP, APM_MUL,
                        Apt[MAT_DIM], APM_DUP, APM_MUL, 
                                        APM_ADD, NULL ) ;
    apmCalc( Rnorm_two, *Apt2, APM_DUP, APM_MUL,
                        Apt2[MAT_DIM], APM_DUP, APM_MUL, 
                                        APM_ADD, NULL ) ;
    apmMultiply( Rnorm_prod, Rnorm_one, Rnorm_two ) ;
    if( apmCompare( Rnorm_prod, zero ) == 1 ) {
        apmMultiply( Rx, Rarea, Rarea ) ;
        apmDivide( Rsin_sq, precision, (APM) NULL, Rx, Rnorm_prod ) ;

        if( apmCompare( Rsin_sq, Rsmall_sinsq ) == -1 ) {
            Rm_sign( Rsign, Rarea ) ;

            if( apmCompare( Rnorm_two, Rnorm_one ) != 1 ) {
                    apmCalc( Rx, Rcthet, *Apt2, APM_MUL,
                     Rsign, Rsthet, Apt2[MAT_DIM], APM_MUL, APM_MUL,
                     APM_SUB, NULL ) ;
                    apmCalc( Ry, Rsthet, *Apt2, Rsign, APM_MUL, APM_MUL,
                     Rcthet, Apt2[MAT_DIM], APM_MUL, 
                     APM_ADD, NULL ) ;

                    apmAssign( *Apt2, Rx ) ;
                    apmAssign( Apt2[MAT_DIM], Ry ) ;
            }
            else {
                apmCalc( Rsign, Rsign, APM_NEG, NULL ) ;
                    apmCalc( Rx, Rcthet, *Apt, APM_MUL,
                     Rsign, Rsthet, Apt[MAT_DIM], APM_MUL, APM_MUL,
                     APM_SUB, NULL ) ;
                    apmCalc( Ry, Rsthet, *Apt, Rsign, APM_MUL, APM_MUL,
                     Rcthet, Apt[MAT_DIM], APM_MUL, 
                     APM_ADD, NULL ) ;

                    apmAssign( *Apt, Rx ) ;
                    apmAssign( Apt[MAT_DIM], Ry ) ;
            }
        }
    }
}
\end{verbatim}

%% file: appendices/code/cr_rows.tex
\begin{verbatim} 
# include <stdio.h>
# include <math.h>
# include "apm.h"
# include "apmSpecial.h"
# include "converse.h"
# include "bounding.h"
# include "rows.h"

int          isNewPrism ;

APM          cr_scratch ;
APM          RBmat[MAT_SZ], Rconst_mat[DF_SQ], Rcopy[4 * DF_SQ] ;
APM          *Rcopy_rows[TWO_DF] ;
APM          RBu_rows[DEG_FREE], RBv_rows[DEG_FREE] ;
APM          Rbd_star, Rgd_star, Rstar, RPvp_star ;
APM          Rcenter_err[MAT_DIM] ;
APM          Rup_rows[DEG_FREE], Ruu_rows[DEG_FREE], Ruv_rows[DEG_FREE] ;
APM          Rvp_rows[DEG_FREE], Rvu_rows[DEG_FREE], Rvv_rows[DEG_FREE] ;

double    Bmat[MAT_SZ], const_mat[DF_SQ], copy[4 * DF_SQ] ;
double          *copy_rows[TWO_DF] ;
double    Bu_rows[DEG_FREE], Bv_rows[DEG_FREE] ;
double    bd_star, gd_star, star, Pvp_star ;
double          center_err[MAT_DIM] ;
double    up_rows[DEG_FREE], uu_rows[DEG_FREE], uv_rows[DEG_FREE] ;
double    vp_rows[DEG_FREE], vu_rows[DEG_FREE], vv_rows[DEG_FREE] ;

Bdd_dbl          *cr_factors[NUM_FACTS] ;
Bdd_term  cr_terms[NUM_TERMS] ;
Bdd_expr  beta_prod ;

Bdd_apm          *Rcr_factors[NUM_FACTS] ;
Bapm_term Rcr_terms[NUM_TERMS] ;
Bapm_expr Rbeta_prod ;
/* ++++++++++++++++++++++++++++++ */

init_crRows()
/*
        Set up the expressions and terms as described in 
        appendix B.
*/
{
    int              j, k ;
    APM              *Rcpt ;
    double    *cpt ;
    Bdd_dbl   **dpt ;
    Bdd_apm   **apt ;

/*
        Initialize a batch of APM's.
*/
    for(j=0 ; j < DEG_FREE ; j++ ) {
        Rvp_rows[j] = apmNew( BASE ) ;
        Rup_rows[j] = apmNew( BASE ) ;
        Ruu_rows[j] = apmNew( BASE ) ;
        Ruv_rows[j] = apmNew( BASE ) ;
        Rvu_rows[j] = apmNew( BASE ) ;
        Rvv_rows[j] = apmNew( BASE ) ;
        RBu_rows[j] = apmNew( BASE ) ;
        RBv_rows[j] = apmNew( BASE ) ;
    }

    Rstar = apmNew( BASE ) ;
    Rgd_star = apmNew( BASE ) ;
    Rbd_star = apmNew( BASE ) ;
    RPvp_star = apmNew( BASE ) ;
    cr_scratch = apmNew( BASE ) ;
    for( j=0 ; j < MAT_SZ ; j++ )  {
        Bmat[j] = 0.0 ;
        RBmat[j] = apmNew( BASE ) ;
    }

    for( j=0 ; j < DF_SQ ; j++ ) 
        Rconst_mat[j] = apmNew( BASE ) ;

    for( j=0 ; j < (4 * DF_SQ) ; j++ )
        Rcopy[j] = apmNew( BASE ) ;

    for( j=0 ; j < MAT_DIM ; j++ )
        Rcenter_err[j] = apmNew( BASE ) ;

    cpt = copy ;
    Rcpt = Rcopy ;
    for( j=0 ; j < TWO_DF ; j++ ) {
        copy_rows[j] = cpt ;
        Rcopy_rows[j] = Rcpt ;

        cpt += TWO_DF ;
        Rcpt += TWO_DF ;
    }

/*
        Set the number of terms in the bounded expressions
*/

    beta_prod.nterms = Rbeta_prod.nterms = 3 ;

/*        
        Assign terms
*/

    beta_prod.terms = cr_terms ;
    Rbeta_prod.terms = Rcr_terms ;

/*
        Set nfactors.
*/
    
    Rbeta_prod.terms[0].nfactors = beta_prod.terms[0].nfactors = 1  ;
    Rbeta_prod.terms[1].nfactors = beta_prod.terms[1].nfactors = 1  ;
    Rbeta_prod.terms[2].nfactors = beta_prod.terms[2].nfactors = 1  ;

/*
        Assign factors. 
*/

    dpt = cr_factors ;
    apt = Rcr_factors ;
    for( k=0 ; k < beta_prod.nterms ; k++ ) {
        beta_prod.terms[k].factors = dpt ;
        Rbeta_prod.terms[k].factors = apt ;

        dpt += beta_prod.terms[k].nfactors ;
        apt += Rbeta_prod.terms[k].nfactors ;
    }

/*
        Set up those of the "bound" attributes which are 
        bounded APM's.
*/

    newBapm( Rbeta_prod.bound, BASE ) ;
    for( j=0 ; j < NUM_TERMS ; j++ ) {
        newBapm( Rcr_terms[j].bound, BASE ) ;
    }

/*
        Set up the terms and expressions. 
*/

/* beta_prod */
 
    Rbeta_prod.const = apmNew( BASE ) ;
        Rbeta_prod.terms[0].coef = apmNew( BASE ) ;

            beta_prod.terms[0].factors[0] = &a_sin.bound ;
            Rbeta_prod.terms[0].factors[0] = &Ra_sin.bound ;

        Rbeta_prod.terms[1].coef = apmNew( BASE ) ;

            beta_prod.terms[1].factors[0] = &c_sin.bound ;
            Rbeta_prod.terms[1].factors[0] = &Rc_sin.bound ;

        Rbeta_prod.terms[2].coef = apmNew( BASE ) ;

            beta_prod.terms[2].factors[0] = &b_sin.bound ;
            Rbeta_prod.terms[2].factors[0] = &Rb_sin.bound ;
}
/* +++++++++++++++++++++++++++++++++++++++++ */

Rcr_rows( Rw, Amat, Deriv, Priz ) 

APM        *Rw, *Amat, *Deriv ;
RPrism  *Priz ;
/*
    Obtain bounds on the sums of the absolute values of
    the entries in the rows of 
       -1
    [A]   * Deriv * Pmat,

    put the results in w.
*/
{
    int                j ;
    APM                *end_row, *end_mat, *Pmat, *inv_pt ;
    APM                *p1pt, *p2pt, *b1pt, *b2pt, *wu_pt, *wv_pt ;

    Pmat = Priz->matrix ;
    Rset_inverse( Amat ) ;

/*
        Do up some row sums for the inverse; these
        are used to calculate center_err[].
*/
    b1pt = RBmat + STAID_LEN + N_PARMS ;
    b2pt = b1pt + MAT_DIM * DEG_FREE ;
    for( j=0  ; j < DEG_FREE ; j++ ) {
        apmAssign( RBu_rows[j], zero ) ;
        apmAssign( RBv_rows[j], zero ) ;

        for( end_row = b1pt + TWO_DF ; b1pt < end_row ; ) {
            apmCalc( RBu_rows[j], RBu_rows[j], *b1pt++,
                                  APM_ABS, APM_ADD, NULL ) ;
            apmCalc( RBv_rows[j], RBv_rows[j], *b2pt++,
                                  APM_ABS, APM_ADD, NULL ) ;
        }
    }

/*
        Call functions which calculate upper bound on the 
        sums of the elements of various matrices.
        Before any bounding of matrices, one must invoke
        global_bounds( Pmat ) to set such global variables,
        as cos_one, and sin_sum.  This is done in Rtry_prism.
*/
    Rbeta_dif_star( Rbd_star, Deriv ) ;
    Rgamdif_star( Rgd_star, Deriv, Pmat ) ;

/* 
        Calculate bounds on the sums of the absolute values
        of the elements in various blocks.
*/
    /* up & vp blocks */

    apmAssignLong( RPvp_star, 0L, 0, BASE ) ;
    p1pt = Pmat + STAID_LEN + (MAT_DIM * DEG_FREE) ;
    end_mat = p1pt + (DEG_FREE * MAT_DIM) ;
    for( ; p1pt < end_mat ; p1pt += TWO_DF ) {
        for( end_row = p1pt + N_PARMS ; p1pt < end_row ; p1pt++ ) 
            apmCalc( RPvp_star, RPvp_star, *p1pt, APM_ABS,
                                                  APM_ADD, NULL ) ;
    }

    apmCalc( Rstar, Rgd_star, Rbd_star, RPvp_star, 
                             APM_MUL, APM_ADD, NULL ) ;
    b1pt = RBmat + STAID_LEN + N_PARMS + DEG_FREE ;
    b2pt = RBmat + STAID_LEN + N_PARMS + DEG_FREE + (MAT_DIM * DEG_FREE) ;
    for( j=0 ; j < DEG_FREE ; j++ ) {
        apmAssignLong( Rup_rows[j], 0L, 0, BASE ) ;
        apmAssignLong( Rvp_rows[j], 0L, 0, BASE ) ; 
        for( end_row = b1pt + DEG_FREE ; b1pt < end_row ; 
                                         b1pt++, b2pt++ ) {
            apmCalc( Rup_rows[j], Rup_rows[j], *b1pt, APM_ABS,
                                                APM_ADD, NULL ) ;
            apmCalc( Rvp_rows[j], Rvp_rows[j], *b2pt, APM_ABS,
                                                APM_ADD, NULL ) ;
        }

        apmCalc( Rup_rows[j], Rup_rows[j], Rstar, APM_MUL, NULL ) ;
        apmCalc( Rvp_rows[j], Rvp_rows[j], Rstar, APM_MUL, NULL ) ;

        b1pt += N_PARMS + DEG_FREE ;
        b2pt += N_PARMS + DEG_FREE ;
    }

/*
        Do the remaining blocks - those that actually arise 
        from the derivatives of the (u,v) -> (u',v') part of
        the map.  This section uses the mighty bound_rows(),
        which may be found below.
*/

        /* (u,u) block :
                B(u,u) * P(v,u) + B(u,v) * { beta * P(v,u) -
                                                    P(u,u) }
        */

    p1pt = Pmat + STAID_LEN + (DEG_FREE * MAT_DIM) + N_PARMS ;
    p2pt = Pmat + STAID_LEN + N_PARMS ;

    b1pt = RBmat + STAID_LEN + N_PARMS ;
    b2pt = RBmat + STAID_LEN + N_PARMS + DEG_FREE ;
    Rbound_rows( Ruu_rows, b1pt, p1pt, b2pt, p2pt ) ;

        /* (u,v) block :
                B(u,u) * P(v,v) + B(u,v) * { beta * P(v,v) -
                                                    P(u,v) }

        */

    p1pt = Pmat + STAID_LEN + (DEG_FREE*MAT_DIM) + N_PARMS + DEG_FREE ;
    p2pt = Pmat + STAID_LEN + N_PARMS + DEG_FREE ;

/*         The same parts of RBmat as used to find uu_rows. */
    Rbound_rows( Ruv_rows, b1pt, p1pt, b2pt, p2pt ) ;

        /* (v,u) block :
                B(v,u) * P(v,u) + B(v,v) * { beta * P(v,u) -
                                                    P(u,u) }
        */

    p1pt = Pmat + STAID_LEN + (DEG_FREE*MAT_DIM) + N_PARMS ;
    p2pt = Pmat + STAID_LEN + N_PARMS ;

    b1pt = RBmat + STAID_LEN + (DEG_FREE*MAT_DIM) + N_PARMS ;
    b2pt = RBmat + STAID_LEN + (DEG_FREE*MAT_DIM) + N_PARMS + DEG_FREE ;
    Rbound_rows( Rvu_rows, b1pt, p1pt, b2pt, p2pt ) ;

        /* (v,v) block :
                B(v,u) * P(v,v) + B(v,v) * { beta * P(v,v) -
                                                    P(u,v) }
        */

    p1pt = Pmat + STAID_LEN + (DEG_FREE*MAT_DIM) + N_PARMS + DEG_FREE ;
    p2pt = Pmat + STAID_LEN + N_PARMS + DEG_FREE ;

/*        Same parts of RBmat as are used to find vu_rows. */
    Rbound_rows( Rvv_rows, b1pt, p1pt, b2pt, p2pt ) ;

/*
        Get the contibutions to Rw[] that arise from
        errors in the computation of the image of the
        prism's center.
*/
    for( j=0 ; j < DEG_FREE ; j++ ) {
        center_err[j+N_PARMS] = Bu_rows[j] * DBL_ERR ;
        center_err[j+N_PARMS+DEG_FREE] = Bv_rows[j] * DBL_ERR ;
        apmMultiply( Rcenter_err[j+N_PARMS], RBu_rows[j], max_error ) ;
        apmMultiply( Rcenter_err[j+N_PARMS+DEG_FREE], RBu_rows[j], 
                                                      max_error ) ;
    }

/*
        Compute the components of w[].
*/
    wu_pt = &Rw[N_PARMS] ;
    wv_pt = &Rw[N_PARMS + DEG_FREE] ;
    for( j=0 ; j < DEG_FREE ; j++, wu_pt++, wv_pt++ ) {
        apmCalc( *wu_pt, Rup_rows[j], Ruu_rows[j], Ruv_rows[j], max_error,
                         APM_ADD, APM_ADD, APM_ADD, NULL ) ;
        apmCalc( *wv_pt, Rvp_rows[j], Rvu_rows[j], Rvv_rows[j], max_error,
                         APM_ADD, APM_ADD, APM_ADD, NULL ) ;
    }

    /*
        Include errors due to miscalculation of
        prism's center.
    */

    for( j= N_PARMS ; j < MAT_DIM ; j++ ) 
        apmCalc( Rw[j], Rw[j], Rcenter_err[j], APM_ADD, NULL ) ;

    return ;
}
/* +++++++++++++++++++++++++++++++ */

Rbound_rows( rows, first_b, first_p, second_b, second_p ) 

APM  *rows, *first_b, *second_b, *first_p, *second_p ;
{
/*
        Obtain upper bounds on the sums of the absolute 
        values of rows of matricies given by expressions
        like:
            B1 * S1 + B2 * ( [beta] * S1 - S2 ).

        Expressions like these arise in cr_rows() above.
        The idea is to cast these rows as bounded expressions
        and then use the usual machinery to find their limits.
*/
    int            j, k ;
    APM     *bpt_a, *bpt_b, *ppt_a, *ppt_b, *end_row, *cpt ;

/*
        Evaluate the constant part of the matrix expression.
        It's :
            (B1 + 2.0 * B2) * S1  -  B2 * S2 
*/
    cpt = Rconst_mat ;
    for( j=0 ; j < DEG_FREE ; j++ ) {

        bpt_a = first_b + j * MAT_DIM ;
        bpt_b = second_b + j * MAT_DIM ;
        for( k=0 ; k < DEG_FREE ; k++ ) {
            apmAssignLong( *cpt, 0L, 0, BASE ) ;

            ppt_a = first_p + k ;
            ppt_b = second_p + k ;
            for( end_row = bpt_a + DEG_FREE ; bpt_a < end_row ; ) {
                apmCalc( *cpt, *cpt, *bpt_a, 
                                     *bpt_b, two, APM_MUL, 
                                     APM_ADD,
                                     *ppt_a, APM_MUL,
                                     *bpt_b, *ppt_b,
                                     APM_MUL, APM_SUB, 
                                     APM_ADD, NULL ) ;

                bpt_a++, bpt_b++ ;
                ppt_a += MAT_DIM ;
                ppt_b += MAT_DIM ;
            }

            bpt_a -= DEG_FREE ;
            bpt_b -= DEG_FREE ;
            cpt++ ;
        }
    }

    cpt = Rconst_mat ;
    for( j=0 ; j < DEG_FREE ; j++ ) {
        apmAssignLong( rows[j], 0L, 0, BASE ) ;

        bpt_a = second_b + j * MAT_DIM ;
        bpt_b = bpt_a + 1 ;
        for( k=0 ; k < DEG_FREE ; k++ ) {
            ppt_a = first_p + k ;
            ppt_b = ppt_a + MAT_DIM ;

/*   a * sin( v[0] ) term */
            apmMultiply( cr_scratch, *bpt_a, *ppt_a ) ;
            apmNegate( Rbeta_prod.terms[0].coef, cr_scratch ) ;

/*   c * sin( v[0] + v[1] )  term */
            apmCalc( cr_scratch, *bpt_a, *bpt_b, APM_ADD,
                                 *ppt_a, *ppt_b, APM_ADD,
                                        APM_MUL, NULL ) ;
            apmNegate( Rbeta_prod.terms[1].coef, cr_scratch ) ;

/*   b * sin( v[0] + v[1] )  term */
            apmMultiply( cr_scratch, *bpt_b, *ppt_b ) ;
            apmNegate( Rbeta_prod.terms[2].coef, cr_scratch ) ;

            apmAssign( Rbeta_prod.const, *cpt++ ) ;
            Rbound_expr( &Rbeta_prod ) ;

            RmaxAbs( cr_scratch, Rbeta_prod.bound.ub, 
                                 Rbeta_prod.bound.lb ) ;
            apmCalc( rows[j], rows[j], cr_scratch, APM_ADD, NULL ) ;
        }
    }
}
/* ++++++++++++++++++++++++++++++ */

Rset_inverse( mat ) 

APM   *mat ;
{
    APM   *end_row, *end_block, *end_col ;
    APM   *ipt_a, *ipt_b, *ipt_c, *ipt_set, *mpt_a, *mpt_b ;

    if( isNewPrism == YES ) { 
        end_block = RBmat + N_PARMS * (MAT_DIM + 1) ;
        for( ipt_a=RBmat, mpt_a=mat ; ipt_a < end_block ; ) {
            apmDivide( *ipt_a, precision, (APM)NULL, one, *mpt_a ) ;

            mpt_a += MAT_DIM + 1 ;
            ipt_a += MAT_DIM + 1 ;
        }

        isNewPrism = NO ;
    }
    
    Rinvert_corner( mat ) ;

/* 
    Set the (u,p) part of the inverse.
*/
    ipt_a = RBmat + STAID_LEN + N_PARMS ;
    ipt_b = RBmat + STAID_LEN + N_PARMS + DEG_FREE ;

    ipt_set = RBmat + STAID_LEN  ;
    end_block = ipt_set + (MAT_DIM * DEG_FREE) ;
    for( ; ipt_set < end_block ; ipt_set += TWO_DF ) {
        ipt_c = RBmat ;

        mpt_a = mat + STAID_LEN ;
        mpt_b = mat + STAID_LEN + (DEG_FREE * MAT_DIM) ;

        end_row = ipt_set + N_PARMS ;
        for( ; ipt_set < end_row ; ipt_set++ ) {
            apmAssignLong( *ipt_set, 0L, 0, BASE ) ; 

            end_col = mpt_a + (DEG_FREE * MAT_DIM) ;
            for( ; mpt_a < end_col ; mpt_a += MAT_DIM ) {
                apmCalc( *ipt_set, *ipt_a, *mpt_a, APM_MUL,
                                   *ipt_b, *mpt_b, APM_MUL,
                                   APM_ADD, APM_NEG, 
                                   *ipt_set, APM_ADD, NULL ) ;

                ipt_a++ ;
                ipt_b++ ;
                mpt_b += MAT_DIM ;
            }
            apmCalc( *ipt_set, *ipt_set, *ipt_c, APM_MUL, NULL ) ;

            ipt_a -= DEG_FREE ;
            ipt_b -= DEG_FREE ;
            ipt_c += MAT_DIM + 1 ;

            mpt_a -= (MAT_DIM * DEG_FREE) - 1 ;
            mpt_b -= (MAT_DIM * DEG_FREE) - 1 ;
        }

        ipt_a += MAT_DIM ;
        ipt_b += MAT_DIM ;
        mpt_a -= DEG_FREE ;
        mpt_b -= DEG_FREE ;
    }

/* 
    Set the (v,p) part of the inverse.
*/
    ipt_a = RBmat + STAID_LEN + N_PARMS + (DEG_FREE * MAT_DIM) ;
    ipt_b = RBmat + STAID_LEN + N_PARMS + (DEG_FREE*MAT_DIM) + DEG_FREE ;

    ipt_set = RBmat + STAID_LEN + (DEG_FREE * MAT_DIM)  ;
    end_block = ipt_set + (MAT_DIM * DEG_FREE) ;
    for( ; ipt_set < end_block ; ipt_set += TWO_DF ) {
        ipt_c = RBmat ;

        mpt_a = mat + STAID_LEN ;
        mpt_b = mat + STAID_LEN + (DEG_FREE * MAT_DIM) ;

        end_row = ipt_set + N_PARMS ;
        for( ; ipt_set < end_row ; ipt_set++ ) {
            apmAssignLong( *ipt_set, 0L, 0, BASE ) ; 

            end_col = mpt_a + (DEG_FREE * MAT_DIM) ;
            for( ; mpt_a < end_col ; mpt_a += MAT_DIM ) {
                apmCalc( *ipt_set, *ipt_a, *mpt_a, APM_MUL,
                                   *ipt_b, *mpt_b, APM_MUL,
                                   APM_ADD, APM_NEG, 
                                   *ipt_set, APM_ADD, NULL ) ;

                ipt_a++ ;
                ipt_b++ ;
                mpt_b += MAT_DIM ;
            }
            apmCalc( *ipt_set, *ipt_set, *ipt_c, APM_MUL, NULL ) ;

            ipt_a -= DEG_FREE ;
            ipt_b -= DEG_FREE ;
            ipt_c += MAT_DIM + 1 ;

            mpt_a -= (MAT_DIM * DEG_FREE) - 1 ;
            mpt_b -= (MAT_DIM * DEG_FREE) - 1 ;
        }

        ipt_a += MAT_DIM ;
        ipt_b += MAT_DIM ;
        mpt_a -= DEG_FREE ;
        mpt_b -= DEG_FREE ;
    }
}
/* +++++++++++++++++++++ */

Rinvert_corner( mat )

APM  *mat ;
{
/*
        Set up matrices to prepare 'em for use by Rgauss().
        Note that we use the matirx called const_mat[].  
        At the times this function is called const_mat[] 
        doesn't contain anything important.
*/
    int             j ;
    APM      *end_row, *mpt, *bpt, *cpt ;

/*
        Copy the matrix.
*/
    mpt = mat + STAID_LEN + N_PARMS ;
    for( j=0 ; j < TWO_DF ; j++ ) {

        cpt = Rcopy_rows[j] ;
        end_row = mpt + TWO_DF ;
        while( mpt < end_row ) 
            apmAssign( *cpt++, *mpt++ ) ;

        mpt += N_PARMS ;
    }

/* 
        Do the inversion.
*/
    Rgauss( Rcopy_rows ) ;

/*
        Copy the answer.
*/

    bpt = RBmat + STAID_LEN + N_PARMS ;
    for( j=0 ; j < TWO_DF ; j++ ) {

        cpt = Rcopy_rows[j] ;
        end_row = bpt + TWO_DF ;
        while( bpt < end_row ) 
            apmAssign( *bpt++, *cpt++ ) ;

        bpt += N_PARMS ;
    }
}
\end{verbatim}

%% file: appendices/code/ff_fatten.tex
\begin{verbatim} 
# include  <stdio.h>
# include  <math.h>
# include  "apm.h"
# include  "apmSpecial.h"
# include  "converse.h"

/* ++++++++++++++++++++++++++++++++++++++++++ */

Rfxed_form( Amat, Deriv, Prizmat )

APM  *Amat, *Deriv, *Prizmat ;
/*
     Prepares the matrix called "A" in my notes.  Eventually we want to
     have A = DF*Priz, but early in a calculation, when Priz is singular,
     we want to fatten A up by requiring it to have a certain fixed form.
     In the inerest of speed we have coded the calculations below in 
     terms of pointers.  Our hope is that the resulting function will
     scream along at ultrasonic speed.  Unfortunately it is quite
     unreadable.
*/
{
    APM                  *Aend, *Aend2, *Dend, *Pend, *Pend2  ;
    register APM  *Apt, *Apt2, *Dpt, *Ppt, *Ppt2 ;

/*
        Copy the few terms which appear in the top rows of Amat.
*/
    Aend = Amat + N_PARMS * (MAT_DIM + 1) ;
    for( Apt = Amat, Ppt = Prizmat ; Apt < Aend ; Apt += (MAT_DIM + 1 ),
                                                  Ppt += (MAT_DIM + 1 ) ) 
        apmAssign( *Apt, *Ppt ) ;

/*
        Clear out those parts of Amat which change from iteration to 
        iteration. 
*/

    Aend = Amat + MAT_SZ ;
    for( Apt = Amat + STAID_LEN ; Apt < Aend ; Apt++ )
        apmAssignLong( *Apt, 0L, 0, 0 ) ;

/*
        Set the (u,p) part of A
        It's equal to the (v,p) part of Prizmat.
*/
    
    Aend = Amat + STAID_LEN + (DEG_FREE * MAT_DIM) - TWO_DF ;
    Ppt = Prizmat + STAID_LEN + (DEG_FREE * MAT_DIM) ;
    for( Apt = Amat + STAID_LEN ; Apt < Aend ; Apt += TWO_DF ) {
        for( Pend = Ppt + N_PARMS ; Ppt < Pend ; Ppt++, Apt++ )
            apmCalc( *Apt, *Apt, *Ppt, APM_ADD, NULL ) ;

        Ppt += TWO_DF ;
    }
/*
        Set the (v,p) part - three terms.
*/
        /*  First term - equal to Deriv(v,p) * Prizmat(p,p)  */

    Dpt = Deriv + STAID_LEN + (DEG_FREE * MAT_DIM) ;
    Apt = Amat + STAID_LEN + (DEG_FREE * MAT_DIM) ;
    for( Aend = Apt + (DEG_FREE*MAT_DIM) ; Apt < Aend ; Apt += TWO_DF ) {
        Ppt = Prizmat ;
        for( Dend = Dpt + N_PARMS ; Dpt < Dend ; Dpt++ ) {
            apmMultiply( *Apt++, *Dpt, *Ppt ) ;
            Ppt += MAT_DIM + 1 ;
        }

        Dpt += TWO_DF ;
    }

        /*  Second term - equal to negative Prizmat(u,p) */

    Ppt = Prizmat + STAID_LEN ;
    Apt = Amat + STAID_LEN + (DEG_FREE * MAT_DIM) ;
    for( Pend = Ppt + (DEG_FREE * MAT_DIM) ; Ppt < Pend ; Ppt += TWO_DF ) {
        for( Aend = Apt + N_PARMS ; Apt < Aend ; Apt++, Ppt++ )
            apmCalc( *Apt, *Apt, *Ppt, APM_SUB, NULL ) ;

        Apt += TWO_DF ;
    }

        /*  Third term - equal to Deriv(v,v) * Prizmat(v,p)  */

    Dpt = Deriv + STAID_LEN + (DEG_FREE * (MAT_DIM + 1)) + N_PARMS ;
    Dend = Deriv + MAT_SZ ;
    Apt = Amat + STAID_LEN + (DEG_FREE * MAT_DIM) ;
    while( Dpt < Dend ) {
        Ppt = Prizmat + STAID_LEN + (DEG_FREE * MAT_DIM) ;
        Pend = Prizmat + MAT_SZ - TWO_DF ;
        while( Ppt < Pend ) {
            Aend = Apt + N_PARMS ;
            while( Apt < Aend ) {
                apmCalc( *Apt, *Dpt, *Ppt, APM_MUL, *Apt, APM_ADD, NULL ) ;
                Apt++ ;
                Ppt++ ;
            }

            Dpt++ ;
            Ppt += TWO_DF ;
            Apt -= N_PARMS ;
        }

        Dpt += N_PARMS + DEG_FREE ;
        Apt += MAT_DIM ;
    }
/*
        (u,v) part
        equals Priz(v,u) + Priz(v,v)
*/

    Apt = Amat + STAID_LEN + N_PARMS + DEG_FREE ;
    Aend = Amat + STAID_LEN + (DEG_FREE * MAT_DIM) ;
    Ppt = Prizmat + STAID_LEN + (DEG_FREE * MAT_DIM) + N_PARMS ;
    Ppt2 = Ppt + DEG_FREE ;
    while( Apt < Aend ) {
        Pend = Ppt + DEG_FREE ;
        while( Ppt < Pend )  {
            apmCalc( *Apt, *Ppt, *Ppt2, APM_ADD, *Apt, APM_ADD, NULL ) ;
            Apt++ ;
            Ppt++ ;
            Ppt2++ ;
        }

        Apt  += N_PARMS + DEG_FREE ;
        Ppt  += N_PARMS + DEG_FREE ;
        Ppt2 += N_PARMS + DEG_FREE ;
    }
/*
        The (v,u) part
        equal to Deriv(v,v) * { Priz(v,u) + Priz(v,v) },
        which also equals Deriv(v, v) * A(u,v)
*/

    Apt = Amat + STAID_LEN + (DEG_FREE * MAT_DIM) + N_PARMS ;
    Dpt = Deriv + STAID_LEN + (DEG_FREE * MAT_DIM) + N_PARMS + DEG_FREE ;
    Dend = Deriv + MAT_SZ ;
    while( Dpt < Dend ) {
        Apt2 = Amat + STAID_LEN + N_PARMS + DEG_FREE ;
        Aend2 = Apt2 + (DEG_FREE * MAT_DIM) ;
        while( Apt2 < Aend2 ) {
            Aend = Apt + DEG_FREE ;
            while( Apt < Aend ) {
                apmCalc( *Apt, *Apt, *Dpt, *Apt2, APM_MUL, APM_ADD, NULL ) ;
                Apt++ ;
                Apt2++ ;
            }

            Dpt++ ;
            Apt -= DEG_FREE ;
            Apt2 += DEG_FREE + N_PARMS ;
        }

        Apt += MAT_DIM ;
        Dpt += N_PARMS + DEG_FREE ;
    }

/*
        (v,v) part - equals Deriv(v,v) * Priz(v,v) - Priz(u,v) 
*/

        /* First term */
    Apt = Amat + STAID_LEN + (DEG_FREE * MAT_DIM) + N_PARMS + DEG_FREE ;
    Dpt = Deriv + STAID_LEN + (DEG_FREE * MAT_DIM) + N_PARMS + DEG_FREE ;
    Dend = Deriv + MAT_SZ ;
    while( Dpt < Dend ) {
        Ppt = Prizmat + STAID_LEN + (DEG_FREE * MAT_DIM) + N_PARMS + DEG_FREE ;
        Pend = Prizmat + MAT_SZ ;
        while( Ppt < Pend ) {
            Aend = Apt + DEG_FREE ;
            while( Apt < Aend ) {
                apmCalc( *Apt, *Apt, *Dpt, *Ppt, APM_MUL, APM_ADD, NULL ) ;

                Apt++ ;
                Ppt++ ;
            }

            Dpt++ ;
            Apt -= DEG_FREE ;
            Ppt += DEG_FREE + N_PARMS ;
        }

        Apt += MAT_DIM ;
        Dpt += N_PARMS + DEG_FREE ;
    }
        /* Second term */
    Apt = Amat + STAID_LEN + (DEG_FREE * MAT_DIM) + N_PARMS + DEG_FREE ;
    Ppt = Prizmat + STAID_LEN + N_PARMS + DEG_FREE ;
    Pend = Ppt + (MAT_DIM * DEG_FREE) ;
    while( Ppt < Pend ) {
        Aend = Apt + DEG_FREE ;
        while( Apt < Aend ) {
            apmCalc( *Apt, *Apt, *Ppt, APM_SUB, NULL ) ;

            Apt++ ;
            Ppt++ ;
        }

        Ppt += N_PARMS + DEG_FREE ;
        Apt += N_PARMS + DEG_FREE ;
    }
}
\end{verbatim}

%% file: appendices/code/ff_rows.tex
\begin{verbatim} 
# include <stdio.h>
# include <math.h>
# include "apm.h"
# include "apmSpecial.h"
# include "converse.h"
# include "bounding.h"
# include "rows.h"

APM          Rerr_star ;
APM          ff_scratch ;
APM          Rcenter_err[MAT_DIM] ;
APM          Rdet_vu, Rdet_uv, Rstar ;
APM          RAvv_star, RAuvInv_star ;
APM          Rb_star, Rbd_star, Rgd_star ;
APM          RPvv_star, RPvp_star, RPvu_star ;

double          beta_star() ;
double          center_err[MAT_DIM] ;
Bdd_dbl          *ff_factors[NUM_FACTS] ;
Bdd_term  ff_terms[NUM_TERMS] ;
Bdd_expr  beta[3] ;

Bdd_apm          *Rff_factors[NUM_FACTS] ;
Bapm_term Rff_terms[NUM_TERMS] ;
Bapm_expr Rbeta[3] ;
/* ++++++++++++++++++++++++++++++ */

init_ffRows()
/*
        Set up the expressions and terms as described in my notes
        from 11/14.
*/
{
    int              j, k ;
    Bdd_dbl   **dpt ;
    Bdd_apm   **apt ;
    Bdd_term  *tpt ;
    Bapm_term *Rtpt ;

/*
    Set up some APM's to be used to hold intermediate 
    results.
*/
    Rstar = apmNew( BASE ) ;
    Rdet_uv = apmNew( BASE ) ;
    Rdet_vu = apmNew( BASE ) ;
    Rb_star = apmNew( BASE ) ;
    Rbd_star = apmNew( BASE ) ;
    Rgd_star = apmNew( BASE ) ;
    Rerr_star = apmNew( BASE ) ;
    RAvv_star = apmNew( BASE ) ;
    RPvv_star = apmNew( BASE ) ;
    RPvp_star = apmNew( BASE ) ;
    RPvu_star = apmNew( BASE ) ;
    ff_scratch = apmNew( BASE ) ;
    RAuvInv_star = apmNew( BASE ) ;

    for( j = 0 ; j < MAT_DIM ; j++ )
        Rcenter_err[j] = apmNew( BASE ) ;

/*
        Set the number of terms in the bounded expressions
*/

    beta[0].nterms = Rbeta[0].nterms = 2 ;
    beta[1].nterms = Rbeta[1].nterms = 1 ;
    beta[2].nterms = Rbeta[2].nterms = 2 ;

/*        
        Assign terms
*/

    tpt = ff_terms ;
    Rtpt = Rff_terms ;
    for( j=0 ; j < 3 ; j++ ) {
        beta[j].terms = tpt ;
        Rbeta[j].terms = Rtpt ;
        tpt += beta[j].nterms ;
        Rtpt += Rbeta[j].nterms ;
    }

/*
        Set nfactors.
*/
    
    Rbeta[0].terms[0].nfactors = beta[0].terms[0].nfactors = 1  ;
    Rbeta[0].terms[1].nfactors = beta[0].terms[1].nfactors = 1  ;
    Rbeta[1].terms[0].nfactors = beta[1].terms[0].nfactors = 1  ;
    Rbeta[2].terms[0].nfactors = beta[2].terms[0].nfactors = 1  ;
    Rbeta[2].terms[1].nfactors = beta[2].terms[1].nfactors = 1  ;

/*
        Assign factors. 
*/

    dpt = ff_factors ;
    apt = Rff_factors ;
    for( j=0 ; j < 3 ; j++ ) {
    /* 
        beta
    */
        for( k=0 ; k < beta[j].nterms ; k++ ) {
            beta[j].terms[k].factors = dpt ;
            Rbeta[j].terms[k].factors = apt ;

            dpt += beta[j].terms[k].nfactors ;
            apt += Rbeta[j].terms[k].nfactors ;
        }
    }

/*
        Set up those of the "bound" attributes which are 
        bounded APM's.
*/

    for( j=0 ; j < NUM_TERMS ; j++ ) {
        newBapm( Rff_terms[j].bound, BASE ) ;
    }

    for( j=0 ; j < 3 ; j++ ) {
        newBapm( Rbeta[j].bound, BASE ) ;
    }

/*
        Set up the terms and expressions. 
*/

/* beta */
 
        /* beta[0] = 2.0 - a * sin(v[0]) - c * sin(v[0] + v[1])  */
    beta[0].const = 2.0, Rbeta[0].const = two ;
        beta[0].terms[0].coef = -1.0 ;
        Rbeta[0].terms[0].coef = neg_one ;

            beta[0].terms[0].factors[0] = &a_sin.bound ;
            Rbeta[0].terms[0].factors[0] = &Ra_sin.bound ;

        beta[0].terms[1].coef = -1.0 ;
        Rbeta[0].terms[1].coef = neg_one ;

            beta[0].terms[1].factors[0] = &c_sin.bound;
            Rbeta[0].terms[1].factors[0] = &Rc_sin.bound;

        /* beta[1] = - 2.0 * c * sin( v[0] + v[1] )        */
    beta[1].const = 0.0, Rbeta[1].const = zero ;
        beta[1].terms[0].coef = -2.0 ;
        Rbeta[1].terms[0].coef = neg_two ;

            beta[1].terms[0].factors[0] = &c_sin.bound;
            Rbeta[1].terms[0].factors[0] = &Rc_sin.bound;

        /* beta[2] = 2.0 - b * sin(v[1]) - c * sin(v[1] + v[0])  */
    beta[2].const = 2.0, Rbeta[2].const = two ;
        beta[2].terms[0].coef = -1.0 ;
        Rbeta[2].terms[0].coef = neg_one ;

            beta[2].terms[0].factors[0] = &b_sin.bound;
            Rbeta[2].terms[0].factors[0] = &Rb_sin.bound;

        beta[2].terms[1].coef = -1.0 ;
        Rbeta[2].terms[1].coef = neg_one ;

            beta[2].terms[1].factors[0] = &c_sin.bound ;
            Rbeta[2].terms[1].factors[0] = &Rc_sin.bound ;
}
/* +++++++++++++++++++++++++++++++++ */

Rff_rows( w, Amat, Deriv, Priz ) 

APM         *w, *Amat, *Deriv ;
RPrism   *Priz ;
/*
    Obtain bounds on the sums of the absolute values of
    the entries in the rows of 
       -1
    [A]   * Deriv * Pmat,

    put the results in w.
*/
{
    APM   *apt, *mpt, *end_row, *end_mat, *Pmat ;

/*
        Check that A(u,v) is invertible.  If not, die.
*/
    Pmat = Priz->matrix ;

    apt = Amat + STAID_LEN + N_PARMS + DEG_FREE ;
    apmMultiply( Rdet_uv, *apt, *(apt + MAT_DIM + 1) ) ;
    apt++ ;
    apmCalc( Rdet_uv, Rdet_uv, *apt, *(apt + MAT_DIM -1), 
                                APM_MUL, APM_SUB, NULL ) ;
    apmAbsoluteValue( ff_scratch, Rdet_uv ) ;

    if( apmCompare( ff_scratch, max_error ) != 1 ) {
        fprintf( stderr, 
                 "The determinant of A(u,v) is too small. Died. \n" ) ;
        fprintf( stderr, "\t %.12e \n", apmtodbl( ff_scratch ) ) ;
        cease() ;
    }
/*
        Call functions which calculate upper bound on the 
        sums of the elements of various matrices.
        Before any bounding of matrices, one must invoke
        global_bounds( Pmat ) to set such global variables,
        as cos_one, and sin_sum.  It is called in Rtry_prism().
*/
    Rbeta_star( Rb_star ) ;
    Rbeta_dif_star( Rbd_star, Deriv ) ;
    Rgamdif_star( Rgd_star, Deriv, Pmat ) ;

/*
        Find sums of the absolute values of the entries
        of Pmat(v,v), Pmat(v,u), and Pmat(v,p)
*/
    end_mat = Pmat + MAT_SZ ;

    apmAssign( RPvv_star, zero ) ; 
    mpt = Pmat + STAID_LEN + (DEG_FREE * MAT_DIM) + N_PARMS + DEG_FREE ;
    for( ; mpt < end_mat ; mpt += (N_PARMS + DEG_FREE) ) {
        for( end_row = mpt + DEG_FREE ; mpt < end_row ; mpt++ ) {
            apmCalc( RPvv_star, RPvv_star, *mpt, APM_ABS, 
                                                   APM_ADD, NULL ) ;
        }
    }

    apmAssign( RPvu_star, zero ) ;
    mpt = Pmat + STAID_LEN + (DEG_FREE * MAT_DIM) + N_PARMS ;
    for( ; mpt < end_mat ; mpt += (N_PARMS + DEG_FREE) ) {
        for( end_row = mpt + DEG_FREE ; mpt < end_row ; mpt++ ) {
            apmCalc( RPvu_star, RPvu_star, *mpt, APM_ABS,
                                                   APM_ADD, NULL ) ;
        }
    }

    apmAssign( RPvp_star, zero ) ; 
    mpt = Pmat + STAID_LEN + (DEG_FREE * MAT_DIM) ;
    for( ; mpt < end_mat ; mpt += TWO_DF ) {
        for( end_row = mpt + N_PARMS ; mpt < end_row ; mpt++ ) {
            apmCalc( RPvp_star, RPvp_star, *mpt, APM_ABS,
                                                   APM_ADD, NULL ) ;
        }
    }

    apmAssign( RAvv_star, RSmBlock_err ) ;
    mpt = Amat + STAID_LEN + DEG_FREE * MAT_DIM + DEG_FREE + N_PARMS ;
    for( ; mpt < end_mat ; mpt += TWO_DF ) {
        for( end_row = mpt + N_PARMS ; mpt < end_row ; mpt++ ) {
            apmCalc( RAvv_star, RAvv_star, *mpt, 
                                APM_ABS, APM_ADD, NULL ) ;
        }
    }

    apmAssign( RAuvInv_star, RSmBlock_err ) ;
    mpt = Amat + STAID_LEN + N_PARMS + DEG_FREE ;
    for( ; mpt < end_mat ; mpt += TWO_DF ) {
        for( end_row = mpt + N_PARMS ; mpt < end_row ; mpt++ ) {
            apmCalc( RAuvInv_star, RAuvInv_star, *mpt, 
                                   APM_ABS, APM_ADD, NULL ) ;
        }
    }

    apmDivide( ff_scratch, precision, (APM) NULL, 
                           RAuvInv_star, Rdet_uv ) ;
    apmAssign( RAuvInv_star, ff_scratch ) ;

/*
        Check that A(v,u) is invertible.  If not, die.
        If it is, set the harder-to-compute elements of w.
*/

    apt = Amat + STAID_LEN + N_PARMS + (DEG_FREE * MAT_DIM) ;
    apmMultiply( Rdet_vu, *apt, *(apt + MAT_DIM + 1) ) ;
    apt++ ;
    apmCalc( Rdet_vu, Rdet_vu, *apt, *(apt + MAT_DIM - 1),
                                     APM_MUL, APM_SUB, NULL ) ;
    apmAbsoluteValue( ff_scratch, Rdet_vu ) ;

    if( apmCompare( ff_scratch, max_error ) != 1 ) {
        fprintf( stderr, 
                 "The determinant of A(v,u) is too small. Died. \n") ;
        fprintf( stderr, "\t %.12e \n", apmtodbl( ff_scratch ) ) ;
        cease() ;
    }
/*
        Note that the sums below seem to contain some misplaced
        elements of Amat.  These are to be thought of as elements
        of A(v,u) inverse.
*/
    else {
        apmCalc( w[3], Amat[MAT_SZ-DEG_FREE-1], APM_ABS, 
                       Amat[STAID_LEN+(DEG_FREE*MAT_DIM)+N_PARMS+1],
                       APM_ABS, max_error, APM_ADD, APM_ADD, NULL ) ;

        apmCalc( w[4], Amat[MAT_SZ-TWO_DF], APM_ABS, 
                       Amat[STAID_LEN+(DEG_FREE*MAT_DIM)+N_PARMS],
                       APM_ABS, max_error, APM_ADD, APM_ADD, NULL ) ;

        apmCalc( Rerr_star, RAvv_star, RAuvInv_star, APM_MUL, 
                                       one, APM_ADD, NULL );
        apmCalc( Rcenter_err[3], w[3], Rerr_star, max_error,
                                       APM_MUL, APM_MUL, NULL ) ;
        apmCalc( Rcenter_err[4], w[4], Rerr_star, max_error, 
                                       APM_MUL, APM_MUL, NULL ) ;
        apmMultiply( Rcenter_err[5], RAuvInv_star, max_error ) ;
        apmAssign( Rcenter_err[6], Rcenter_err[5] ) ;

        apmCalc( Rstar, RPvp_star, RPvv_star, APM_ADD,
                        Rbd_star, APM_MUL,
                        Rb_star, RPvu_star, APM_MUL,
                        Rgd_star, APM_ADD, APM_ADD, NULL ) ;

        apmCalc( ff_scratch, Rcenter_err[3], Rstar, w[3], 
                                 APM_MUL, APM_ADD, NULL ) ;
        apmDivide( w[3], precision, (APM) NULL, ff_scratch, Rdet_vu ) ;
        apmCalc( ff_scratch, Rcenter_err[4], Rstar, w[4], 
                                 APM_MUL, APM_ADD, NULL ) ;
        apmDivide( w[4], precision, (APM) NULL, ff_scratch, Rdet_vu ) ;
        apmAdd( w[5], one, Rcenter_err[5] ) ;
        apmAdd( w[6], one, Rcenter_err[6] ) ;
    }

    return ;
}
/* +++++++++++++++++++++++++++++++ */

Rbeta_star( answer )

APM        answer ;
{
    Rbound_expr( &Rbeta[0] ) ;
    Rbound_expr( &Rbeta[1] ) ;
    Rbound_expr( &Rbeta[2] ) ;

    RmaxAbs( answer, Rbeta[0].bound.ub, Rbeta[0].bound.lb )  ;
    RmaxAbs( Rrow_abs[0], Rbeta[1].bound.ub, Rbeta[1].bound.lb ) ;
    RmaxAbs( Rrow_abs[1], Rbeta[2].bound.ub, Rbeta[2].bound.lb ) ;

    apmCalc( answer, answer, Rrow_abs[0], Rrow_abs[1], 
                                        APM_ADD, APM_ADD, NULL ) ; 
}
\end{verbatim}

%% file: appendices/code/gauss.tex
\begin{verbatim} 
# include <stdio.h>
# include <math.h>
# include "apm.h"
# include "apmSpecial.h"
# include "converse.h"

                                                apmAssign(y, t) )
/*
        The Numerical Recipes Gauss-Jordan matrix inverter as adaptaed
        for a converse KAM code.
        I have removed the dimension arguments n and m and replaced 
        them with TWO_DF and 1.  I have also changed all the floats 
        into doubles and replaced some automatically allocated
        arrays with arrays of fixed dimension.  Finally, I have 
        replaced the error handling code with some of my own.

        Rgauss, the rigorous version, also does a host of checks to
        guarantee that the inverse it produces, when multiplied by 
        the original matrix, a, gives something equal to the 
        identity to the accuracy specified by the global variable,
        "precision".
*/

int        extra_dp, last_inv_dp ;
int        inv_depth ;        /* Used to make sure that we don't keep trying
                           to invert singular matrices by using 
                           ever increasing precision.
                                                                        */

APM        a_abs, Rbig, Rdum, Rpivinv, Rtemp ;
APM        Rrow_max, Rcol_max, Rmat_min, Rmat_max ;
APM        *Rmat[TWO_DF], Rmat_block[4*DF_SQ] ;
APM        Rdiv_err, Rrow_err, Rinv_err, Rtotal_err, Rpiv_err ;
/* ++++++++++++++++++++++++++++++++ */

initGauss()
{
    int   j, k ;
    APM   *mpt ;

    inv_depth = 0 ;
    extra_dp = 0 ;

    Rbig = apmNew( BASE ) ;
    Rdum = apmNew( BASE ) ;
    a_abs = apmNew( BASE ) ;
    Rtemp = apmNew( BASE ) ;
    Rpivinv = apmNew( BASE ) ;
    Rinv_err = apmNew( BASE ) ;
    Rrow_err = apmNew( BASE ) ;
    Rpiv_err = apmNew( BASE ) ;
    Rdiv_err = apmNew( BASE ) ;
    Rrow_max = apmNew( BASE ) ;
    Rcol_max = apmNew( BASE ) ;
    Rmat_min = apmNew( BASE ) ;
    Rmat_max = apmNew( BASE ) ;
    Rtotal_err = apmNew( BASE ) ;

    mpt = Rmat_block ;
    for( j=0 ; j < TWO_DF ; j++ ) {
        Rmat[j] = mpt ;
        for( k=0 ; k < TWO_DF ; k++ ) 
            *mpt++ = apmNew( BASE ) ;
    }
}
/* ++++++++++++++++++++++++++++ */

Rgauss( a )

APM  **a ;
{
        int indxc[TWO_DF],indxr[TWO_DF],ipiv[TWO_DF];
        int i,icol,irow,j,k,l,ll;
        int inv_dp, err_dp ;

        if( ++inv_depth > MAX_RECUR ) {
            fprintf( stderr, "Singular matrix in Rgauss. Died. \n" ) ;
            cease() ;
        }

        for( j=0 ; j < TWO_DF ; j++ ) {
            ipiv[j] = 0 ;
            indxr[j] = 0 ;
            indxc[j] = 0 ;
        }
/*
            If this is the attempt to invert a,
            copy the matrix in case of a loss of precision. 
            Also, choose
            the precision to which to do the inversion calculations.
*/
        if( inv_depth == 1 ) {
            copyRmat( Rmat, a ) ;
            inv_dp =  choosePrecis( a ) ;
        }
        else {
            if( extra_dp == 0 ) 
                inv_dp = last_inv_dp + DFLT_XDP ;
            else
                inv_dp = last_inv_dp + extra_dp ;
        }
        last_inv_dp = inv_dp ;                

/*
        Initialize the error propagation stuff.
*/
        apmAssignLong( Rdiv_err, 1L, -inv_dp, BASE ) ;
        apmAssignLong( Rinv_err, 0L, 0, BASE ) ;
        apmAssign( Rpiv_err, Rinv_err ) ;

        for (i=0;i<TWO_DF;i++) {
            apmAssignLong( Rbig, 0L, 0, BASE ) ;
            for (j=0;j<TWO_DF;j++) {
                if (ipiv[j] != 1) {
                    for (k=0;k<TWO_DF;k++) {
                        if (ipiv[k] == 0) {
                            apmAbsoluteValue( a_abs, a[j][k] ) ;
                            if( apmCompare(a_abs, Rbig) != -1 ) {
                                apmAssign( Rbig, a_abs ) ;
                                irow=j;
                                icol=k;
                            }
                        } 
                        else if (ipiv[k] > 1) {
                            fprintf( stderr,
                                "Singular matrix in gauss. Died.\n" ) ;
                            cease() ;
                        }
                    }
                }
            }

            ++(ipiv[icol]);
            if(irow != icol) {
                for (l=0;l<TWO_DF;l++) 
                    Rm_swap(a[irow][l],a[icol][l],Rtemp) ;
            }

            indxr[i]=irow;
            indxc[i]=icol;
/*
                Check that the pivot interval does not
                contain zero.  If it does, restart the
                calculation and carry more decimal places.
*/
            apmCalc( Rtemp, a[icol][icol], APM_ABS,
                            Rinv_err, APM_SUB, NULL ) ;
            if( apmCompare( Rtemp, zero ) != 1 ) { 
                copyRmat( a, Rmat ) ;
                Rgauss( a ) ;
                return ;
            }
/*
                Get the new pivot error.  It is here that we face
                the possibility of catastrophic loss of precision.
*/
            apmDivide( Rpiv_err, inv_dp, (APM)NULL, Rinv_err, Rtemp ) ;
            apmCalc( Rpiv_err, Rpiv_err, Rdiv_err, Rdiv_err, 
                                         APM_ADD, APM_ADD, NULL ) ;
            apmDivide(Rpivinv,inv_dp,(APM)NULL,one,a[icol][icol]) ;
            apmAssignLong( a[icol][icol], 1L, 0, BASE ) ;
          
            apmAssignLong( Rrow_max, 0L, 0, BASE ) ;
            for (l=0;l<TWO_DF;l++) { 
                if( l != icol ) {
                    apmAbsoluteValue( Rtemp, a[icol][l] ) ;
                    if( apmCompare( Rtemp, Rrow_max ) < 0 )
                        apmAssign( Rrow_max, Rtemp ) ;
                }

                apmCalc(a[icol][l], a[icol][l], Rpivinv,APM_MUL,NULL) ;
            }

/*
                Get a bound on the size of the errors in the elements
                of the pivot row.
*/
            apmCalc( Rrow_err,  Rinv_err, Rpivinv, APM_MUL,
                                Rrow_max, Rinv_err, APM_ADD,
                                Rpiv_err, APM_MUL, APM_ADD, NULL ) ;

            apmAssignLong( Rcol_max, 0L, 0, BASE ) ;
            for (ll=0;ll<TWO_DF;ll++) {
                if (ll != icol) {
                    apmAssign( Rdum, a[ll][icol] ) ;
                    apmAbsoluteValue( Rtemp, Rdum ) ;
                    if( apmCompare( Rtemp, Rcol_max ) == 1 )
                        apmAssign( Rcol_max, Rtemp ) ;

                    apmAssignLong( a[ll][icol], 0L, 0, BASE ) ;
                    for (l=0;l<TWO_DF;l++) 
                        apmCalc( a[ll][l], a[ll][l], a[icol][l], Rdum,
                                           APM_MUL, APM_SUB, NULL ) ;
                }
            }

/*
                Calculate the new upper bound on errors in the matrix.
*/
            apmCalc( Rinv_err,  Rrow_max, Rrow_err, APM_ADD, 
                                Rinv_err, APM_MUL,
                                Rcol_max, Rrow_err, APM_MUL,
                                Rinv_err, APM_ADD, 
                                APM_ADD, APM_ADD, NULL ) ;
/*
        Add an extra Rdiv_err to Rinv_err and truncate everything.
        This will probably speed the calculation considerably.
*/
            apmCalc( Rinv_err, Rinv_err, Rdiv_err, APM_ADD, NULL ) ;

            apmTruncate( Rinv_err, inv_dp ) ;
            for( l = 0 ; l < TWO_DF ; l++ )
                    for( ll=0 ; ll < TWO_DF ; ll++ ) 
                    apmTruncate( a[l][ll], inv_dp ) ;
        }

        for (l=(TWO_DF-1);l>=0;l--) {
            if (indxr[l] != indxc[l])
                for (k=0;k<TWO_DF;k++)
                    Rm_swap(a[k][indxr[l]],a[k][indxc[l]],Rtemp);
        }
/*
        Check the overall size of the error.
        If it is too big, set extra_dp and try again.
*/
        err_dp = -(apmLogBd( Rinv_err ) + OOM_DF) ;
        if( err_dp < precision ) { 
                extra_dp = precision - err_dp + 2 ;
                copyRmat( a, Rmat ) ;
                Rgauss( a ) ;
                return ;
        }

/*
        Tidy up.
        If we reach this line, all is well, the inversion is
        good to the desired precision, so all we want to do is 
        restore the recurrsive variables to their initial state.
*/
    inv_depth = 0 ;
    extra_dp = 0 ;
    return ;
}
/* +++++++++++++++++++++++++++++++++ */

copyRmat( copy, mat ) 

APM   **copy, **mat ;
{
    int  j, k ;

    for( j=0 ; j < TWO_DF ; j++ ) 
        for( k=0 ; k < TWO_DF ; k++ )
            apmAssign( copy[j][k], mat[j][k] ) ;
}
/* ++++++++++++++++++++++++++++++++++ */

choosePrecis( mat )

APM  **mat ;
{
    APM  *mpt, *end_mat ;
    int  oom_min, oom_max, oom_err, oom_twos ;

/*
        Find the minimum and maximum entries of the matrix.
        If none of the entries has absolute value bigger than 
        one, use one as the maximum; this ensures that the 
        resulting inverse will have entries good to at least
        "precision" decimal places.
*/
    mpt = mat[0] ;
    apmAssignLong( Rmat_min, 0L, 0, BASE ) ;
    apmAssignLong( Rmat_max, 1L, 0, BASE ) ;

    for( end_mat = mpt + (TWO_DF*TWO_DF) ; mpt < end_mat ; mpt++ ) {
        apmAbsoluteValue( Rtemp, *mpt ) ;
        if( apmCompare( Rmat_min, Rtemp ) > 0 )
                apmAssign( Rmat_min, Rtemp ) ;
        else if( apmCompare( Rmat_max, Rtemp ) < 0 )
            apmAssign( Rmat_max, Rtemp ) ;
    }

/*
        Do a basic estimate of the number of digits one must carry 
        to get an answer whose precision is as good as the code
        requires.
            First find the orders of magnitude ("oom"'s) of various things.
*/
    oom_max = apmLogBd( Rmat_max ) ;
    oom_twos = (TWO_DF / 3) ;

    oom_err = oom_twos + OOM_DF + (2 * TWO_DF + 1) * abs( oom_max )  ;

    if( oom_err < 0 )
        return( precision ) ;
    else
        return( precision + oom_err ) ;
}
\end{verbatim}

%% file: biblio.tex
%
\newcommand{\ETDS}{{\em Ergodic Theory and Dynamical Systems \/} }
\newcommand{\CMP}{{\em Communications in Mathematical Physics \/}}
\newcommand{\RMP}{{\emReviews of Modern Physics \/}}
\newcommand{\JMP}{{\em Journal of Mathematical Physics \/}}
\newcommand{\Vol}[1]{{\bf {#1}}} 

%
%
	\newpage
	\renewcommand{\baselinestretch}{1.0}
	\large

%% file: thesis.bbl
\normalsize

\begin{thebibliography}{FrednEthel00}

\bibitem[Arn64]{Arnold:diff} V. I. Arnold,
				``Instability of Dynamical Systems
				  with Several Degrees of Freedom,''
				  {\em Soviet Mathematics-Doklady}
				  \Vol{5}, 581-585 (1964).

\bibitem[Arn78]{Arnold:Math-Meth} V.~I.~Arnold,
				{\em Mathematical Methods of
				Classical Physics}, 
				(Springer-Verlag, New York, 1978).

\bibitem[Aub83a]{Aub:ref} S. Aubry,
				``The twist map, the extended 
				  Frenkel-Kontorova model and the
				  devil's staircase,'' 
				  {\em Physica \/} \Vol{7D},
				  240-258 (1983).

\bibitem[Aub83b]{Aub:cond} S.~Aubry, 
				``Devil's staircase and order
				  without periodicity in classical
				  condensed matter,'' 
				  {\em J. Physique \/}
				  \Vol{44}, 147-162 (1983).

\bibitem[Bang87]{Bang:min}  V.~Bangert
				``Minimal Geodesics,'' preprint (1987).

\bibitem[BGGS80]{BenStr:ex} G.~Benettin, L.~Galgani, A.~Giorgilli and
			    J-M.~Strelcyn,
				``Lyapunov Characteristic Exponents for
				  Smooth Dynamical Systems and for
				  Hamiltonian Systems; a Method for
				  Computing all of Them. Part 2:
				  Numerical Application,''
				  {\em Meccanica \/}
				  \Vol{15}, 21-30 (1980).

\bibitem[Birk22]{Birk:circ} G.D. Birkhoff,
				``Surface transformations and their
				  dynamical  applications,''
				  {\em Acta Mathematica \/}
				  \Vol{43}, 1-119 (1922);
				  reprinted in {\em Collected
				  Mathematical Papers}, vol. II.
				  Amer. Math. Soc.: New York, 1950,
				  pp. 111-229.

\bibitem[Birk27]{Birk:orb} G.D. Birkhoff,
				``On the periodic motions of dynamical 
				  systems,''  {\em Acta Mathematica \/}
				  \Vol{50}, 359-379 (1927).

\bibitem[Bost86]{Bost:KAM} J. Bost,
				``Tores invariants des syst\`{e}ms
				  dynamiques Hamiltoniens,''
				  {\em Asterisque} \Vol{133-134},
				  113-157 (1986).

\bibitem[CC88]{Alessandra} A. Celletti and L. Chierchia,
				``Construction of Analytic KAM
				  Surfaces and Effective Stability 
				  Bounds,''\CMP \Vol{118},
				  119-161 (1988).

\bibitem[CMP87]{CMP:ext} Q. Chen, J.D. Meiss and I.C. Percival,
				``Orbit extension method for finding
				  unstable orbits,'' {\em Physica \/}
				  \Vol{29D}, 143-154 (1987).

\bibitem[Chkv79]{Chkv:Osc} B. Chirikov,
				``A Universal Instability of 
				  Many-Dimensional Oscillator Systems,''
				  {\em Physics Reports \/} \Vol{52} \#5, 
				  263-379 (1979).

\bibitem[FPU55]{Fermi} E.~Fermi, J.~Pasta and S.~Ulam,
				``Studies of Non Linear Problems,''
				  Los Alamos Report LA-1940, May 1955;
				  reprinted in E.~Fermi, {\em Collected
				  Works,} University of Chicago Press,
				  Chicago, (1965), Volume 2, pgs. 978-988.

\bibitem[Fro71]{Fro:int} C.~Froeschl\'{e},
				``On the number of isolating integrals
				  in systems with three degrees of freedom,''
				  {\em Astrophys. Space Sci. \/} \Vol{14}, 
				  110-117 (1971).

\bibitem[Fro72]{Fro:num1} C.~Froeschl\'{e},
				``Numerical Study of a Four-Dimensional 
				  Mapping,''  {\em Astron. \& Astrophys. \/}
				  \Vol{16}, 172-189 (1972).

\bibitem[Fro73]{Fro:num2} C.~Froeschl\'{e} and J.P.~Scheideker,
				``Numerical Study of a Four-Dimensional
				  Mapping,'' { \em Astron. \& Astrophys. \/}
				  \Vol{22}, 431-436 (1973).

\bibitem[Grn79]{Grn:meth} J.M. Greene,
				``A method for determining a stochastic
				  transition,'' \JMP  \Vol{20} \#6, 
				  1183-1201 (1979).

\bibitem[Hed32]{Hed:ex} G.A. Hedlund,
				``Geodesics on a two-dimensional Riemannian
				  manifold with periodic coefficients,''
				  {\em Annals of Mathematics \/} \Vol{33},
				  719-739 (1932).

\bibitem[Herm88]{Herm:non} Michael~R.~Herman,
				``Existence et Non Existence de Tores
				  Inavriants par des Diffeomorphismes
				  Symplectiques,'' Preprint (1988).

\bibitem[Herm83]{Herm:big} Michael~R.~Herman,
				``Sur les courbes invariantes par les
  				  diffeomorphismes de l'anneau, Vol. 1,"
				  {\em Asterisque \/} \Vol{103-104}, (1983).

\bibitem[KnBg85]{Bagley:Diff} K. Kaneko and R. Bagley,
				``Arnold Diffusion, Ergodicity and 
				  Intermittency in a Coupled Standard 
				  Mapping,'' {\em Physics Letters}
				  \Vol{110A} \#9, 435-440, (1985).

\bibitem[Kat82]{Kat:rem} A. Katok,
				``Remarks on Birkhoff and Mather twist
				  map theorems,'' 
				  \ETDS  \Vol{2}, 185-194 (1982). 

\bibitem[Kat88]{Kat:min} A. Katok,
				``Minimal Orbits for Small Perturbations
				  of Completely Integrable Hamiltonian 
				  Systems,'' Preprint (1988).

\bibitem[KB87]{KB:Birk} A. Katok and D. Bernstien,
				``Birkhoff periodic orbits for small 
				  perturbations of completely integrable
				  Hamiltonian systems with 
				  convex Hamiltonians,''
				  {\em Inventiones mathematicae \/} \Vol{88},
				  225-241 (1987).


\bibitem[Khin64]{Khin:num} A.Ya. Khinchin,
				{\em Continued Fractions,\/}
				(University of Chicago Press, Chicago, 
				1964).

\bibitem[KimOst86]{KO}    S. Kim and S. Ostlund,
				``Simultaneous rational approximations in
				  the study of dynamical systems,''
				  {\em Physical Review A \/} \Vol{34} \#4,
				  3426-3434 (1986).

\bibitem[KM88]{KM:orbs} Hyung-tae~Kook and James~D.~Meiss, 
				``Periodic Orbits for Reversible,
				  Symplectic Mappings,'' (1988),
				  to appear in {\em Physica \/} \Vol{D}.

\bibitem[LR88]{LlR:small}  Rafael de la Llave and David Rana, 
				``Accurate Strategies for Small
				  Divisor Problems,'' preprint (1988).

\bibitem[McK88]{MacKay:crit} R.S.~MacKay,
				``A criterion for non-existence of 
				  invariant tori for Hamiltonian 
				  systems,'' (1988),
				  to appear in {\em Physica \/} \Vol{D}.

\bibitem[MMP84]{MMP:tran} R.S.~MacKay, J.D.~Meiss and I.C.~Percival,
				``Transport in Hamiltonian systems,''
				  {\em Physica \/} \Vol{13D}, 55-81 (1984).

\bibitem[MMS89]{MMS:Conv} R.S.~MacKay, J.D.~Meiss and J.~Stark,
				``Converse KAM Theory for Symplectic
				  Twist Maps,'' Preprint (1989).

\bibitem[MP85]{MP:Conv} R.S.~MacKay and I.C.~Percival,
				``Converse KAM : Theory and Practice,"
				  \CMP  \Vol{98}, 469-512 (1985).

\bibitem[Ma82a]{Ma:top} J. Mather,
				``Existence of quasi-periodic orbits for
				  twist maps of the annulus,''
				  {\em Topology \/} \Vol{21} \#4, 
				  457-467 (1982).

\bibitem[Ma82b]{Ma:balls} J. Mather,
				``Glancing billiards,''
				  \ETDS  \Vol{2}, 397-403 (1982).

\bibitem[Ma84]{Ma:circ} J. Mather,
				``Non-existence of invariant circles,''
				  \ETDS  \Vol{4}, 301-311 (1984).

\bibitem[Ma86]{Ma:crit} J. Mather,
				``A criterion for the non-existence 
				  of invariant circles,''
				  {\em Math. Publ. IHES.} \Vol{63}, 
				  153-204 (1986).

\bibitem[Max77]{Maxwell} J.~C.~Maxwell,
				{\em Matter and Motion}, (1877).
				Reprinted by The MacMillan Co.,
				New York, 1920.

\bibitem[MP87]{MP:Newt} B.~Metsel and I.C. Percival,
				``Newton method for highly unstable 
				  orbits,'' {\em Physica \/} 
				  \Vol{24D}, 172-178 (1987).

\bibitem[Moser73]{Moser:book} J.~Moser,
				{\em Stable and Random Motions in 
				Dynamical Systems with Special 
				Emphasis on Celestial Mechanics}, 
				(Princeton University Press,
				Princeton, New Jersey, 1973).
\bibitem[Nekh71]{Nek:diff} N. N. Nekhoroshev
				``Behaviour of Hamiltonian systems 
				  close to integrable,'' {\em
				  Functional Analysis and Applications\/}
				  \Vol{5}, 338-339 (1971).

\bibitem[Osc68]{Osc}  V.I.Oseledec,
				``A Multiplicative Ergodic Theorem:
				  Lyapunov Characteristic Numbers 
				  for Dynamical Systems,'' {\em
				  Trans. Moscow Math. Soc.} \Vol{19},
				  197-231 (1968).

\bibitem[PFTV86]{NumRecp} W.H.~Price, B.P.~Flannery, S.A.~Teukolsky,
			  W.T.~Vetterling, {\em Numerical Recipes,\/}
			  (Cambridge University Press, Cambridge, 1987).

\bibitem[Rana87]{Rana:thesis} D. Rana,
				``Proof of Accurate Upper and Lower
				  Bounds to Stability Domains in 
				  Small Denominator Problems,''
				  PhD thesis, Princeton (1987).

\bibitem[Rob78]{Rob:num} J. Roberts,
				{\em Elementary Number Theory, A Problem
				     Oriented Approach, \/}
				     (MIT Press, Cambridge, Massachusettes,
				      1978).

\bibitem[Smale65]{Smale:horse} S.~Smale,
				``Diffeomorphisms with many periodic 
				  points,'' in 
				  S.~S.~Cairns, ed., 
				  {\em Differential and 
				  Combinatorial Topology}, 
				 (Princeton University Press,
				  Princeton, New Jersey, 1965).

\bibitem[Smale80]{Smale:time} S.~Smale,
				{\em The Mathematics of Time,}
				(Springer-Verlag, New York, 1980).

\bibitem[Strk88]{Stark:ex} J.~Stark, 
				``An Exhaustive Criterion for the
				  Non-Existence of invariant Circles
				  for Area-Preserving Twist Maps,''
	  			  \CMP  \Vol{117}, 177-189 (1988).

\bibitem[Ttch39]{Tit:real}  	E.C.~Titchmarsh,
			    	{\em The Theory of Functions},
			    	(Oxford University Press, Oxford, 1939).

\bibitem[Wig88]{Wiggins:book} S.~Wiggins,
				{\em Global Bifurcations and Chaos,}
				(Springer-Verlag, NewYork, 1988).

\bibitem[Wilb87]{Wilb:circ} J. Wilbrink,
				``Erratic Behavior of Invariant Circles in
				  Standard-like Mappings,'' {\em Physica \/}
				  \Vol{26D}, 358-368 (1987).

\end{thebibliography}
